\documentclass[a4,10pt]{article}
\usepackage{amsmath,amscd,amssymb}
\usepackage{mathrsfs}

\newcommand{\nbiga}{\mathcal{A}}
\newcommand{\nbigb}{\mathcal{B}}
\newcommand{\nbigc}{\mathcal{C}}
\newcommand{\nbigd}{\mathcal{D}}
\newcommand{\nbige}{\mathcal{E}}
\newcommand{\nbigf}{\mathcal{F}}
\newcommand{\nbigg}{\mathcal{G}}
\newcommand{\nbigh}{\mathcal{H}}
\newcommand{\nbigi}{\mathcal{I}}
\newcommand{\nbigj}{\mathcal{J}}

\newcommand{\nbigl}{\mathcal{L}}
\newcommand{\nbigm}{\mathcal{M}}
\newcommand{\nbign}{\mathcal{N}}
\newcommand{\nbigo}{\mathcal{O}}
\newcommand{\nbigp}{\mathcal{P}}

\newcommand{\nbigs}{\mathcal{S}}
\newcommand{\nbigt}{\mathcal{T}}
\newcommand{\nbigu}{\mathcal{U}}
\newcommand{\nbigv}{\mathcal{V}}

\newcommand{\nbigy}{\mathcal{Y}}
\newcommand{\nbigz}{\mathcal{Z}}

\newcommand{\proj}{\mathbb{P}}
\newcommand{\seisuu}{{\mathbb Z}}
\newcommand{\rnum}{{\mathbb Q}}

\newcommand{\cnum}{{\mathbb C}}
\newcommand{\real}{{\mathbb R}}

\newcommand{\EE}{\mathbb{E}}

\newcommand{\KK}{\mathbb{K}}

\newcommand{\gbiga}{\mathfrak A}
\newcommand{\gbigb}{\mathfrak B}

\newcommand{\gbigf}{\mathfrak F}
\newcommand{\gbigg}{\mathfrak G}
\newcommand{\gbigh}{\mathfrak H}
\newcommand{\gbigi}{\mathfrak I}
\newcommand{\gbigj}{\mathfrak J}

\newcommand{\gbign}{\mathfrak N}

\newcommand{\gbigp}{\mathfrak P}
\newcommand{\gbigq}{\mathfrak Q}
\newcommand{\gbigr}{\mathfrak R}

\newcommand{\gbigw}{\mathfrak W}

\newcommand{\gbigz}{\mathfrak Z}

\newcommand{\gminia}{\mathfrak a}
\newcommand{\gminib}{\mathfrak b}
\newcommand{\gminic}{\mathfrak c}

\newcommand{\gminik}{\mathfrak k}

\newcommand{\gminim}{\mathfrak m}

\newcommand{\gminip}{\mathfrak p}

\newcommand{\gminiy}{\mathfrak y}

\newcommand{\veceta}{{\boldsymbol \eta}}

\newcommand{\vecv}{{\boldsymbol v}}

\newcommand{\vecw}{{\boldsymbol w}}

\newcommand{\veca}{{\boldsymbol a}}
\newcommand{\vecb}{{\boldsymbol b}}

\newcommand{\vecdelta}{{\boldsymbol \delta}}
\newcommand{\vecs}{{\boldsymbol s}}

\newcommand{\vecm}{{\boldsymbol m}}

\newcommand{\vecomega}{{\boldsymbol \omega}}
\newcommand{\vecx}{{\boldsymbol x}}

\newcommand{\veczeta}{{\boldsymbol \zeta}}
\newcommand{\vecz}{{\boldsymbol z}}

\newcommand{\vecU}{{\boldsymbol U}}

\newcommand{\vecS}{{\boldsymbol S}}

\newcommand{\lrarr}{\longrightarrow}

\newcommand{\pf}{{\bf Proof}\hspace{.1in}}
\newcommand{\qed}{\mbox{\rule{1.2mm}{3mm}}}

\def\Hom{\mathop{\rm Hom}\nolimits}

\def\Cok{\mathop{\rm Cok}\nolimits}

\def\Image{\mathop{\rm Im}\nolimits}

\def\Re{\mathop{\rm Re}\nolimits}

\def\Gr{\mathop{\rm Gr}\nolimits}

\def\rank{\mathop{\rm rank}\nolimits}

\def\Ker{\mathop{\rm Ker}\nolimits}

\def\Gr{\mathop{\rm Gr}\nolimits}

\def\ord{\mathop{\rm ord}\nolimits}

\def\id{\mathop{\rm id}\nolimits}

\def\gcd{\mathop{\rm g.c.d.}\nolimits}

\def\Irr{\mathop{\rm Irr}\nolimits}

\newcommand{\del}{\partial}

\newcommand{\nhom}{{\mathcal Hom}}

\newcommand{\sankaku}{\triangle}

\newcommand{\zetabar}{\overline{\zeta}}

\newcommand{\Abar}{\overline{A}}

\newcommand{\tildepsi}{\widetilde{\psi}}
\newcommand{\psitilde}{\tildepsi}

\newcommand{\closedopen}[2]{[#1,#2[}
\newcommand{\openclosed}[2]{]#1,#2]}
\newcommand{\openopen}[2]{]#1,#2[}
\newcommand{\closedclosed}[2]{[#1,#2]}

\newcommand{\rhotilde}{\widetilde{\rho}}

\newcommand{\Vhat}{\widehat{V}}
\newcommand{\nablahat}{\widehat{\nabla}}
\newcommand{\Vtilde}{\widetilde{V}}
\newcommand{\nablatilde}{\widetilde{\nabla}}

\newcommand{\wbar}{\overline{w}}

\newcommand{\pitilde}{\widetilde{\pi}}
\newcommand{\Phitilde}{\widetilde{\Phi}}

\newcommand{\ftilde}{\widetilde{f}}

\newcommand{\Ftilde}{\widetilde{F}}

\newcommand{\nbigltilde}{\widetilde{\nbigl}}

\newcommand{\Rtilde}{\widetilde{R}}

\newcommand{\Fhat}{\widehat{F}}

\newcommand{\Rhat}{\widehat{R}}

\newcommand{\Sbar}{\overline{S}}

\newcommand{\Xtilde}{\widetilde{X}}
\newcommand{\vecy}{\boldsymbol y}

\newcommand{\Ltilde}{\widetilde{L}}

\def\ord{\mathop{\rm ord}\nolimits}

\def\op{\mathop{\rm op}\nolimits}
\def\Mod{\mathop{\rm Mod}\nolimits}

\def\Hol{\mathop{\rm Hol}\nolimits}

\def\Sing{\mathop{\rm Sing}\nolimits}

\def\DR{\mathop{\rm DR}\nolimits}

\def\hol{\mathop{\rm hol}\nolimits}

\def\Pro{\mathop{\rm Pro}\nolimits}

\def\Mero{\mathop{\rm Mero}\nolimits}
\def\mero{\mathop{\rm mero}\nolimits}
\def\Bl{\mathop{\rm Bl}\nolimits}
\def\Quad{\mathop{\rm Quad}\nolimits}
\def\sm{\mathop{\rm sm}\nolimits}
\def\RNC{\mathop{\rm RNC}\nolimits}
\def\Sub{\mathop{\rm Sub}\nolimits}
\def\MFV{\mathop{\rm MFV}\nolimits}
\def\Accum{\mathop{\rm Accum}\nolimits}
\def\Crit{\mathop{\rm Crit}\nolimits}
\def\Sol{\mathop{\rm Sol}\nolimits}

\def\Turn{\mathop{\rm Turn}\nolimits}

\newcommand{\Ubar}{\overline{U}}

\newcommand{\Ibar}{\overline{I}}

\newcommand{\psihat}{\widehat{\psi}}

\newcommand{\gtilde}{\widetilde{g}}
\newcommand{\Ztilde}{\widetilde{Z}}
\newcommand{\nbigttilde}{\widetilde{\nbigt}}

\newcommand{\varphitilde}{\widetilde{\varphi}}

\newcommand{\Stilde}{\widetilde{S}}

\newcommand{\Atilde}{\widetilde{A}}

\newcommand{\Deltatilde}{\widetilde{\Delta}}

\newcommand{\Ytilde}{\widetilde{Y}}
\newcommand{\Ybar}{\overline{Y}}

\newcommand{\nbigsbar}{\overline{\nbigs}}

\newcommand{\vecnbigi}{{\boldsymbol \nbigi}}

\newcommand{\Lambdatilde}{\widetilde{\Lambda}}

\newcommand{\nrhom}{R{\mathcal Hom}}
\newcommand{\DDD}{\boldsymbol D}

\newcommand{\Ktilde}{\widetilde{K}}

\newcommand{\Htilde}{\widetilde{H}}

\newcommand{\nutilde}{\widetilde{\nu}}

\newcommand{\gminiatilde}{\widetilde{\gminia}}

\newcommand{\cnumtilde}{\widetilde{\cnum}}

\newcommand{\phitilde}{\widetilde{\phi}}

\newcommand{\Vbar}{\overline{V}}

\newcommand{\gammatilde}{\widetilde{\gamma}}

\newcommand{\Btilde}{\widetilde{B}}

\newcommand{\Bbar}{\overline{B}}

\newcommand{\nbigcbar}{\overline{\nbigc}}

\newcommand{\vecX}{\boldsymbol X}
\newcommand{\vecnbigf}{\boldsymbol \nbigf}

\newcommand{\vecphi}{\boldsymbol \phi}
\newcommand{\vecnbigc}{\boldsymbol \nbigc}

\newcommand{\realc}{\real{\textrm -c}}
\newcommand{\Ecat}{{\sf E}}
\newcommand{\Dcat}{{\sf D}}

\newcommand{\nihom}{{\mathcal{I}hom}}
\newcommand{\vecXtilde}{\boldsymbol \Xtilde}
\newcommand{\vecDelta}{\boldsymbol \Delta}

\newcommand{\realbar}{\overline{\real}}

\newcommand{\vecY}{\boldsymbol Y}
\newcommand{\Subbar}{\overline{\Sub}}

\newcommand{\Khat}{\widehat{K}}
\newcommand{\kappahat}{\widehat{\kappa}}

\newcommand{\pibar}{\overline{\pi}}

\newcommand{\indlim}{``\varinjlim\!\!\textrm{''}}
\newcommand{\vecZ}{\boldsymbol Z}

\newcommand{\Gammabar}{\overline{\Gamma}}
\newcommand{\nbigjbar}{\overline{\nbigj}}
\newcommand{\IC}{\mathrm{I}\cnum}
\newcommand{\II}{\mathbb{I}}
\newcommand{\sfa}{{\sf a}}

\newcommand{\gbiggtilde}{\widetilde{\gbigg}}
\newcommand{\IA}{\mathrm{I}\nbiga}
\newcommand{\vecII}{\pmb{\II}}
\newcommand{\gbigqbar}{\overline{\gbigq}}

\newtheorem{thm}{Theorem}[section]
\newtheorem{cor}[thm]{Corollary}

\newtheorem{rem}[thm]{Remark}
\newtheorem{lem}[thm]{Lemma}
\newtheorem{prop}[thm]{Proposition}
\newtheorem{df}[thm]{Definition}

\newtheorem{condition}[thm]{Condition}
\newtheorem{assumption}[thm]{Assumption}

\newtheorem{notation}[thm]{Notation}

\begin{document}

\title{Curve test for enhanced ind-sheaves 
and holonomic $D$-modules}
\author{Takuro Mochizuki}
\date{}
\maketitle

\begin{abstract}{
Recently, the Riemann-Hilbert correspondence
was generalized
to the context of general holonomic $D$-modules
by A. D'Agnolo and M. Kashiwara.
Namely, they proved that
their enhanced de Rham functor
induces a fully faithful embedding of
the derived category of cohomologically holonomic 
complexes of $D$-modules
into the derived category of 
complexes of real constructible enhanced ind-sheaves.

In this paper, we study a condition
when a complex of real constructible enhanced ind-sheaves $K$
is induced by a cohomologically holonomic complex of $D$-modules.
We characterize such $K$ in terms 
of the restriction of $K$ to holomorphic curves.

\vspace{.1in}
}

\end{abstract}

\section{Introduction}

\subsection{Main result}

In \cite{DAgnolo-Kashiwara1},
D'Agnolo and Kashiwara established
the Riemann-Hilbert correspondence
for holonomic $\nbigd$-modules,
by generalizing the classical Riemann-Hilbert correspondence
\cite{Kashiwara-Riemann-Hilbert, Mebkhout-Riemann-Hilbert, Mebkhout-Riemann-Hilbert2}
between complexes of regular holonomic $\nbigd$-modules
and cohomologically $\cnum$-constructible complexes.
They introduced the concept of
$\real$-constructible enhanced ind-sheaves
on the basis of the theory of ind-sheaves
\cite{Kashiwara-Schapira-ind-sheaves, Kashiwara-Schapira2}.
For any complex manifold $X$,
they constructed the de Rham functor 
$\DR^{\Ecat}_X$
from the derived category of 
cohomologically holonomic complexes of $\nbigd_X$-modules
$\Dcat^b_{\hol}(\nbigd_X)$
to the derived category of $\real$-constructible
enhanced ind-sheaves $\Ecat^b_{\realc}(\IC_X)$,
and they proved that 
$\DR^{\Ecat}_X$ is fully faithful
and compatible with
the $6$-operations and the duality.
They also gave the reconstruction
of a holonomic complex $\nbigm^{\bullet}$
from its solution complex
$\Sol^{\Ecat}(\nbigm^{\bullet})
\in \Ecat^b_{\realc}(\IC_X)$.
They efficiently used 
the study of the formal structure and the asymptotic analysis
of meromorphic flat bundles
(\cite{kedlaya}, \cite{kedlaya2},
 \cite{majima}, \cite{mochi11}, 
 \cite{Mochizuki-wild}, \cite{sabbah4}).
In \cite{DAgnolo-Kashiwara2},
they also introduced the natural perversity condition
for $\Ecat^b_{\realc}(\IC_X)$,
and they proved that $\DR^{\Ecat}_X$ is exact
with respect to the natural $t$-structure of $\Dcat^b_{\hol}(\nbigd_X)$
and the $t$-structure of $\Ecat^b_{\realc}(\IC_X)$
with respect to the perversity condition.

We may still ask an interesting question.
Let $\Ecat^b_{\nbigd}(\IC_X)$
denote the essential image of $\DR^{\Ecat}_X$.
It is natural to ask 
a condition for an object $K\in \Ecat^b_{\realc}(\IC_X)$
to be contained in
$\Ecat^b_{\nbigd}(\IC_X)$.
In the regular singular case,
it is given by the cohomological
 $\cnum$-constructibility condition.
As far as the author knows,
such a clear condition has not yet been given
in the enhanced case.

\vspace{.1in}

In this paper, we study ``a curve test''.
We consider the full subcategory
$\Ecat^b_{\sankaku}(\IC_X)
\subset
 \Ecat^b_{\realc}(\IC_X)$
determined by the following condition
for objects $K\in\Ecat^b_{\realc}(\IC_X)$.
\begin{itemize}
\item
 Set $\Delta:=\{|z|<1\}$.
 Let $\varphi:\Delta\lrarr X$ be any holomorphic map.
 Then,
 $\Ecat\varphi^{-1}(K)
 \in \Ecat^b_{\nbigd}(\IC_{\Delta})$.
\end{itemize}
By the compatibility of 
the de Rham functors $\DR^{\Ecat}$
and $6$-operations,
$\Ecat^b_{\nbigd}(\IC_X)$
is clearly a full subcategory of
$\Ecat^b_{\sankaku}(\IC_X)$.
The following is the main theorem
of this paper.

\begin{thm}[Theorem \ref{thm;16.7.11.1}]
\label{thm;18.11.16.10}
$\Ecat^b_{\sankaku}(\IC_X)$
is equal to
$\Ecat^b_{\nbigd}(\IC_X)$.
\end{thm}

\begin{rem}
Recently, Kuwagaki {\rm\cite{Kuwagaki}} introduced
another approach to the irregular Riemann-Hilbert correspondence.
\hfill\qed
\end{rem}

\subsection{Meromorphic flat connections and enhanced ind-sheaves}

A holonomic $\nbigd$-module
can be locally described as 
the gluing of meromorphic flat connections
given on subvarieties.
Hence, it is a key step to study 
such a characterization 
for meromorphic flat connections
(Theorem \ref{thm;16.5.11.20}).

Let $X$ be an $n$-dimensional complex manifold.
Let $H$ be a normal crossing hypersurface of $X$.
Let $\vecX(H)$ denote the bordered space $(X\setminus H,X)$
in the sense of \cite{DAgnolo-Kashiwara1, DAgnolo-Kashiwara2}.
(A bordered space is defined to be
a pair $(M,\check{M})$ of a good topological space 
$\check{M}$ and an open subset $M$,
where a topological space is called good 
if it is Hausdorff, locally compact, 
countable at infinity and has finite flabby dimension.
See \cite[Page 86]{DAgnolo-Kashiwara1}.)
Suppose that
an object $K\in \Ecat^b_{\realc}(\IC_{\vecX(H)})$
satisfies the following condition.
\begin{itemize}
\item
 $K_{|X\setminus H}$ is induced by a local system on $X\setminus H$.
\item
 Let $\varphi:\Delta\lrarr X$ be any holomorphic map
 such that
 $\varphi(\Delta\setminus\{0\})\subset X\setminus H$
and $\varphi(0)\in H$.
Then,
 $\Ecat\varphi^{-1}(K)$
 comes from a meromorphic flat bundle 
 on $(\Delta,0)$.
\end{itemize}
We would like to prove that there exists
a meromorphic flat connection
$(V,\nabla)$ on $(X,H)$ with an isomorphism
$\DR^{\Ecat}_{\vecX(H)}(V)[-n]\simeq K$
in $\Ecat^b_{\realc}(\IC_{\vecX(H)})$.
Once we obtain a meromorphic flat connection
$(V,\nabla)$  on $(X,H)$
such that 
$\Ecat\varphi^{-1}(K)
\simeq
 \DR^{\Ecat}_{\vecDelta(0)}\varphi^{\ast}(V,\nabla)[-1]$
in a natural way
for any $\varphi:\Delta\lrarr X$ as above,
then it is not so difficult to prove that
$K\simeq
 \DR^{\Ecat}_{\vecX(H)}(V,\nabla)[-n]$
(Proposition \ref{prop;16.6.23.20}).
Therefore, we would like to construct
such a meromorphic connection $V$.

\vspace{.1in}

We explain a brief outline for the construction of such $V$
in the essential case $n=2$ (Proposition \ref{prop;16.5.10.10}).
Let us introduce an auxiliary condition.
(See \S\ref{subsection;16.7.25.30}
for more details.)
Suppose that we are given
a {\em good} set of ramified irregular values
 $\nbigi_P$ at a smooth point $P\in H$
 and a multiplicity function
 $\gminim_P:\nbigi_P\lrarr \seisuu_{\geq 0}$
 such that the following holds.
\begin{description}
\item[(Condition \ref{condition;18.11.26.3})]
For any holomorphic map $\varphi:\Delta\lrarr X$
such that
$\varphi(\Delta\setminus\{0\})\subset X\setminus H$
and $\varphi(0)$ is close to $P$ in $H$,
$\Irr(\Ecat\varphi^{-1}K)
\simeq
 \varphi^{\ast}\nbigi_P$
holds,
and the multiplicity functions are the same.
\end{description}
Then, we can construct a good meromorphic flat bundle
$(V,\nabla)$ on a neighbourhood $X_P$ of $P$ in $X$,
with an isomorphism
$\DR^{\Ecat}_{\vecX_P(H\cap X_P)}(V)[-2]
\simeq K$
(Proposition \ref{prop;16.7.25.40}).
There exists a similar condition at cross points of $H$.
Although Condition \ref{condition;18.11.26.3}
is not always satisfied even if $K$ comes from a meromorphic flat connection,
it is our strategy to modify 
$X$ in a birational way so that
Condition \ref{condition;18.11.26.3} is satisfied.

There exists a stratification
$X\setminus H=\coprod\nbigc$
by locally closed subanalytic subsets
such that
there exist subanalytic functions
$h^{\nbigc}_i$ $(i=1,\ldots,r)$ on $(\nbigc,X)$
and an isomorphism
$\pi^{-1}(\cnum_{\nbigc})
\otimes K
\simeq
 \bigoplus
 \cnum^{\Ecat}\overset{+}{\otimes}
 \cnum_{t\geq h^{\nbigc}_i}$.
Here, $r$ is the rank of the local system $K_{|X\setminus H}$.
By using such a local description of $K$,
we can prove that there exists
a $0$-dimensional closed subanalytic subset
$Z\subset H$
such that 
any $P\in H\setminus Z$ is a smooth point of $H$,
and that
Condition \ref{condition;18.11.26.3}
is satisfied for $K$ at $P$.
As a result, we obtain a meromorphic flat bundle
$(V,\nabla)$ on $(X\setminus Z,H\setminus Z)$
with an isomorphism
$\DR^{\Ecat}(V)[-2]
\simeq
 K_{|X\setminus Z}$
(Proposition \ref{prop;16.6.26.10}).

If $Z$ is not empty,
we define $\psi_1:X_1\lrarr X$ 
as the complex blowing up 
at the points of $Z$,
and we set $H_1:=\psi_1^{-1}(H)$.
As in the previous stage,
there exists the $0$-dimensional closed subanalytic subset
$Z_1\subset \psi_1^{-1}(Z)\subset H_1$
such that 
there exists a meromorphic flat bundle
$(V_1,\nabla)$ 
on $(X_1\setminus Z_1,H_1\setminus Z_1)$
with an isomorphism
$\DR^{\Ecat}(V_1,\nabla)[-2]
\simeq \Ecat\psi_1^{-1}(K)_{|X_1\setminus Z_1}$.
If $Z_1$ is not empty,
we define again 
$\psi_2:X_2\lrarr X_1$ 
as the complex blowing up at $Z_1$.
By continuing the procedure successively,
we obtain the following sequence:
\[
 \cdots \lrarr
 (X_{\ell},H_{\ell})
\stackrel{\psi_{\ell}}{\lrarr}
 (X_{\ell-1},H_{\ell-1})
\stackrel{\psi_{\ell-1}}{\lrarr}
\cdots
\lrarr
 (X_1,H_1)
\stackrel{\psi_1}{\lrarr}
 (X,H).
\]
Here, 
for each $j$,
there exists a $0$-dimensional closed subanalytic subset
$Z_j\subset H_j$
such that 
Condition \ref{condition;18.11.26.3}
is satisfied for
$\Ecat(\psi_1\circ\cdots\psi_{j})^{-1}K$
at $P\in H_j\setminus Z_j$,
and 
$\psi_{j+1}$ are the blowings up
along $Z_j$.
Moreover, $Z_{j+1}$ are contained in
$(\psi_{j+1})^{-1}(Z_j)$.
We would like to prove that this construction will stop 
after finite steps,
i.e.,
$Z_{\ell}$ can be empty if $\ell$ is large.
Note that this type of issue also appeared in the study of
Sabbah's conjecture,
i.e., the higher dimensional generalization of 
the Hukuhara-Levelt-Turrittin theorem.
(See \cite{kedlaya, kedlaya2, mochi6, Mochizuki-wild, sabbah4}.)

For the above problem,
it is important to study the pull back of subanalytic functions
by the composition of $\psi_p$.
(See \cite{Bierstone-Milman} and
\cite{Hironaka-Introduction-real-analytic}
for the general theory of subanalytic sets.)
Let $\nbigc$ be an open strata in the stratification.
Let $\varpi_j:\Xtilde_j(H_j)\lrarr X_j$
be the oriented real blowing up of $X_j$ along $H_j$.
Let $Q_j$ be a point of $\varpi_j^{-1}(Z_j)$,
and let $\nbigu_{Q_j}$ be a small neighbourhood 
of $Q_j$ in $\Xtilde_j(H_j)$.
Let $\Phi_j:\Xtilde_j(H_j)\lrarr X$
be the morphism induced by 
$\varpi_j$ and the composition of $\psi_p$.
If $\Phi_j(\nbigu_{Q_j})\cap \nbigc\neq\emptyset$,
we are lead to study the functions 
$\Phi_j^{\ast}(h^{\nbigc}_i)$
and 
$\Phi_j^{\ast}(h^{\nbigc}_i-h^{\nbigc}_k)$
on $\Phi_j^{-1}(\nbigc)\cap\nbigu_{Q_j}$,
and to prove that
$\Phi_j^{\ast}(h^{\nbigc}_i)$ 
and 
$\Phi_j^{\ast}(h^{\nbigc}_i-h^{\nbigc}_k)$ 
have nice properties.
However, it is not easy to deal with 
such functions directly.
Thus, we are invited to study
a lifting problem with respect to
sequences of local real blowings up,
which we shall explain in the next subsection.

\subsection{Rectilinearization and a lifting problem}

Let $M$ be an $m$-dimensional real analytic manifold.
Let $U$ be an open subanalytic subset in $M$.
Let $f$ be a real analytic function on $U$
which is subanalytic on $(U,M)$.
It is not easy to study directly the behaviour of $f$
around boundary points of $U$.
The fundamental and useful tools are rectilinearization theorems
developed by Hironaka, Bierstone-Milman and Parusi\'{n}ski.
According to the rectilinearization theorems
for subanalytic subsets
\cite{Bierstone-Milman, Hironaka-Introduction-real-analytic}
and subanalytic functions \cite{Parusinski},
there exists a locally finite family of
real analytic maps
$\phi_{\alpha}:W_{\alpha}\lrarr M$ $(\alpha\in\Lambda)$
such that the following holds.
\begin{itemize}
\item
Each $W_{\alpha}$ is equipped with
a real analytic coordinate system
$(x_{\alpha,1},\ldots,x_{\alpha,m})$,
which induces
$W_{\alpha}\simeq\real^m$.
\item
Each $\phi_{\alpha}$ is the composition
of local real blowings up.
\item
$\phi_{\alpha}^{-1}(U)\subset W_{\alpha}\simeq\real^m$ is the union of
the quadrants contained in $U$.
\item
The restriction of $\phi_{\alpha}^{\ast}(f)$
to each $m$-dimensional quadrant 
$\gbigq\subset\phi_{\alpha}^{-1}(U)$ 
is constantly $0$,
or expressed as 
$a_{\alpha,\gbigq}
\prod_{i=1}^{m}|x_{\alpha,i}|^{r_{\alpha,i}}$
for some $(r_{\alpha,i})\in\rnum^m$
and a nowhere vanishing ramified analytic function $a_{\alpha,\gbigq}$
on the closure of $\gbigq$.
\end{itemize}
(See \S\ref{subsection;18.11.16.2} for more precise.)

We would like to apply the rectilinearization theorem
to the study of the above functions
$\Phi_j^{\ast}(h_i^{\nbigc})$
and 
$\Phi_j^{\ast}(h_i^{\nbigc}-h_k^{\nbigc})$.
There exists a rectilinearization
$\phi_{\alpha}:W_{\alpha}\lrarr \Xtilde(H)$ $(\alpha\in\Lambda)$
for $h^{\nbigc}_i$
and $h^{\nbigc}_i-h^{\nbigc}_k$.
If there exists a real analytic map
$\Psi_j:\nbigu_{Q_j}\lrarr W_{\alpha_0}$
for some $\alpha_0\in\Lambda$
such that 
$\phi_{\alpha_0}\circ\Psi_j$
is equal to the restriction of $\Phi_j$ to $\nbigu_{Q_j}$,
then 
we can deduce that
$\Phi_j^{\ast}(h_i^{\nbigc})
=\Psi_j^{\ast}\bigl(
 \phi_{\alpha_0}^{\ast}(h_i^{\nbigc})
 \bigr)$
and 
$\Phi_j^{\ast}(h_i^{\nbigc}-h^{\nbigc}_k)
=\Psi_j^{\ast}\bigl(
 \phi_{\alpha_0}^{\ast}(h_i^{\nbigc}-h^{\nbigc}_k)
 \bigr)$
have nice properties 
because 
$\phi_{\alpha_0}^{\ast}(h_i^{\nbigc})$
and 
$\phi_{\alpha_0}^{\ast}(h_i^{\nbigc}-h^{\nbigc}_k)$
have nice properties.
Thus, we are lead to study
the lifting problem
whether there exists such a real analytic map $\Psi_j$.
Indeed, it is one of the main issues studied
in this paper.
We shall study the problem in 
\S\ref{section;16.9.4.1}--\ref{section;16.9.4.3}
by using the technical preliminaries in \S\ref{section;16.9.4.2}.
There are sophisticated studies 
about infinite sequences of complex blowing up of surfaces
(for example, see \cite{Favre-Jonsson}).
However, for the purpose of the above lifting problem,
at this moment,
it looks more convenient to adapt a more naive approach
by considering the families of curves around the centers of blowings up
and their limit curve.

\paragraph{Acknowledgement}

I am grateful to Masaki Kashiwara
for asking this question and for many discussions.
This study grew from my effort to understand
the interesting works of
Andrea D'Agnolo,
Kashiwara and Pierre Schapira
on enhanced ind-sheaves and
the generalized Riemann-Hilbert correspondence.
I also thank Giovanni Morando
who first attracted my attention to 
the theory of ind-sheaves.
I thank Claude Sabbah for his kindness
and for discussions on many occasions.
I am grateful to 
Akira Ishii and Yoshifumi Tsuchimoto
for their constant encouragement.
I heartily thank the referees 
for their careful and patient readings
and for their numerous helpful comments
to improve and clarify this manuscript.

This work was partially supported by 
the Grant-in-Aid for Scientific Research (S) (No. 17H06127),
the Grant-in-Aid for Scientific Research (S) (No. 16H06335),
the Grant-in-Aid for Scientific Research (S) (No. 24224001),
and the Grant-in-Aid for Scientific Research (C) (No. 15K04843), 
Japan Society for the Promotion of Science.

\part{Preliminaries}

\section{Subanalytic geometry}

We give preliminaries from subanalytic geometry.
We mention 
\cite{Bierstone-Milman}
as a useful general reference for the subject.
In \S\ref{subsection;18.11.16.40},
we prepare some notation.
In \S\ref{subsection;18.11.16.2},
we recall the useful existence theorems of
rectilinearizations for subanalytic subsets
(Theorem \ref{thm;16.4.13.1})
and subanalytic functions
(Theorem \ref{thm;16.7.20.2}).
We also state some consequences of the theorems
in a way convenient to this study.
In \S\ref{subsection;18.11.16.41},
we collect some technical results
which explain that subanalytic maps and subanalytic functions
have nice properties outside of small subsets.
In \S\ref{subsection;18.11.16.42},
complementary statements are collected.

\subsection{Some notation}
\label{subsection;18.11.16.40}

For any $\epsilon>0$,
we set $\II_{\epsilon}:=\{0\leq t<\epsilon\}$
and $\II^{\circ}_{\epsilon}:=\II_{\epsilon}\setminus\{0\}$.
If $\epsilon=1$, we use $\II$ and $\II^{\circ}$
instead of $\II_{1}$ and $\II_1^{\circ}$,
respectively.
More generally, for any $a<b$,
intervals are denoted as
$\closedclosed{a}{b}:=\{a\leq t\leq b\}$,
$\openopen{a}{b}:=\{a<t<b\}$,
$\openclosed{a}{b}:=\{a<t\leq b\}$
and 
$\closedopen{a}{b}:=\{a\leq t< b\}$.

We set $\realbar:=\real\cup\{\pm\infty\}$.
We regard it as real analytic manifold with boundary 
by the coordinate $(\pm t)^{-1}$ around $\pm\infty$.

Let $A$ be any subanalytic subset in a real analytic manifold $M$.
For each $k\in\seisuu_{\geq 0}$,
let $A_k^{\sm}$ denote the set of
the $k$-dimensional smooth points of $A$.
According to \cite[Theorem 7.2]{Bierstone-Milman},
$A_k^{\sm}$ is a subanalytic open subset in $A$.
We set
$A^{\sm}=\bigcup_{k\geq 0} A_k^{\sm}$.
The set $\Sing(A):=A\setminus A^{\sm}$
is called the singular locus of $A$,
which is a closed subset of $A$.
Let $\Abar$ denote the closure of $A$ in $M$,
and $\del A:=\Abar\setminus A$.

Let $N$ be a real analytic manifold.
Let $B$ be a subanalytic subset of $N$.
A map $F:A\lrarr B$ is called subanalytic on $(A,M)$ if
the graph is subanalytic in $M\times N$.
For subanalytic subsets 
$A_i\subset A$ $(i\in\Lambda)$ and $B_i\subset B$ $(i\in\Lambda)$,
a subanalytic map 
$F:(\{A_i\,|\,i\in\Lambda\},A)\lrarr 
 (\{B_i\,|\,i\in\Lambda\},B)$
means a subanalytic map $F:A\lrarr B$
such that $F(A_i)\subset B_i$ $(i\in\Lambda)$.
If $A$ is a real analytic manifold,
a real analytic map 
$F:(\{A_i\,|\,i\in\Lambda\},A)\lrarr 
 (\{B_i\,|\,i\in\Lambda\},B)$
means a real analytic map $F:A\lrarr N$
such that $F(A_i)\subset B_i$ and $F(A)\subset B$.

A function $f:A\lrarr \real$ is called a subanalytic function
on $(A,M)$
if the graph of $f$ is a subanalytic subset of
$M\times \proj^1(\real)$.
Note that a function $f:A\lrarr\real$ is regarded as
as a map $A\lrarr\proj^1(\real)$ such that $f(A)\subset\real$.
If moreover $f:A\lrarr\real$ is continuous,
then $f$ is called continuous subanalytic on $(A,M)$.
Similarly, if moreover $f:A\lrarr\real$ is real analytic,
then $f$ is called real analytic and subharmonic on $(A,M)$.

Let $P$ be a point of $\del A$
at which $\Abar$ is a real analytic submanifold of $M$
with corner,
i.e.,  there exists a neighbourhood
$(\nbigu,t_1,\ldots,t_n)$ of $M$ around $P$
such that $\Abar\cap\nbigu=\bigcup_{i=1}^{\ell}\{t_i\geq 0\}$
for some $\ell$.
A real analytic function $f$
defined on an open subset $U\subset M$
is called ramified real analytic around $P\in U$
if $f$ is expressed as a convergent power series
\[
 f=\sum_{j_1\geq -N}
 \cdots
 \sum_{j_{\ell}\geq -N}
 a_{j_1,\ldots,j_{\ell}}(t_{\ell+1},\ldots,t_n)
 \prod_{i=1}^{\ell}(t_i^{1/e})^{j_i}
\]
for an integer $N$, a positive integer $e$,
and real analytic functions $a_{j_1,\ldots,j_{\ell}}$.

\subsection{Rectilinearization}
\label{subsection;18.11.16.2}

\subsubsection{Rectilinearization of subanalytic subsets}
\label{subsection;16.7.22.1}

A subset $\gbigq\subset\real^n$ is called
a quadrant if there exists a decomposition
 $\{1,\ldots,n\}=I_0\sqcup I_+\sqcup I_-$
 for which
\[
 \gbigq=\bigl\{
 (x_1,\ldots,x_n)\,\big|\,
 x_i=0\,\,(i\in I_0),\,\,
 x_i> 0\,\,(i\in I_+),\,\,
 x_i< 0\,\,(i\in I_-)
 \bigr\}.
\]
\begin{df}
A subset $S\subset \real^n$ is called rectilinearized
if $S$ is the union of the quadrants contained in $S$.
\hfill\qed
\end{df}

Let $B$ be any subanalytic subset in
an $n$-dimensional real analytic manifold $M$.
A rectilinearization of $B$ is a locally finite family of 
real analytic maps
$\phi_{\alpha}:W_{\alpha}\lrarr M$ $(\alpha\in\Lambda)$
with the following property.
\begin{itemize}
\item
 Each $W_{\alpha}$ is equipped with
 a real analytic coordinate system $(x_1,\ldots,x_n)$
 which induces an isomorphism
 $W_{\alpha}\simeq\real^n$.
\item
 Each $\phi_{\alpha}$ is the composition of
 a finite sequence of local real blowings up 
 along smooth real analytic centers.
 Namely,
 $\phi_{\alpha}$ is factorized as follows:
\begin{equation}
\label{eq;16.4.13.2}
 W_{\alpha}=W_{\alpha}^{(k(\alpha))}
 \stackrel{\phi_{\alpha}^{(k(\alpha))}}\lrarr 
 W_{\alpha}^{(k(\alpha)-1)}
\stackrel{\phi_{\alpha}^{(k(\alpha)-1)}}\lrarr
\cdots
\stackrel{\phi_{\alpha}^{(2)}}{\lrarr}
 W_{\alpha}^{(1)}
\stackrel{\phi_{\alpha}^{(1)}}{\lrarr}
 W_{\alpha}^{(0)}
=M.
\end{equation}
Moreover, for each $\ell<k(\alpha)$,
there exist a subanalytic open subset
$U_{\alpha}^{(\ell)}\subset
 W_{\alpha}^{(\ell)}$
and a closed real analytic submanifold
$C_{\alpha}^{(\ell)}\subset U_{\alpha}^{(\ell)}$
so that
$\phi_{\alpha}^{(\ell+1)}$
is the real blowing up of $U_{\alpha}^{(\ell)}$
along $C_{\alpha}^{(\ell)}$.
We also impose that
$C_{\alpha}^{(\ell)}$ is subanalytic
in $W_{\alpha}^{(\ell)}$.
Note that
$C_{\alpha}^{(\ell)}$ can be empty.
\item
There exist compact subsets $K_{\alpha}\subset W_{\alpha}$ 
 such that
 $\bigcup_{\alpha\in \Lambda}
 \phi_{\alpha}(K_{\alpha})=M$.
\item
 For each $\alpha$,
 $\phi_{\alpha}^{-1}(B)$ 
 is rectilinearized
with respect to the coordinate system.
\end{itemize}
Let us recall the following fundamental theorem
due to Hironaka \cite{Hironaka-Introduction-real-analytic}.
(See also \cite{Bierstone-Milman}
and \cite{Parusinski}.)

\begin{thm}
\label{thm;16.4.13.1}
For any subanalytic subset $B$ in $M$,
there exists a rectilinearization of $B$.
\hfill\qed
\end{thm}

\subsubsection{Ramified normal crossing functions on quadrants}

Recall that a real analytic function $g$ on $\real^n$
is called normal crossing
if $g=\prod_{i=1}^nx_i^{m_i}\times g_0$,
where $m_i$ are non-negative integers,
and $g_0$ is a nowhere vanishing real analytic function
on $\real^n$.

We introduce a variant of the concept for functions on quadrants.
Let $\gbigq$ be any $n$-dimensional quadrant in $\real^n$.
A real analytic function $f$ on $\gbigq$
is called normal crossing
if $f=\prod_{i=1}^nx_i^{m_i}\times g$,
where $g$ is a nowhere vanishing real analytic function on 
a neighbourhood of the closure $\overline{\gbigq}$,
and $(m_1,\ldots,m_n)\in \seisuu_{\geq 0}^n$.

\vspace{.1in}

Let $\gbigq_0\subset\real^n$ denote
the $n$-dimensional quadrant defined as
$\{x_i>0\,\,(i=1,\ldots,n)\}$.
Let $\{1,\ldots,n\}=I_+\sqcup I_-$ be a decomposition.
Let $\gbigq(I_+,I_-)$ be the $n$-dimensional quadrant in $\real^n$
defined as
$\gbigq(I_+,I_-)=\{x_i>0\,\,(i\in I_+),\,\,x_i<0\,\,(i\in I_-)\}$.
For any tuple of positive integers
 $(\rho_1,\ldots,\rho_n)$,
let $F:\overline{\gbigq_0}\lrarr \overline{\gbigq(I_+,I_-)}$
be the homeomorphism defined by
$F(y_1,\ldots,y_n)=
 (\epsilon_1y_1^{\rho_1},\ldots,\epsilon_ny_n^{\rho_n})$,
 where $\epsilon_i=1$ $(i\in I_+)$ and 
 $\epsilon_i=-1$ $(i\in I_-)$.
Such a homeomorphism is called
a ramified real analytic isomorphism
of $\overline{\gbigq_0}$ and 
$\overline{\gbigq(I_+,I_-)}$.

A real analytic function $f$ on 
an $n$-dimensional quadrant $\gbigq$ is called
ramified normal crossing
if there exists a ramified real analytic isomorphism
$F:\overline{\gbigq_0}\lrarr \overline{\gbigq}$
such that $F^{\ast}(f)$ is normal crossing.
There exists the continuous extension
$\overline{f}$ of $f$ on $\overline{\gbigq}$.
For any quadrant $\gbigq'$ contained in $\overline{\gbigq}$,
if the restriction $\overline{f}_{|\gbigq'}$ is not constantly $0$,
then $\overline{f}_{|\gbigq'}$ is 
a ramified normal crossing function on $\gbigq'$,
where we regard $\gbigq'$ as a quadrant 
in $\real^{\dim \gbigq'}$.

Let $\RNC_+(\gbigq)$ denote the set of 
ramified normal crossing functions on $\gbigq$.
We set
$\RNC_-(\gbigq):=\bigl\{
 1/f\,\big|\, f\in\RNC_+(\gbigq)
 \bigr\}$
and 
$\RNC(\gbigq):=
\RNC_+(\gbigq)\cup\RNC_-(\gbigq)\sqcup\{0\}$.

\vspace{.1in}
Let $\real^r_{\vecy}=\{(y_1,\ldots,y_r)\in\real^r\}$
and $\real^n_{\vecx}=\{(x_1,\ldots,x_n)\in\real^n\}$.

\begin{lem}
\label{lem;16.7.20.1}
Let $\phi:\real^r_{\vecy}\lrarr\real^n_{\vecx}$
be a real analytic map such that 
each $\phi^{\ast}(x_i)$ $(i=1,\ldots,n)$ is normal crossing
or constantly $0$.
Let $\gbigq_1$  be a quadrant of $\real^n_{\vecx}$.
Then, the following holds.
\begin{itemize}
\item
$\phi^{-1}(\gbigq_1)$ is rectilinearized.
\item
 Let $f\in\RNC(\gbigq_1)$,
 where we naturally regard $\gbigq_1$ as a quadrant
 in $\real^{\dim \gbigq_1}$.
 Let $\gbigq_2$ be any quadrant  of $\real^r_{\vecy}$
 such that $\phi(\gbigq_2)\subset \gbigq_1$.
 Then, 
 $\phi^{\ast}(f)_{|\gbigq_2}$ is contained in
 $\RNC(\gbigq_2)$.
\end{itemize}
\end{lem}
\pf
By the assumption,
one of the following holds
for each $i=1,\ldots,n$;
(i) $\phi^{\ast}(x_i)=0$,
(ii) $\phi^{\ast}(x_i)=
 a_i\cdot\prod_{k=1}^ry_k^{m(i)_k}$
where $a_i$ is nowhere vanishing on $\real^r_{\vecy}$,
and $m(i)_k\in\seisuu_{\geq 0}$.

For any quadrant $\gbigq$ in $\real^r_{\vecy}$,
one of the following holds;
(i) $\prod_{k=1}^ry_k^{m(i)_k}$
are constantly $0$ on $\gbigq$,
(ii) $\prod_{k=1}^ry_k^{m(i)_k}$ is positive on $\gbigq$,
(iii) $\prod_{k=1}^ry_k^{m(i)_k}$ is negative on $\gbigq$.
Hence,
one of the following holds;
(i) $\phi(\gbigq)\subset \gbigq_1$,
(ii) $\phi(\gbigq)\cap \gbigq_1=\emptyset$.
Then, we obtain the first claim.

Let us study the second claim.
It is enough to consider the case
$f\in\RNC_+(\gbigq_1)$.
We may assume 
that $\gbigq_1$ is 
$\bigcap_{i=1}^{\ell}\{x_i>0\}
\cap
 \bigcap_{i=\ell+1}^n\{x_i=0\}$,
and that
$\gbigq_2$ is 
$\bigcap_{k=1}^{p}\{y_k>0\}
\cap
 \bigcap_{k=p+1}^r\{y_k=0\}$.
For any $i=1,\ldots,\ell$,
because
$\phi^{\ast}(x_i)> 0$
on $\gbigq_2$,
we obtain $a_i>0$ on $\gbigq_2$
and $m(i)_k=0$ for $k=p+1,\ldots,r$.

We naturally regard 
$\real_{\vecx}^{\ell}
=\{(x_1,\ldots,x_{\ell})\}$
as a subspace of $\real_{\vecx}^n$.
Then, $\gbigq_1$ is 
a quadrant in $\real^{\ell}_{\vecx}$.
Similarly,
we naturally regard 
$\real_{\vecy}^{p}
=\{(y_1,\ldots,y_p)\}$
as a subspace of $\real_{\vecy}^r$,
and 
$\gbigq_2$ is a quadrant of $\real_{\vecy}^p$.

There exists a positive integer $\rho$
with the following property.
\begin{itemize}
\item
We define the morphism
 $\Phi:\real_{\vecx}^{\ell}\lrarr \real_{\vecx}^{\ell}$ 
by 
$\Phi(x_1,\ldots,x_{\ell})=(x_1^{\rho},\ldots,x_{\ell}^{\rho})$.
Let $\Phi_{\gbigq_1}:\gbigq_1\lrarr \gbigq_1$
denote the induced map.
Then,
$\Phi_{\gbigq_1}^{\ast}(f)$ is a real analytic normal crossing function
on $\gbigq_1$.
\end{itemize}
For $i=1,\ldots,\ell$,
we obtain
$\phi^{\ast}(x_i^{1/\rho})_{|\real_{\vecy}^{p}}
=a^{1/\rho}_{i|\real_{\vecy}^{p}}\prod_{k=1}^p y_k^{m(i)_k/\rho}$.
We define $\Psi:\real_{\vecy}^{p}\lrarr\real_{\vecy}^p$ by
$\Psi(y_1,\ldots,y_p)=(y_1^{\rho},\ldots,y_{p}^{\rho})$.
Then,
$\Psi^{\ast}\phi^{\ast}(x_i^{1/\rho})$
are the restriction of 
normal crossing real analytic functions
on $\real^p_{\vecy}$.
Hence, there exists a real analytic map
$F:\real^p_{\vecy}\lrarr \real^{\ell}_{\vecx}$
such that
$F(\gbigq_2)\subset \gbigq_1$,
and $\Phi\circ F=\phi\circ \Psi$.
It is easy to see that $F^{\ast}(x_i)$ 
$(i=1,\ldots,\ell)$ are normal crossing.
Let $\Psi_{\gbigq_2}$ and $F_{\gbigq_2}$
denote the restriction of $\Psi$ and $F$ to $\gbigq_2$,
respectively.
Then, 
$F_{\gbigq_2}^{\ast}\Phi_{\gbigq_1}^{\ast}(f_{|\gbigq_1})
=\Psi_{\gbigq_2}^{\ast}\phi^{\ast}(f_{|\gbigq_1})$
is a real analytic normal crossing function on $\gbigq_2$.
Thus, we obtain the second claim.
\hfill\qed

\vspace{.1in}

The following lemma is a consequence of
\cite[Theorem 4.4]{Bierstone-Milman}.
\begin{lem}
\label{lem;16.7.20.3}
Let $\gbigq_1$ be any  quadrant  of $\real^n$.
Let $f$ be an element of $\RNC(\gbigq_1)$,
and let $g_1,\ldots,g_m$ be
real analytic functions on $\real^n$
such that each $g_i$ is not constantly $0$.
Then, there exists a rectilinearization 
$\{(W_{\alpha},\phi_{\alpha})\,|\,\alpha\in\Lambda\}$
of $\gbigq_1$ 
with the following property.
\begin{itemize}
\item
$\phi_{\alpha}^{\ast}(f)_{|\gbigq}\in\RNC(\gbigq)$
for any quadrant $\gbigq$
contained in $\phi_{\alpha}^{-1}(\gbigq_1)$.
\item
$\phi_{\alpha}^{\ast}(g_j)$
are normal crossing functions
on $W_{\alpha}$.
\end{itemize}
\end{lem}
\pf
Set $h:=\prod_{i=1}^n x_i\times \prod_{j=1}^mg_j$.
By \cite[Theorem 4.4]{Bierstone-Milman},
there exists a locally finite family of real analytic maps
$\phi_{\alpha}:W_{\alpha}\lrarr \real^n$
$(\alpha\in\Lambda)$
with the following property.
\begin{itemize}
\item
Each $\phi_{\alpha}$ is the composition of
a finite sequence of local real blowings up.
\item
Each $W_{\alpha}$ is equipped with
a coordinate system $(y_1,\ldots,y_n)$
which induces an isomorphism
$W_{\alpha}\simeq\real^n$,
and 
$\phi_{\alpha}^{-1}(h)$
is normal crossing on $W_{\alpha}$.
\item
There exist compact subsets $K_{\alpha}\subset W_{\alpha}$
such that
$\real^n=\bigcup \phi(K_{\alpha})$.
\end{itemize}
Because 
$\phi_{\alpha}^{-1}(h)$
is normal crossing,
$\phi_{\alpha}^{\ast}(x_i)$
and $\phi_{\alpha}^{\ast}(g_j)$
are normal crossing.
By Lemma \ref{lem;16.7.20.1},
$\{(W_{\alpha},\phi_{\alpha})\,|\,\alpha\in\Lambda\}$
is a rectilinearization of 
the quadrant $\gbigq_1$ of $\real^n$.
We also obtain from Lemma \ref{lem;16.7.20.1} that 
$\phi_{\alpha}^{\ast}(f)_{|\gbigq}\in\RNC(\gbigq)$
for any quadrant $\gbigq\subset\phi_{\alpha}^{-1}(\gbigq_1)$.
\hfill\qed

\vspace{.1in}
We give a remark on rectilinearizations of subanalytic subsets.
\begin{prop}
\label{prop;18.1.14.10}
Let $B_1,\ldots,B_m$ be subanalytic subsets
in a real analytic manifold $M$.
There exists a locally finite family of
real analytic maps
$\{(W_{\alpha},\phi_{\alpha})\,|\,\alpha\in\Lambda\}$
which is a rectilinearization of each $B_i$.
\end{prop}
\pf
We use an induction on $m$.
Suppose that we have already obtained
a locally finite family of real analytic maps
$\{(W_{\alpha},\phi_{\alpha})\,|\,\alpha\in\Lambda\}$
which is a rectilinearization of $B_i$ $(i=1,\ldots,m-1)$.
There exists subanalytic compact subsets
$K_{\alpha}\subset W_{\alpha}$
such that
$\bigcup_{\alpha}\phi_{\alpha}(K_{\alpha})=M$.
By Theorem \ref{thm;16.4.13.1},
there exist rectilinearizations
$\{(V_{\alpha,\beta},\psi_{\alpha,\beta})
 \,|\,\beta\in\Gamma(\alpha)\}$
of $\phi_{\alpha}^{-1}(B_m)$.
Each $W_{\alpha}$ is equipped with
the coordinate system $(x^{\alpha}_1,\ldots,x^{\alpha}_n)$.
Each $V_{\alpha,\beta}$ is equipped with
the coordinate system
$(y^{\alpha,\beta}_1,\ldots,y^{\alpha,\beta}_n)$.
We set 
$h_{\alpha,\beta}:=
 \prod_{j=1}^n y^{\alpha,\beta}_j\cdot
 \prod_{i=1}^n \psi_{\alpha,\beta}^{\ast}x^{\alpha}_i$.
Let $L_{\alpha,\beta}$
be any subanalytic compact subsets of $V_{\alpha,\beta}$
such that
$\bigcup_{\beta}\psi_{\alpha,\beta}(L_{\alpha,\beta})$
contains a neighbourhood of
$K_{\alpha}$.
Applying \cite[Theorem 4.4]{Bierstone-Milman}
to $h_{\alpha,\beta}$, 
we also obtain a locally finite family of 
real analytic maps
$\lambda_{\alpha,\beta,\gamma}:
 Y_{\alpha,\beta,\gamma}
\lrarr
 V_{\alpha,\beta}$
$(\lambda_{\alpha,\beta,\gamma}\in\Upsilon(\alpha,\beta))$
such that
(i) each $\lambda_{\alpha,\beta,\gamma}$
is the composition of a finite sequence of
local real blowings up,
(ii) $\lambda_{\alpha,\beta,\gamma}^{\ast}h_{\alpha,\beta}$
are normal crossing,
(iii) there exist compact subsets 
$M_{\alpha,\beta,\gamma}\subset Y_{\alpha,\beta,\gamma}$
such that 
$\bigcup_{\gamma}
 \lambda_{\alpha,\beta,\gamma}(M_{\alpha,\beta,\gamma})$
contains a neighbourhood of $L_{\alpha,\beta}$.
Then, 
as in the proof of Lemma \ref{lem;16.7.20.3},
we can observe that
each $\lambda_{\alpha,\beta,\gamma}^{-1}(B_j)$
is the union of quadrants
contained in 
$\lambda_{\alpha,\beta,\gamma}^{-1}(B_j)$.

For each $\alpha\in\Lambda$,
we define the finite subset
\[
 \Gamma_1(\alpha):=\bigl\{
 \beta\in\Gamma(\alpha)\,\big|\,
 \psi_{\alpha,\beta}(V_{\alpha,\beta})
 \cap K_{\alpha}\neq\emptyset
 \bigr\}. 
\]
For each $\beta\in\Gamma_1(\alpha)$,
we define the finite subset
\[
\Upsilon_1(\alpha,\beta):=
 \bigl\{
 \gamma\in\Upsilon(\alpha,\beta)\,\big|\,
 \lambda_{\alpha,\beta,\gamma}(Y_{\alpha,\beta,\gamma})
 \cap L_{\alpha,\beta}\neq\emptyset
 \bigr\}.
\]
Let $\Lambdatilde$
denote the set of
$(\alpha,\beta,\gamma)$,
where
$\alpha\in\Lambda$,
$\beta\in\Gamma_1(\alpha)$
and $\gamma\in\Upsilon(\alpha,\beta)$.
For each 
$(\alpha,\beta,\gamma)\in\Lambdatilde$,
we set
$\phitilde_{\alpha,\beta,\gamma}:=
 \phi_{\alpha}\circ\psi_{\alpha,\beta}
 \circ\lambda_{\alpha,\beta,\gamma}$.
Then, the locally finite family of real analytic maps
$(Y_{\alpha,\beta,\gamma},\phitilde_{\alpha,\beta,\gamma})$
$((\alpha,\beta,\gamma)\in\Lambdatilde)$
has the desired property.
\hfill\qed

\vspace{.1in}
Similarly, we can obtain the following refinement.

\begin{prop}
Let $B_i$ $(i\in S)$ be a locally finite family of subanalytic subsets
in a real analytic manifold $M$.
There exists a locally finite family of
real analytic maps
$\{(W_{\alpha},\phi_{\alpha})\,|\,\alpha\in\Lambda\}$
which is a rectilinearization of each $B_i$ $(i\in S)$.
\hfill\qed
\end{prop}

\subsubsection{Rectilinearization of subanalytic functions}

The following is the rectilinearization theorem 
for subanalytic functions due to Parusi\'{n}ski
\cite[Theorem 2.7]{Parusinski}.
(See also \cite[Theorem 3.4]{Kurdyka-Panuescu}.)
It is also fundamental in this paper.

\begin{thm}
\label{thm;16.7.20.2}
Let $U$ be a subanalytic open subset of $\real^n$,
and let $f:U\lrarr \real$ be a continuous subanalytic function
on $(U,\real^n)$.
Then, there exists a rectilinearization 
$\{(W_{\alpha},\phi_{\alpha})\,|\,\alpha\in\Lambda\}$
of $U$ with the following property.
\begin{itemize}
\item
 Let $\gbigq$ be any quadrant
 contained in $\phi_{\alpha}^{-1}(U)$.
 Then, 
 $\phi_{\alpha}^{\ast}(f)_{|\gbigq}\in\RNC(\gbigq)$.
\hfill\qed
\end{itemize}
\end{thm}

For any subanalytic open subset $U\subset\real^n$,
let $\Quad_n(U,\real^n)$ denote
the set of the $n$-dimensional quadrants 
contained in $U$.

\begin{cor}
Let $U$ be a rectilinearized open set in $\real^n$.
Let $f_1,\ldots,f_{\ell}$ be continuous subanalytic functions on 
$(U,\real^n)$
such that $f_{i|\gbigq}\in\RNC(\gbigq)$
for each $\gbigq\in \Quad_n(U,\real^n)$.
Let $h_1,\ldots,h_{p}$ be real analytic functions on $\real^n$.
Let $g_1,\ldots,g_m:U\lrarr\real$ be 
continuous subanalytic functions on $(U,\real^n)$.
Let $Y_1,\ldots,Y_k$ denote rectilinearized subsets of $\real^n$.
Then, there exists a rectilinearization
$\{(W_{\alpha},\phi_{\alpha})\,|\,\alpha\in\Lambda\}$
for $U$ and $Y_i$ $(i=1,\ldots,k)$
such that the following holds:
\begin{itemize}
\item
 $\phi_{\alpha}^{\ast}(f_i)_{|\gbigq},
\phi_{\alpha}^{\ast}(g_j)_{|\gbigq}
 \in \RNC(\gbigq)$
 for each $\gbigq\in \Quad_n(\phi_{\alpha}^{-1}(U),W_{\alpha})$,
 where we regard 
$W_{\alpha}\simeq\real^n$
 by the coordinate system.
\item
 Each
 $\phi_{\alpha}^{\ast}(h_j)$
 is normal crossing on $W_{\alpha}\simeq\real^n$,
 or constantly $0$.
\end{itemize}
\end{cor}
\pf
Let us study the case $m=1$.
Let $(x_1,\ldots,x_n)$ be the global coordinate system of $\real^n$.
By Theorem \ref{thm;16.7.20.2}
and Lemma \ref{lem;16.7.20.3},
there exists a rectilinearization
$\{(W_{\alpha},\phi_{\alpha})\,|\,\alpha\in\Lambda\}$
for $U$ such that
(i) $\phi_{\alpha}^{\ast}(x_i)$ 
are normal crossing on $W_{\alpha}$,
(ii) $\phi_{\alpha}^{\ast}(h_j)$
are constantly $0$ or normal crossing on $W_{\alpha}$,
(iii)
$\phi_{\alpha}^{\ast}(g_1)\in\RNC(\gbigq)$
 for each $\gbigq\in\Quad_n(\phi_{\alpha}^{-1}(U),W_{\alpha})$.
We also obtain that 
$\phi_{\alpha}^{\ast}(f_i)_{|\gbigq}
 \in \RNC(\gbigq)$
from Lemma \ref{lem;16.7.20.3}.
By Lemma \ref{lem;16.7.20.1},
$\phi_{\alpha}^{-1}(Y_i)$ are also rectilinearized.
Thus, we are done in the case $m=1$.
We can prove the claim for general $m$
by an easy induction.
\hfill\qed

\begin{cor}
\label{cor;18.11.10.1}
Let $U$ be a subanalytic open subset in $\real^n$.
Let $Y_1,\ldots,Y_k$ be subanalytic subsets in $\real^n$.
Let $g_1,\ldots,g_N:U\lrarr\real$ be functions
which are continuous subanalytic on $(U,\real^n)$.
Let $h_1,\ldots,h_{\ell}$ be real analytic functions on $\real^n$.
Then, there exists a rectilinearization
$\{(W_{\alpha},\phi_{\alpha})\,|\,\alpha\in\Lambda\}$
of $U$ and $Y_i$ $(i=1,\ldots,k)$
such that 
$\phi_{\alpha}^{\ast}(g_i)_{|\gbigq}\in \RNC(\gbigq)$
for any $\gbigq\in \Quad_n(\phi_{\alpha}^{-1}(U),W_{\alpha})$,
and that
each $\phi_{\alpha}^{\ast}(h_j)$ is
normal crossing on $W_{\alpha}$
or constantly $0$.
\hfill\qed
\end{cor}

\subsection{Fibrations and ramified analyticity}
\label{subsection;18.11.16.41}

\subsubsection{Singularity of fibrations}

Let $M$ be a real analytic manifold.
Let $Y$ be a relatively compact subanalytic subset in 
$M\times \real^n$
with $\dim Y=k$.
Let $\Ybar^{\sm}_k$ denote the closure of $Y^{\sm}_k$
in $Y$.
Set $Y':=Y\setminus Y^{\sm}_k$.
Note that $\dim Y'\leq k-1$,
that $Y'$ is closed in $Y$,
and that 
$\Sing(\Ybar^{\sm}_k)$
is contained in $Y'$.

Let $\pi:M\times \real^n\lrarr\real^n$ be the projection.
Let $\ell:=\dim \pi(Y)$.
The following lemma is standard.
(For example, see a stronger theorem in
\cite[Part I \S1.7]{Goresky-MacPherson}
\footnote{The author thanks one of the referees
who informed a good reference.}.)
We state it in a way convenient to us,
and include a proof for the convenience of the readers.

\begin{lem}
\label{lem;16.1.22.1}
There exist closed subanalytic subsets
$\gbigg\subset \Ybar^{\sm}_k$
and $\gbigw\subset \pi(Y)$
with the following property.
\begin{itemize}
\item
 $\dim \gbigw<\ell$ and $\dim \gbigg<k$.
\item
  $\Sing\pi(Y)\subset\gbigw$.
\item
The induced map
 $Y^{\sm}_k\setminus
 (\gbigg\cup \pi^{-1}(\gbigw))
\lrarr
 \pi(Y)\setminus \gbigw$
is submersive.
\item
 For any point $P\in \pi(Y)\setminus \gbigw$,
the inequality
$\dim\bigl(
 (\gbigg\cup Y')\cap\pi^{-1}(P)
 \bigr)\leq k-1-\ell$ holds.
\end{itemize}
In particular, 
for any point $P\in\pi(Y)\setminus \gbigw$,
the inequality
$\dim\bigl(
 Y\cap\pi^{-1}(P)
 \bigr)\leq k-\ell$ holds.

If $k=\ell$, we may assume $\gbigg=\emptyset$.
\end{lem}
\pf
We use an induction on $\dim Y$.
The claim is clear in the case $\dim Y=0$.
We assume that we have already known the claim
in the case $\dim Y<k$.

Let $\Ytilde_k$ denote the closure of $Y_k^{\sm}$
in $M\times\real^n$,
which is compact.
There exists a $k$-dimensional real analytic compact manifold 
$\Xtilde_k$
with a real analytic map
$\rhotilde:\Xtilde_k\lrarr M\times\real^n$
such that $\rhotilde(\Xtilde_k)=\Ytilde_k$.
We set $X_k:=\rhotilde^{-1}(\Ybar^{\sm}_k)$,
and $\rho:=\rhotilde_{|X_k}$.
There exists 
a closed subanalytic subset $\gbigg\subset \Ybar^{\sm}_k$
such that
$\Sing\Ybar^{\sm}_k\subset\gbigg$
and that the induced map
$\rho^{-1}(\Ybar^{\sm}_k\setminus \gbigg)\lrarr
 \Ybar^{\sm}_k\setminus \gbigg$
is a local diffeomorphism.

We set $\pitilde:=\pi\circ\rho:X_k\lrarr\real^n$.
By the construction,
$\dim \pitilde(X_k)=\dim\pi(\Ybar_{k}^{\sm})$ holds.
Suppose 
$\dim\pitilde(X_k)=\ell$.
Let $H_1$ denote the closed real analytic subset
of the points $x\in X_k$
such that $\rank(T_x\pitilde)<\ell$.
The set of the critical values of
$X_k\setminus \pitilde^{-1}\Sing\pi(Y)
\lrarr
 \pi(Y)\setminus \Sing\pi(Y)$
is equal to
$\pitilde(H_1)\setminus\Sing\pi(Y)$.
It is subanalytic and its measure is $0$,
and hence
$\dim\bigl(
 \pitilde(H_1)\setminus\Sing\pi(Y)
 \bigr)\leq\ell-1$.
Let $\gbigw_1$ denote the closure of 
$\pitilde(H_1)\cup\Sing\pi(Y)$
in $\pi(Y)$.
Then,
$X_k\setminus \pitilde^{-1}(\gbigw_1)
\lrarr
 \pi(Y)\setminus \gbigw_1$
is submersive.
Because
$X_k\setminus
 \bigl(
 \pitilde^{-1}(\gbigw_1)\cup\rho^{-1}(\gbigg)
 \bigr)
\lrarr
 \Ybar^{\sm}_k\setminus(\pi^{-1}(\gbigw_1)\cup \gbigg)$
is a local diffeomorphism,
the morphism
$\Ybar^{\sm}_k\setminus
 (\pi^{-1}(\gbigw_1)\cup \gbigg)
\lrarr
 \pi(Y)\setminus \gbigw_1$
is submersive.

If $\dim\pi(\Ybar^{\sm}_k)<\ell$,
let $\gbigw_1$ denote the closure of $\pi(\Ybar^{\sm}_k)$
in $\pi(Y)$,
and we set $\gbigg:=\emptyset$.
Then, we clearly obtain that
$\Ybar^{\sm}_k\setminus
 (\pi^{-1}(\gbigw_1)\cup \gbigg)
\lrarr
 \pi(Y)\setminus \gbigw_1$
is submersive
because 
$\Ybar^{\sm}_k\setminus
 (\pi^{-1}(\gbigw_1)\cup \gbigg)
=\emptyset$.

Note that $\dim(\gbigg\cup Y')\leq k-1$.
Suppose that $\dim \pi(\gbigg\cup Y')=\ell$.
By using the hypothesis of the induction,
there exists a closed subanalytic subset
$\gbigw_2\subset \pi(\gbigg\cup Y')$
such that 
(i) $\dim \gbigw_2<\ell$,
(ii) $\dim \bigl(
 (\gbigg\cup Y')\cap\pi^{-1}(P)\bigr)\leq k-\ell-1$
for any $P\in\pi(\gbigg\cup Y')\setminus \gbigw_2$.
If $\dim \pi(\gbigg\cup Y')<\ell$,
we set $\gbigw_2:=\pi(\gbigg\cup Y')$.
Then, we clearly obtain
$\dim \bigl(
 (\gbigg\cup Y')\cap\pi^{-1}(P)\bigr)\leq k-\ell-1$
for any $P\in\pi(\gbigg\cup Y')\setminus \gbigw_2=\emptyset$.

Let $\gbigw$ denote the closure of the union of
$\gbigw_1$ and $\gbigw_2$.
Then, $\gbigg$ and $\gbigw$ have the desired property.

If $k=\ell$, 
we obtain $\pi(\gbigg)\subset \gbigw$
because 
$\dim\bigl(
 (\gbigg\cup Y')\cap\pi^{-1}(P)
 \bigr)\leq -1$ for any $P\in \pi(Y)\setminus \gbigw$.
Hence, we may assume $\gbigg=\emptyset$.
\hfill\qed

\subsubsection{Ramified real analyticity at boundary points}
\label{subsection;16.7.25.1}

Let $U$ be a subanalytic open subset in 
an $n$-dimensional real analytic manifold $M$.
Let $f$ be a real analytic function on $U$
which is subanalytic on $(U,M)$.
Let $\del U$ denote the boundary of $U$,
which is $(n-1)$-dimensional.

\begin{lem}
\label{lem;15.12.28.10}
There exists a closed $(n-2)$-dimensional
subset $Z\subset \del U$
with the following property.
\begin{itemize}
\item
 $Z\supset\Sing(\del U)$.
\item
 $f$ is ramified real analytic around any point
 $P\in \del U\setminus Z$.
\end{itemize}
\end{lem}
\pf
There exists a rectilinearization 
$(W_{\alpha},\phi_{\alpha})$ $(\alpha\in\Lambda)$
of $(U,f)$ and $\del U$.
There exist compact subsets
$K_{\alpha}\subset W_{\alpha}$ $(\alpha\in\Lambda)$
such that
$\bigcup\phi_{\alpha}(K_{\alpha})=\real^n$.
Let $T^2_{\alpha}$ denote the union of
the quadrants $\gbigq$ in $W_{\alpha}$ such that
(i) $\gbigq\subset\phi_{\alpha}^{-1}(\del U)$,
(ii) $\dim \gbigq\leq n-2$.
Each $\phi_{\alpha}$ is factorized
as in (\ref{eq;16.4.13.2}) of \S\ref{subsection;16.7.22.1}.
Let $C_{\alpha}^{(\ell)}$ denote the centers of local blowings up.
Let $\psi_{\alpha}^{(\ell)}:W_{\alpha}^{(\ell)}\lrarr \real^n$
denote the induced map.
There exist the closed subanalytic subsets
$Z_{\alpha}^{(\ell)}:=
 \phi_{\alpha}(K_{\alpha})
\cap
 \psi_{\alpha}^{(\ell)}(
 C_{\alpha}^{(\ell)})$.
We set
\[
 Z:=
\Bigl(
 \bigcup_{\alpha}
 \phi_{\alpha}(K_{\alpha}\cap T^2_{\alpha})
\Bigr)
\cup
\Bigl(
 \del U\cap 
\bigcup_{\alpha,\ell} Z_{\alpha}^{(\ell)}
 \Bigr)\cup
\Sing(\del U).
\]

Let $P$ be any point of $\del U\setminus Z$.
There exists $\alpha\in\Lambda$
such that $P\in\phi_{\alpha}(K_{\alpha})$.
By our choice of $Z$,
there exists a neighbourhood $M_P$ of $P$ in $\real^n$
such that 
$\phi_{\alpha}$ induces a diffeomorphism
$\phi_{\alpha}^{-1}(M_P)\simeq M_P$.
Let $(s_1,\ldots,s_n)$ be the coordinate system
of $W_{\alpha}$.
Because $\del U$ is smooth around $P$,
we may assume that
$\phi_{\alpha}^{-1}(M_P\cap U)
=\{s_1>0\}\cap \phi_{\alpha}^{-1}(M_P)$.
Note that 
$\phi_{\alpha}^{\ast}(f)$ is real analytic
on $\{s_1>0\}\cap \phi_{\alpha}^{-1}(M_P)$,
and that
$\phi_{\alpha}^{\ast}(f)
 \in \RNC(\gbigq)$ for any 
$\gbigq\in\Quad_n(\phi_{\alpha}^{-1}(U),W_{\alpha})$.
Hence,
$\phi_{\alpha}^{\ast}(f)$
is expressed as
$\sum a_j(s_2,\ldots,s_n)s_1^{j/\rho}$
on $\phi_{\alpha}^{-1}(M_P\cap U)$.
\hfill\qed

\subsubsection{Ramified analyticity at corners}
\label{subsection;18.11.16.50}

The results in 
\S\ref{subsection;18.11.16.50}--\ref{subsection;18.11.16.51}
are preliminary for the study in \S\ref{section;18.11.15.20}.

Let us consider the case
$M=M_0\times\real$,
and $U$ is contained in
$M_0\times\{\tau>0\}$,
where $\tau$ denotes the standard coordinate on $\real$.
Here, $M_0$ is an $(n-1)$-dimensional
real analytic manifold.
We set $V:=\Ubar\cap (M_0\times\{0\})$.
Let $V^{\circ}$ be the set of interior points
of $V$ as a subset of $M_0\times\{0\}$.
Let $\overline{V^{\circ}}$
denote the closure of $V^{\circ}$
in $M_0\times\{0\}$.
Let $\del(V^{\circ})$ denote the boundary 
of $V^{\circ}$ as a subset of $M_0\times\{0\}$.
Let $f$ be a continuous subanalytic function on $(U,M)$.
We obtain the following lemma
as a direct consequence of Lemma \ref{lem;15.12.28.10}.

\begin{lem}
\label{lem;15.12.29.1}
There exists a closed subanalytic subset
$Z_0\subset \overline{V^{\circ}}$
such that the following holds:
\begin{itemize}
\item
 $\dim Z_0\leq n-2$
 and 
 $Z_0\supset \del V^{\circ}$.
\item
 For any $P\in V^{\circ}\setminus Z_0$,
 there exists a neighbourhood
 $M_P$ in $M$
 such that
 $U\cap M_P=M_P\cap\{\tau>0\}=:U_P$
 and that 
 $f_{|U_P}$ is expressed as
 $\sum \alpha_j(x_1,\ldots,x_{n-1})\tau^{j/\rho}$,
 where $\rho$ is a positive integer,
 $(x_1,\ldots,x_{n-1})$ is a real analytic coordinate system
 on a neighbourhood of $P$ in $M_P\cap V^{\circ}$,
 and $\alpha_j$ are real analytic functions.

\hfill\qed
\end{itemize}
\end{lem}

Let us study the behaviour of $f$
around any sufficiently general point of $Z_0$.

\begin{prop}
\label{prop;15.12.28.1}
There exists a closed subanalytic subset 
$Z_1\subset Z_0$
with the following property:
\begin{itemize}
\item
 $\dim Z_1\leq n-3$
and
 $Z_1\supset \Sing(\del V^{\circ})\cup\Sing(Z_0)$.
\item
 Let $P$ be any point of $\del V^{\circ}\setminus Z_1$.
 Let $(\nbigu_P,y_1,\ldots,y_{n-1})$ 
 be any real analytic coordinate neighbourhood of $P$
 in $M_0\times\{0\}$
 such that $\nbigu_P\cap V^{\circ}=\{y_{n-1}>0\}$.
 Then, there exist a positive integer $m>0$
 and a positive number $C>0$
 such that 
\[
 \nbigb:=
 \bigl\{
 (y_1,\ldots,y_{n-1},\tau)\,\big|\,
 0<\tau<Cy_{n-1}^m
 \bigr\}
\subset U.
\]
Moreover,
$f_{|\nbigb}$ is real analytic with respect to
$(y_1,\ldots,y_{n-2},y_{n-1}^{1/\rho},(y_{n-1}^{-m}\tau)^{1/\rho})$
for a positive integer $\rho$,
i.e.,
$f_{|\nbigb}$ is expressed as a convergent power series
\[
 f_{|\nbigb}=
 \sum_{i\geq -N_1}
 \sum_{j\geq -N_2}
 A_{ij}(y_1,\ldots,y_{n-2})
 \cdot
 y_{n-1}^{i/\rho}
 \cdot
 (y_{n-1}^{-m}\tau)^{j/\rho},
\]
where $A_{ij}$ are real analytic functions of
$(y_1,\ldots,y_{n-2})$.
\item
 Let $P$ be any point of 
 $(V^{\circ}\cap Z_0)\setminus Z_1$.
 Let $(\nbigu_P,y_1,\ldots,y_{n-1})$ be
 any real analytic coordinate neighbourhood of $P$
 in $M_0\times\{0\}$
 such that 
 $\nbigu_P\cap (V^{\circ}\setminus Z_0)
 =\{y_{n-1}\neq 0\}$.
Then,
there exist
a positive integer $m>0$
 and a positive number $C>0$
 such that 
\[
 \nbigb_+:=
 \bigl\{
 (y_1,\ldots,y_{n-1},\tau)\,\big|\,
 y_{n-1}>0,\,\,
 0<\tau<Cy_{n-1}^m
 \bigr\} 
\subset U,
\]
\[
 \nbigb_-:=
 \bigl\{
 (y_1,\ldots,y_{n-1},\tau)\,\big|\,
 y_{n-1}<0,\,\,
 0<\tau<C(-y_{n-1})^m
 \bigr\} 
\subset U.
\]
Moreover,
$f$ is real analytic with respect to
$\bigl(y_1,\ldots,y_{n-2},
 (\pm y_{n-1})^{1/\rho},
 ((\pm y_{n-1})^{-m}\tau)^{1/\rho}
 \bigr)$
on $\nbigb_{\pm}$
for a positive integer $\rho$,
i.e.,
$f_{|\nbigb_{\pm}}$ is expressed as a convergent power series
\[
 f_{|\nbigb_{\pm}}=
 \sum_{i\geq -N_1}
 \sum_{j\geq -N_2}
 A_{\pm,ij}(y_1,\ldots,y_{n-2})
 \cdot
 (\pm y_{n-1})^{i/\rho}
 \cdot
 ((\pm y_{n-1})^{-m}\tau)^{j/\rho},
\]
\end{itemize}
where $A_{\pm,ij}$ are real analytic functions of
$(y_1,\ldots,y_{n-2})$.
If $U$ is relatively compact in $M$,
the positive numbers $\rho$ and $m$
in the second and third properties
are bounded.
\end{prop}
\pf
It is enough to study 
the connected components $\nbigc$
of 
$U\cap
 \Bigl(
 ((M_0\times\{0\})\setminus Z_0)\times\real_{\geq 0}
 \Bigr)$
and the restriction of $f$ to $\nbigc$.
Note that the third property
of $(U,f)$
follows from the second properties
of $(\nbigc,f_{|\nbigc})$
for any connected components $\nbigc$.
Hence,
we may assume
$Z_0=\del V^{\circ}$
from the beginning.

By Corollary \ref{cor;18.11.10.1},
there exists a rectilinearization
$(W_{\alpha},\phi_{\alpha})$ $(\alpha\in\Lambda)$
of $(U,f)$,
$V^{\circ}$ and $\del V^{\circ}$
such that 
$\phi_{\alpha}^{\ast}(\tau)$ are normal crossing
for any $\alpha\in\Lambda$.
There exist compact subsets $K_{\alpha}\subset W_{\alpha}$
such that $M=\bigcup \phi_{\alpha}(K_{\alpha})$.
For each $\alpha$,
let $T^3_{\alpha}$ denote the union of
the quadrants $\gbigq$ in $W_{\alpha}$
such that
(i) $\gbigq\cap K_{\alpha}\subset\phi_{\alpha}^{-1}(\del V^{\circ})$,
(ii) $\dim \phi_{\alpha}(\gbigq)\leq n-3$.

For each $\alpha$,
there exists a factorization of $\phi_{\alpha}$
into a finite sequence of local real blowings up
as in (\ref{eq;16.4.13.2}) in \S\ref{subsection;16.7.22.1}.
Let $C_{\alpha}^{(\ell)}$ denote the centers of 
the local real blowings up.
Note that 
$\dim C_{\alpha}^{(\ell)}\leq n-2$.
Let $\psi_{\alpha}^{(\ell)}:W_{\alpha}^{(\ell)}\lrarr \real^n$
be the induced map.
Let us look at the strict transforms
$(\overline{V^{\circ}})^{(\ell)}_{\alpha}
\subset 
 W_{\alpha}^{(\ell)}$
of $\overline{V^{\circ}}$
with respect to $\psi_{\alpha}^{(\ell)}$.
Let 
$\nbiga^{(\ell)}_{\alpha}$
denote the set of
$(n-2)$-dimensional smooth points of
$C_{\alpha}^{(\ell)}\cap 
 \del(\overline{V^{\circ}})^{(\ell)}_{\alpha}$.
Set
$\nbigb^{(\ell)}_{\alpha}:=
C_{\alpha}^{(\ell)}\cap 
 \del(\overline{V^{\circ}})^{(\ell)}_{\alpha}
\setminus
 \nbiga^{(\ell)}_{\alpha}$.
By the construction,
$\dim \nbigb^{(\ell)}_{\alpha}\leq n-3$ holds.
Let $Z_{\alpha}^{(\ell)}$
denote the closure of
$\psi_{\alpha}^{(\ell)}(\nbigb^{(\ell)}_{\alpha})
 \cap
 \phi_{\alpha}(K_{\alpha})$.
Let $Z_1$ denote the union of
$\Sing(\del\overline{V^{\circ}})
=\Sing(Z_0)$,
$\phi_{\alpha}(K_{\alpha}\cap T^3_{\alpha})$
for all $\alpha$,
and 
$Z_{\alpha}^{(\ell)}$ for all of $\alpha$ and $\ell$.
By construction,
$Z_1$ is a subanalytic subset with
$\dim Z_1\leq n-3$.

Let us study the second property.
Let $P$ be any point of 
$\del V^{\circ}\setminus Z_1$.
Let $(\nbigu_P,y_1,\ldots,y_{n-1})$ be
any real analytic coordinate neighbourhood
of $P$ in $M_0$
such that
$V^{\circ}\cap\nbigu_P=\{y_{n-1}>0\}$.
By shrinking $\nbigu_P$,
we may assume that
$\nbigu_P=\{(y_1,\ldots,y_{n-1})\,|\,|y_i|<\epsilon\}$
for some $\epsilon>0$.
Let $M_P$ be a neighbourhood of
$P$ in $M=M_0\times\real$
defined as
$M_P=\nbigu_P\times\{|\tau|<\epsilon\}$,
which is the product of 
$M_{1,P}:=\{(y_1,\ldots,y_{n-2})\,|\,|y_i|<\epsilon\}$
and 
$M_{2,P}:=\{(y_{n-1},\tau)\,|\,|y_{n-1}|<\epsilon,|\tau|<\epsilon\}$.

As a preparation for the further argument,
we construct a sequence of blowings up inductively.
Set $Y^{(0)}:=\real^2=\{(y_{n-1},\tau)\}$,
$H^{(0)}:=\{\tau=0\}$,
and $Q^{(0)}:=(0,0)$.
Let $\kappa_1:Y^{(1)}\lrarr Y^{(0)}$
be the real blowing up at $Q^{(0)}$.
Let $H^{(1)}$ denote the strict transform of $H^{(0)}$,
and set $Q^{(1)}:=\kappa_1^{-1}(Q^{(0)})\cap H^{(1)}$.
Suppose that we have already constructed a sequence
\[
 Y^{(i)}\stackrel{\kappa_i}{\lrarr}
 Y^{(i-1)}\stackrel{\kappa_{i-1}}{\lrarr}
 \cdots
 \stackrel{\kappa_1}{\lrarr}
 Y^{(0)}
\]
with points $Q^{(j)}\in Y^{(j)}$
such that $\kappa_{j+1}$
are the real blowing up at $Q^{(j)}$,
and $H^{(j)}\subset Y^{(j)}$
which is the strict transform of $H^{(0)}$.
Then, we define $\kappa_{i+1}:Y^{(i+1)}\lrarr Y^{(i)}$
as the blowing up at $Q^{(i)}$.
We also define $H^{(i+1)}$ as the strict transform of
$H^{(i)}$ with respect to $\kappa_{i+1}$,
and we put $Q^{(i+1)}:=\kappa_{i+1}^{-1}(Q^{(i)})\cap H^{(i+1)}$.
Thus, the inductive construction can proceed.
For each $i$,
$(y_{n-1},\tau y_{n-1}^{-i})$
is a natural local coordinate system around $Q^{(i)}$.
Let $\nu_i$ denote the naturally induced morphism
$M_{1,P}\times Y^{(i)}\lrarr M_{1,P}\times Y^{(0)}$.

Let $\gamma:(\II^{\circ}_{\delta},0,\II_{\delta})\lrarr 
 (V^{\circ},P,\Vbar^{\circ})$
be any real analytic path 
such that 
$\gamma(\II_{\delta}^{\circ})\setminus\{0\})$
does not intersect with 
the set of the critical values of 
$\phi_{\alpha}$ for any $\alpha\in\Lambda$.
After making $\delta$ smaller,
there exists $\alpha_0\in\Lambda$
such that
$\gamma(\II_{\delta})\subset
 \phi_{\alpha_0}(K_{\alpha_0})$.
We obtain the path
$\gammatilde_{\alpha_0}$
to $W_{\alpha_0}$
such that 
$\phi_{\alpha_0}\circ\gammatilde_{\alpha_0}
=\gamma$.

\begin{lem}
\label{lem;18.11.23.1}
There exist
a positive integer $k$,
an open neighbourhood 
$\nbigv_k$
of the image of $\gammatilde_{\alpha_0}$
in $W_{\alpha_0}$,
an open subset
$\nbigu_k\subset
 M_{1,P}\times Y^{(k)}$,
and a real analytic isomorphism
$\mu:
 \nbigu_k\lrarr
 \nbigv_k$
such that
$\nu_k=\phi_{\alpha_0}\circ \mu$
on $\nbigu_k$.
\end{lem}
\pf
Let $M_P'$ be any small neighbourhood 
of  $P$ in $M_P$
such that $M_P'\cap Z_1=\emptyset$.
Because of the existence of 
the lift $\gammatilde_{\alpha_0}$,
for each $\ell$
the strict transform of $V^{\circ}$
with respect to $\psi^{(\ell)}$
is non-empty.
By the construction of $Z_1$,
for each $\ell$,
if $C_{\alpha_0}^{(\ell-1)}\cap
 (\psi_{\alpha_0}^{(\ell-1)})^{-1}(M_P')$
is non-empty,
it is equal to
$\del(\overline{V^{\circ}})^{(\ell-1)}_{\alpha_0}
\cap
 (\psi_{\alpha_0}^{(\ell-1)})^{-1}(M_P')$.
We set
$k:=\#\bigl\{
 \ell\,\big|\,
 C_{\alpha_0}^{(\ell-1)}\cap
 (\psi_{\alpha_0}^{(\ell-1)})^{-1}(M_P')\neq\emptyset
 \bigr\}$.
Then, the claim of the lemma is clear.
\hfill\qed

\vspace{.1in}

After making $\epsilon$ smaller,
we may assume that 
$M_{1,P}\times\{(y_{n-1},\tau)\in M_{2,P}\,|\,y_{n-1}>0,\tau=0\}$
is contained in
$\phi_{\alpha_0}(W_{\alpha_0})$.
Let $(s_1,\ldots,s_n)$ be the coordinate system
of $W_{\alpha_0}$.
Because $\phi_{\alpha_0}^{-1}(V^{\circ})$
is rectilinearized,
and because of the description of
$\psi_{\alpha_0}$ around $P$
in Lemma \ref{lem;18.11.23.1},
we may assume that
$\phi_{\alpha_0}^{-1}(V^{\circ}\cap M_P)
=\{s_n=0,s_{n-1}>0\}\cap
 \phi_{\alpha_0}^{-1}(M_P)$.
We define the real analytic functions
$b_1:=\phi_{\alpha_0}^{\ast}(\tau y_{n-1}^{-k})$
and 
$b_2:=\phi_{\alpha_0}^{\ast}(y_{n-1})$
on $\phi_{\alpha_0}^{-1}(M_P)$.
We regard $s_n$ as a real analytic function
on $\phi_{\alpha_0}^{-1}(M_P)$.
Both of the differentials $db_1$ and $ds_n$ 
are nowhere vanishing on 
$\phi_{\alpha_0}^{-1}(M_P)$.
Note that 
$\phi_{\alpha_0}^{-1}(V^{\circ}\cap M_P)
\subset b_1^{-1}(0)\cap s_n^{-1}(0)$.
Then, we obtain that 
$A_0:=s_n/b_1$ is a real analytic function on
$\phi_{\alpha_0}^{-1}(M_P)$,
which is nowhere vanishing.
By replacing $s_n$ with $-s_n$ if necessary,
we may assume that $A_0>0$.
Note that
$\phi_{\alpha_0}^{-1}(M_P\cap V^{\circ})
=\{s_n=0\}\cap \{b_2>0\}
 \cap
 \phi_{\alpha_0}^{-1}(M_P)$.
The derivatives $db_2$
and $ds_{n-1}$ are nowhere vanishing
on $\{s_n=0\}\cap\phi_{\alpha_0}^{-1}(M_P)$.
Hence,
on $\{s_n=0\}\cap\phi_{\alpha_0}^{-1}(M_P)$,
$s_{n-1}/b_2$ 
is real analytic and positive.
Hence, after shrinking $M_P$,
$s_{n-1}$ is expressed as
$A_1b_2+B b_1$
on $\phi_{\alpha_0}^{-1}(M_P)$,
where $A_1$ is a nowhere vanishing
real analytic function on $\phi_{\alpha_0}^{-1}(M_P)$,
and $B$ is a real analytic function on
$\phi_{\alpha_0}^{-1}(M_P)$.

Because
$\phi_{\alpha_0}^{-1}(V^{\circ}\cap M_P)$
contains
$\phi_{\alpha_0}^{-1}(M_P)
\cap \{s_n=0,s_{n-1}>0\}$,
and because
$\phi_{\alpha_0}^{-1}(V^{\circ})$
and $\phi_{\alpha_0}^{-1}(U)$
are rectilinearized,
we obtain that
$\phi_{\alpha_0}^{-1}(U\cap M_P)$
contains
$\phi_{\alpha_0}^{-1}(M_P)
\cap \{s_n>0,s_{n-1}>0\}$.
Note that
$T^3_{\alpha_0}\cap \phi_{\alpha_0}^{-1}(M_P)=\emptyset$.
Hence, 
$s_i$ $(i<n-1)$ are nowhere vanishing
on $\phi_{\alpha_0}^{-1}(M_P)$.

Let $m$ be any integer strictly larger than $k+1$.
If $C>0$  is sufficiently small,
$\phi_{\alpha_0}^{-1}(\nbigb\cap M_P)$ is contained in
$\{s_n>0,s_{n-1}>0\}\subset W_{\alpha_0}$.
Moreover,
on $\phi_{\alpha_0}^{-1}(\nbigb\cap M_P)$,
we obtain that
$s_{n}^{j/\rho}=
 A_0^{j/\rho}\phi_{\alpha_0}^{\ast}(\tau y_{n-1}^{-k})^{j/\rho}
=A_0^{j/\rho}\phi_{\alpha_0}^{\ast}(\tau y_{n-1}^{-m})^{j/\rho}
 \phi_{\alpha_0}^{\ast}(y_{n-1})^{(m-k)/\rho}$,
and that
\[
 s_{n-1}^{j/\rho}=
 A_1^{j/\rho}\phi_{\alpha_0}^{\ast}(y_{n-1})^{j/\rho}
 \bigl(1+BA_1^{-1}\phi_{\alpha_0}^{\ast}(\tau y_{n-1}^{-m})
 \phi_{\alpha_0}^{\ast}(y_{n-1})^{m-k-1}
  \bigr)^{j/\rho}.
\]
Because
$f$ is expressed as
$\sum c_{j_1,j_2}(s_1,\ldots,s_{n-2})
 s_{n}^{j_1/\rho}s_{n-1}^{j_2/\rho}$
on $\phi_{\alpha_0}^{-1}(M_P)
\cap \{s_n>0,s_{n-1}>0\}$
for some real analytic functions
$c_{j_1,j_2}$ of $(s_1,\ldots,s_{n-2})$,
the second property is clear.
\hfill\qed

\subsubsection{Fibrations and ramified real analyticity}
\label{subsection;18.11.16.51}

We reformulate the results in \S\ref{subsection;18.11.16.50}
in a way convenient for the study in
\S\ref{section;18.11.15.20}.

Let $M_1$ be an $(n-2)$-dimensional real analytic manifold.
Set $A:=\real_{\geq 0}\times M_1$
and $B:=S^1\times A$.
Note that
$\del B=S^1\times\del A$.
Let $U$ be an open subset in 
$B\setminus\del B$ which is subanalytic in $B$.
Let $f$ be a real analytic and subanalytic on $(U,B)$.
Note that $\dim U=n$.
We set
$V:=\overline{U}\cap \del B$.
Let $V^{\circ}$ denote the set of the interior points of
$V$ in $\del B$.
If $V^{\circ}\neq\emptyset$,
then $\dim V^{\circ}=n-1$.
There exist closed subanalytic subsets
$Z_1\subset Z_0\subset \overline{V^{\circ}}$
as in Lemma \ref{lem;15.12.29.1} and 
Proposition \ref{prop;15.12.28.1}.

Let $q:B\lrarr A$ be the projection.
We obtain $q(\overline{V^{\circ}})\subset \del A$.
Clearly, $\dim q(\overline{V^{\circ}})\leq \dim \del A=n-2$ holds.
Note that
$q(Z_{1})\subset q(Z_0)\subset q(\overline{V^{\circ}})$,
and the inequalities
$\dim q(Z_1)\leq n-3$
and $\dim q(Z_0)\leq n-2$.
We obtain the following lemma as a consequence of
Lemma \ref{lem;16.1.22.1}
\begin{lem}
There exists a closed subanalytic subset
$\gbigw_0\subset q(\overline{V^{\circ}})$ with the following property.
\begin{itemize}
\item
 $\dim \gbigw_0\leq n-3$
and  $\del q(\overline{V^{\circ}})\subset \gbigw_0$.
\item
 $Z_0\setminus q^{-1}(\gbigw_0)$
 is horizontal with respect to $q$,
 i.e.,
 $\dim \bigl(q^{-1}(P)\cap Z_0\bigr)=0$
 for any $P\in q(\overline{V^{\circ}})\setminus \gbigw_0$.
\hfill\qed
\end{itemize}
\end{lem}

Because
$\dim q^{-1}(\gbigw_0)\leq n-2$,
we may assume that
$Z_0$ contains
$Y_0:=q^{-1}(\gbigw_0)\cap \overline{V^{\circ}}$,
by enlarging $Z_0$ and $Z_1$.
Let $Y_1$ denote the closure of
$Z_0\setminus Y_0$ in $Z_0$.
Note that $\dim(Y_0\cap Y_1)\leq n-3$.
We may assume that $Z_1$ contains
$Y_0\cap Y_1$
by enlarging $Z_1$.

\begin{lem}
\label{lem;15.12.30.1}
There exist a closed subanalytic subset
$\gbigw_1\subset \gbigw_0$
and a closed subanalytic subset 
$\gbigg\subset \overline{(Y_0)^{\sm}_{n-2}}$
with the following property.
\begin{itemize}
\item
 $\dim \gbigw_1\leq n-4$
 and $\dim \gbigg\leq n-3$.
\item
 $\gbigw_1\supset\Sing(\gbigw_0)$
 and $\gbigg\supset\Sing\overline{(Y_0)_{n-2}^{\sm}}$.
\item
 $\bigl(
 Z_1\cup \gbigg
 \bigr)\setminus q^{-1}(\gbigw_1)$ is horizontal
 with respect to $q$.
\item
 $(Y_0)_{n-2}^{\sm}\setminus
 (\gbigg\cup q^{-1}(\gbigw_1))\lrarr \gbigw_0\setminus \gbigw_1$
 is submersive.
\end{itemize}
\end{lem}
\pf
It follows from Lemma \ref{lem;16.1.22.1}.
\hfill\qed

\vspace{.1in}

Let $r$ be the coordinate of $\real_{\geq 0}$,
and $\theta$ be the local coordinate of $S^1=\real/2\pi\seisuu$.
Let $(x_1,\ldots,x_{n-2})$ be the coordinate system of $\real^{n-2}$.
We set
$\gbigr_0:=Z_0$
and 
$\gbigr_1:=Z_1\cup \gbigg$.
We summarize the property of
$\gbigr_i$ and $\gbigw_i$ $(i=0,1)$.
\begin{prop}
\label{prop;16.1.27.10}
The following holds.
\begin{itemize}
\item
$\dim \gbigr_0\leq n-2$,
$\dim \gbigr_1\leq n-3$,
$\dim \gbigw_0\leq n-3$
and $\dim \gbigw_1\leq n-4$.
\item
$\gbigr_1$ is contained in $q^{-1}(\gbigw_0)$.

\item
 $\gbigr_0\setminus q^{-1}(\gbigw_0)$
 and
 $\gbigr_1\setminus q^{-1}(\gbigw_1)$
 are horizontal with respect to $q$.
\item
 Let $P$ be any point of $q(\overline{V^{\circ}})\setminus \gbigw_0$
 such that
 $\dim (q^{-1}(P)\cap \overline{V^{\circ}})=1$.
 Let $Q$ be any interior point of
 $q^{-1}(P)\cap \overline{V^{\circ}}\setminus \gbigr_0
 \subset q^{-1}(P)$.
 Then, there exists a neighbourhood 
 $\nbigu$ of $Q$
 in $\Ubar$
 and a positive integer $\rho$
 such that the restriction of $f$ to $\nbigu$
 is expressed as
\[
 f=\sum_{j\geq -N_1}
 \alpha_j(\theta,x_1,\ldots,x_{n-2})\cdot r^{j/\rho}.
\]
Here, $\alpha_j$ are real analytic functions.
\item
 Let $P$ be any point of 
$\gbigw_0\setminus (\gbigw_1\cup\del q(\overline{V^{\circ}}))$
 such that $\dim (q^{-1}(P)\cap \overline{V^{\circ}})=1$.
 Let $Q$ be any interior point of
 $q^{-1}(P)\cap \overline{V^{\circ}}\setminus \gbigr_1
 \subset q^{-1}(P)$.
There exists
a real analytic coordinate neighbourhood
$(\nbign;y_1,\ldots,y_{n-2})$
 around $P$ in $\real^{n-2}$,
 such that $\gbigw_0\cap\nbigu=\{y_{n-2}=0\}$.
There exist real numbers
 $\theta_1<\theta_2$ such that
the interval
 $\closedclosed{\theta_1}{\theta_2}$
 is a small neighbourhood of $Q$
 in $q^{-1}(P)\cap \overline{V^{\circ}}$
and that
$\closedclosed{\theta_1}{\theta_2}
\cap 
 (q^{-1}(P)\cap
 \gbigr_1)=\emptyset$.
Moreover,
 there exist 
 a positive integer $\rho$,
 a positive integer $m$,
and a positive number $C>0$,
such that
\[
 \nbigu_{\pm}=
 \bigl\{(y_1,\ldots,y_{n-2},\theta,r)\,\big|\,
 (y_1,\ldots,y_{n-2})\in\nbign,\,\,\,
 \theta_1<\theta<\theta_2,\,\,\,
 0<\pm y_{n-2},\,\,\,
 0<r<C(\pm y_{n-2})^{m}\,\,
 \bigr\}
\subset
 U
\]
and that the restriction of $f$ to $\nbigu_{\pm}$
are expressed as
\[
 f=\sum_{i\geq -N_1}\sum_{j\geq -N_2}
 \alpha_{\pm,i,j}(y_1,\ldots,y_{n-3},\theta)\cdot
 y_{n-2}^{i/\rho}\cdot
 \bigl((\pm y_{n-2}\bigr)^{-m}r)^{j/\rho}.
\]
Here $\alpha_{\pm,i,j}$ are real analytic functions.
\item
 Let $P$ be any point of $\del q(\overline{V^{\circ}})\setminus \gbigw_1$
 such that $\dim (q^{-1}(P)\cap \overline{V^{\circ}})=1$.
 Let $Q$ be any interior point of
 $q^{-1}(P)\cap \overline{V^{\circ}}\setminus\gbigr_1
 \subset q^{-1}(P)$.
There exists a real analytic coordinate neighbourhood
 $(\nbign;y_1,\ldots,y_{n-2})$ around $P$
in $\real^{n-2}$
such that $q(\overline{V^{\circ}})\cap\nbign=\{y_{n-2}\geq 0\}$.
There exist
 real numbers $\theta_1<\theta_2$
such that
$\closedclosed{\theta_1}{\theta_2}$
is a neighbourhood of $Q$
in $q^{-1}(P)\cap \overline{V^{\circ}}$
such that
$\closedclosed{\theta_1}{\theta_2}
\cap
 (q^{-1}(P)\cap\gbigr_1)
=\emptyset$.
Moreover, there exist
 a positive integer $\rho$,
 a positive integer $m$,
and a positive number $C>0$,
such that
\[
 \nbigu=
 \bigl\{(y_1,\ldots,y_{n-2},\theta,r)\,\big|\,
 (y_1,\ldots,y_{n-2})\in\nbign,\,\,\,
 0<y_{n-2},\,\,\,
 \theta_1<\theta<\theta_2,\,\,\,
 0<r<C y_{n-2}^{m}\,\,
 \bigr\}
\subset
 U
\]
and that the restriction of $f$ to $\nbigu$
is expressed as
\[
 f=\sum_{i\geq -N_1}\sum_{j\geq -N_2}
 \alpha_{i,j}(y_1,\ldots,y_{n-3},\theta)\cdot
 y_{n-2}^{i/\rho}\cdot
 (y_{n-2}^{-m}r)^{j/\rho}.
\]
Here, $\alpha_{i,j}$ are real analytic functions.
\end{itemize}
If $U$ is relatively compact,
the numbers $\rho$ and $m$ are bounded.
\hfill\qed
\end{prop}

\subsection{Complements}
\label{subsection;18.11.16.42}

\subsubsection{Dominant subanalytic functions}

This subsection is preliminary for
the study in \S\ref{subsection;18.11.16.100}.
Let $H$ be a closed subanalytic subset of 
a real analytic manifold $M$.
\begin{lem}
\label{lem;16.7.23.1}
There exists a continuous subanalytic function $\chi_H$
on $M$ such that
$\chi_H^{-1}(0)=H$.
\end{lem}
\pf
Let $(\nbigu_{\lambda},x^{\lambda}_1,\ldots,x^{\lambda}_n)$
$(\lambda\in\Lambda)$
be a locally finite covering
of $M$ by coordinate neighbourhoods.
Let $d_{\nbigu_{\lambda}}$
denote the distance on $\nbigu_{\lambda}$
induced by the coordinate system
$(x^{\lambda}_1,\ldots,x^{\lambda}_n)$
and the standard Euclidean distance on $\real^n$.

Let $\nbigv_{\mu}$ $(\mu\in \Gamma)$ 
be a refined locally finite covering of $M$
such that 
each $\nbigv_{\mu}$ is a relatively compact subanalytic
subset in $\nbigu_{\lambda(\mu)}$
for some $\lambda(\mu)\in\Lambda$.
Let $C_{\mu}$ be the subanalytic subset
of $\nbigu_{\lambda(\mu)}$
obtained as the union of
$H\cap\nbigu_{\lambda(\mu)}$
and 
$\nbigu_{\lambda(\mu)}\setminus
 \nbigv_{\mu}$.
We define the map
 $\chi_{\mu}:\nbigu_{\lambda(\mu)}\lrarr \real_{\geq 0}$
by 
\[
 \chi_{\mu}(P):=d_{\nbigu_{\lambda(\mu)}}(P,C_{\mu}).
\]
It is a continuous subanalytic function on $\nbigu_{\lambda(\mu)}$
(see \cite[Remark 3.11]{Bierstone-Milman}).
The support is contained in $\nbigv_{\mu}$.
We extend it to a continuous subanalytic function on $M$
by setting $\chi_{\mu}(P)=0$ for any $P\not\in \nbigu_{\lambda(\mu)}$.
We set 
$\chi_H(P):=\max_{\mu\in \Gamma}\chi_{\mu}(P)$.
Then, it has the desired property.
\hfill\qed

\vspace{.1in}
We fix any continuous subanalytic function $\chi_H$
on $M$ such that $\chi_H^{-1}(0)=H$.
Note that 
$\frac{1}{\chi_H}$ is a subanalytic function
on $(M\setminus H,M)$.

\begin{lem}
\label{lem;16.7.23.2}
Let $U$ be a relatively compact subanalytic open subset of $M$
such that $U\cap H=\emptyset$.
Let $f$ be a continuous subanalytic function on $(U,M)$
with the following property.
\begin{itemize}
\item
 For any relatively compact subset $V\subset M\setminus H$,
 $|f_{|V\cap U}|$ is bounded.
\end{itemize}
Then, there exist positive constants $N$ and $C$
such that 
$|f|\leq C(\frac{1}{\chi_H})^{N}_{|U}$
on $U$.
\end{lem}
\pf
There exists a rectilinearization
$\{(W_{\alpha},\phi_{\alpha})\,|\,\alpha\in\Lambda\}$
for $f$ and $\chi_H$.
It is easy to see that
for each $\alpha\in\Lambda$
there exist positive constants
$C_{\alpha}$ and $N_{\alpha}$
such that
$|\phi_{\alpha}^{\ast}(f)|
\leq
 C_{\alpha}\phi_{\alpha}^{\ast}(\frac{1}{\chi_{H}})^{N_{\alpha}}$.
By using the relative compactness of $U$,
we obtain the claim of the lemma.
\hfill\qed

\subsubsection{Lift of maps at boundary}

This subsection is a preliminary for the study
in \S\ref{subsection;18.11.16.101}.
Let $F:N\lrarr M$ be a real analytic map
of real analytic manifolds.
Let $U$ be a subanalytic relatively compact open subset of $N$.
Let $W$ be a real analytic manifold
equipped with an isomorphism
$W\simeq\real^m$.
Let $(y_1,\ldots,y_m)$ be the induced coordinate system on $W$.
Let $\phi:W\lrarr M$ be a real analytic map.
Let $g:U\lrarr W$ be a real analytic map
such that
(i) $F_{|U}=\phi\circ g$,
(ii) $g(U)$ is relatively compact subset of $W$,
(iii) $g$ is subanalytic on $(U,N)$.

\begin{lem}
\label{lem;16.7.26.1}
There exist a closed subanalytic subset
$Z\subset \del U$
with $\dim Z\leq \dim N-2$
such that 
(i) $\Sing(\del U)\subset Z$,
(ii) there exists a continuous subanalytic map
$\bar{g}:\Ubar\setminus Z\lrarr W$
such that $\bar{g}_{|U}=g$.
\end{lem}
\pf
We obtain the subanalytic functions
$g^{\ast}(y_i)$ on $(U,N)$.
They are bounded
because $g(U)$ is a relatively compact subset of $W$.
By Lemma \ref{lem;15.12.28.10},
there exists a closed subanalytic subset
$Z\subset \del U$
with $\dim Z\leq \dim N-2$
such that
$g^{\ast}(y_i)$ are ramified real analytic
around any point of $\del U\setminus Z$.
Then, the claim is clear.
\hfill\qed

\vspace{.1in}
Let $\Crit(\phi)$ denote the set of the critical values of $\phi$.

\begin{lem}
\label{lem;16.7.26.20}
Suppose that $\phi$ 
is obtained as the composition of 
a finite sequence of local real blowings up,
and that 
$\dim(U\cap F^{-1}(\Crit(\phi)))<\dim N$.
Then, the map $\overline{g}$ in Lemma {\rm\ref{lem;16.7.26.1}}
is real analytic on $\Ubar\setminus Z$.
Here, we regard $\Ubar\setminus Z$
as a real analytic manifold with smooth boundary.
\end{lem}
\pf
Let $P$ be any point of $\del U\setminus Z$.
 Let $(N_P;x_1,\ldots,x_n)$ 
 be a coordinate neighbourhood
of $N$ around $P$
such that $N_P\cap\del U=\{x_n=0\}$
and $N_P\cap U=\{x_n>0\}$.
Then,
$\overline{g}^{\ast}(y_j)$
are expressed as convergent power series
$\sum_{k\geq 0} a_{j,k}(x_1,\ldots,x_{n-1})x_n^{k/\rho}$
where $\rho$ is a positive integer
and $a_{j,k}$ are real analytic.
Note that 
$\dim\overline{(F^{-1}(\Crit(\phi))\cap U)}\cap \del U<\dim(\del U)-1$.
Let $Q=(x_1,x_2,\ldots,x_{n-1},0)$
be any point of 
$N_P\setminus \overline{F^{-1}(\Crit(\phi))\cap U}$.
We consider the path
$\gamma_{Q}:\II\lrarr N_P$
defined by 
$\gamma_Q(t)=(x_1,\ldots,x_{n-1},t)$.
Because $\phi$ is assumed to be the composition of
a finite sequence of local real blowings up,
and because the image of $F\circ\gamma_Q$
is contained in the image of $W$,
we can observe that 
there exists a real analytic map
$\gammatilde_Q:\II\lrarr W$
such that
$\phi\circ\gammatilde_Q=F\circ\gamma_Q$.
It implies that
$\gamma_Q^{\ast}\overline{g}^{\ast}(y_j)$
are real analytic functions on $\II$.
Hence, we obtain that $\rho=1$,
i.e.,
$\overline{g}$ is real analytic around $P$.
\hfill\qed

\subsubsection{Boundedness}

This subsection is a preliminary for the study
in \S\ref{subsection;18.11.16.102}.
Let $M$ be a real analytic manifold.
Let $\pi:M\times \II\lrarr M$ be the projection.
Let $B$ be any subanalytic subset in $M\times \II$
such that  $\dim B=\dim M+1$.
For any $x\in M$,
let $B_x:=B\cap(\{x\}\times \II)$.
Let $\overline{B_x}$ denote the closure of $B_x$.

\begin{lem}
\label{lem;16.8.12.10}
Suppose that
$\overline{B_x}$ does not contain $(x,1)$
for any $x\in M$.
Then, there exists a subanalytic closed subset 
$Z\subset M$ 
such that the following holds.
\begin{itemize}
\item $\dim Z<\dim M$.
\item
The closure of
$B\cap((M\setminus Z)\times \II)$
in $(M\setminus Z)\times \II$
does not intersect with
$(M\setminus Z)\times\{1\}$.
\end{itemize}
\end{lem}
\pf
Let $\overline{B}$ denote the closure of $B$
in $M\times \II$.
Let $B^{\circ}$ denote the interior part of $B$.
We set $R:=\Bbar\setminus B^{\circ}$.
Because $\dim R\leq \dim M$,
there exists a closed subanalytic subset
$Z_1\subset M$
such that 
(i) $\dim Z_1<\dim M$,
(ii) each connected component of
 $M\setminus Z_1$ is simply connected,
(iii) the induced map 
$R\setminus\pi^{-1}(Z_1)\lrarr M\setminus Z_1$ 
 is proper and a local diffeomorphism.
On each connected component $\nbigc$
of $M\setminus Z_1$,
there exist subanalytic functions
$h^{\nbigc}_i$ on $(\nbigc,M)$
such that 
$\pi^{-1}(\nbigc)\cap R$ 
is the union of the graph of $h^{\nbigc}_i$.
We may assume
$h^{\nbigc}_i<h^{\nbigc}_{i+1}$.
Then, $\Bbar\cap \pi^{-1}(\nbigc)$
is a union of sets of the form
$\{(x,t)\,|\,x\in\nbigc,\,\,h_i(x)\leq t\leq h_{i+1}(x)\}$.
Moreover,
$\pi^{-1}(\nbigc)\cap(\Bbar\setminus B)$
is relatively $0$-dimensional over $\nbigc$.
Then, the claim of the lemma is clear.
\hfill\qed

\vspace{.1in}

Let $Y_i$ $(i=1,2)$ be real analytic manifolds.
Let $\nbiga$ be a relatively compact subanalytic subset
of $Y_1\times Y_2$.
Let $f$ be a subanalytic function
on $(\nbiga,Y_1\times Y_2)$.
Let $q:Y_1\times Y_2\lrarr Y_2$ denote the projection.

\begin{lem}
\label{lem;16.8.12.11}
Suppose that 
$f_{|\nbiga\cap q^{-1}(x)}$ is bounded
for any $x\in Y_2$.
Then, there exists a closed subanalytic subset
$Z\subset Y_2$
such that
(i) $\dim Z<\dim Y_2$,
(ii) any $x\in Y_2\setminus Z$
 has a neighbourhood $U_x$
 in $Y_2\setminus Z$
 on which $f_{|\nbiga\cap q^{-1}(U_x)}$
 is bounded.
\end{lem}
\pf
Let $\Gamma_f\subset Y_1\times Y_2\times \realbar$
denote the graph of $f$.
Let $B_f$ denote the image of $\Gamma_f$
by the projection
$Y_1\times Y_2\times\realbar
\lrarr Y_2\times\realbar$.
It is a subanalytic subset.
Then, we obtain the claim of this lemma
from Lemma \ref{lem;16.8.12.10}.

\hfill\qed

\section{Preliminaries for 
infinite sequences of blowings up of mixed type}
\label{section;16.9.4.2}

We give some technical preliminaries
to study infinite sequences of complex blowings up
in Part \ref{part;18.11.16.31}.
The results in
\S\ref{subsection;16.5.31.30}--\ref{subsection;16.5.31.31}
will be used in \S\ref{subsection;18.11.13.1},
and the results in \S\ref{subsection;16.5.31.32}
will be used in \S\ref{subsection;18.11.13.2}.

\subsection{Vanishing}
\label{subsection;16.5.31.30}

\subsubsection{Statements}
\label{subsection;16.5.31.20}

Let $\kappahat$ be a positive real number.
Let $\nbigs$ be an infinite subset
in $\rnum\cap\openopen{0}{\kappahat}$
satisfying the following condition.
\begin{itemize}
\item
 For any $\gminiy\in\openopen{0}{\kappahat}$,
 the intersection
 $\nbigs_{\gminiy}:=\nbigs\cap\openopen{0}{\gminiy}$
 is finite.
\end{itemize}
We set
$T_m(\nbigs):=
 \seisuu_{\geq 0}\times \nbigs^m$ $(m\geq 1)$
and 
$T_0(\nbigs)=\seisuu_{\geq 0}$.
We put
$T_+(\nbigs):=
 \coprod_{m\geq 1}
 T_m(\nbigs)$
and $T(\nbigs):=T_0(\nbigs)\cup T_+(\nbigs)$.
For any element $\vecs=(s_1,\ldots,s_m)\in\nbigs^m$,
the number $m$ is denoted by $|\vecs|$.

Let $e$ be a positive integer.
For any $L\geq 0$,
we set
\[
 T_+(\nbigs,L):=
 \Bigl\{
 (i,s_1,\ldots,s_m)\in
 T_+(\nbigs)\,\Big|\,
 \frac{i}{e}+\sum_{p=1}^m s_p=L
 \Bigr\},
\]
and $T_0(\nbigs,L):=\{i\in T_0(\nbigs)\,|\,i/e=L\}$.
We set
$T(\nbigs,L):=T_+(\nbigs,L)\cup T_0(\nbigs,L)$.

Let $\nbigv_j$ $(j=1,2,\ldots)$ be a sequence of
neighbourhoods of $0$ in $\cnum$
such that $\nbigv_j\supset\nbigv_{j+1}$ for any $j$.
Let $\nbigo(\nbigv_j)$ be the space of
holomorphic functions on $\nbigv_j$.
Let $\sfa_j:\real_{>0}\lrarr \nbigo(\nbigv_j)$
$(j=1,2,\ldots)$
be maps with the following property:
\begin{itemize}
\item
$\sfa_j(\gminiy)=0$ unless $\gminiy\in \nbigs$.
\item
For each $\gminiy\in\nbigs$,
there exists $j_0(\gminiy)$
 such that 
 $\sfa_{j_1}(\gminiy)=\sfa_{j_2}(\gminiy)$
 for any $j_1,j_2\geq j_0(\gminiy)$.
Moreover,
 if $j\geq j_0(\gminiy)$,
then
 $\sfa_j(\gminiy)$ are constant functions.
\end{itemize}
We set 
$\sfa_{\infty}(\gminiy):=\lim_{j\to\infty}\sfa_j(\gminiy)$
in the stalk $\nbigo_{\cnum,0}$
of the sheaf of holomorphic functions at $0$.
Because $\sfa_{\infty}(\gminiy)$
are the germs of constant functions,
we may also regard $\sfa_{\infty}(\gminiy)$
as complex numbers.
We further assume the following.
\begin{itemize}
\item
 $\#\{\gminiy\,|\,\sfa_{\infty}(\gminiy)\neq 0\}=\infty$.
\end{itemize}

For any $\gminiy\in\real_{>0}$,
we define the real analytic functions
$R_j(\gminiy):\real\times\nbigv_j\lrarr \real$
and 
$I_j(\gminiy):\real\times\nbigv_j\lrarr \real$
$(j=1,2,\ldots)$
by 
\[
 R_j(\gminiy)(\phi):
=\Re\Bigl(
 \sfa_j(\gminiy)e^{\sqrt{-1}\gminiy\phi}
 \Bigr),
\quad
 I_j(\gminiy)(\phi):
=\Image\Bigl(
 \sfa_j(\gminiy)e^{\sqrt{-1}\gminiy\phi}
 \Bigr).
\]
For any point $\phi_0\in\real$,
the induced germs of real analytic functions at 
$(\phi_0,0)\in\real\times\cnum$
are also denoted by the same notation.
For each $\gminiy$, there exists $j_0$ such that
$R_{j_1}(\gminiy)=R_{j_2}(\gminiy)$
and $I_{j_1}(\gminiy)=I_{j_2}(\gminiy)$ for any $j_1,j_2\geq j_0$
as germs of real analytic functions at $(\phi_0,0)$.
We set
$R_{\infty}(\gminiy):=
 \lim_{j\to\infty}R_j(\gminiy)$
and 
$I_{\infty}(\gminiy):=
 \lim_{j\to\infty}I_j(\gminiy)$
in the space of germs of real analytic functions at $(\phi_0,0)$.
We define
the real analytic functions
$R_{\infty}(\gminiy),I_{\infty}(\gminiy):\real\lrarr\real$
by 
\[
 R_{\infty}(\gminiy)(\phi):
=\Re\Bigl(
 \sfa_{\infty}(\gminiy)e^{\sqrt{-1}\gminiy\phi}
 \Bigr),
\quad
 I_{\infty}(\gminiy)(\phi):
=\Image\Bigl(
 \sfa_{\infty}(\gminiy)e^{\sqrt{-1}\gminiy\phi}
 \Bigr).
\]
We may also regard them as real analytic functions on $\real$,
or germs of real analytic functions at $\phi_0\in\real$.

Let $\nbigi$ be an interval.
Let $f_{i,c_1,c_2}:\nbigi\lrarr \real$
$((i,c_1,c_2)\in\seisuu_{\geq 0}^3)$
be a family of real analytic functions.
We assume that $f_{0,0,0}$ is constantly $0$.

Let $L$ be a non-negative real number.
For any $(\ell_1,\ell_2,m)\in\seisuu_{\geq 0}^3$
satisfying $m=\ell_1+\ell_2$,
we obtain the following functions 
on $\nbigi\times \nbigv_j$ $(j=1,2,\ldots)$:
\begin{multline}
 A^{(m)}_{j}(L,(\ell_1,\ell_2)):=
 \sum_{(i,\veczeta)\in T_+(\nbigs,L)}
 \sum_{\substack{(c_1,c_2)\in\seisuu_{\geq 0}^2\\
 c_1+c_2=|\veczeta|+m\\
 c_i\geq \ell_i}}
 f_{i,c_1,c_2}
 \frac{c_1!c_2!}{(c_1-\ell_1)!(c_2-\ell_2)!}
 \prod_{q=1}^{c_1-\ell_1}R_{j}(\zeta_q)
 \prod_{q=c_1-\ell_1+1}^{|\veczeta|}I_{j}(\zeta_q)
\\
+\sum_{i\in T_0(\nbigs,L)}
 \ell_1!\ell_2!
 f_{i,\ell_1,\ell_2}.
\end{multline}
Here, we use the convention
$\prod_{q=c}^{c-1}a_q:=1$.
We fix any point $\phi_0\in\nbigi$,
and the induced germs of real analytic functions
at $(\phi_0,0)\in\nbigi\times\cnum$
are also denoted by
$A_j^{(m)}(L,(\ell_1,\ell_2))$.

For each pair of $L$ and $m$,
the set 
$\bigl\{
 (i,\veczeta,c_1,c_2)\,\big|\,
 (i,\veczeta)\in T_+(\nbigs,L),\,\,
 c_1+c_2=|\veczeta|+m,\,\,c_i\geq \ell_i
 \bigr\}$
is finite.
Hence, there exists $j_0$ depending on $(L,m)$
such that 
$A_{j_1}^{(m)}(L,\ell_1,\ell_2)
=A_{j_2}^{(m)}(L,\ell_1,\ell_2)$
for any $j_1,j_2\geq j_0$
as germs of real analytic functions
at $(\phi_0,0)\in\nbigi\times\cnum$.
We set
$A_{\infty}^{(m)}(L,\ell_1,\ell_2):=
 \lim_{j\to\infty}A^{(m)}_{j}(L,\ell_1,\ell_2)$
in the space of germs of real analytic functions 
at $(\phi_0,0)$.
We obtain
\begin{multline}
 A^{(m)}_{\infty}(L,(\ell_1,\ell_2))=
 \sum_{(i,\veczeta)\in T_+(\nbigs,L)}
 \sum_{\substack{(c_1,c_2)\in\seisuu_{\geq 0}^2\\
 c_1+c_2=|\veczeta|+m\\c_i\geq \ell_i}}
 f_{i,c_1,c_2}
 \frac{c_1!c_2!}{(c_1-\ell_1)!(c_2-\ell_2)!}
 \prod_{q=1}^{c_1-\ell_1}R_{\infty}(\zeta_q)
 \prod_{q=c_1-\ell_1+1}^{|\veczeta|}I_{\infty}(\zeta_q)
\\
+\sum_{i\in T_0(\nbigs,L)}
 \ell_1!\ell_2!
 f_{i,\ell_1,\ell_2}.
\end{multline}
We may also regard them
as real analytic functions on $\real$.
When $(m,\ell_1,\ell_2)=(0,0,0)$,
we use the notation
$A^{(0)}_{j}(L)$
and 
$A^{(0)}_{\infty}(L)$,
instead of 
$A^{(0)}_{j}(L,(0,0))$
and 
$A^{(0)}_{\infty}(L,(0,0))$.
We shall prove the following propositions.

\begin{thm}
\label{thm;16.4.6.11}
Let $L_1$ be the real number such that 
$A^{(0)}_{\infty}(L)=0$ for any $L<L_1$.
Then, the following claims hold:
\begin{itemize}
\item
 $A_{\infty}^{(m)}(L-m\kappahat,(\ell_1,\ell_2))=0$ holds
 for any $L\leq L_1$,
 and any $(\ell_1,\ell_2)\in\seisuu_{\geq 0}^2$
 satisfying $\ell_1+\ell_2=m\geq 1$.
\item
There exists $j_0$ such that
  $A_{j}^{(m)}(L-m\kappahat,(\ell_1,\ell_2))=0$
 hold
 for any $L\leq L_1$,
 any $j\geq j_0$,
 and any $(\ell_1,\ell_2)\in\seisuu_{\geq 0}^2$
 satisfying $\ell_1+\ell_2=m\geq 1$.
 Moreover, 
 $A_{j}^{(0)}(L_1)=A_{\infty}^{(0)}(L_1)$
 and $A_j^{(0)}(L)=0$ $(L<L_1)$ hold
 for any $j\geq j_0$.
\end{itemize}
\end{thm}

\begin{thm}
\label{thm;16.4.6.12}
Suppose that $A^{(0)}_{\infty}(L)=0$ hold
for any $L\geq 0$.
Then, $f_{i,c_1,c_2}=0$ hold
for any $(i,c_1,c_2)$.
\end{thm}

\subsubsection{Preliminary}

As a preliminary, we consider the following condition
for subsets $S\subset\real_{\geq 0}$.
\begin{condition}
\label{condition;18.11.23.100}
 For any $a\in\real_{\geq 0}$,
 there exists $\epsilon>0$
 such that
 $S\cap \openopen{a}{a+\epsilon}$
 is empty.
\end{condition}

\begin{lem}
\label{lem;18.1.13.1}
Let $S$ be a subset in $\real_{\geq 0}$
satisfying Condition {\rm\ref{condition;18.11.23.100}}.
Let $s_i$ $(i=1,2,\ldots)$
be a decreasing sequence in $S$.
Then, there exists $i_0$
such that $s_i=s_{i_0}$ for any $i\geq i_0$.
\end{lem}
\pf
We set $s_{\infty}:=\underset{i\to\infty}{\lim}s_i$.
There exists $\epsilon>0$ such that
$\openopen{s_{\infty}}{s_{\infty}+\epsilon}\cap S
=\emptyset$.
Then, we obtain
$s_i=s_{\infty}$ for any sufficiently large $i$.
\hfill\qed

\vspace{.1in}

For any subset $S\subset\real_{\geq 0}$,
let $\Accum(S)$ denote the set of
accumulation points in $S$,
i.e.,
$P\in \real_{\geq 0}$ is contained in
$\Accum(S)$
if and only if
$P$ is contained in the closure of
$S\setminus\{P\}$.

\begin{lem}
\label{lem;16.6.1.11}
Let $S_i$ $(i=1,2)$ be subsets
in $\real_{\geq 0}$
satisfying Condition {\rm\ref{condition;18.11.23.100}}.
Then,
$S_1\cup S_2$ satisfies 
Condition {\rm\ref{condition;18.11.23.100}}.
Moreover,
$\Accum(S_1\cup S_2)=
\Accum(S_1)\cup \Accum(S_2)$
holds.
\end{lem}
\pf
For any $a\in\real_{\geq 0}$,
there exist $\epsilon_i>0$
such that
$\openopen{a}{a+\epsilon_i}\cap S_i=\emptyset$.
Set $\epsilon_0:=\min\{\epsilon_1,\epsilon_2\}$.
Then, we obtain
$\openopen{a}{a+\epsilon_0}
\cap(S_1\cup S_2)=\emptyset$.
Hence, $S_1\cup S_2$ also satisfies 
Condition \ref{condition;18.11.23.100}.

Clearly
$\Accum(S_1\cup S_2)\supset
\Accum(S_1)\cup \Accum(S_2)$
holds.
Let $s\in \Accum(S_1\cup S_2)$.
There exists a sequence
$s_i\in (S_1\cup S_2)\setminus\{s\}$
such that $\lim s_i=s$.
One of $S_1$ or $S_2$ contains
infinite subsequence of $s_i$.
Hence, $s\in \Accum(S_1)\cup\Accum(S_2)$.
\hfill\qed

\begin{lem}
\label{lem;16.6.1.1}
Let $S_i\subset \real_{\geq 0}$ $(i=1,2,\ldots)$
be subsets satisfying 
Condition {\rm \ref{condition;18.11.23.100}}.
Suppose that 
\[
 \lim_{i\to\infty}
 \inf\bigl(S_i\setminus\{0\}\bigr)
=\infty.
\]
Then, 
$\bigcup_{i\geq 1}S_i$ also satisfies 
Condition {\rm\ref{condition;18.11.23.100}}.
Moreover,
$\Accum\Bigl(
 \bigcup_{i\geq 1}S_i
 \Bigr)=\bigcup_{i\geq 1}\Accum(S_i)$.
\end{lem}
\pf
For any any $a\in\real_{\geq 0}$ and $\epsilon_0>0$,
there exists $i_0$ such that
$\openopen{a}{a+\epsilon_0}\cap S_i=\emptyset$
for any $i\geq i_0$.
We obtain
$\openopen{a}{a+\epsilon_0}
\cap
\Bigl(
\bigcup_{i\geq 1}S_i
\Bigr)
=
\openopen{a}{a+\epsilon_0}
\cap
\Bigl(
\bigcup_{i=1}^{i_0-1}S_i
\Bigr)$.
Because 
$\bigcup_{i=1}^{i_0-1}S_i$
satisfies 
Condition \ref{condition;18.11.23.100}.
we obtain the first claim of the lemma.
Clearly,
$\Accum\Bigl(
 \bigcup_{i\geq 1}S_i
 \Bigr)
\supset
 \bigcup_{i\geq 1}
 \Accum(S_i)$ holds.
For any $s\in \Accum\Bigl(\bigcup_{i\geq 1}S_i\Bigr)$
and $\epsilon>0$,
there exists $i_0$ such that 
$\openopen{s-\epsilon}{s+\epsilon}
 \cap\Accum\Bigl(\bigcup_{i\geq 1}S_i\Bigr)
=\openopen{s-\epsilon}{s+\epsilon}
 \cap\Accum\Bigl(\bigcup_{i=1}^{i_0}S_i\Bigr)$.
Then, we obtain the second claim
by using the argument in the proof of 
Lemma \ref{lem;16.6.1.11}.
\hfill\qed

\vspace{.1in}
The following lemma is clear.

\begin{lem}
\label{lem;18.1.13.2}
Let $S$ be a subset in $\real_{\geq 0}$
satisfying
Condition {\rm\ref{condition;18.11.23.100}}.
Then, any subset $S'\subset S$
also satisfies
Condition {\rm\ref{condition;18.11.23.100}}.
\hfill\qed
\end{lem}

\vspace{.1in}

Let $S_i\subset\real_{\geq 0}$ $(i=1,2)$
be subsets satisfying the following conditions.
\begin{itemize}
\item
There exists  $\kappahat>0$
 such that
 $S_1$ is contained in 
 $\openopen{0}{\kappahat}$,
 and that $S_1$ is an infinite set.
 Moreover,
 for any $\gminiy<\kappahat$,
 the intersection
 $S_1\cap \openclosed{0}{\gminiy}$
 is finite.
\item
 $S_2$ satisfies
Condition \ref{condition;18.11.23.100}.
\end{itemize}

We set
$S_1+S_2:=
 \bigl\{t\in\real_{\geq 0}\,\big|\,
 \exists(a_1,a_2)\in S_1\times S_2,\,\,
 a_1+a_2=t
 \bigr\}$.

\begin{lem}
\label{lem;16.4.4.10}
$S_1+S_2$ also satisfies 
Condition {\rm\ref{condition;18.11.23.100}}.
Moreover,
for any $c\in S_1+S_2$,
the set
$\bigl\{(c_1,c_2)\in S_1\times S_2\,\big|\,
 c_1+c_2=c\bigr\}$
is finite.
\end{lem}
\pf
Let $a$ be any real number.
Let $t_i$ be a decreasing sequence in $S_1+S_2$
such that  $\lim t_i=a$.
There exist $(b_{i,1},b_{i,2})\in S_1\times S_2$
such that $t_i=b_{i,1}+b_{i,2}$.
By taking a sub-sequence,
we may assume
$b_{i,1}\leq b_{i+1,1}$ for any $i$.
Then,
$b_{i,2}=t_{i}-b_{i,1}$ is decreasing.
Hence, there exists $i_0$
such that
$b_{i,2}=b_{i_0,2}$ for any $i\geq i_0$.
Because $t_i$ is decreasing
and $b_{i,1}$ is increasing,
there exists  $i_1$ such that 
$t_{i}=t_{i_1}$
and $b_{i,1}=b_{i_1,1}$ 
for any $i\geq i_1$,
i.e.,
$t_i=a$ for any $i\geq i_1$.
Hence, we obtain the first claim.

Let $c\in S_1+S_2$.
If the set 
$\bigl\{(c_1,c_2)\in S_1\times S_2\,\big|\,
 c_1+c_2=c\bigr\}$
is infinite,
there exists an infinite sequence
$(c_{i,1},c_{i,2})$ in $S_1\times S_2$
such that
(i) $c_{i,1}+c_{i,2}=c$,
(ii) $c_{i,1}$ is increasing.
Then, $c_{i,2}$ is decreasing,
and hence there exists $i_0$
such that $c_{i,2}=c_{i_0,2}$ for any $i\geq i_0$.
It contradicts the assumption that
the sequence $(c_{i,1},c_{i,2})$ is infinite.
\hfill\qed

\vspace{.1in}

Let $\Sbar_i$ $(i=1,2)$ denote the closure of
the above sets $S_i$ in $\real_{\geq 0}$.
Let $\Accum'(S_1+ S_2)$
denote the image of the map
$\bigl(
\Accum(S_1)\times \Sbar_2
\bigr)
\cup
\bigl(
 \Sbar_1\times\Accum(S_2)\bigr)
\lrarr
\real_{\geq 0}$
defined by
$(c_1,c_2)\longmapsto c_1+c_2$.
Note that
$\Accum(S_1)=\{\kappahat\}$.

\begin{lem}
\label{lem;16.6.1.10}
$\Accum(S_1+ S_2)
=\Accum'(S_1+S_2)$
holds.
\end{lem}
\pf
Clearly,
$\Accum(S_1+S_2)
\supset\Accum'(S_1+ S_2)$
holds.
Let us prove 
$\Accum(S_1+S_2)
\subset\Accum'(S_1+S_2)$.
Let $c$ be any element of $\Accum(S_1+S_2)$.
There exists an increasing sequence $c_i$
in $(S_1+S_2)\setminus\{c_i\}$
such that
$\underset{i\to\infty}\lim c_i=c$.
There exist
$(c_i^{(1)},c_i^{(2)})\in S_1\times S_2$
such that
$c_i^{(1)}+c_i^{(2)}=c_i$.
By taking a sub-sequence,
we may assume that
$c_i^{(1)}$ is increasing.
Suppose that there exists $i_0$ such that
$c_i^{(1)}=c_{i_0}^{(1)}$ for any $i\geq i_0$.
Then, the sequence $c_i^{(2)}=c_i-c_i^{(1)}$
is convergent to $c-c_{i_0}^{(1)}$.
Hence,
we obtain
$c-c_{i_0}^{(1)}\in\Accum(S_2)$,
and 
$c\in \Accum'(S_1\times S_2)$.
Suppose that
$\underset{i\to\infty}{\lim}c_i^{(1)}=\kappahat$.
Then, the sequence $c_i^{(2)}=c_i-c_i^{(1)}$
is convergent to 
$c-\kappahat$.
Hence,
we obtain $c-\kappahat\in\Sbar_2$,
and $c\in \Accum'(S_1\times S_2)$.
\hfill\qed

\vspace{.1in}
Let $S_1$ be the above set.
For $\ell\geq 1$,
let $\sum^{\ell}S_1$
denote the image of the map
$S_1^{\ell}\lrarr\real_{\geq 0}$
defined by
$(s_1,\ldots,s_{\ell})\longmapsto
 \sum_{i=1}^{\ell} s_i$.
We formally set $\sum^0S_1:=\{0\}$.

\begin{cor}
\label{cor;16.6.1.2}
The set $\sum^{\ell}S_1$
satisfies 
Condition {\rm\ref{condition;18.11.23.100}}.
Moreover, the following holds:
\[
 \Accum\bigl(\sum^{\ell}S_1\bigr)
=\bigcup_{1\leq m\leq \ell}
\Bigl(
 \{m\kappahat\}
+\sum^{\ell-m}S_1
\Bigr).
\]
\end{cor}
\pf
We obtain the first claim from Lemma \ref{lem;16.4.4.10}
and an induction.
Let us study the second claim.
We use an induction on $\ell$.
The claim is clear in the cases $\ell=0,1$.
According to Lemma \ref{lem;16.6.1.10},
$\Accum(\sum^{\ell}S_1)
=\bigl(
 \{\kappahat\}+\overline{\sum^{\ell-1}S_1}
 \bigr)
\cup
 \bigl(
 \overline{S_1}+\Accum(\sum^{\ell-1}S_1)
 \bigr)$
holds.
By the assumption of the induction,
$\Accum(\sum^{\ell-1}S_1)
=\bigcup_{1\leq m\leq \ell-1}
 \Bigl(
 \{m\kappahat\}+\sum^{\ell-1-m}S_1
 \Bigr)$
holds.
Note $\overline{S_1}=S_1\sqcup\{\kappahat\}$.
Hence, we obtain
\begin{multline}
 \overline{S_1}+\Accum\Bigl(\sum^{\ell-1}S_1\Bigr)
=\bigcup_{1\leq m\leq\ell-1}
 \Bigl(
 \{(m+1)\kappahat\}+\sum^{\ell-1-m}S_1
 \Bigr)
\cup
 \bigcup_{1\leq m\leq \ell-1}
 \Bigl(
 \{m\kappahat\}+\sum^{\ell-m}S_1
 \Bigr)
 \\
=\bigcup_{1\leq m\leq \ell}
\Bigl(
 \{m\kappahat\}+\sum^{\ell-m}S_1
 \Bigr).
\end{multline}
Because
$\{\kappahat\}+\overline{\sum^{\ell-1}S_1}
=\bigl(
 \{\kappahat\}+
 \sum^{\ell-1}S_1
 \bigr)\cup
 \bigl(
 \{\kappahat\}+
 \Accum\bigl(\sum^{\ell-1}S_1\bigr)
 \bigr)$,
we obtain the claim in the case $\ell$.
\hfill\qed

\begin{cor}
\label{cor;16.6.1.3}
Let $e$ be a positive integer.
The set 
$\bigcup_{i\geq 0}
 \bigcup_{\ell\geq 0}
\bigl(
\{i/e\}+\sum^{\ell}S_1
 \bigr)$
satisfies 
Condition {\rm\ref{condition;18.11.23.100}}.
Moreover, the following holds:
\[
\Accum\Bigl(
 \bigcup_{i\geq 0}
 \bigcup_{\ell\geq 0}
\bigl(
\{i/e\}+\sum^{\ell}S_1
 \bigr)
\Bigr)
=\bigcup_{i\geq 0}
 \bigcup_{m\geq 1}
 \bigcup_{\ell\geq 0}
 \Bigl(
 \{i/e+m\kappahat\}
+\sum^{\ell}S_1
\Bigr).
\]
\end{cor}
\pf
Let $\beta_0>0$ be the infimum of $S_1$.
The infimum of
$\{i/e\}+\sum^{\ell}S_1$
is $i/e+\ell\beta_0$.
The claim follows from
Lemma \ref{lem;16.6.1.1}
and Corollary \ref{cor;16.6.1.2}.
\hfill\qed

\subsubsection{Accumulation}

We return to the situation in \S\ref{subsection;16.5.31.20}.
We set
$\nbigt(\nbigs):=\bigl\{
 i/e+\sum\zeta_p\,\big|\,
 (i,\veczeta)\in T_+(\nbigs)
 \bigr\}
\cup \frac{1}{e}\seisuu_{\geq 0}$.
According to Corollary \ref{cor;16.6.1.3},
the set $\nbigt(\nbigs)$
satisfies 
Condition \ref{condition;18.11.23.100}.
and the following holds:
\[
  \Accum\bigl(\nbigt(\nbigs)\bigr)
=\bigcup_{m\geq 1}
 \Bigl(
 \{m\kappahat\}
+\nbigt(\nbigs)
 \Bigr).
\]
Similarly,
we put
$\nbigt(\nbigs_{\gminiy}):=\bigl\{
 i/e+\sum_{j=1}^m\zeta_j\,\big|\,
 (i,\veczeta)\in T_+(\nbigs_{\gminiy})
 \bigr\}
\cup \frac{1}{e}\seisuu_{\geq 0}$
for any $\gminiy\in\openopen{0}{\kappahat}$.
The set $\nbigt(\nbigs_{\gminiy})$
is discrete.

Fix $L\in\real_{\geq 0}$.
There exists $\delta>0$
such that the following holds
for any $m\geq 0$:
\[
\openopen{L-m\kappahat}{L-m\kappahat+\delta}
\cap
\nbigt(\nbigs)
=\emptyset.
\]
Because the set
$\bigcup_{m\geq 0}
 T(\nbigs,L-m\kappahat)$ is finite,
there exists
$\gminiy_0\in\openopen{\kappahat-\delta}{\kappahat}$
such that
\[
 \bigcup_{m\geq 0}
 T(\nbigs,L-m\kappahat)
=\bigcup_{m\geq 0}
 T(\nbigs_{\gminiy_0},L-m\kappahat).
\]

\begin{lem}
\label{lem;16.5.31.10}
For any $\gminiy>\gminiy_0$,
there exists $\epsilon>0$
such that the following holds:
\begin{itemize}
\item
If $(i,\zeta_1,\ldots,\zeta_q)\in T(\nbigs)$
satisfies
$L-m\kappahat-\epsilon<i/e+\sum\zeta_p<L-m\kappahat$ 
for some $m\geq 0$,
then $\max\{\zeta_p\}>\gminiy$ holds.
\end{itemize}
\end{lem}
\pf
Because
the set $\nbigt(\nbigs_{\gminiy})$ is discrete,
there exists $\epsilon>0$
satisfying
$\epsilon<
 L-m\kappahat
-\max\{\nu\in\nbigt(\nbigs_{\gminiy})\,|\,\nu<L-m\kappahat\}$
for any $m\geq 0$.
Then, the condition is satisfied.
\hfill\qed

\begin{lem}
\label{lem;16.6.1.20}
If $(i,\zeta_1,\ldots,\zeta_q)\in T(\nbigs)$
 satisfies
 $i/e+\sum\zeta_p=L-m\kappahat$ for some $m\geq 0$,
then $\max\{\zeta_p\}<\gminiy_0$.
\end{lem}
\pf
Suppose that $(i,\zeta_1,\ldots,\zeta_q)\in T(\nbigs)$
 satisfies
 $i/e+\sum\zeta_p=L-m\kappahat$ for some $m\geq 0$.
We assume that $\zeta_1\geq \gminiy_0$,
and we shall derive a contradiction.
Because
$0<\kappahat-\zeta_1\leq \kappahat-\gminiy_0<\delta$,
we obtain 
$L-(m+1)\kappahat<
 i/e+\sum_{p=2}^q\zeta_p
 < L-(m+1)\kappahat+\delta$.
Such $(i,\zeta_2,\ldots,\zeta_q)$
cannot exist by our choice of $\delta$.
\hfill\qed

\begin{lem}
\label{lem;16.4.5.2}
For any $\gminiy\in\openopen{\gminiy_0}{\kappahat}$
there exists $\epsilon>0$
such that the following holds:
\begin{itemize}
\item
If $(i,\zeta_1,\ldots,\zeta_q)\in T(\nbigs)$
 satisfies
 (i) $L-\epsilon<i/e+\sum\zeta_p<L$,
 (ii) $\zeta_1\geq\cdots\geq \zeta_q$,
then  there exists a unique positive integer
 $\ell$ 
 for which the following holds:
\[
\zeta_1\geq\cdots\geq \zeta_{\ell}>\gminiy,
\quad
\zeta_{\ell+1}+\cdots+\zeta_q=L-\ell\kappahat.
\]
\end{itemize}
Note that 
$\zeta_{p}<\gminiy_0$ $(p=\ell+1,\ldots,q)$
hold by Lemma {\rm\ref{lem;16.6.1.20}}.
\end{lem}
\pf
We set
$\gamma(L):=
 \#\Bigl\{
 m\geq 1\,\Big|\,
 L-m\kappahat
 \in \overline{\nbigt(\nbigs)}
\Bigr\}$.
We use an induction on $\gamma(L)$.

Suppose $\gamma(L)=0$.
For any $\gminiy>\gminiy_0$,
let $\epsilon>0$ be as in Lemma \ref{lem;16.5.31.10}.
Moreover, if $\epsilon$ is sufficiently small,
then
$\nbigt(\nbigs)\cap
 \openopen{L-\kappahat-\epsilon}{L-\kappahat+\delta}
=\emptyset$ holds.
If $(i,\zeta_1,\ldots,\zeta_q)$ satisfies
$L-\epsilon<i/e+\sum\zeta_j<L$
and $\zeta_1\geq\cdots\geq\zeta_q$,
then $\zeta_1>\gminiy$ holds,
which implies
\[
 L-\kappahat-\epsilon
<i/e+\sum_{j\geq 2}\zeta_j
<L-\kappahat+(\kappahat-\gminiy)
<L-\kappahat+\delta.
\]
Hence, there does not exist such 
$(i,\zeta_1,\ldots,\zeta_q)$.
Thus, we are done in the case $\gamma(L)=0$.

Suppose that we have already proved
the case $\gamma(L)\leq n$,
and we shall prove the claim 
in the case $\gamma(L)=n+1$.
Note that
$\gamma(L)=n+1$
implies
$\gamma(L-\kappahat)=n$.
Indeed, because $n+1\geq 1$,
there exists $m\geq 1$ such that
$L-m\kappahat\in \overline{\nbigt(S)}$.
Then, it is easy to observe
$L-\kappahat=L-m\kappahat+(m-1)\kappahat
 \in\overline{\nbigt(S)}$.
Hence, $\gamma(L-\kappahat)=\gamma(L)-1=n$.

Let $\gminiy>\gminiy_0$.
By the assumption of the induction,
there exists a positive constant $\epsilon_1$ 
for the pair of $L-\kappahat$ and $\gminiy$
satisfying the property in the statement of Lemma \ref{lem;16.4.5.2}.
There also exists a positive constant $\epsilon$
for the pair of $L-\kappahat$ and $\gminiy$
as in Lemma \ref{lem;16.5.31.10}.
Let $\epsilon_2$ be a real number satisfying
$0<\epsilon_2<\min\{\epsilon_1,\epsilon\}$.
Suppose that 
$(i,\zeta_1,\ldots,\zeta_q)$
satisfies
$L-\epsilon_2<i/e+\sum\zeta_j<L$
and $\zeta_1\geq\cdots\geq\zeta_q$.
Because $\kappahat>\zeta_1>\gminiy$,
the following holds:
\[
 L-\kappahat-\epsilon_2
<i/e+\sum_{j\geq 2}\zeta_j
<(L-\kappahat)+(\kappahat-\zeta_1)
<L-\kappahat+\delta.
\]
By our choice of $\delta$,
we obtain
$i/e+\sum_{j\geq 2}\zeta_j
\leq L-\kappahat$.
It implies either one of the following:
\begin{description}
\item[(i)]
$i/e+\sum_{j\geq 2}\zeta_j=L-\kappahat$.
\item[(ii)]
$L-\kappahat-\epsilon_2
<i/e+\sum_{j\geq 2}\zeta_j<L-\kappahat$.
\end{description}
If (i) occurs, we are done by Lemma \ref{lem;16.6.1.20}.
If (ii) occurs, we may apply the hypothesis of the induction.
\hfill\qed

\vspace{.1in}

Let $L\in \real_{\geq 0}$.
Suppose that
the set 
$\{m\in\seisuu_{\geq 0}\,|\,L-m\kappahat
 \in\overline{\nbigt(\nbigs)}\}$
is not empty.
Let $m_0$ be the maximum element
of the set.
\begin{lem}
\label{lem;16.6.1.30}
$L-m_0\kappahat$ is an isolated point in $\nbigt(\nbigs)$.
\end{lem}
\pf
If $L-m_0\kappahat$ is an accumulation point of
$\nbigt(\nbigs)$,
it is contained in
$\bigcup_{m\geq 1}
 \bigl(\{m\kappahat\}+\nbigt(\nbigs)\bigr)$,
i.e.,
$L-m_0\kappahat
=m\kappahat+s$
for some $m\geq 1$
and $s\in\nbigt(\nbigs)$.
We obtain
$L-(m_0+m)\kappahat=s\in\nbigt(\nbigs)$,
which contradicts the choice of $m_0$.
Hence, $L-m_0\kappahat$ is an isolated point
of $\nbigt(\nbigs)$.
\hfill\qed

\subsubsection{A formula}

For any $L\geq 0$, $\gminiy\in\openopen{0}{\kappahat}$
and $d\in\seisuu_{>0}$,
we set 
$\nbigu(d,\gminiy,L):=\bigl\{
 (\zeta_1,\ldots,\zeta_d)\in \nbigs^d\,\big|\,
 \sum_{i=1}^d\zeta_i=L,\,\,
 \zeta_i>\gminiy
 \bigr\}$.
For any $(d_1,d_2,d)\in\seisuu_{\geq 0}^3$ such that 
$d_1+d_2=d$,
and for any $j=1,\ldots,\infty$,
we set
\[
 B_j(L,\gminiy,(d_1,d_2)):=
 \sum_{\veczeta\in \nbigu(d,\gminiy,L)}
 \frac{1}{d_1!d_2!}
 \prod_{q=1}^{d_1}R_j(\zeta_q)
 \prod_{q=d_1+1}^{d}I_j(\zeta_q).
\]

\begin{cor}
\label{cor;16.4.5.10}
For any $L\in\real_{\geq 0}$,
there exist $\epsilon_0>0$ and $0<\gminiy_0<\kappahat$
such that the following holds
for any $L'$ with $L-\epsilon_0<L'<L$:
\begin{equation}
\label{eq;16.4.5.1}
 A_{j}^{(m)}(L',(b_1,b_2))=\\
 \sum_{d\geq 1}
 \sum_{\substack{d_1+d_2=d\\d_i\geq 0}}
 A_j^{(m+d)}\bigl(
 L-d\kappahat,
 (d_1,d_2)+(b_1,b_2)
 \bigr)
\cdot
 B_j(L'-(L-d\kappahat),\gminiy_0,(d_1,d_2)).
\end{equation}
\end{cor}
\pf
It follows from
Lemma \ref{lem;16.6.1.20}
and Lemma \ref{lem;16.4.5.2}.
\hfill\qed

\subsubsection{Vanishing at $\infty$}

We fix $L>0$.
\begin{lem}
\label{lem;16.6.1.40}
Suppose that
$A^{(0)}_{\infty}(L')=0$ for any $L'< L$.
Then,
$A^{(m)}_{\infty}(L'-m\kappahat,(b_1,b_2))=0$
holds for any $L'\leq L$ and 
$(b_1,b_2,m)\in\seisuu_{\geq 0}^2\times\seisuu_{\geq 1}$
satisfying $b_1+b_2=m$.
\end{lem}
\pf
We set
$\nbigttilde(\nbigs):=
 \bigcup_{m\geq 0}
 \bigl(
 \{m\kappahat\}+\nbigt(\nbigs)
 \bigr)$.
The set $\nbigttilde(\nbigs)$
satisfies 
Condition \ref{condition;18.11.23.100}.

If $L'\not\in \nbigttilde(\nbigs)$,
then 
$A^{(m)}_{\infty}(L'-m\kappahat,(b_1,b_2))=0$
holds for any $m\geq 0$.
Hence, it is enough to prove the vanishing
$A^{(m)}_{\infty}(L'-m\kappahat,(b_1,b_2))=0$ $(m\geq 1)$
under the following assumption:
\begin{description}
\item[(Q)]
$A^{(m)}_{\infty}(L''-m\kappahat,(b_1,b_2))=0$ $(m\geq 0)$
for any $L''<L'$.
\end{description}

It is enough to consider the case that
the set
$\bigl\{
 m\in\seisuu_{\geq 1}\,\big|\,
 L'-m\kappahat\in\overline{\nbigt(\nbigs)}
 \bigr\}$
is not empty.
Let $m_0$ be the maximum element of the set.
By the choice of $m_0$,
$A^{(m)}_{\infty}(L'-m\kappahat,(b_1,b_2))=0$
holds for any $m>m_0$
and for any $(b_1,b_2)\in\seisuu_{\geq 0}^2$
satisfying $b_1+b_2=m$.
Let us prove the vanishing for any $m$
by a descending induction.

According to Lemma \ref{lem;16.6.1.30},
$L'-m_0\kappahat$ is an isolated point
of $\nbigt(\nbigs)$.
By the assumption {\bf(Q)},
the following holds
for any $\gminiy\in\nbigs$
and for any $(b_1,b_2)\in\seisuu_{\geq 0}^2$
with $b_1+b_2=m_0-1$:
\[
A_{\infty}^{(m_0-1)}(L'-m_0\kappahat+\gminiy,(b_1,b_2))=
A_{\infty}^{(m_0-1)}(L'+\gminiy-\kappahat-(m_0-1)\kappahat,(b_1,b_2))
=0.
\]
If $\gminiy$ is sufficiently close to $\kappahat$,
we obtain the following equality 
from Corollary \ref{cor;16.4.5.10}:
\[
0=A_{\infty}^{(m_0-1)}(L'-m_0\kappahat+\gminiy,(b_1,b_2))
=\sum_{c_1+c_2=1}
A_{\infty}^{(m_0)}(L'-m_0\kappahat,(c_1,c_2)+(b_1,b_2))
 \cdot 
 B_{\infty}(\gminiy,\gminiy_0,(c_1,c_2)).
\]
Then, 
by varying $\gminiy$,
we obtain 
$A_{\infty}^{(m_0)}(L'-m_0\kappahat,(d_1,d_2))=0$
for any $(d_1,d_2)\in\seisuu^2_{\geq 0}$
satisfying $d_1+d_2=m_0$.

Suppose that we have already known
the vanishing 
$A_{\infty}^{(m)}(L'-m\kappahat,(b_1,b_2))=0$
for any $(b_1,b_2)\in\seisuu^2_{\geq 0}$
satisfying $b_1+b_2=m$
for any $m>m_1\geq 1$,
and let us prove the vanishing in the case of $m_1$.
By the assumption {\bf(Q)},
$A_{\infty}^{(m_1-1)}(L'-m_1\kappahat+\gminiy,(b_1,b_2))=0$
holds for any $\gminiy\in\nbigs$
and for any $(b_1,b_2)\in\seisuu_{\geq 0}^2$ 
satisfying $b_1+b_2=m_1-1$.
By Corollary \ref{cor;16.4.5.10},
the following equality holds:
\begin{multline}
 0=A_{\infty}^{(m_1-1)}(L'-m_1\kappahat+\gminiy,(b_1,b_2))=
 \\
 \sum_{d\geq 0}
 \sum_{c_1+c_2=d+1}
 A_{\infty}^{(m_1+d)}(L'-(m_1+d)\kappahat,(c_1,c_2)+(b_1,b_2))
 \cdot
 B_{\infty}(\gminiy+d\kappahat,\gminiy_0,(c_1,c_2)).
\end{multline}
By using the assumption of the induction,
we can rewrite it as follows:
\[
 0=\sum_{c_1+c_2=1}
 A_{\infty}^{(m_1)}(L'-m_1\kappahat,(c_1,c_2)+(b_1,b_2))
 \cdot
 B_{\infty}(\gminiy,\gminiy_0,(c_1,c_2)).
\]
Hence, we obtain 
$A_{\infty}^{(m_1)}(L'-m_1\kappahat,(d_1,d_2))=0$
for any $(d_1,d_2)\in\seisuu_{\geq 0}^2$
satisfying $d_1+d_2=m_1$.
Thus, we obtain Lemma \ref{lem;16.6.1.40}.
\hfill\qed

\subsubsection{Proof of Theorem \ref{thm;16.4.6.11}}

Let us prove Theorem \ref{thm;16.4.6.11}.
We have already known the first claim
by Lemma \ref{lem;16.6.1.40}.

\begin{lem}
\label{lem;16.4.6.21}
For any $L\leq L_1$,
there exists $j_0(L)\in\seisuu_{\geq 1}$
and $\epsilon(L)>0$
such that 
$A^{(m)}_{j}(L'-m\kappahat,(c_1,c_2))=0$ holds
for any 
$L'\in
 \openopen{L-\epsilon(L)}{L+\epsilon(L)}
 \cap
 \closedclosed{0}{L_1}$,
any 
$j\geq j_0(L)$,
and any
$(c_1,c_2,m)\in\seisuu_{\geq 0}^2\times\seisuu_{\geq 1}$
satisfying $c_1+c_2=m$.
Moreover,
$A^{(0)}_j(L')=A^{(0)}_{\infty}(L')$ holds
for any 
$L'\in
 \openopen{L-\epsilon(L)}{L+\epsilon(L)}
 \cap
 \closedclosed{0}{L_1}$,
and any 
$j\geq j_0(L)$.
\end{lem}
\pf
It is enough to study
the claim for 
$L'\in
 \openopen{L-\epsilon(L)}{L}
 \cap
 \closedclosed{0}{L_1}$.
There exists
$\gminiy_0(L)$ as in Corollary \ref{cor;16.4.5.10}.
By Corollary \ref{cor;16.4.5.10},
if $\epsilon(L)$ is sufficiently small,
we obtain the following:
\[
 A_{j}^{(m)}(L',(b_1,b_2))=\\
 \sum_{d\geq 1}
 \sum_{c_1+c_2=d}
 A_j^{(m+d)}\bigl(L-d\kappahat,
 (c_1,c_2)+(b_1,b_2)
 \bigr)
\cdot
 B_{j}(L'-L+d\kappahat,\gminiy_0(L),(c_1,c_2)).
\]
If $j_0(L)$ is sufficiently large,
$A_j^{(m+d)}\bigl(L-d\kappahat,
 (c_1,c_2)+(b_1,b_2)
 \bigr)=0$
holds for any $j\geq j_0(L)$.
Then, we obtain
the claim of Lemma \ref{lem;16.4.6.21}.
\hfill\qed

\vspace{.1in}
Then, by using the compactness 
of the interval $[0,L_1]$,
we obtain the claim of Theorem \ref{thm;16.4.6.11}.
\hfill\qed

\subsubsection{Proof of Theorem \ref{thm;16.4.6.12}}

We assume $f\neq 0$,
and we shall deduce a contradiction.
We set $\beta_0:=\min\nbigs>0$.
Let $L_0$ be the minimum of the non-empty set
$\bigl\{
 i/e+(j+k)\beta_0\,\big|\,
 f_{i,j,k}\neq 0
 \bigr\}$.
We set
$T(L_0):=\bigl\{(i,j,k)\in\seisuu_{\geq 0}^3
 \,\big|\,f_{i,j,k}\neq 0,\,\,\,i/e+(j+k)\beta_0=L_0
 \bigr\}$.
We set 
$\ell_0:=\max\{j+k\,|\,(i,j,k)\in T(L_0)\}$
and $i_0:=e(L_0-\ell_0\beta_0)$.

\begin{lem}
$A_{\infty}^{(\ell_0)}(L_0-\ell_0\beta_0,(\ell_1,\ell_2))
=\ell_1!\ell_2!f_{i_0,\ell_1,\ell_2}$
holds for any $(\ell_1,\ell_2)$
with $\ell_1+\ell_2=\ell_0$.
\end{lem}
\pf
By definition, 
the following equality holds:
\begin{multline}
 \label{eq;16.3.8.30}
 A_{\infty}^{(\ell_0)}(L_0-\ell_0\beta_0,(\ell_1,\ell_2))
=\ell_1!\ell_2!f_{i_0,\ell_1,\ell_2}
 \\
+\sum_{(i,\veczeta)\in T_+(\nbigs,L_0-\ell_0\beta_0)}
 \sum_{j+k=|\veczeta|+\ell_0}
 f_{i,j,k}(\theta)
 \frac{j!k!}{(j-\ell_1)!(k-\ell_2)!}
 \prod_{q=1}^{j-\ell_1}R_{\infty}(\zeta_q)
 \prod_{q=j-\ell_1+1}^{|\veczeta|}
 I_{\infty}(\zeta_q).
\end{multline}
Let us observe that 
the second term in the right hand side of (\ref{eq;16.3.8.30})
does not appear.
Suppose that the summand 
for $(i,\zeta_1,\ldots,\zeta_M,j,k)$ is non-trivial,
and we shall derive a contradiction.
There exists $M>0$ such that $f_{i,j,k}\neq 0$
for some $(i,j,k)$ satisfying the relations
$i/e+\sum_{q=1}^M\zeta_q+\ell_0\beta_0=L_0$
and $j+k=M+\ell_0$.
Because $\beta_0\leq \zeta_q$,
we obtain
\begin{equation}
\label{eq;16.6.2.2}
 i/e+(j+k)\beta_0= i/e+M\beta_0+\ell_0\beta_0
\leq
 i/e+\sum_{q=1}^M\zeta_q+\ell_0\beta_0=L_0.
\end{equation}
By our choice of $L_0$,
the inequality in (\ref{eq;16.6.2.2}) should be an equality,
and 
we obtain $\zeta_1=\cdots=\zeta_q=\beta_0$.
By our choice of $\ell_0$,
we obtain $M=0$,
which contradicts $M>0$.
Hence, the second term
in the right hand side of (\ref{eq;16.3.8.30}) is $0$.
\hfill\qed

\vspace{.1in}

We obtain
$A_{\infty}^{(m)}(L,(a,b))=0$ for any $L$
and $a+b=m$
by the assumption and Lemma \ref{lem;16.6.1.40}.
In particular,
we obtain
$A_{\infty}^{(\ell_0)}(L_0-\ell_0\beta_0,(\ell_1,\ell_2))=0$
for any $(\ell_1,\ell_2)$
such that $\ell_1+\ell_2=\ell_0$,
i.e.,
$f_{i_0,\ell_1,\ell_2}=0$
for any $(\ell_1,\ell_2)$
such that $\ell_1+\ell_2=\ell_0$.
But, it contradicts our choice of
$\ell_0$.
Hence, we obtain Theorem \ref{thm;16.4.6.12}.
\hfill\qed

\subsection{Pull back of ramified real analytic functions}
\label{subsection;16.5.31.31}

\subsubsection{Setting}

Let $m_k$ $(k=1,2\ldots)$ 
be an increasing sequence of integers
such that $m_k\to\infty$.
Let $\kappa_k\in \frac{1}{m_k}\seisuu_{>0}$
be a strictly increasing sequence.
We set $\kappahat:=\lim_{k\to\infty}\kappa_k
 \in \real_{\geq 0}\cup\{\infty\}$.
We assume $\kappahat<\infty$,
which implies either one of the following.
\begin{description}
\item[Case 1]
 There exists a prime $\gminip_0$
 such that
 $-\ord_{\gminip_0}(\kappa_i)\to\infty$.
\item[Case 2]
 There exists a sequence of primes
 $\gminip_i$ $(i=1,2,\ldots)$
 with $\lim_{i\to\infty}\gminip_i=\infty$
 such that 
 $\ord_{\gminip_i}(\kappa_i)<0$
 and
 $\ord_{\gminip_i}(\kappa_j)\geq 0$ $(\forall j<i)$.
\end{description}

Let $\nbigp^{(k)}(t,\gminia)
=\sum_{\gminiy>0}
 \nbigp^{(k)}_{\gminiy}(\gminia)\,t^{\gminiy}
\in \cnum[\![t^{1/m_k},\gminia]\!]$
$(k=1,2,\ldots)$
be a family satisfying the following conditions:
\begin{itemize}
\item
$\nbigp^{(k)}(t,\gminia)$ are convergent power series
of $(t^{1/m_k},\gminia)$.
\item
$\nbigp^{(k)}_{\kappa_k}(0)\neq 0$.
\item
If $\gminiy<\kappa_k$,
then
$\nbigp^{(k)}_{\gminiy}(\gminia)$
are independent of $\gminia$.
We denote $\nbigp^{(k)}_{\gminiy}(\gminia)$
just by $\nbigp^{(k)}_{\gminiy}$
for such $\gminiy$.
\item
$\nbigp^{(k)}_{\gminiy}=\nbigp^{(k')}_{\gminiy}$ holds
for $\gminiy<\min\{\kappa_k,\kappa_{k'}\}$,
and 
if $k<k'$ then
$\nbigp^{(k)}_{\kappa_k}(0)
=\nbigp^{(k')}_{\kappa_k}$ holds.
\item
For any $\gminiy<\kappa_k$ such that
$\nbigp^{(k)}_{\gminiy}\neq 0$,
 $\ord_{\gminip_0}(\gminiy)>\ord_{\gminip_0}(\kappa_k)$
holds in {\bf (Case 1)},
or $\ord_{\gminip_k}(\gminiy)\geq 0$ holds
in {\bf (Case 2)}.
\end{itemize}
For any $\gminiy<\kappahat$,
we put 
$\nbigp_{\gminiy}:=\lim_{k\to\infty}
 \nbigp_{\gminiy}^{(k)}$.
Indeed,
$\nbigp_{\gminiy}=\nbigp_{\gminiy}^{(k)}$ holds
for any $\gminiy<\kappa_k$.
We set
\[
 \nbigs:=
\bigcup_{k}
 \bigl\{
 0<\gminiy<\kappahat\,\big|\,
 \nbigp_{\gminiy}^{(k)}\neq 0
 \bigr\}.
\]
Note that $\bigl\{
 0<\gminiy<\kappahat\,\big|\,
 \nbigp_{\gminiy}\neq 0
 \bigr\}
\subset\nbigs$,
and that $\nbigs$ is discrete in
$\{0\leq \gminiy<\kappahat\}$.

\vspace{.1in}

Let $\phi_1^{(k)}>0$ $(k=1,2,\ldots,)$
be a decreasing sequence.
We set $\nbigi^{(k)}:=
 \bigl\{
 -\phi_1^{(k)}<\phi<\phi_1^{(k)}\bigr\}$.
Let $r_k>0$ be a decreasing sequence of real numbers,
and let $\nbigu^{(k)}$ be a decreasing sequence of
neighbourhoods of $0$ in $\cnum_{\gminia}$
such that 
$\nbigp^{(k)}(t,\gminia)$
are absolutely convergent 
on $(t,\gminia)\in \closedopen{0}{r_k}\times\nbigu^{(k)}$.
We set
$J^{(k)}:=
 \closedopen{0}{r_k}\times\nbigi^{(k)}\times\nbigu^{(k)}
=\{(t,\phi,\gminia)\}$.
We define the maps
$\varphi^{(k)}:
 J^{(k)}\lrarr \real_{\geq 0}\times\real\times\real^2$
by
\[
 \sfa^{(k)}(t,\phi,\gminia)
=\Bigl(
 t,\phi,
 \Re\bigl(\nbigp^{(k)}(te^{\sqrt{-1}\phi},\gminia)\bigr),
 \Image\bigl(\nbigp^{(k)}(te^{\sqrt{-1}\phi},\gminia)\bigr)
 \Bigr).
\]

\subsubsection{Pull back of ramified real analytic functions}

Let $f$ be a ramified real analytic function on 
a neighbourhood $Y$ of $(0,0,0,0)$
in $\real_{\geq 0}\times\real\times\real^2$
with the following power series expansion,
where $e$ is a positive integer:
\[
 f=\sum_{i,j,k\geq 0}
 f_{i,j,k}(\theta)r^{i/e}x^jy^k.
\]
We assume that $f_{0,0,0}(\theta)=0$ for any $\theta$,
and that $f$ is not constantly $0$.

We set
$f^{(k)}(\phi,t,\gminia):=
 (\sfa^{(k)})^{\ast}(f)(\phi,t,\gminia)$,
which has the expansion
$f^{(k)}(\phi,t,\gminia)
=\sum f^{(k)}_{L}(\phi,\gminia)t^{L}$.
We may regard
$f^{(k)}_L(\phi,\gminia)$
as germs of real analytic functions
at $(0,0)\in\real\times\cnum$.

\vspace{.1in}

For any $\gminiy\in\nbigs$,
we set
$R(\gminiy):=\Re(\nbigp_{\gminiy}e^{\sqrt{-1}\gminiy\phi})$
and
$I(\gminiy):=\Image(\nbigp_{\gminiy}e^{\sqrt{-1}\gminiy\phi})$.
For any $L\geq 0$,
we obtain the following real analytic functions
of $\phi\in\real$:
\[
 A(L):=
 \sum_{(i,\veczeta)\in T_+(\nbigs,L)}
 \sum_{\substack{
 (c_1,c_2)\in\seisuu_{\geq 0}^2\\
 c_1+c_2=|\veczeta|}}
 f_{i,c_1,c_2}(\phi)
 \prod_{q=1}^{c_1}R(\zeta_q)
 \prod_{q=c_1+1}^{|\veczeta|}I(\zeta_q)
+\sum_{i\in T_0(\nbigs,L)}
 f_{i,0,0}(\phi).
\]
We may naturally regard
$R(\gminiy)$, $I(\gminiy)$ and $A(L)$
as germs of real analytic functions
at $(0,0)\in\real\times\cnum$.

\begin{lem}
\label{lem;16.6.1.50}
There exists $L_1$
such that $A(L)=0$ for any $L<L_1$
and $A(L_1)\neq 0$
as a germ of real analytic functions at $(0,0)\in\real\times\cnum$,
\end{lem}
\pf
Because $f$ is assumed to be non-constant,
the set
$B:=\{L\in\real_{>0}\,|\,A(L)\neq 0\}$ is non-empty
by Theorem \ref{thm;16.4.6.12}.
By Lemma \ref{lem;18.1.13.1} and 
Lemma \ref{lem;18.1.13.2},
there exists the minimum of $B$,
which is the desired $L_1$.
\hfill\qed

\begin{thm}
\label{thm;16.3.8.20}
Let $L_1$ be as in Lemma {\rm\ref{lem;16.6.1.50}}.
There exists $k_1$ 
such that the following holds
for any $k\geq k_1$:
\begin{itemize}
\item
$f^{(k)}_L=0$ for any $L<L_1$
 as germs of real analytic functions.
\item
 $f^{(k)}_{L_1}
=A(L_1)$ as a germ  of real analytic functions.
In particular, $f^{(k)}_{L_1}$ is independent of $\gminia$.
\end{itemize}
\end{thm}
\pf
For any subset $S\subset\real_{\geq 0}$,
we put 
$T_+(S):=
 \coprod_{m\geq 1}
 \seisuu_{\geq 0}
 \times
 S^m$
and $T_0(S):=\seisuu_{\geq 0}$.
For any $L\geq 0$,
we set
\[
 T_+(S,L):=\bigl\{
 (i,s_1,\ldots,s_m)\in T_0(S)\,\big|\,
 \frac{i}{e}+\sum_{j=1}^m s_j=L
 \bigr\},
\]
and 
$T_0(S,L):=\bigl\{
 i\in T_0(S)\,\big|\,i/e=L
 \bigr\}$.
We define
$T(S)=T_0(S)\sqcup T_+(S)$
and 
$T(S,L)=T_0(S,L)\sqcup T_+(S,L)$.
For any element $\vecs=(s_1,\ldots,s_m)\in S^m$,
the number $m$ is denoted by $|\vecs|$.

We set
$\nbigs(\nbigp^{(k)}):=
 \bigl\{
 \gminiy\in\real\,\big|\,
 \nbigp^{(k)}_{\gminiy}\neq 0
 \bigr\}$.
For any $\gminiy\in\nbigs(\nbigp^{(k)})$,
we set
$R^{(k)}(\gminiy):=
 \Re(\nbigp^{(k)}_{\gminiy}(\gminia)e^{\sqrt{-1}\gminiy\phi})$
and
$I^{(k)}(\gminiy):=
 \Image(\nbigp^{(k)}_{\gminiy}(\gminia)e^{\sqrt{-1}\gminiy\phi})$.
We obtain the following description:
\[
 f^{(k)}_L=
 \sum_{(i,\veczeta)\in T_+(\nbigs(\nbigp^{(k)}),L)}
 \sum_{\substack{(c_1,c_2)\in\seisuu_{\geq 0}^2\\
    c_1+c_2=|\veczeta|}}
 f_{i,c_1,c_2}
 \prod_{q=1}^{c_1}
 R^{(k)}(\zeta_q)
 \prod_{q=c_1+1}^{|\veczeta|}
 I^{(k)}(\zeta_q)
\\
+\sum_{i\in T_0(\nbigs(\nbigp^{(k)}),L)}
 f_{i,0,0}.
\]

We set
$\nbigs_{<\kappahat}(\nbigp^{(k)}):=
 \bigl\{
 \gminiy\in\real\,\big|\,
 \nbigp^{(k)}_{\gminiy}\neq 0,\,\,
 \gminiy<\kappahat
 \bigr\}$.
For $(m,\ell_1,\ell_2)\in\seisuu_{\geq 0}^3$
such that $m=\ell_1+\ell_2$,
we put
\begin{multline}
 A^{(m)}_{<\kappahat}(\nbigp^{(k)},L,(\ell_1,\ell_2)):=
 \\
 \sum_{(i,\veczeta)\in T_+(\nbigs_{<\kappahat}(\nbigp^{(k)}),L)}
 \sum_{\substack{(c_1,c_2)\in\seisuu_{\geq 0}^2\\
    c_1+c_2=|\veczeta|+m}}
 f_{i,c_1,c_2}(\phi)
 \frac{c_1!c_2!}{(c_1-\ell_1)!(c_2-\ell_2)!}
 \prod_{q=1}^{c_1-\ell_1}
 R^{(k)}(\zeta_q)
 \prod_{q=c_1-\ell_1+1}^{|\veczeta|}
 I^{(k)}(\zeta_q)
\\
+\sum_{i\in T_0(\nbigs_{<\kappahat}(\nbigp^{(k)}),L)}
 \ell_1!\ell_2!
 f_{i,\ell_1,\ell_2}(\phi).
\end{multline}

For any $\delta> 0$ and $(b_1,b_2)\in\seisuu_{\geq 0}^2$,
we set
\[
 \nbigu^{(k)}(\delta,b_1,b_2):=
 \Bigl\{
 (\zeta_1,\ldots,\zeta_{b_1+b_2})
 \in\nbigs(\nbigp^{(k)})^{b_1+b_2}
 \,\Big|\,
 \zeta_i\geq \kappahat,\,\,
 \sum\zeta_i=\delta
 \Bigr\}.
\]
If $\nbigu^{(k)}(\delta,b_1,b_2)\neq\emptyset$,
we define
\[
 B^{(k)}(\delta,b_1,b_2):=
 \sum_{\veczeta\in\nbigu^{(k)}(\delta,b_1,b_2)}
 \frac{1}{b_1!b_2!}
 \prod_{q=1}^{b_1}R^{(k)}(\zeta_q)
 \prod_{q=b_1+1}^{b_1+b_2}I^{(k)}(\zeta_q).
\]
If $\nbigu^{(k)}(\delta,b_1,b_2)=\emptyset$,
we set $B^{(k)}(\delta,b_1,b_2):=0$.
We also formally set
$B^{(k)}(0,0,0):=1$
and 
$B^{(k)}(0,b_1,b_2):=0$ for 
$(b_1,b_2)\in\seisuu_{\geq 0}^2\setminus\{(0,0)\}$.
Note that for each $k$
the set $\{\delta\,|\,\nbigu^{(k)}(\delta,b_1,b_2)\neq\emptyset\}$
is discrete.
For any $L\leq L_1$,
we obtain the following equality:
\[
 f^{(k)}_L
=\sum_{d\geq 0}
\sum_{\delta\geq 0}
 \sum_{\ell_1+\ell_2=d}
 A^{(d)}_{<\kappahat}
 \bigl(\nbigp^{(k)},L-\delta,(\ell_1,\ell_2)\bigr)
\cdot
 B^{(k)}(\delta,\ell_1,\ell_2).
\]
By Theorem \ref{thm;16.4.6.11},
there exists $k_0$ such that 
$A^{(d)}_{<\kappahat}
 \bigl(\nbigp^{(k)},L-\sum\zeta_j,(\ell_1,\ell_2)\bigr)=0$
for any $L\leq L_1$,
any $k\geq k_0$
and any $(\ell_1,\ell_2)\in\seisuu^2_{\geq 0}$
satisfying $\ell_1+\ell_2=d\geq 1$.
We obtain $f^{(k)}_L=0$
for any $L<L_1$
and any $k\geq k_0$.
We also obtain
$f^{(k)}_{L_1}=A(L_1)\neq 0$
for any $k\geq k_0$.
Hence, we obtain Theorem \ref{thm;16.3.8.20}.
\hfill\qed

\begin{cor}
\label{cor;16.5.6.20}
There exists a discrete subset
$Z\subset\real$ and 
$k_0\in\seisuu_{>0}$ such that 
the following holds
for any $k\geq k_0$:
\begin{itemize}
\item
The orders
$\ord_tf^{(k)}(\phi,t,\gminia)$
for $(\phi,\gminia)\in (\nbigi^{(k)}\setminus Z)\times\nbigu^{(k)}$
are constant
and independent of $k$.
Note that the functions 
$f^{(k)}(\phi,t,\gminia)$ are not constant.
\hfill\qed
\end{itemize}
\end{cor}

\subsection{Formal paths}
\label{subsection;16.5.31.32}

\subsubsection{Ringed spaces}

Let $Y$ be a real analytic manifold.
Let $\nbigo_Y^{\real}$ denote 
the sheaf of real analytic functions on $Y$.
For any open subset $U$ of $Y$,
let $\nbign'_Y(U)$ be the Novikov type ring over $\nbigo^{\real}_Y(U)$,
i.e.,
\[
 \nbign'_Y(U):=\Bigl\{
 \sum_{i=0}^{\infty}
 b_i t^{\gminiy_i}\,\Big|\,
 b_i\in\nbigo^{\real}_Y(U),\,\,
 \gminiy_i\in\rnum_{\geq 0},\,\,
 \gminiy_i<\gminiy_{i+1},\,\,
 \lim_{i\to\infty}\gminiy_i=\infty
 \Bigr\}.
\]
A presheaf on $Y$ is defined
by the correspondence $U\longmapsto \nbign'_Y(U)$.
Let $\nbign_Y$ denote the associated sheaf on $Y$.
If $U$ is connected,
then $\nbign_Y(U)=\nbign'_Y(U)$ holds.
We formally denote sections $s$ of $\nbign_Y$ on $U$
as $\sum a_{\gminiy}t^{\gminiy}$.
There exists the natural morphism of pre-sheaves
$\nbign_Y'\lrarr \nbigo^{\real}_Y$
induced by
$\sum_{\gminiy}a_{\gminiy}t^{\gminiy}
 \longmapsto
 a_0$.
It induces
a morphism of sheaves
$\nbign_Y\lrarr\nbigo^{\real}_Y$.

For any element 
$s=\sum a_{\gminiy}t^{\gminiy}$
in the stalk $\nbign_{Y,P}$,
we put
$\ord_{P,t}(s):=\min\{\gminiy\,|\,a_{\gminiy}\neq 0\}$
if $s\neq 0$,
and 
$\ord_{P,t}(s):=\infty$
if $s=0$.
Here, $a_{\gminiy}\neq 0$
means that $a_{\gminiy}$ is not constantly $0$
on a neighbourhood of $P$.

Let $P$ be any point of $Y$.
Let $\iota_P:\{P\}\lrarr Y$ denote the inclusion.
We obtain the natural morphism
$\iota_P^{-1}\nbign_Y\lrarr\nbign_P$.
For any section $s$ of $\nbign_Y$
on an open set $U\subset Y$,
and for any $P\in U$,
let $s_{P}$ and $s_{|P}$
denote the induced elements
of $\nbign_{Y,P}$ and $\nbign_P$,
respectively.
Note that $\ord_{P,t}(s)\leq \ord_{P,t}(s_{|P})$.

Let $s$ be an element of $\nbign_{Y,P}$.
There exists a small neighbourhood $Y_P$ of $P$
such that $s=\sum a_{\gminiy}t^{\gminiy}$,
where $a_{\gminiy}\in\nbigo^{\real}_{Y}(Y_P)$.
We say that $s$ is convergent
if the following holds.
\begin{itemize}
\item
There exists $m\in\seisuu_{>0}$
such that
$a_{\gminiy}=0$ unless $m\gminiy\in\seisuu$.
Moreover,
$s$ comes from a real analytic function on
$Y_P\times I_{t^{1/m}}$.
\end{itemize}

\subsubsection{Morphisms of ringed spaces}

Let $M$ be any real analytic manifold.
An $\nbign_Y$-path in $M$
is a morphism of
ringed spaces
$F:(Y,\nbign_Y)\lrarr (M,\nbigo^{\real}_M)$,
i.e.,
it consists of 
a real analytic map $F_0:Y\lrarr M$
and a homomorphism of sheaves of algebras
$F^{\ast}:
 F_0^{-1}\nbigo^{\real}_M\lrarr \nbign_Y$
such that the composite with $\nbign_Y\lrarr\nbigo^{\real}_Y$
is equal to the natural map 
$F_0^{\ast}:\nbigo^{\real}_M
 \lrarr\nbigo_Y^{\real}$.

Let $\iota:Y_1\subset Y$ be any real analytic submanifold.
Any $\nbign_Y$-path $F$
naturally induces the $\nbign_{Y_1}$-path
$F_{|Y_1}:=F\circ\iota$.

Let $F,F'$ be $\nbign_Y$-paths in $M$.
Let $P\in Y$.
Suppose that $F_0(P)=F'_0(P)$.
Then, 
we set
$\ord_{P,t}(F,F')
=\min_i\ord_t \bigl(
F^{\ast}(x_i)-(F')^{\ast}(x_i)
 \bigr)$,
where $(x_1,\ldots,x_n)$
is a real analytic coordinate system
around $F_0(P)$
such that $x_i(F_0(P))=0$.
It is easy to see that the number 
$\ord_{P,t}(F,F')$
is independent of the choice of
the coordinate system.
If $F_0(P)\neq F'_0(P)$,
then we set
$\ord_{P,t}(F,F')=0$.
Let $Y_1\subset Y$ be any real analytic manifold
such that $P\in Y_1$.
Note 
$\ord_{P,t}(F,F')
\leq
\ord_{P,t}(F_{|Y_1},F'_{|Y_1})$.

\vspace{.1in}
Suppose that $M$ is equipped with a global coordinate system
$(x_1,\ldots,x_n)$.
We also assume that $Y$ is connected.
Then, $\nbign_Y$-paths in $M$
are equivalent to 
tuples $(f_1,\ldots,f_n)\in\nbign_Y(Y)^n$
by the correspondence
$F\longmapsto f_i=F^{\ast}(x_i)$
$(i=1,\ldots,n)$.

\subsubsection{Real blowings up}
\label{subsection;16.6.3.1}

Let $C$ be a real analytic submanifold in $M$.
We say that 
an $\nbign_Y$-path $F$
factors through $C$
if it factors through $(C,\nbigo^{\real}_C)$
as a morphism of ringed spaces.

Suppose that 
an $\nbign_Y$-path $F$
does not factor through $C$.
Let $P\in Y$.
Let $(x_1,\ldots,x_n)$ be
a real analytic coordinate system
around $F_0(P)$
such that
$C=\{x_1=\cdots=x_{\ell}=0\}$.
If $F_0(P)\in C$,
we set
\[
 \ord_{P,t}(F,C):=
 \min_{1\leq i\leq \ell}
 \ord_{P,t}F^{\ast}(x_i).
\]
It is easy to see that 
$\ord_{P,t}(F,C)$
is independent of the choice of
the coordinate system.
Note that for any $P\in Y$
the $\nbign_P$-path
$\varphi_{|P}$ of $M$ is induced,
and $\ord_{P,t}(F_{|P},C)
 \geq
 \ord_{P,t}(F,C)$ holds.

\begin{lem}
\label{lem;16.4.10.1}
Let $F$ be an $\nbign_Y$-path in $M$
which does not factor through $C$.
Let $P$ be a point of $M$.
Suppose 
$\ord_{P,t}(F,C)
=\ord_{P,t}(F_{|P},C)$.
Let $\Bl_CM$ denote the real blowing up of $M$
along $C$.
Then, the following holds.
\begin{itemize}
\item
There exist a neighbourhood  $Y_1$ of $P$ in $Y$
and a $\nbign_{Y_1}$-path $\Ftilde$
in $\Bl_CM$
such that 
$p\circ\Ftilde
=F_{|Y_1}$,
where $p:\Bl_C(M)\lrarr M$ denotes the projection.
\item
If there exist another neighbourhood $Y_2$ of $P$
and $\nbign_{Y_2}$-path $\varphitilde'$
in $\Bl_{C}M$ 
such that 
$p\circ \Ftilde'
=F_{|Y_2}$,
then
$\Ftilde_{|Y_1\cap Y_2}
=\Ftilde'_{|Y_1\cap Y_2}$ holds.
\end{itemize}
\end{lem}
\pf
By shrinking $Y$ and $M$,
we may assume that $M$ is equipped with
a coordinate system $(x_1,\ldots,x_n)$
such that $C=\{x_1=\cdots=x_{\ell}=0\}$.
The $\nbign_Y$-path $F$
is expressed as
$(F_1,\ldots,F_n)
 \in \nbign_Y(Y)^n$.
Each $F_i$ has the expansion
$F_{i}=
 \sum_{\gminiy\in\real}F_{i,\gminiy}t^{\gminiy}$.
We may assume that 
$\ord_{Q,t}(F,C)=\ord_{Q,t}(F_{1})=:\gminiy_0$
for any $Q$,
and that $F_{1,\gminiy_0}$ is nowhere vanishing.

Recall
$\Bl_C(M)=\bigl\{
 \bigl(
 (x_1,\ldots,x_n),[y_1:\cdots:y_{\ell}]
 \bigr)\,\big|\,
 x_iy_j-x_jy_i=0,\,(1\leq i,j\leq\ell)
 \bigr\}$.
Let $U_k\subset \Bl_C(M)$ $(k=1,\ldots,\ell)$
denote the open subsets 
determined by the conditions $y_k\neq 0$.
Let $(u^{(k)}_1,\ldots,u^{(k)}_n)$ 
be the coordinate system of $U_k$
induced by
$u^{(k)}_j=x_j/x_k$ if $1\leq j\leq \ell$ and $j\neq k$,
by
$u^{(k)}_j=x_j$ if $j=k$ or $j>\ell$.
We define the $\nbign_Y$-path $\Ftilde$ in $U_1$
as follows 
with respect to $(u^{(1)}_1,\ldots,u^{(1)}_n)$:
\[
 \Bigl(
 F_1,
 F_{2}/F_1,\ldots,F_{\ell}/F_1,
 F_{\ell+1},\ldots,F_n
 \Bigr).
\]
Clearly,  $p\circ\Ftilde=F$ holds.

Let $\Ftilde'$ be an $\nbign_Y$-path
in $U_k$ such that
$p\circ\Ftilde'=F$.
It is expressed as
$(\Ftilde'_1,\ldots,\Ftilde'_n)$
with respect to
$(u^{(k)}_1,\ldots,u^{(k)}_n)$.
Note that
$F_1=\Ftilde'_1\cdot\Ftilde_k'$
and $F_k=\Ftilde_k'$.
By comparing the orders
of $F_1$ and $F_k$,
we obtain that
$\Ftilde'_1$ is invertible.
Then, it is easy to see that 
$\Ftilde'$ is 
equal to
$\Ftilde$ as
an $\nbign_Y$-path in $\Bl_C(M)$.
\hfill\qed

\vspace{.1in}

Let $M$, $F$ and $P\in Y$ be as in 
Lemma {\rm\ref{lem;16.4.10.1}}.
Let $F'$ be another $\nbign_Y$-path in $M$
with the following property.
\begin{itemize}
\item
 The underlying map
 $F_0,F'_0:Y\lrarr M$ are the same
 on a neighbourhood of $P$.
\item
$\ord_{P,t}(F',F)>\ord_{P,t}(F,C)$.
\end{itemize}
Then, we can easily check the following
by a direct computation.
\begin{lem}
$F'$ does not factor through $C$,
and 
$\ord_{P,t}(F',C)=\ord_{P,t}(F'_{|P},C)$.
\hfill\qed
\end{lem}

If $Y_1$ is a small neighbourhood  of $P$ in $Y$,
then there exist
the $\nbign_{Y_1}$-paths
$\Ftilde$ and $\Ftilde'$ in $\Bl_CM$
such that 
$p\circ\Ftilde=F_{|Y_1}$
and 
$p\circ\Ftilde'=F'_{|Y_1}$.
We can check the following by a direct computation.
\begin{lem}
The underlying real analytic maps
$\Ftilde_0,\Ftilde'_0:
 Y_1\lrarr M$
are the same.
Moreover,
\[
 \ord_{Q,t}(\Ftilde,\Ftilde')
=\ord_{Q,t}(F,F')-\ord_{Q,t}(F,C)
\]
holds at any $Q\in Y_1$.
\hfill\qed
\end{lem}

\subsubsection{A complement on convergence}

Set $I:=\{0\leq \theta\leq 1\}$ and $J:=\{0\leq r\leq 1\}$.
Let 
$\iota:M\lrarr X$
and
$p:X\lrarr I\times J$
be real analytic maps.
We set $f:=p\circ \iota$.
We assume that
(i) $\dim M=2$,
(ii) $M$ is connected,
(iii)
$I\times\{0\}\subsetneq f(M)$,
(iv)
$\dim f^{-1}(I\times\{0\})=1$.

Let $F^{(0)}$
be the  $\nbign_I$-path in $I\times J$
induced by the identity map
$I\simeq I\times\{0\}$
and the correspondence $r\longmapsto t$.

\begin{lem}
\label{lem;16.8.28.20}
Let $F$ be an $\nbign_I$-path in $X$
which factors through $M$.
We assume the following.
\begin{itemize}
\item
$p\circ F$ is equal to
the $\nbign_I$-path $F^{(0)}$.
\end{itemize}
Then, there exists
a $0$-dimensional subset $Z\subset I$
such that 
$F$ is convergent at any $P\in I\setminus Z$
in the following sense.
\begin{itemize}
\item
 Let $(X_{P},x_1,\ldots,x_n)$ 
 be a real analytic coordinate neighbourhood of $X$
 around $F_0(P)$.
 Let  $(F_1,\ldots,F_n)$
 be the description of $F$
 with respect to the coordinate system.
 Then, $F_i$ are convergent. 
\end{itemize}
\end{lem}
\pf
There exists a $0$-dimensional real analytic subset
$Z_0\subset f^{-1}(I\times\{0\})$ such that
the following holds.
\begin{itemize}
\item
$f^{-1}(I\times\{0\})\setminus Z_0$
is a smooth submanifold in $M$.
\item
 Let $P_1$ be any point of 
 $f^{-1}(I\times\{0\})\setminus Z_0$.
 Let $(U;x,y)$ be any real analytic coordinate neighbourhood 
 around $P_1$
 such that $x^{-1}(0)=U\cap f^{-1}(I\times\{0\})$.
 Then, for an appropriate choice of $y$,
 $U\lrarr I\times J$ is expressed as
 $(y,\sum_{j=e}^{\infty} a_j(y)x^j)$
 such that
 $a_e$ is nowhere vanishing on $\{x=0\}$.
\end{itemize}
Then, the claim of the lemma is clear.
\hfill\qed

\section{Meromorphic flat bundles}

We recall some basic results
on meromorphic flat bundles
as a preliminary for the study in Part \ref{part;18.11.16.30}.

\subsection{One dimensional case}
\label{subsection;18.11.26.1}

Let $X:=\{z\in\cnum\,|\,|z|<1\}$ and $H:=\{0\}$.
For any positive integer $e$,
set $X^{(e)}:=\{\zeta_e\in\cnum\,|\,|\zeta_e|<1\}$
and $H^{(e)}:=\{0\}$.
We define the ramified covering
$\varphi_e:(X^{(e)},H^{(e)})\lrarr (X,H)$
by $\varphi_e(\zeta_e)=\zeta_e^e$.
Let $(V,\nabla)$ be a meromorphic flat bundle
on $(X,H)$.
According to the Hukuhara-Levelt-Turrittin theorem,
there exist a positive integer $e$,
a set
$\Irr(V,\nabla)\subset \nbigo_{X^{(e)}}(\ast H^{(e)})/\nbigo_{X^{(e)}}$
and a decomposition
\[
 \varphi_e^{\ast}
 (V,\nabla)_{|\widehat{H^{(e)}}}
=\bigoplus_{\gminia\in\Irr(V,\nabla)}
 (\Vhat_{\gminia},\nabla_{\gminia})
\]
such that
$\nabla_{\gminia}-d\gminiatilde\id_{\Vhat_{\gminia}}$
are regular singular,
where $\gminiatilde\in\nbigo_{X^{(e)}}(\ast H^{(e)})$
are representatives of $\gminia$.
We assume that $\rank\Vhat_{\gminia}\neq 0$
for any $\gminia\in\Irr(V,\nabla)$,
then $\Irr(V,\nabla)$ is uniquely determined.
The set $\Irr(V,\nabla)$ is invariant
under the action of the Galois group
of the ramified covering $X^{(e)}\lrarr X$.
In this paper, it is called the set of irregular values of $(V,\nabla)$.
The ranks $\rank \Vhat_{\gminia}$ are called
the multiplicity of $\gminia$.
They induce a map
$\rank:\Irr(V,\nabla)\lrarr \seisuu_{>0}$,
called the multiplicity function.
Set $\Rhat:=\varinjlim\cnum[\![z^{1/e}]\!]$
and $\Khat:=\varinjlim\cnum(\!(z^{1/e})\!)$.
We may naturally regard
$\Irr(V,\nabla)$
as a subset in $\Khat/\Rhat$.
If $\Irr(V,\nabla)$ is contained in
$\cnum(\!(z)\!)/\cnum[\![z]\!]$,
then $(V,\nabla)$ is called unramified.
We shall often use the natural bijection
$\zeta_e^{-1}\cnum[\zeta_e^{-1}]
\simeq
 \cnum(\!(\zeta_e)\!)/\cnum[\![\zeta_e]\!]$.

\vspace{.1in}

Let $\varpi:\Xtilde(H)\lrarr X$ be the oriented real blowing up.
A $C^{\infty}$-function $f$ on
an open subset $\nbigu\subset\Xtilde(H)$ is called holomorphic
if $f_{|\nbigu\setminus \varpi^{-1}(H)}$ is holomorphic.
Let $\nbigo_{\Xtilde(H)}$ denote the sheaf of holomorphic
functions on $\Xtilde(H)$.
For any $\nbigo_X$-module $\nbigm$,
let $\varpi^{\ast}\nbigm:=
 \varpi^{-1}\nbigm\otimes_{\varpi^{-1}\nbigo_X}
 \nbigo_{\Xtilde(H)}$.
For simplicity, we assume that
$(V,\nabla)$ is unramified,
i.e.,
$\Irr(V,\nabla)\subset z^{-1}\cnum[z^{-1}]$.
According to a classical theory,
for any $Q\in\varpi^{-1}(H)$,
there exists a neighbourhood $\nbigu_Q$ of $Q$
in $\Xtilde(H)$
and a decomposition
\begin{equation}
 \label{eq;16.2.7.1}
 \varpi^{\ast}(V,\nabla)_{|\nbigu_Q}
\simeq
\bigoplus_{\gminia\in\Irr(V,\nabla)}
 (V_{\gminia,\nbigu_Q},\nabla_{\gminia,\nbigu_Q})
\end{equation}
such that
$(V_{\gminia,\nbigu_Q},
 \nabla_{\gminia,\nbigu_Q})_{|
 \widehat{\varpi^{-1}(H)\cap\nbigu_Q}}
=\varpi^{\ast}(\Vhat_{\gminia},\nabla_{\gminia})$.
Such a decomposition is not uniquely determined.
We set
$\nbigf_{\gminia}^{\nbigu_Q}(\varpi^{\ast}V)
=\bigoplus_{\gminib\leq_{Q}\gminia}
 V_{\gminib,\nbigu_Q}$.
Here, we define the partial order $\leq_Q$
on $\Irr(V,\nabla)$ by setting
$\gminia\leq_{Q}\gminib$
if $-\Re(\gminia)\leq -\Re(\gminib)$
on $\nbigu_Q\setminus\varpi^{-1}(H)$,
which is independent of 
the choice of any sufficiently small $\nbigu_Q$.
Thus, we obtain the well defined filtration
$\nbigf^{\nbigu_Q}$ 
indexed by
$(\Irr(\nabla),\leq_Q)$.
Such a filtration $\nbigf^{\nbigu_Q}$
is called Stokes filtration.

\vspace{.1in}

Let $L$ be the local system on $X\setminus H$
obtained as the sheaf of flat sections of $(V,\nabla)$.
Let $\Ltilde$ be the local system on $\Xtilde(H)$
induced by $L$.
Let $\Ltilde_Q$ denote the stalk of $\Ltilde$ at $Q$.
Suppose that $(V,\nabla)$ is unramified.
Then, we obtain the filtration $\nbigf^Q$ of $\Ltilde_Q$
indexed by 
the partially ordered set $(\Irr(V,\nabla),\leq_Q)$
induced by $\nbigf^{\nbigu_Q}$.
We choose a frame $\vecv$ of $V$,
and let $h_V$ be a metric determined by
$h_V(v_i,v_j)=\delta_{i,j}$.
Then, the filtration $\nbigf^Q$
is characterized by the following condition.
\begin{itemize}
\item
 $s\in\nbigf^{Q}_{\gminic}(\Ltilde_Q)$
 if and only if
 $|s\cdot\exp(\gminic)|_{h_V}=O(|z|^{-C})$ 
for some $C>0$.
\end{itemize}
Set $\Gr^{\nbigf^Q}_{\gminic}(\Ltilde_Q):=
 \nbigf^{Q}_{\gminic}(\Ltilde_Q)\big/
 \nbigf^Q_{<\gminic}(\Ltilde_Q)$.
We may regard $\nbigf^Q$
as a filtration indexed by
$(z^{-1}\cnum[z^{-1}],\leq_Q)$.
There exists a decomposition
$\Ltilde_Q=\bigoplus G_{Q,\gminia}$
such that
$\nbigf^Q_{\gminia}(\Ltilde_Q)
=\bigoplus_{\gminib\leq_Q\gminia}G_{Q,\gminib}$.

If $Q'$ is sufficiently close to $Q$,
then the identity on $\Irr(V,\nabla)$
induces a morphism of partially ordered sets
$(\Irr(V,\nabla),\leq_Q)
\lrarr
 (\Irr(V,\nabla),\leq_{Q'})$.
There exists the isomorphism
$\Ltilde_{Q}\simeq \Ltilde_{Q'}$
induced by the parallel transport.
Then, any splitting of $\nbigf^Q$
induces a splitting of $\nbigf^{Q'}$.
In this sense,
the family 
$\{\nbigf^Q\,|\,Q\in \varpi^{-1}(H)\}$
satisfies a compatibility condition.
It is called the Stokes structure of $(V,\nabla)$.

For any $\gminia,\gminib\in\Irr(V,\nabla)$
with $\gminia\neq\gminib$,
we obtain 
a $C^{\infty}$-function
$F_{\gminia,\gminib}:=
 |z|^{-\ord(\gminia-\gminib)}(\gminia-\gminib)$
on $\Xtilde(H)$.
Let 
$Z(\gminia,\gminib):=
 \bigl\{
 Q\in\varpi^{-1}(H)\,\big|\,
 F_{\gminia,\gminib}(Q)=0
 \bigr\}$.
The following lemma is well known.
\begin{lem}
\label{lem;16.2.7.2}
Let $I$ be an open interval in $\varpi^{-1}(H)$.
Suppose that
$|Z(\gminia,\gminib)\cap I|\leq 1$
for any pair $\gminia,\gminib\in\Irr(V,\nabla)$
with $\gminia\neq\gminib$.
Then, there exists a decomposition
$L_{|I}=\bigoplus_{\gminia\in\Irr(V,\nabla)}
 G_{I,\gminia}$
which is compatible with the filtrations
$\nbigf^Q$ $(Q\in I)$.
We set $\nbigf^I_{\gminia}:=
 \bigoplus_{\gminib\leq_I\gminia}G_{I,\gminib}$
which is independent of the choice of such a splitting.
\hfill\qed
\end{lem}

Let $I=\openopen{\theta_1}{\theta_2}$ 
be as in Lemma {\rm\ref{lem;16.2.7.2}}.
Let $\Ltilde_I$ denote the restriction of $\Ltilde$ to $I$.
We choose $\epsilon>0$
such that
$\openopen{\theta_1}{\theta_1+\epsilon}$
and 
$\openopen{\theta_2-\epsilon}{\theta_2}$
do not intersect with any of $Z(\gminia,\gminib)$.
Let $\nbigh(\Ltilde_{|I})$ denote the space of global sections.
Set $Q_1=\theta_1+\epsilon$
and $Q_2=\theta_2-\epsilon$.
We naturally identify 
$\nbigh(\Ltilde_I)$
with $\Ltilde_{Q_i}$ $(i=1,2)$.

\begin{lem}
\label{lem;16.2.7.10}
Let $\nbigh(\Ltilde_I)=
 \bigoplus_{\gminia\in\Irr(V,\nabla)}
 \nbigg^{(1)}_{\gminib}$
be a decomposition such that
$\nbigf^{Q_i}_{\gminia}\nbigh(\Ltilde_I)
=\bigoplus_{\gminib\leq_{Q_i}\gminia}
 \nbigg^{(1)}_{\gminib}$
for $i=1,2$.
Then, 
$\nbigf^{Q}_{\gminia}\nbigh(\Ltilde_I)
=\bigoplus_{\gminib\leq_{Q}\gminia}
 \nbigg^{(1)}_{\gminib}$
holds for any $Q\in I$.
In other words,
the filtrations
$\nbigf^Q$ $(Q\in I)$
are determined by
$\nbigf^{Q_i}$ $(i=1,2)$.
\end{lem}
\pf
There exists a decomposition
$\bigoplus\nbigg_{\gminia}$
as in Lemma \ref{lem;16.2.7.2}.
Set 
$T(\gminia,Q_1,Q_2):=\bigl\{
 \gminib\,\big|\,
 \gminib\leq_{Q_i}\gminia\,(i=1,2)
 \bigr\}$.
Then,
$\nbigg^{(1)}_{\gminia}$
is contained in
$\bigoplus_{\gminib\in T(\gminia,Q_1,Q_2)}
 \nbigg_{\gminib}$.
Note that
for any $\gminib\in T(\gminia,Q_1,Q_2)$
and $Q\in I$,
$\gminib\leq_Q\gminia$ holds.
Hence, 
for any $\gminia\leq_Q\gminic$,
we obtain
$\nbigg^{(1)}_{\gminia}
\subset
\nbigf^{Q}_{\gminic}$.
Then, we can easily deduce the claim of the lemma.
\hfill\qed

\subsection{Higher dimensional case}

Let $X$ be a complex manifold
with a simply normal crossing hypersurface $H$.
Let $P$ be any point of $H$.
Let $(X_P,z_1,\ldots,z_n)$ be a holomorphic coordinate neighbourhood
around $P$ such that
$H_P:=H\cap X_P=\bigcup_{i=1}^{\ell}\{z_i=0\}$.
For any $\vecm\in\seisuu^{\ell}$,
set $\vecz^{\vecm}:=\prod_{i=1}^{\ell}z_i^{m_i}$.

Let $f\in\nbigo_X(\ast H)_P$.
Suppose that there exists an element
$\vecm\in\seisuu_{\leq 0}^{\ell}\setminus\{(0,\ldots,0)\}$
such that 
the function $\vecz^{-\vecm}f$ is holomorphic at $P$
and $(\vecz^{-\vecm}f)(P)\neq 0$.
Then, we say that $f$ has order $\vecm$,
and $\vecm$ is denoted by $\ord(f)$.
We also define $\ord(f):=(0,\ldots,0)$
for any $f\in\nbigo_{X,P}$.
Otherwise,
we say that $\ord(f)$ does not exist.

Let $f\in\nbigo_X(\ast H)_P/\nbigo_{X,P}$.
Let $\ftilde\in\nbigo_X(\ast H)_P$ 
be any lift of $f$.
We say that $f$ has order $\vecm$
if $\ftilde$ has order $\vecm$.
Note that the condition is independent of
the choice of a lift $\ftilde$.

Let $\leq_{\seisuu^{\ell}}$
denote the partial order on $\seisuu^{\ell}$
defined by
$\vecm\leq_{\seisuu^{\ell}}\vecm'
\stackrel{\rm}{\Longleftrightarrow}
 m_i\leq m_i'$.

A subset $\nbigi\subset \nbigo_X(\ast H)_P/\nbigo_{X,P}$
is called a good set of irregular values
if the following holds.
\begin{itemize}
\item
For any $\gminia\in\nbigi$,
$\ord(\gminia)$ exists.
The set 
$\{\ord(\gminia)\,|\,\gminia,\gminib\in\nbigi\}$
is totally ordered 
with respect to $\leq_{\seisuu^{\ell}}$.
\item
For any $\gminia,\gminib\in\nbigi$,
$\ord(\gminia-\gminib)$ exists.
The set
$\{\ord(\gminia-\gminib)\,|\,\gminia,\gminib\in\nbigi\}$
is totally ordered with respect to $\leq_{\seisuu^{\ell}}$.
\end{itemize}

A meromorphic flat connection means
a coherent reflexive $\nbigo_X(\ast H)$-module
with an integrable connection $(V,\nabla)$.
If $V$ is a locally free $\nbigo_X(\ast H)$-module
then $(V,\nabla)$ is called a meromorphic flat bundle.

Let $(V,\nabla)$ be a meromorphic flat bundle on $(X,H)$.
We say that $(V,\nabla)$ is unramifiedly good  at $P$
if the following holds.
\begin{itemize}
\item
 There exist a good set of irregular values
 $\Irr(V,\nabla,P)$ at $P$
 and a decomposition
\[
 (V,\nabla)_{|\widehat{P}}
=\bigoplus_{\gminia\in\Irr(V,\nabla,P)}
 (\Vhat_{\gminia},\nablahat_{\gminia}),
\]
such that
$\nablahat_{\gminia}-d\gminiatilde\id_{\Vhat_{\gminia}}$
are regular singular.
Here, $\gminiatilde\in\nbigo_{X}(\ast H)_P$
are lifts of $\gminia$.
\end{itemize}
We say that $(V,\nabla)$ is good at $P$
if the following holds.
\begin{itemize}
\item
 Let $(X_P,z_1,\ldots,z_n)$ be 
 any small holomorphic coordinate neighbourhood of $P$
 such that
 $H_P:=X_P\cap H=\bigcup_{i=1}^{\ell}\{z_i=0\}$.
 For any positive integer $e$,
 let $X_P^{(e)}$ denote a small neighbourhood
 of $(0,\ldots,0)$ in 
 $\cnum^n=\{(\zeta_{e,1},\ldots,\zeta_{e,n})\}$,
 and let
 $\varphi_e:X_P^{(e)}\lrarr X_P$
 denote the map defined by
\[
  \varphi_e(\zeta_{e,1},\ldots,\zeta_{e,n})
=(\zeta_{e,1}^e,\ldots,\zeta_{e,\ell}^e,\zeta_{e,\ell+1},\ldots,\zeta_{e,n}).
\]
Then, $\varphi_e^{\ast}(V,\nabla)_{|X_P}$
is unramifiedly good
for an appropriate positive integer $e$.
\end{itemize}
We say that $(V,\nabla)$ is (unramifiedly) good 
if it is (unramifiedly) good at any point of $P$.

A meromorphic flat bundle is not necessarily good.
The following fundamental theorem is due to Kedlaya
\cite{kedlaya,kedlaya2}.
(See \cite{mochi6,Mochizuki-wild} for the algebraic case.)

\begin{thm}
Let $(V,\nabla)$ be any meromorphic flat connection
on $(X,H)$.
For any $P\in H$,
there exist a small neighbourhood $X_P$
and a projective morphism of complex manifolds
$\psi_P:\check{X}_P\lrarr X_P$
such that 
(i) $\check{H}_P:=\psi_P^{-1}(H)$ is normal crossing,
(ii) $\check{X}_P\setminus \check{H}_P\simeq X_P\setminus H$,
(iii) $\psi_P^{\ast}(V,\nabla)$ is good.
\hfill\qed
\end{thm}

Let $(V,\nabla)$ be 
an unramifiedly good meromorphic flat connection
on $(X,H)$.
Let $\varpi:\Xtilde(H)\lrarr X$
be the oriented real blowing up of $X$ along $H$.
Let $L$ be the local system on $X\setminus H$
which is the sheaf of flat sections of $(V,\nabla)_{|X\setminus H}$.
Let $\Ltilde$ be the local system on $\Xtilde(H)$.
Let $Q\in \Xtilde(H)$.
As in the one dimensional case,
we obtain the filtration
$\nbigf^Q(\Ltilde_Q)$ 
indexed by
the ordered set
$\bigl(
 \Irr(V,\nabla,\varpi(Q)),\leq_Q
 \bigr)$.
The filtration is called the Stokes filtration.
There exists a splitting
$\Ltilde_Q=\bigoplus G_{Q,\gminia}$
such that
$\nbigf^Q_{\gminia}(\Ltilde_Q)
=\bigoplus_{\gminib\leq_Q\gminia}
 G_{Q,\gminib}$.
If $Q'$ is sufficiently close to $Q$,
there exists a natural map
$\bigl(
 \Irr(V,\nabla,\varpi(Q)),\leq_Q
 \bigr)
\lrarr
 \bigl(
 \Irr(V,\nabla,\varpi(Q')),\leq_{Q'}
 \bigr)$,
which is order preserving.
Any splitting of $\nbigf^Q$
induces a splitting of $\nbigf^{Q'}$.
The family of Stokes filtrations
$\bigl\{
 \nbigf^Q\,\big|\,
 Q\in\varpi^{-1}(H)
 \bigr\}$
is called the Stokes structure 
associated to $(V,\nabla)$.
The following is proved 
in \cite[Corollary 4.3.3]{Mochizuki-wild}.

\begin{thm}
Let $\vecnbigi$ be a good system of irregular values
on $(X,H)$.
The above construction induces
the functorial equivalence between 
unramifiedly good meromorphic flat bundles
over $\vecnbigi$
and 
local systems with Stokes structure
over $\vecnbigi$.
\hfill\qed
\end{thm}

\subsubsection{Extension} 

Set $\Delta:=\{z\in\cnum\,|\,|z|<1\}$.
Let $H$ be a complex manifold.
We set $X:=\Delta\times H$.
We naturally identify $H$ with $\{0\}\times H$.
Let $\vecnbigi$ be a good system of
irregular values on $(X,H)$.
Let $H_0$ be a complex submanifold of $H$.
We set $X_0:=\Delta\times H_0$.
By restricting $\vecnbigi$ to 
$(X_0,H_0)$,
we obtain a good system of
irregular values $\vecnbigi_0$ on $(X_0,H_0)$.
The following is
a special case of 
\cite[Corollary 4.4.4]{Mochizuki-wild}.
\begin{thm}
\label{thm;17.12.27.1}
Suppose that
the inclusion $H_0\lrarr H$ induces an isomorphism
of the fundamental groups.
Then, the restriction induces
a functorial bijective correspondence
between
local systems on $\Xtilde(H)$ with Stokes structure
over $\vecnbigi$
and local systems on $\Xtilde_0(H_0)$ with Stokes structure
over $\vecnbigi_0$.
\hfill\qed
\end{thm}

Let us observe a variant.
Set $X:=\Delta^2$, $H_i:=\{z_i=0\}$ $(i=1,2)$
and $H:=H_1\cup H_2$.
Let $\Hol(X)$ denote the space of holomorphic functions
on $X$.
Let $\Mero(X,H)$ be the space of meromorphic functions
on $(X,H)$.
Similarly,
$\Mero(X,H_i)$ be the space of meromorphic functions
on $(X,H_i)$.

Let $\nbigi$ be a good set of irregular values on $(X,H)$.
Note that $\nbigi\subset\Mero(X,H)/\Hol(X)$.
By exchanging $z_1$ and $z_2$ if necessary,
we may assume the following.
\begin{itemize}
\item
 $\nbigi\lrarr \Mero(X,H)/\Mero(X,H_2)$ is injective.
\end{itemize}

Let $\MFV(X,H,\vecnbigi)$
be the category of unramifiedly good 
meromorphic flat bundles on $(X,H,\vecnbigi)$.
Let $\vecnbigi_1$ be the good system of irregular values
on $(X\setminus H_2,H_1\setminus H_2)$
induced by $\vecnbigi$.
Let $\MFV(X\setminus H_2,H_1\setminus H_2,\vecnbigi_1)$
be the category of unramifiedly good 
meromorphic flat bundles on 
$(X\setminus H_2,H_1\setminus H_2,\vecnbigi_1)$.

\begin{prop}
\label{prop;16.2.7.21}
The restriction induces 
an equivalence of the categories
$\MFV(X,H,\vecnbigi)
\lrarr
\MFV(X\setminus H_2,H_1\setminus H_2,\vecnbigi_1)$.
\end{prop}
\pf
We set $X^{\circ}:=X\setminus H_2$
and $H_1^{\circ}:=H_1\setminus H_2$.
Let $\Ltilde$ be a local system on $\Xtilde(H)$.
Suppose that 
$\Ltilde_{|\Xtilde^{\circ}(H_1^{\circ})}$
is equipped with a Stokes structure
$\vecnbigf$ over $\vecnbigi_1$.
For any $\gminia,\gminib\in\nbigi$
with $\gminia\neq\gminib$,
we set
$F_{\gminia,\gminib}:=|\vecz^{-\vecm}|\Re(\gminia-\gminib)$
which naturally induces a $C^{\infty}$-function
on $\Xtilde(H)$.
Let $Z(\gminia,\gminib)=F_{\gminia,\gminib}^{-1}(0)$.
We use the polar coordinate system
$(r_1,\theta_1,r_2,\theta_2)$
determined by
$z_i=r_i\exp(\sqrt{-1}\theta_i)$.
Let $Q=(0,\theta_1^{(0)},0,\theta_2^{(0)})$
be any point of $\varpi^{-1}(0,0)$.
There exists an interval
$\closedclosed{\theta_1^{(1)}}{\theta_1^{(2)}}$
such that
$\theta_1^{(1)}<\theta_1^{(0)}<\theta_1^{(2)}$
and 
\[
 \bigl\{
 (0,\theta_1,0,\theta^{(0)}_2)\,\big|\,
 \theta_1^{(1)}\leq \theta_1\leq\theta_1^{(2)}
 \bigr\}
\cap
\Bigl(
 \bigcup_{\substack{
 \gminia,\gminib\in\nbigi\\
 \gminia\neq\gminib}}
 Z(\gminia,\gminib)
\Bigr)
\subset \{Q\}.
\]
There exists an interval
$\closedclosed{\theta_2^{(1)}}{\theta_2^{(2)}}$
such that
 $\theta_2^{(1)}<\theta_2^{(0)}<\theta_2^{(2)}$
and
\[
  \bigl\{
 (0,\theta^{(i)}_1,0,\theta_2)\,\big|\,
 \theta_2^{(1)}\leq \theta_2\leq\theta_2^{(2)}
 \bigr\}
\cap
\Bigl(
 \bigcup_{\substack{
 \gminia,\gminib\in\nbigi\\
 \gminia\neq\gminib}}
 Z(\gminia,\gminib)
\Bigr)
=\emptyset
\quad (i=1,2).
\]
There exist small positive numbers,
$\delta$ and $\epsilon$ such that 
the following set is empty:
 \[
  \bigl\{
 (0,\theta_1,r_2,\theta_2)\,\big|\,
 0\leq r_2\leq\delta,\,\,
 \theta_2^{(1)}\leq \theta_2\leq\theta_2^{(2)},
 \,\,\,
 \bigl(
 \theta_1^{(1)}\leq \theta_1\leq \theta_1^{(1)}+\epsilon,
 \mbox{\rm or }
  \theta_1^{(2)}-\epsilon\leq \theta_1\leq \theta_1^{(2)}
 \bigr)
 \bigr\}
\cap
\Bigl(
 \bigcup_{\substack{
 \gminia,\gminib\in\nbigi\\
 \gminia\neq\gminib}}
 Z(\gminia,\gminib)
\Bigr).
\]
For $0<r_2\leq \delta$
and 
$\theta^{(1)}_2\leq \theta_2\leq \theta^{(2)}_2$,
set $Q_1(r_2,\theta_2)=(0,\theta^{(1)}_1+\epsilon,r_2,\theta_2)$
and 
$Q_2(r_2,\theta_2)=
 (0,\theta^{(2)}_1-\epsilon,r_2,\theta_2)$.
Let $\nbigh$ denote the space of sections of
$\nbigltilde$
on 
\[
 W=
 \bigl\{
 (0,\theta_1,r_2,\theta_2)\,\big|\,
 0< r_2\leq\delta,\,\,
 \theta_1^{(1)}\leq \theta_1\leq\theta_1^{(2)},
 \theta_2^{(1)}\leq \theta_2\leq \theta_2^{(2)}
 \bigr\}.
\]
There exists the natural isomorphism
between $\nbigh$  and
$\nbigltilde_{Q}$ for any $Q\in W$,
with which we identify them.
Let 
$\nbigf^{Q_i(r_2,\theta_2)}$ 
$(i=1,2)$
denote the filtrations 
on $\nbigh$
induced by the Stokes filtrations
at $Q_i(r_2,\theta_2)$.
There exists a decomposition
\begin{equation}
 \label{eq;16.2.7.11}
 \nbigh=\bigoplus_{\gminia\in\nbigi}
 \nbigg_{\gminia}
\end{equation}
which induces a splitting 
to each filtration
$\nbigf^{Q_i(r_2,\theta_2)}$,
independently from $(r_2,\theta_2)$.
By Lemma \ref{lem;16.2.7.10},
the decomposition (\ref{eq;16.2.7.11})
induces a splitting of each filtration
$\nbigf^{Q'}$ $(Q'\in W)$.

Note that 
$\gminia\leq_Q\gminib$
if and only if
$\gminia\leq_{Q'}\gminib$
for any $Q'\in W$.
We define $\nbigf^Q$ 
by
\[
 \nbigf^Q_{\gminia}=
 \bigoplus_{\gminib\leq_Q\gminia}
 \nbigg_{\gminib}.
\]
It is independent of the choice of 
a decomposition (\ref{eq;16.2.7.11}).
It is easy to see that
the family $\{\nbigf^Q\,|\,Q\in\varpi^{-1}(0,0)\}$
induces a Stokes structure of $\nbigltilde$ along 
$\varpi^{-1}(0,0)$.
Hence, there exists a meromorphic flat bundle
$(V_1,\nabla)$ on a neighbourhood of $(0,0)$
whose Stokes structure along $\varpi^{-1}(0,0)$
is induced by $\{\nbigf^Q\,|\,Q\in\varpi^{-1}(0,0)\}$.
It is easy to see that 
the Stokes filtrations of $(V_1,\nabla)$
at $Q'\in\varpi^{-1}(H_1\setminus H_2)$
are the same as $\nbigf^{Q'}$ for $(V,\nabla)$.
Hence, we obtain a Stokes structure 
of $\nbigltilde$
whose restriction to $X\setminus H_2$
is equal to the prescribed one.
The claim of the proposition follows.
\hfill\qed

\section{Enhanced ind-sheaves}

We give some preliminaries from
the theory of enhanced ind-sheaves
for the study in Part \ref{part;18.11.16.30}.
We use the notation and the terminology in 
\cite{DAgnolo-Kashiwara1, DAgnolo-Kashiwara2,
Kashiwara-Schapira, Kashiwara-Schapira-ind-sheaves,
Kashiwara-Schapira2}.
In \S\ref{subsection;18.11.16.110},
we recall exactness of some basic functors.
In \S\ref{subsection;18.11.16.111},
we study filtrations associated to stably free enhanced ind-sheaves
in some special cases.
(See \cite[Definition 3.3.6]{DAgnolo-Kashiwara2}
for stably free enhanced ind-sheaves.)
In \S\ref{subsection;18.11.16.112},
we observe that prolongations of local systems
to $\real$-constructible enhanced ind-sheaves
are characterized by the prolongations of
their restriction to paths.
In \S\ref{subsection;18.11.16.113},
we give a sufficient condition
for the existence of a global filtration on a local system
which induces a prescribed enhanced ind-sheaf.

\subsection{Preliminary}
\label{subsection;18.11.16.110}

\subsubsection{Ind-sheaves}

Let $M$ be any good topological space.
Let $H\in \IC_{M}$.
For any relatively compact open subset $U$,
there exists the natural morphisms
$\cnum_U\otimes H\lrarr H$.
\begin{lem}
\label{lem;16.7.24.1}
$H$ is isomorphic to
$\varinjlim_U(\cnum_U\otimes H)$
in $\IC_M$.
\end{lem}
\pf
There exists a small filtrant category $I$
and a functor $\alpha:I\lrarr \Mod_c(\cnum_M)$
so that $H=\indlim\alpha$,
where $\Mod_c(\cnum_M)$
denote the category of $\cnum_M$-modules
whose supports are compact.
Recall that
$\cnum_U\otimes H
=\indlim\bigl(\cnum_U\otimes \alpha\bigr)$
by definition \cite{Kashiwara-Schapira-ind-sheaves}.
For any $G\in \Mod_c(\cnum_{M})$,
if $U$ is sufficiently large,
the induced morphisms
$\Hom(G,\cnum_U\otimes \alpha_i)
\lrarr
 \Hom(G,\alpha_i)$ are isomorphisms
for any $i\in I$.
Hence,
$\varinjlim_i\Hom(G,\cnum_U\otimes \alpha_i)
\lrarr
\varinjlim_i\Hom(G,\alpha_i)$ is an isomorphism.
We regard
$\cnum_U\otimes H$
and $H$ as objects
in $\Mod_c(\cnum_M)^{\wedge}$.
By the above,
for any $G\in\Mod_c(\cnum_M)$,
if $U$ is sufficient large,
$(\cnum_U\otimes H)(G)$
is naturally isomorphic to $H(G)$.
It implies that 
$\varinjlim_U(\cnum_U\otimes H)
\lrarr H$ is an isomorphism
in $\Mod_c(\cnum_M)^{\wedge}$.
Because the functor
$\IC_M\lrarr \Mod_c(\cnum_M)^{\wedge}$
commutes with inductive limits
\cite[Theorem 6.1.8]{Kashiwara-Schapira2},
we obtain
$\varinjlim_U(\cnum_U\otimes H)\simeq H$
in $\IC_M$.
\hfill\qed

\vspace{.1in}
Let $\nbigc$ be any locally closed relatively compact subset of $M$.
\begin{lem}
\label{lem;16.7.23.12}
The functor
$H\longmapsto \cnum_{\nbigc}\otimes H$
is an exact functor
on $\IC_M$.
\end{lem}
\pf
Let $0\lrarr H^{(1)}\lrarr H^{(2)}\lrarr H^{(3)}\lrarr 0$
be an exact sequence in $\IC_M$.
There exist a filtrant small category $I$
and an exact sequence of functors 
$0\lrarr \alpha^{(1)}
 \lrarr 
 \alpha^{(2)}
 \lrarr 
 \alpha^{(3)}
 \lrarr 0$
from $I$ to $\Mod_c(\cnum_M)$
such that 
$0\lrarr 
 \indlim \alpha^{(1)}
 \lrarr
 \indlim \alpha^{(2)}
 \lrarr 
 \indlim  \alpha^{(3)}
 \lrarr 0$
is $0\lrarr H^{(1)}\lrarr H^{(2)}\lrarr H^{(3)}\lrarr 0$.
(See \cite[Theorem 1.3.1]{Kashiwara-Schapira-ind-sheaves}
and the remark right after the theorem.)
By \cite[Proposition 2.3.6, Proposition 2.3.10]{Kashiwara-Schapira},
$0\lrarr 
 \cnum_{\nbigc}\otimes \alpha^{(1)}_i
 \lrarr
 \cnum_{\nbigc}\otimes \alpha^{(2)}_i
 \lrarr 
 \cnum_{\nbigc}\otimes \alpha^{(3)}_i\lrarr 0$
are exact for any $i\in I$,
we obtain that
$0\lrarr
 \cnum_{\nbigc}\otimes H^{(1)}
\lrarr
 \cnum_{\nbigc}\otimes H^{(2)}
\lrarr
 \cnum_{\nbigc}\otimes H^{(3)}
\lrarr 0$
is exact.
\hfill\qed

\vspace{.1in}
Let $\nbigz$ be a closed subset of $\nbigc$.
We set $\nbigc_0:=\nbigc\setminus \nbigz$.
There exists the naturally defined exact sequence
$0\lrarr \cnum_{\nbigc_0}\lrarr \cnum_{\nbigc}\lrarr \cnum_{\nbigz}\lrarr 0$
in $\Mod_c(\cnum_M)$.
For any $H\in \IC_M$,
we obtain the following morphisms
in $\IC_M$:
\begin{equation}
\label{eq;16.7.23.11}
0\lrarr
 \cnum_{\nbigc_0}\otimes H
\lrarr 
 \cnum_{\nbigc}\otimes H
\lrarr
 \cnum_{\nbigz}\otimes H
\lrarr 0.
\end{equation}

\begin{lem}
\label{lem;16.7.23.13}
The sequence {\rm(\ref{eq;16.7.23.11})}
is exact.
\end{lem}
\pf
There exists a functor $\alpha$
from a small filtrant category $I$
to $\Mod_c(\cnum_M)$
such that $H=\indlim \alpha$.
The sequences
$0\lrarr 
\cnum_{\nbigc_0}\otimes \alpha_i
\lrarr
\cnum_{\nbigc}\otimes \alpha_i
\lrarr
 \cnum_{\nbigz}\otimes \alpha_i
\lrarr 0$
are exact for any $i$.
Hence, we obtain that (\ref{eq;16.7.23.11})
is exact.

\hfill\qed

\subsubsection{Enhanced ind-sheaves}

Following \cite{DAgnolo-Kashiwara1},
we define the bordered spaces
$\real_{\infty}:=(\real,\realbar)$
and 
$\realbar:=(\realbar,\realbar)$.
Let $j:M\times \real_{\infty}\lrarr M\times\realbar$
be the inclusion.
Let $\pi:M\times\real_{\infty}\lrarr M$
and $\pibar:M\times\realbar\lrarr M$
denote the projections.

For any $K\in \Ecat^b(\IC_M)$,
and for any relatively compact open subset
$U$ of $M$,
there exists a natural isomorphism
$L^{\Ecat}\bigl(
 \pi^{-1}(\cnum_U)\otimes K
 \bigr)
\simeq
 \pi^{-1}(\cnum_U)
 \otimes L^{\Ecat}(K)$
in $\Dcat^b(\IC_{\real_{\infty}\times M})$,
which follows from
\cite[Lemma 4.3.1]{DAgnolo-Kashiwara1}
and the construction of $L^{\Ecat}$
in \cite[(4.4.1) and Notation 4.4.5]{DAgnolo-Kashiwara1}.
There also exists a natural isomorphism
$Rj_{!!}\bigl(
 \pi^{-1}(\cnum_U)
 \otimes
 L^{\Ecat}K
 \bigr)
\simeq
 \pibar^{-1}(\cnum_U)\otimes
Rj_{!!}L^{\Ecat}(K)$
in $\Dcat^b(\IC_{\realbar\times M})$,
which follows from 
\cite[Theorem 5.2.7]{Kashiwara-Schapira-ind-sheaves}.
Hence, we obtain an isomorphism
$Rj_{!!}L^{\Ecat}(\pi^{-1}(\cnum_U)\otimes K)
\simeq
 \pibar^{-1}(\cnum_U)\otimes
Rj_{!!}L^{\Ecat}(K)$
in $\Dcat^b(\IC_{\realbar\times M})$.

Let $\nbigh^0\Ecat^b(\IC_M)$
denote the heart of $\Ecat^b(\IC_M)$
with respect to the $t$-structure 
in \cite[\S4.6]{DAgnolo-Kashiwara1}.
Recall that 
an object $K\in \Ecat^b(\IC_M)$
is contained in $\nbigh^0\Ecat^b(\IC_M)$
if and only if
the object
$Rj_{!!}L^{\Ecat}K\in \Dcat^b(\IC_{\realbar\times M})$
is contained in
$\IC_{\realbar\times M}$.
For such $K$,
and for any relatively compact open subset
$U$ of $M$,
$\pi^{-1}(\cnum_U)\otimes K$
is contained in $\nbigh^0\Ecat^b(\IC_M)$
because 
$Rj_{!!}L^{\Ecat}(\pi^{-1}(\cnum_U)\otimes K)
\simeq
 \pibar^{-1}(\cnum_U)\otimes
Rj_{!!}L^{\Ecat}(K)$.
For inclusions $U_1\subset U_2$,
there exists the natural morphism
$\pi^{-1}(\cnum_{U_1})\otimes K
\lrarr
\pi^{-1}(\cnum_{U_2})\otimes K$.

\begin{lem}
\label{lem;16.7.14.1}
For any object $K$ of $\nbigh^0\Ecat^b(\IC_M)$,
there exists the natural isomorphism
$K\simeq
 \varinjlim_{U}
 \bigl(
 \pi^{-1}(\cnum_U)\otimes K
 \bigr)$
in $\nbigh^0\Ecat^b(\IC_M)$.
\end{lem}
\pf
By the definition of the $t$-structures,
there exists the fully faithful embedding
of the abelian categories
$\nbigh^0\Ecat^b(\IC_M)
\lrarr \IC_{\realbar\times M}$
induced by
$G\longmapsto
 Rj_{!!}L^{\Ecat}(G)$.
By Lemma \ref{lem;16.7.24.1},
there exist the natural isomorphisms
$Rj_{!!}L^{\Ecat}(K)\simeq
 \varinjlim_U
 \bigl(
 \pibar^{-1}(\cnum_U)
 \otimes
 Rj_{!!}L^{\Ecat}(K)
 \bigr)
\simeq
 \varinjlim_U
 Rj_{!!}L^{\Ecat}(\pi^{-1}(\cnum_U)\otimes K)$.
Hence, we obtain the natural isomorphism
$K\simeq
 \varinjlim_{U}\bigl(
 \pi^{-1}(\cnum_U)\otimes K
\bigr)$
in $\nbigh^0\Ecat^b(\IC_M)$.
\hfill\qed

\vspace{.1in}

Let $\nbigc$ be any relatively compact locally closed subset of $M$.
We obtain the following from 
Lemma \ref{lem;16.7.23.12}.
\begin{lem}
\label{lem;16.7.23.21}
The functor $K\longmapsto \pi^{-1}(\cnum_{\nbigc})\otimes K$
on $\Ecat^b(\cnum_M)$
induces an exact functor
on $\nbigh^0\Ecat^b(\cnum_M)$.
\hfill\qed
\end{lem}

Let $\nbigz$ be a closed subset in $\nbigc$.
Set $\nbigc_0:=\nbigc\setminus \nbigz$.
We obtain the following from 
Lemma \ref{lem;16.7.23.13}.
\begin{lem}
\label{lem;16.7.23.20}
$0\lrarr
 \pi^{-1}(\cnum_{\nbigc_0})\otimes K
\lrarr
 \pi^{-1}(\cnum_{\nbigc})\otimes K
\lrarr
 \pi^{-1}(\cnum_{\nbigz})\otimes K
\lrarr 0$
is an exact sequence 
in $\nbigh^0\Ecat^b(\IC_M)$.

\hfill\qed
\end{lem}

\subsection{Stably free enhanced ind-sheaves}
\label{subsection;18.11.16.111}

\subsubsection{Enhanced ind-sheaves associated to 
global subanalytic functions}

Suppose that $M$ is a subanalytic space.
(See \cite[Definition 2.3.1]{DAgnolo-Kashiwara1}
for subanalytic spaces.)
Let $\nbigc\subset M$ be any locally closed subanalytic subset.
Let $g$ be a continuous subanalytic function on $(\nbigc,M)$.
Set $Z(g):=\{(x,t)\in \nbigc\times\real\,|\,t\geq g(x)\}\subset M\times\real$.
Let $\cnum_{t\geq g}$ denote 
the $\real$-constructible sheaf $\cnum_{Z(g)}$
on $M\times\real$,
which can be extended to an $\real$-constructible sheaf
on $M\times\realbar$.
We obtain $\cnum_M^{\Ecat}\overset{+}{\otimes}\cnum_{t\geq g}$
in $\nbigh^0\Ecat^b(\IC_M)$.
(See \cite{DAgnolo-Kashiwara1}.)

\begin{lem}
Let $L$ be any local system on $M$.
Let $g_i$ $(i=1,2)$ be continuous subanalytic functions on $M$
such that $g_1(x)\leq g_2(x)$ for any $x\in M$.
The natural morphism
$\kappa:\pi^{-1}(L)\otimes
 \bigl(
 \cnum_M^{\Ecat}\overset{+}{\otimes}
 \cnum_{t\geq g_1}
 \bigr)
\lrarr 
 \pi^{-1}(L)\otimes
 \bigl(
\cnum_M^{\Ecat}\overset{+}{\otimes}
 \cnum_{t\geq g_2}
 \bigr)$
is an isomorphism
in $\Ecat^b(\IC_M)$.
\end{lem}
\pf
By Lemma \ref{lem;16.7.14.1},
it is enough to prove that
the natural morphism
\begin{equation}
\label{eq;16.7.14.11}
 \pi^{-1}(\cnum_{U}\otimes L)
\otimes
 \bigl(\cnum_M^{\Ecat}\overset{+}{\otimes}\cnum_{t\geq g_1}\bigr)
\lrarr
 \pi^{-1}(\cnum_U\otimes L)
\otimes
 \bigl(\cnum_M^{\Ecat}\overset{+}{\otimes}\cnum_{t\geq g_2}\bigr) 
\end{equation}
is an isomorphism
for any relatively compact open subanalytic subset $U\subset M$.
By \cite[Lemma 4.3.1]{DAgnolo-Kashiwara1},
there exist the following natural isomorphisms:
\[
 \pi^{-1}(\cnum_U\otimes L)
\otimes
 (\cnum_M^{\Ecat}\overset{+}{\otimes}\cnum_{t\geq g_i})
\simeq
 \cnum_M^{\Ecat}\overset{+}{\otimes}
\bigl(
\pi^{-1}(L_{|U})\otimes
\cnum_{t\geq g_{i|U}}
\bigr).
\]
By \cite[Remark 4.1.3]{DAgnolo-Kashiwara1},
there exist natural isomorphisms
$\cnum_{t\geq a}
 \overset{+}{\otimes}
 \cnum_{t\geq g_{i|U}}
\simeq
\cnum_{t\geq g_{i|U}+a}$
for any $a\in\real$.
By \cite[Proposition 4.7.9]{DAgnolo-Kashiwara1},
there exist the following natural isomorphisms:
\begin{multline}
\label{eq;16.7.14.10}
 \Hom_{\Ecat^b(\IC_M)}
 \Bigl(
 \cnum_M^{\Ecat}\overset{+}{\otimes}
\bigl(
\pi^{-1}(L_{|U})\otimes
\cnum_{t\geq g_{i|U}}
\bigr),
 \cnum_M^{\Ecat}\overset{+}{\otimes}
\bigl(
\pi^{-1}(L_{|U})\otimes
\cnum_{t\geq g_{j|U}}
\bigr)
 \Bigr)
 \\
\simeq
 \varinjlim_{a\to\infty}
 \Hom_{\Mod(\cnum_{M\times\real})}
 \bigl(
\pi^{-1}(L_{|U})\otimes
\cnum_{t\geq g_{i|U}},
\pi^{-1}(L_{|U})\otimes
\cnum_{t\geq g_{j|U}+a}
 \bigr).
\end{multline}
Let $b$ be any large number such that
$g_2(x)\leq g_1(x)+b$ for any $x\in U$.
There exists the natural morphism
$\pi^{-1}(L_{|U})\otimes
 \cnum_{t\geq g_{2|U}}
\lrarr
 \pi^{-1}(L_{|U})\otimes
 \cnum_{t\geq g_{1|U}+b}$.
By (\ref{eq;16.7.14.10}),
it induces a morphism
\[
 \kappa_1:
 \cnum_M^{\Ecat}
 \overset{+}{\otimes}
 \bigl(
 \pi^{-1}(L_{|U})\otimes
 \cnum_{t\geq g_{2|U}}
 \bigr)
\lrarr
 \cnum_M^{\Ecat}
 \overset{+}{\otimes}
 \bigl(
 \pi^{-1}(L_{|U})\otimes
 \cnum_{t\geq g_{1|U}}
 \bigr)
\]
in $\Ecat^b(\IC_M)$.
The composite
$\kappa_1\circ\kappa$
is induced by the natural morphism
$\pi^{-1}(L_{|U})\otimes
 \cnum_{t\geq g_{1|U}}
\lrarr
 \pi^{-1}(L_{|U})\otimes
 \cnum_{t\geq g_{1|U}+b}$.
Hence, we can easily observe that
$\kappa_1\circ\kappa$
is equal to the identity of 
$\cnum_M^{\Ecat}
 \overset{+}{\otimes}
 \bigl(
 \pi^{-1}(L_{|U})\otimes
 \cnum_{t\geq g_{1|U}}
 \bigr)$.
Similarly, we can observe that 
$\kappa\circ\kappa_1$ is the identity.
Hence, the morphism (\ref{eq;16.7.14.11})
is an isomorphism.
\hfill\qed

\vspace{.1in}
For any continuous subanalytic function $g$ on $M$,
we define the continuous subanalytic function $g_-$ on $M$
by $g_-(x):=\min\{0,g(x)\}$.
There exist the following natural morphisms
\[
\begin{CD}
 \pi^{-1}(L)\otimes
 \cnum_M^{\Ecat}
\simeq
 \pi^{-1}(L)\otimes
\bigl(
 \cnum_M^{\Ecat}
 \overset{+}{\otimes}
 \cnum_{t\geq 0}
\bigr)
@<{c_1}<<
 \pi^{-1}(L)\otimes
\bigl(
 \cnum_M^{\Ecat}
 \overset{+}{\otimes}
 \cnum_{t\geq g_-}
\bigr)
@>{c_2}>>
 \pi^{-1}(L)\otimes
\bigl(
 \cnum_M^{\Ecat}
 \overset{+}{\otimes}
 \cnum_{t\geq g}
\bigr).
\end{CD}
\]
\begin{cor}
\label{cor;16.7.14.20}
The morphisms $c_i$ are isomorphisms
in $\Ecat^b(\IC_M)$.
In this sense,
there exists a canonical isomorphism
$\pi^{-1}(L)\otimes
 \cnum_M^{\Ecat}
\simeq
\pi^{-1}(L)\otimes
 \bigl(
 \cnum_M^{\Ecat}
 \overset{+}{\otimes}
 \cnum_{t\geq g}
 \bigr)$
for any continuous subanalytic function $g$
on $M$.
\hfill\qed
\end{cor}

\subsubsection{Pre-orders on the sets of continuous subanalytic functions}

Let $\Sub(\nbigc,M)$ denote the set of continuous subanalytic functions on 
$(\nbigc,M)$.
We define the pre-order $\prec$ on $\Sub(\nbigc,M)$
by the condition;
$g\prec f$ holds
if $f-g$ is bounded from above on $U\cap \nbigc$
for any relatively compact subset $U$ in $M$.
We remark that this definition is suitable
for filtrations studied in \S\ref{subsection;16.1.25.20} below.
The pre-order induces the equivalence relation $\sim$
on $\Sub(\nbigc,M)$
defined by
$g\sim f$
$\stackrel{\rm def}{\Longleftrightarrow}$
$g\prec f$ and $f\prec g$.
Let $\Subbar(\nbigc,M)$ denote the quotient set
of $\Sub(\nbigc,M)$
by the above equivalence relation $\sim$.

Let $\nbigc'$ be a locally closed subanalytic subset
in a real analytic manifold $M'$.
Let $F:(\nbigc',M')\lrarr (\nbigc,M)$ be a real analytic map.
The map
$\Sub(\nbigc,M)\lrarr \Sub(\nbigc',M')$
is induced by the pull back
$f\longmapsto F^{\ast}(f)$.
If $f\prec g$ holds in $\Sub(\nbigc,M)$, then
$F^{\ast}f\prec F^{\ast}g$ holds
in $\Sub(\nbigc',M')$.
The map
$\Subbar(\nbigc,M)\lrarr \Subbar(\nbigc',M')$
is naturally induced.

\begin{lem}
\label{lem;16.7.14.40}
For $g_i\in\Sub(\nbigc,M)$ $(i=1,2)$,
$g_1\prec g_2$ holds
if and only if the following condition is satisfied.
\begin{itemize}
\item
Let $\gamma:(\II^{\circ},\II)\lrarr (\nbigc,M)$
be any real analytic map.
Then,
 $\gamma^{\ast}(g_1)
 \prec
 \gamma^{\ast}(g_2)$ holds
 in
 $\Sub(\II^{\circ},\II)$.
\end{itemize}
\end{lem}
\pf
The ``only if'' part is clear.
Let us prove the ``if'' part.
Suppose that
$g_1\not\prec g_2$.
Set $h:=g_2-g_1$.
There exists a compact subset $V\subset M$
such that 
$h$ is not bounded from above on $V\cap \nbigc$.
Let $\Gamma_{h,V\cap \nbigc}$
denote the graph of $h_{|V\cap \nbigc}$.
The closure of $\Gamma_{h,V\cap \nbigc}$
in $M\times\real_{\infty}$
contains a point in $M\times\{\infty\}$.
By the curve selection lemma,
we can choose a real analytic map
$\gamma:(\II^{\circ},\II)\lrarr (\nbigc,M)$
such that 
$\gamma^{\ast}(h)$ is not bounded from above,
i.e.,
$\gamma^{\ast}(g_1)\not\prec
 \gamma^{\ast}(g_2)$.
\hfill\qed

\subsubsection{Spaces of morphisms in some basic cases}
\label{subsection;16.7.14.30}

Let $i_{\nbigc}:\nbigc\lrarr M$ denote the inclusion.
Let $L_i$ $(i=1,2)$ be a local system on $\nbigc$.
We put $L_{i,M}:=i_{\nbigc!}L_i\in\Mod(\cnum_M)$.
Let $g_i$ $(i=1,2)$ be continuous subanalytic functions
on $(\nbigc,M)$.
By Corollary \ref{cor;16.7.14.20},
there exists a natural isomorphism
$\Ecat i_{\nbigc}^{-1}\Bigl(
 \cnum_M^{\Ecat}
\overset{+}{\otimes}
\bigl(
\pi^{-1}(L_{i,M})\otimes
 \cnum_{t\geq g_i}
\bigr)
 \Bigr)
\simeq
 \cnum_{\nbigc}^{\Ecat}\otimes \pi^{-1}(L_i)$
in $\Ecat^b(\IC_{\nbigc})$.
Recall that
the fully faithful embedding
$\Dcat^b(\cnum_{\nbigc})
\lrarr
 \Ecat^b(\IC_{\nbigc})$
is induced by
$G\longmapsto
 \pi^{-1}(G)\otimes\cnum_M^{\Ecat}$.
(See \cite[Proposition 5.1.1]{Kashiwara-Schapira-ind-sheaves}
and \cite[Proposition 4.7.15]{DAgnolo-Kashiwara1}.)
Hence,
the functor $\Ecat i_{\nbigc}^{-1}$
induces the following morphism:
\begin{equation}
 \label{eq;16.7.14.21}
 \Hom_{\Ecat^b(\IC_{M})}
 \Bigl(
\cnum_M^{\Ecat}\overset{+}{\otimes}
\bigl(
 \pi^{-1}(L_{1,M})\otimes\cnum_{t\geq g_1}
\bigr),
\cnum_M^{\Ecat}\overset{+}{\otimes}
\bigl(
 \pi^{-1}(L_{2,M})\otimes\cnum_{t\geq g_2}
\bigr)
 \Bigr)
\lrarr
\Hom_{\Mod(\cnum_{\nbigc})}(L_1,L_2).
\end{equation}

\begin{lem}
\label{lem;16.7.14.20}
The morphism {\rm(\ref{eq;16.7.14.21})}
is injective.
If $g_2\prec g_1$ in $\Sub(\nbigc,M)$,
the morphism is an isomorphism.
If $\nbigc$ is connected
and if $g_2\not\prec g_1$ in $\Sub(\nbigc,M)$,
then we obtain
\[
 \Hom_{\Ecat^b(\IC_{M})}
 \Bigl(
\cnum_M^{\Ecat}\overset{+}{\otimes}
 \bigl(\pi^{-1}(L_{1,M})\otimes\cnum_{t\geq g_1}\bigr),
\cnum_M^{\Ecat}\overset{+}{\otimes}
 \bigl(\pi^{-1}(L_{2,M})\otimes \cnum_{t\geq g_2} \bigr)
 \Bigr)=0.
\]
\end{lem}
\pf
Let $U$ be any relatively compact open subanalytic subset
in $M$.
There exist the following morphisms:
\begin{multline}
\label{eq;16.7.14.22}
\Hom_{\Ecat^b(\IC_{M})}
 \Bigl(
\cnum_M^{\Ecat}\overset{+}{\otimes}
\bigl(
 \pi^{-1}(\cnum_U\otimes L_{1,M})\otimes\cnum_{t\geq g_1}
\bigr),
\cnum_M^{\Ecat}\overset{+}{\otimes}
\bigl(
 \pi^{-1}(L_{2,M})\otimes\cnum_{t\geq g_2}
\bigr)
 \Bigr) \\
\simeq
 \varinjlim_{a\to\infty}
 \Hom_{\Mod(\cnum_{M\times\real})}
 \bigl(
 \pi^{-1}(\cnum_U\otimes L_{1,M})\otimes
 \cnum_{t\geq g_1-a},
 \pi^{-1}(L_{2,M})\otimes
 \cnum_{t\geq g_2}
 \bigr)
\\
\stackrel{\kappa_U}{\lrarr}
 \Hom_{\Mod(\cnum_{M\times\real})}
 \bigl(
 \pi^{-1}(\cnum_U\otimes L_{1,M}),
 \pi^{-1}(L_{2,M})\otimes
 \cnum_{t\geq g_2}
 \bigr)
\\
\simeq
  \Hom_{\Mod(\cnum_{M\times\real})}
 \bigl(
 \pi^{-1}(\cnum_U\otimes L_{1,M}),
 \pi^{-1}(L_{2,M})
 \bigr)
\simeq
 \Hom_{\cnum_{\nbigc}}
 \bigl(
 \cnum_{\nbigc\cap U}\otimes L_1,L_2
 \bigr).
\end{multline}
Here, $\kappa_U$ is injective.
Moreover,
if $g_1-g_2$ is bounded from above
on $\nbigc\cap U$,
then $\kappa_U$ is an isomorphism.
The morphism (\ref{eq;16.7.14.21})
is equal to the projective limit of
(\ref{eq;16.7.14.22}),
where $U$ runs through the set of relatively compact
open subsets of $M$.
Hence, we obtain that 
(\ref{eq;16.7.14.21}) is injective,
and that 
(\ref{eq;16.7.14.21}) is an isomorphism
in the case $g_2\prec g_1$.

Suppose that $\nbigc$ is connected 
and $g_2\not\prec g_1$.
There exists $P\in \nbigcbar$
around which $g_1-g_2$ is not bounded from above.
Let $U$ be a relatively compact subanalytic open neighbourhood
of $P$.
Let $U\cap \nbigc=\coprod\nbigc_i$
denote the decomposition
into connected components.
There exists $i_0$  such that
$g_1-g_2$ is not bounded from above
on $\nbigc_{i_0}$.
We obtain the following:
\begin{multline}
\label{eq;16.7.15.1}
\Hom_{\Ecat^b(\IC_{M})}
 \Bigl(
\cnum_M^{\Ecat}\overset{+}{\otimes}
\bigl(
 \pi^{-1}(\cnum_{\nbigc_{i_0}}\otimes L_{1,M})\otimes\cnum_{t\geq g_1}
\bigr),
\cnum_M^{\Ecat}\overset{+}{\otimes}
\bigl(
 \pi^{-1}(L_{2,M})\otimes\cnum_{t\geq g_2}
\bigr)
 \Bigr)
 \\
\simeq
 \varinjlim_{a\to\infty}
 \Hom_{\Mod(\cnum_{M\times\real})}
 \bigl(
 \pi^{-1}(\cnum_{\nbigc_{i_0}}\otimes L_{1,M})\otimes
 \cnum_{t\geq g_1-a},
 \pi^{-1}(L_{2,M})\otimes
 \cnum_{t\geq g_2}
 \bigr)=0
\end{multline}
Let $s$ be an element of the image of
(\ref{eq;16.7.14.21}).
By (\ref{eq;16.7.15.1}),
the restriction of $s$
to $\nbigc_{i_0}$ is $0$.
Because $\nbigc$ is assumed to be connected,
we obtain $s=0$.
\hfill\qed

\begin{cor}
Suppose that 
$g_1-g_2$ is positive
and that $g_2\not\prec g_1$.
Then, 
$\cnum_M^{\Ecat}\overset{+}{\otimes}
 \cnum_{g_2\leq t<g_1}$
is not isomorphic to $0$
in $\Ecat^b(\IC_M)$.
\end{cor}
\pf
There exists the distinguished triangle
$\cnum_M^{\Ecat}\overset{+}{\otimes}\cnum_{g_2\leq t<g_1}
\lrarr
\cnum_M^{\Ecat}\overset{+}{\otimes}\cnum_{t\geq g_2}
\stackrel{a}\lrarr
 \cnum_M^{\Ecat}\overset{+}{\otimes}\cnum_{t\geq g_1}
\lrarr
\cnum_M^{\Ecat}\overset{+}{\otimes}\cnum_{g_2\leq t<g_1}[1]$.
Because $a$ is not an isomorphism in
$\Ecat^b(\IC_M)$ by Lemma \ref{lem;16.7.14.20},
$\cnum_M^{\Ecat}\overset{+}{\otimes}\cnum_{g_2\leq t<g_1}$
is not isomorphic to $0$.
\hfill\qed

\vspace{.1in}
The restriction
$\Ecat\iota_{\nbigc}^{-1}(K)$
is denoted by
$K_{|\nbigc}$.
If $K=\cnum_M^{\Ecat}\overset{+}{\otimes}
 \cnum_{t\geq g_i}$,
then we obtain
$K_{|U}\simeq \cnum_U^{\Ecat}$
as observed in Corollary \ref{cor;16.7.14.20}.

\begin{cor}
\label{cor;16.6.22.10}
An isomorphism
$\Phi:
(\cnum_M^{\Ecat}\overset{+}{\otimes}
 \cnum_{t\geq g_1})_{|\nbigc}
\simeq
 (\cnum_M^{\Ecat}\overset{+}{\otimes}
 \cnum_{t\geq g_2})_{|\nbigc}$
can be extended to a morphism
$\Phitilde:\cnum_M^{\Ecat}\overset{+}{\otimes}
 \cnum_{t\geq g_1}
\lrarr
 \cnum_M^{\Ecat}\overset{+}{\otimes}
 \cnum_{t\geq g_2}$
if and only if the following holds:
\begin{itemize}
\item
For any real analytic map 
$\gamma:(\II^{\circ},\II)\lrarr (\nbigc,M)$,
the isomorphism
$\gamma^{-1}(\Phi)$
extends to a morphism
$\Ecat\gamma^{-1}\Bigl(
 \cnum_M^{\Ecat}\overset{+}{\otimes}
 \cnum_{t\geq g_1}
 \Bigr)
\lrarr
 \Ecat\gamma^{-1}\Bigl(
 \cnum_M^{\Ecat}\overset{+}{\otimes}
 \cnum_{t\geq g_2}\Bigr)$.
\end{itemize}
If such a morphism $\Phitilde$ exists,
it is unique.
\end{cor}
\pf
The ``only if'' part is clear.
Let us see the ``if'' part.
Under the assumption,
$\gamma^{-1}(g_2)\prec\gamma^{-1}(g_1)$
holds
for any real analytic map 
$\gamma:(\II^{\circ},\II)\lrarr (\nbigc,M)$.
Hence, $g_2\prec g_1$ holds
by Lemma \ref{lem;16.7.14.40}.
\hfill\qed

\begin{lem}
\label{lem;16.1.25.110}
Let $h$, $\alpha$ and $\beta$ be 
continuous subanalytic functions on $(U,M)$.
We assume that $\alpha(P)<\beta(P)$ for any $P\in U$.
Let $L_i$ $(i=1,2)$ be local systems on $U$.
Then, we obtain
\begin{equation}
 \label{eq;16.4.13.22}
 \Hom_{\Ecat^b(\IC_M)}\Bigl(
 \cnum_M^{\Ecat}\overset{+}{\otimes}
 \bigl(
\pi^{-1}(L_{1,M})\otimes
 \cnum_{t\geq h}
 \bigr),
 \,\,
 \cnum_M^{\Ecat}\overset{+}{\otimes}
 \bigl(
 \pi^{-1}(L_{2,M})\otimes
 \cnum_{\alpha\leq t<\beta}
 \bigr)
 \Bigr)
=0.
\end{equation}
\end{lem}
\pf
Let $U$ be a relatively compact subanalytic open subset
in $M$.
There exists the following natural isomorphisms:
\begin{multline}
 \Hom_{\Ecat^b(\IC_M)}\Bigl(
 \cnum_M^{\Ecat}\overset{+}{\otimes}
\bigl(
\pi^{-1}(\cnum_U\otimes L_{1,M})\otimes
 \cnum_{t\geq h}\bigr),
 \cnum_M^{\Ecat}\overset{+}{\otimes}
 \bigl(
 \cnum_{\alpha\leq t<\beta}
\otimes
 \pi^{-1}(L_{2,M})
 \bigr)
 \Bigr)
\simeq
 \\
 \varinjlim_{a\to\infty}
 \Hom_{\Mod(\cnum_{M\times\real})}\Bigl(
 \bigl(
 \pi^{-1}(\cnum_U\otimes L_{1,M})\otimes
 \cnum_{t\geq h-a}\bigr),
 \pi^{-1}(L_{2,M})\otimes
 \cnum_{\alpha\leq t<\beta}
 \Bigr)
=0.
\end{multline}
Then, we obtain (\ref{eq;16.4.13.22})
by taking the projective limit.
\hfill\qed

\vspace{.1in}
Similarly, we obtain the following lemma.

\begin{lem}
Let $\alpha_i$ and $\beta_i$  $(i=1,2)$ be
continuous subanalytic functions on $(\nbigc,M)$.
We assume $\alpha_i(P)<\beta_i(P)$ for any $P\in \nbigc$.
Let $L_i$ $(i=1,2)$ be local systems on $\nbigc$.
Unless $\alpha_2\prec \alpha_1$
and $\beta_2\prec\beta_1$,
we obtain
\begin{equation}
 \label{eq;16.6.22.1}
 \Hom_{\Ecat^b(\IC_M)}\Bigl(
 \cnum_M^{\Ecat}\overset{+}{\otimes}
 \bigl(
 \pi^{-1}(L_{1,M})\otimes
 \cnum_{\alpha_1\leq t<\beta_1}
 \bigr),
 \,\,
 \cnum_M^{\Ecat}\overset{+}{\otimes}
 \bigl(
 \pi^{-1}(L_{2,M})\otimes
 \cnum_{\alpha_2\leq t<\beta_2}
 \bigr)
 \Bigr)
=0.
\end{equation}
\hfill\qed
\end{lem}

\subsubsection{Canonical filtrations}
\label{subsection;16.1.25.20}

Assume that $\nbigc$ is connected.
Let $I=\{g_1,\ldots,g_m\}$ be 
a finite tuple in $\Sub(\nbigc,M)$.
Let $\Ibar$ denote the image of
$I\lrarr \Subbar(\nbigc,M)$.
We obtain the induced order $\prec$ on $\Ibar$.
For each $g_j\in I$,
let $[g_j]$ denote the induced element in $\Ibar$.
We consider an object of the form 
$K=\bigoplus_{i=1}^m
 \cnum_M^{\Ecat}\overset{+}{\otimes}\cnum_{t\geq g_i}$
in $\Ecat^b(\IC_M)$.
It is equipped with the filtration $\nbigf^{\nbigc}$
indexed by $(\Ibar,\prec)$
obtained as follows:
\[
 \nbigf^{\nbigc}_{[g]}(K)
:=\bigoplus_{
\substack{g_i\in I\\
 [g_i]\prec[g]}}
 \cnum_M^{\Ecat}\overset{+}\otimes\cnum_{t\geq g_i}.
\]
We may regard $\nbigf^{\nbigc}$
as a filtration indexed by
$(\Subbar(\nbigc,M),\prec_{\nbigc})$.

\begin{lem}
\label{lem;16.1.26.20}
Let $I_1=\{g_{1},\ldots,g_{m}\}$
and $I_2=\{h_1,\ldots,h_{\ell}\}$
be tuples in $\Sub(\nbigc,M)$.
Suppose that
\[
 K
=\bigoplus_{i=1}^{m}
 \cnum_M^{\Ecat}\overset{+}{\otimes}\cnum_{t\geq g_i}
=\bigoplus_{j=1}^{{\ell}}
 \cnum_M^{\Ecat}\overset{+}{\otimes}\cnum_{t\geq h_j}.
\]
Then, 
$(\Ibar_1,\prec)=(\Ibar_2,\prec)$ holds,
and the induced filtrations on $K$ are the same.
\end{lem}
\pf
Let $a\in \Ibar_1$ be a minimal element.
Suppose that the set $\{a'\in\Ibar_2\,|\,a'\prec a\}$
is empty.
If $g_i\in I_1$ such that $[g_i]=a$,
any morphism
$\cnum^{\Ecat}_M\overset{+}{\otimes}\cnum_{t\geq g_i}
\lrarr
\cnum^{\Ecat}_M\overset{+}{\otimes}\cnum_{t\geq h_j}$
is $0$.
It implies that 
any morphism
$\cnum^{\Ecat}_M\overset{+}{\otimes}\cnum_{t\geq g_i}
\lrarr K$ is $0$,
which arises a contradiction.
Hence, the set 
$\{a'\in\Ibar_2\,|\,a'\prec a\}$
is not empty.
There exists a minimal element $b\in\Ibar_2$
such that $b\prec a$.
There also exists a minimal element $c$
in $\Ibar_1$ such that $c\prec b$,
and hence $c\prec a$.
Because $a$ and $c$ are minimal,
we obtain $a=c=b$.
The composite of the morphisms
\[
 \bigoplus_{j=1}^m
 \cnum_M^{\Ecat}\overset{+}{\otimes}\cnum_{t\geq g_j}
\lrarr
 \bigoplus_{j=1}^{\ell}
 \cnum_M^{\Ecat}\overset{+}{\otimes}\cnum_{t\geq h_j}
\lrarr
 \cnum_M^{\Ecat}\overset{+}{\otimes}
 \Bigl(
 \bigoplus_{j=1}^{\ell}
 \cnum_{t\geq h_j}
 \Bigl/\bigoplus_{[h_j]=b}\cnum_{t\geq h_j}
\Bigr)
\]
induces
\begin{equation}
\label{eq;16.1.24.10}
 \cnum_M^{\Ecat}\overset{+}{\otimes}
 \Bigl(
 \bigoplus_{j=1}^m
 \cnum_{t\geq g_j}
\Big/
  \bigoplus_{[g_j]=b}
 \cnum_{t\geq g_j}
 \Bigr)
\lrarr
  \cnum_M^{\Ecat}\overset{+}{\otimes}
 \Bigl(
 \bigoplus_{j=1}^{\ell}
 \cnum_{t\geq h_j}
 \Bigl/\bigoplus_{[h_j]=b}\cnum_{t\geq h_j}
\Bigr).
\end{equation}
It is easy to see that 
the morphism (\ref{eq;16.1.24.10})
is an isomorphism.
Hence, by an easy induction,
we obtain the claim of the lemma.
\hfill\qed

\begin{df}
If $K=\bigoplus_{g\in I}
 \cnum_M^{\Ecat}\overset{+}{\otimes}\cnum_{t\geq g}$,
we say that the tuple $I$ controls
the growth order of $K$.
\hfill\qed
\end{df}

As observed in Corollary \ref{cor;16.7.14.20},
for such $K$,
the restriction
$K_{|\nbigc}\in \Ecat^b(\IC_{\nbigc})$
is isomorphic to
$(\cnum_{\nbigc}^{\Ecat})^{\oplus |I|}$,
i.e.,
it comes from the free $\cnum_{\nbigc}$-module
which is equipped with 
the induced filtration $\nbigf$
indexed by $(\Ibar,\prec)$.
We can recover $K$ from 
the free $\cnum_{\nbigc}$-module
with the filtration $\nbigf$.

\begin{lem}
\label{lem;16.1.25.32}
Let $I_i$ $(i=1,2)$ be tuples of
continuous subanalytic functions on $(\nbigc,M)$.
We set
$K_i:=
 \bigoplus_{g\in I_i}\cnum_M^{\Ecat}\overset{+}{\otimes}\cnum_{t\geq g}$.
Let $(L_i,\vecnbigf)$ be the underlying filtered $\cnum_{\nbigc}$-free modules.
There exists a natural bijection
\[
 \Hom_{\Ecat^b(\IC_M)}\bigl(
 K_1,K_2
 \bigr)
\simeq
 \bigl\{
f\in\Hom_{\Mod(\cnum_{\nbigc})}(L_1,L_2)
\,\big|\,
 f(\nbigf_aL_1)\subset
 \nbigf_aL_2\,\,
 (\forall a\in\Subbar(\nbigc,M))
 \bigr\}.
\]
In particular,
there exists the natural injection
$\Hom_{\Ecat^b(\IC_M)}\bigl(
 K_1,K_2
 \bigr)
\lrarr
 \Hom_{\Mod(\cnum_{\nbigc})}(L_1,L_2)$.
\end{lem}
\pf
It follows from Lemma \ref{lem;16.7.14.20}.
\hfill\qed

\vspace{.1in}

Let $I_i$ $(i=1,2)$ be tuples in $\Sub(\nbigc,M)$.
Let $\alpha_j<\beta_j$ $(j=1,\ldots,m)$ be continuous subanalytic functions 
on $(\nbigc,M)$.
We set
\[
 K_i:=\bigoplus_{g\in I_i}\cnum_M^{\Ecat}\overset{+}{\otimes}\cnum_{t\geq g},
\quad\quad
 F:=\bigoplus_{j=1}^m
 \cnum_M^{\Ecat}\overset{+}{\otimes}\cnum_{\alpha_j\leq t<\beta_j}.
\]
Note that
$\Hom_{\Ecat^b(\IC_M)}
 \bigl(
 K_i,F
 \bigr)=0$
by Lemma \ref{lem;16.1.25.110}.

\begin{lem}
\label{lem;16.1.25.112}
If $\rho:K_1\simeq K_2\oplus F$ in 
$\Ecat^b(\IC_M)$,
then we obtain $F=0$ and $K_1\simeq K_2$.
In particular,
$\alpha_j\sim\beta_j$ holds for any $j$.
Moreover,
$\Ibar_1=\Ibar_2$ holds
in $\Subbar(\nbigc,M)$,
which is compatible with the multiplicities.
\end{lem}
\pf
Let 
$a:F\lrarr K_2\oplus F$
and 
$b:K_2\oplus F\lrarr F$
denote the natural morphisms
for which $b\circ a=\id_F$ holds.
Because any morphism
$K_1\lrarr F$ is $0$,
we obtain $\id_F=0$,
and hence $F=0$.
\hfill\qed

\subsubsection{Prolongations of morphisms}

Let $I_i$ $(i=1,2)$ be finite tuples in $\Sub(\nbigc,M)$ $(i=1,2)$.
We set
$K_i:=\bigoplus_{g\in I_i}
 \cnum_M^{\Ecat}\overset{+}{\otimes}\cnum_{t\geq g}$.
Let $(L_i,\nbigf^{\nbigc})$ denote
the corresponding local systems with a filtration
on $\nbigc$.
Let $\Phi:L_1\lrarr L_2$ 
be a morphism of local systems on $\nbigc$.
A prolongation of $\Phi$
is a morphism
$\Phitilde:K_1\lrarr K_2$
such that
$\Ecat i_{\nbigc}^{-1}(\Phitilde)=\Phi$.
If such a prolongation exists,
it is unique.

\begin{lem}
\label{lem;16.6.22.11}
There exists a prolongation of $\Phi$
if and only if the following holds.
\begin{itemize}
\item
For any real analytic map  
$\gamma:(\II^{\circ},\II)\lrarr (\nbigc,M)$,
the morphism
$\gamma_{|\II^{\circ}}^{-1}(\Phi):
 \gamma_{|\II^{\circ}}^{-1}L_1
\lrarr
 \gamma_{|\II^{\circ}}^{-1}L_2$
extends to a morphism
$\Ecat\gamma^{-1}(K_1)
\lrarr
\Ecat\gamma^{-1}(K_2)$.
\end{itemize}
\end{lem}
\pf
The ``only if'' part is clear.
The ``if'' part follows from Corollary \ref{cor;16.6.22.10}.
\hfill\qed

\subsubsection{Appendix:
 Canonical filtrations on stably free enhanced ind-sheaves}

We can also consider a canonical filtration
in a more general case.
We set $\Subbar^{\ast}(\nbigc,M):=
 \Subbar(\nbigc,M)\sqcup\{\infty\}$.
We define the order 
$\prec$ on $\Subbar^{\ast}(\nbigc,M)^2$
 by the condition;
$(a,b)\prec(a',b')$
if 
$a\prec a'$
and $b\prec b'$.

We set $\Sub^{\ast}(\nbigc,M):=
 \Sub(\nbigc,M)\sqcup\{\infty\}$.
Let $(\alpha_i,\beta_i)$ $(i=1,\ldots,m)$
be pairs in $\Sub^{\ast}(\nbigc,M)$
such that $\alpha_i<\beta_i$ on $\nbigc$.
We consider any object of the form
$K=\bigoplus_{i=1}^m
 \cnum_M^{\Ecat}\overset{+}{\otimes}
 \cnum_{\alpha_i\leq t<\beta_i}$
in $\Ecat^b(\IC_M)$.
Then, we obtain the filtration
$\nbigf^{\nbigc}$ on $K$ indexed by
$(\Subbar^{\ast}(\nbigc,M)^2,\prec)$
by
\[
 \nbigf^{\nbigc}_{(\alpha,\beta)}(K):=
 \bigoplus_{([\alpha_i],[\beta_i])\prec([\alpha],[\beta])}
  \cnum_M^{\Ecat}\overset{+}{\otimes}
 \cnum_{\alpha_i\leq t<\beta_i}.
\]
It is canonically defined for $K$,
which can be shown
as in the case of Lemma \ref{lem;16.1.26.20}.

\subsection{Prolongations of local systems}
\label{subsection;18.11.16.112}

\subsubsection{Prolongations}

Let $M$ be a real analytic manifold.
Let $\nbigc$ be a locally closed subanalytic subset
in $M$.
The closure of $\nbigc$ in $M$
is denoted by $\nbigcbar$.
The inclusion $\nbigc\lrarr\nbigcbar$
is also denoted by $i_{\nbigc}$.
Let $\vecnbigc$ denote the bordered space $(\nbigc,\nbigcbar)$.
Let $i_{\vecnbigc}:\vecnbigc\lrarr M$ denote 
the naturally defined morphism of the bordered spaces.
We may naturally regard
$\Ecat^b_{\realc}(\IC_{\vecnbigc})$
as a full subcategory of
$\Ecat^b_{\realc}(\IC_M)$
by the fully faithful functor 
$\Ecat i_{\vecnbigc!!}$.

Let $L$ be a local system on $\nbigc$.
An object $K$ of $\Ecat^b_{\realc}(\IC_{\vecnbigc})$
with an isomorphism 
$\iota_K:\Ecat i_{\nbigc}^{-1}K\simeq 
 \cnum_{\nbigc}^{\Ecat}\otimes\pi^{-1}(L)$
in $\Ecat^b(\IC_{\nbigc})$
is called a prolongation of $L$
in $\Ecat^b_{\realc}(\IC_{\vecnbigc})$.
A morphism 
$\varphi:(K_1,\iota_{K_1})\lrarr
 (K_2,\iota_{K_2})$
of prolongations of $L$
is a morphism
$\varphi:K_1\lrarr K_2$
in $\Ecat^b_{\realc}(\IC_{\vecnbigc})$
such that
$\iota_{K_2}\circ\varphi_{|\nbigc}=\iota_{K_1}$.

Let $\vecII^{\circ}$ denote the bordered space
$(\II^{\circ},\II)$.
We consider the following condition
on prolongations $K\in \Ecat^b_{\realc}(\IC_{\vecnbigc})$
of a local system on $\nbigc$:
\begin{itemize}
\item
For any real analytic map 
$\gamma:(\II^{\circ},\II)\lrarr (\nbigc,M)$,
there exist continuous subanalytic functions
$h_1,\ldots,h_m$
on $(\II^{\circ},\II)$
such that
$\Ecat\gamma^{-1}K
\simeq
 \bigoplus_{j=1}^m
 \cnum_{\II}^{\Ecat}\overset{+}{\otimes}
 \cnum_{t\geq h_j}$
in 
$\Ecat^b_{\realc}(\IC_{\vecII^{\circ}})$.
\end{itemize}
Let $\Pro_{tf}(\nbigc,M)$
denote the category of prolongations 
of local systems satisfying the above condition.

\subsubsection{Filtrations by subanalytic subsets}

We recall \cite[Lemma 4.9.9]{DAgnolo-Kashiwara1}.
For any $K\in\Ecat^b_{\realc}(\IC_{M})$,
there exists a filtration
 $M=M^{(0)}\supset M^{(1)}\supset M^{(2)}\supset\cdots$
 by closed subanalytic subsets such that the following holds.
\begin{itemize}
\item
 $M^{(i)}\setminus M^{(i+1)}$
 are subanalytic submanifold of $M$
 of codimension $i$.
\item
 Let $M^{(i)}\setminus M^{(i+1)}=
 \coprod_{j\in\Lambda(i)}\nbigc^{(i)}_j$
 be the decomposition into the connected components.
Then, for any
$m\in\seisuu$,
$i\in\seisuu_{\geq 0}$
and $j\in\Lambda(i)$,
there exist continuous subanalytic functions 
 $g^{(i)}_{j,k,m}$ $(k\in\Gamma_1(i,j,m))$
 and 
 $\psi^{(i)}_{j,\ell,m}<\phi^{(i)}_{j,\ell,m}$ 
 $(\ell\in\Gamma_2(i,j,m))$
 on $(\nbigc^{(i)}_j,M)$
 such that
\[
 \pi^{-1}(\cnum_{\nbigc^{(i)}_j})
 \otimes
 K
\simeq
 \bigoplus_{m\in\seisuu}
 \bigoplus_{k\in\Gamma_1(i,j,m)}
 \cnum_M^{\Ecat}\overset{+}{\otimes}
 \cnum_{t\geq g^{(i)}_{j,k,m}}[m]
\oplus
 \bigoplus_{m\in\seisuu}
 \bigoplus_{\ell\in\Gamma_2(i,j,m)}
 \cnum_M^{\Ecat}\overset{+}{\otimes}
 \cnum_{\psi^{(i)}_{j,\ell,m}\leq t<\phi^{(i)}_{j,\ell,m}}[m].
\]
\end{itemize}
Such a filtration is called a filtration for $K$.

\begin{lem}
\label{lem;16.8.8.1}
For any $K\in\Pro_{tf}(\nbigc,M)$,
there exists a filtration
$\nbigc=\nbigc^{(0)}\supset \nbigc^{(1)}\supset \nbigc^{(2)}\supset\cdots$
such that the following holds.
\begin{itemize}
\item
$\nbigc^{(i)}\setminus \nbigc^{(i+1)}$
are locally closed subanalytic submanifolds of $M$
of codimension $i$ in $M$.
\item
Let 
$\nbigc^{(i)}\setminus \nbigc^{(i+1)}
=\coprod \nbigc^{(i)}_j$
be the decomposition
into connected components.
Then, there exist
 subanalytic functions 
$g^{(i)}_{j,k}$ $(k\in \Gamma(i,j))$
on $(\nbigc^{(i)}_{j},M)$
such that the following holds 
in $\Ecat^b_{\realc}(\IC_M)$:
\[
 \pi^{-1}(\cnum_{\nbigc^{(i)}_{j}})\otimes 
 \Ecat i_{\vecnbigc!!}K
\simeq \bigoplus_{k\in \Gamma(i,j)}
 \cnum_M^{\Ecat}\overset{+}{\otimes}\cnum_{t\geq g^{(i)}_{j,k}}.
\]
\end{itemize}
\end{lem}
\pf
Let $M=M^{(0)}\supset M^{(1)}\supset\cdots$
be a filtration for $\Ecat i_{\vecnbigc!!}K$ as above.
We may regard it as a filtration of $\nbigc$.
Let $M^{(i)}\setminus M^{(i+1)}
 =\coprod_{j\in \Lambda(i)}\nbigc^{(i)}_j$
be the decomposition into connected components.
For any $m\in\seisuu$, $i\in\seisuu_{\geq 0}$
and $j\in\Lambda(i)$,
there exist subanalytic functions
$g^{(i)}_{j,k,m}$ $(k\in \Gamma_1(i,j,m))$
and 
$\psi^{(i)}_{j,\ell,m}<\phi^{(i)}_{j,\ell,m}$ 
$(\ell\in\Gamma_2(i,j,m))$
with an isomorphism
\[
 \pi^{-1}(\cnum_{\nbigc^{(i)}_j})
 \otimes
 \Ecat i_{\vecnbigc!!}K
\simeq
\bigoplus_{m\in\seisuu}
 \bigoplus_{k\in\Gamma_1(i,j,m)}
 \cnum_M^{\Ecat}\overset{+}{\otimes}
 \cnum_{t\geq g^{(i)}_{j,k,m}}[m]
\oplus
 \bigoplus_{m\in\seisuu}
 \bigoplus_{\ell\in\Gamma_2(i,j,m)}
 \cnum_M^{\Ecat}\overset{+}{\otimes}
 \cnum_{\psi^{(i)}_{j,\ell,m}\leq t<\phi^{(i)}_{j,\ell,m}}[m].
\]
Let $\gamma:(\II^{\circ},\II)\lrarr (\nbigc^{(i)}_j,M)$
be any real analytic map
such that
$\gamma$ is an injection.
Then, there exist isomorphisms
\[
 \pi^{-1}(\cnum_{\gamma(\II^{\circ})})
 \otimes
 \Ecat i_{\vecnbigc!!}K
\simeq
 \Ecat\gamma_{!!}
 \Ecat\gamma^{-1}K
\simeq
 \bigoplus_{i=1}^{N}
 \cnum_M^{\Ecat}
 \overset{+}{\otimes}
 \cnum_{t\geq h_i}
\]
for some subanalytic functions $h_i$
on $(\gamma(\II^{\circ}),M)$.
Hence, we obtain
$\Gamma_{1}(i,j,m)=\emptyset$ unless $m=0$,
and 
$\psi^{(i)}_{j,\ell,m}$
and $\phi^{(i)}_{j,\ell,m}$ are mutually bounded
for any $(i,j,\ell,m)$.
\hfill\qed

\vspace{.1in}
Let $\nbigh^0\Ecat^b_{\realc}(\IC_M)$
denote the heart of $\Ecat^b_{\realc}(\IC_M)$
with respect to the $t$-structure
in \cite[\S4.6]{DAgnolo-Kashiwara1}.
(See also \cite[Lemma 4.9.5]{DAgnolo-Kashiwara1}.)

\begin{cor}
$\Ecat i_{\vecnbigc!!}\Pro_{tf}(\nbigc,M)$
is a full subcategory of
$\nbigh^0\Ecat^b_{\realc}(\IC_M)$.
\hfill\qed
\end{cor}

\subsubsection{Description as quotient}
\label{subsection;18.11.16.100}

Let $(K,\iota_K)\in\Pro_{tf}(\nbigc,M)$
be a prolongation of a local system $L$ on $\nbigc$.
Set $L_M:=i_{\nbigc!}L\in\Mod(\cnum_M)$.
For any continuous subanalytic function $G\in\Sub(\nbigc,M)$,
we obtain the following object
in $\Pro_{tf}(\nbigc,M)$:
\[
K_G:=\Ecat i_{\vecnbigc}^{-1}\Bigl( 
\cnum_M^{\Ecat}\overset{+}{\otimes}
 \bigl(
\pi^{-1}L_M\otimes
 \cnum_{t\geq G}
 \bigr)
 \Bigr).
\]

\begin{lem}
\label{lem;16.7.23.10}
Suppose that $\nbigc$ is relatively compact in $M$.
Then, there exist
a continuous subanalytic function $G\in \Sub(\nbigc,M)$
and a morphism
$\rho:K_G \lrarr  K$
in $\Pro_{tf}(\nbigc,M)$,
such that
$\Ecat i_{\vecnbigc!!}(\rho)$ is an epimorphism
in $\nbigh^0\Ecat^b_{\realc}(\IC_{M})$.
Such a morphism
$K_G\lrarr K$ is unique.
\end{lem}
\pf
By \cite[Lemma 4.6.3, Proposition 4.7.2]{DAgnolo-Kashiwara1},
there exists an $\real$-constructible sheaf
$\gbigf$ on
$\vecnbigc\times\real_{\infty}$
such that
$\cnum_{t\geq 0}\overset{+}{\otimes}\gbigf
\simeq\gbigf$
and
$K\simeq
 \cnum_{\nbigcbar}^{\Ecat}\overset{+}{\otimes}
 \gbigf$.
As in the proof of \cite[Lemma 4.9.9]{DAgnolo-Kashiwara1},
there exists a filtration
$\nbigc=\nbigc^{(0)}\supset
 \nbigc^{(1)}\supset
 \nbigc^{(2)}\supset\cdots$
of $\nbigc$ by closed subanalytic subsets
such that the following holds:
\begin{itemize}
\item
$\nbigc^{(i)}\setminus \nbigc^{(i+1)}$
are submanifolds of codimension $i$ in $M$.
\item
Let $\nbigc^{(i)}\setminus \nbigc^{(i+1)}=
 \coprod\nbigc^{(i)}_j$ be the decomposition
into the connected components.
Then, there exist
subanalytic functions
$g^{(i)}_{j,k}$ $(k\in \Gamma(i,j))$
and $\psi^{(i)}_{j,\ell}<\phi^{(i)}_{j,\ell}$
$(\ell\in\Lambda(i,j))$
on connected components $\nbigc^{(i)}_j$,
and isomorphisms
\[
 \pi^{-1}(\cnum_{\nbigc^{(i)}_j})
 \otimes \gbigf
\simeq
 \bigoplus_{k\in\Gamma(i,j)}
 \cnum_{t\geq g^{(i)}_{j,k}}
\oplus
 \bigoplus_{\ell\in\Lambda(i,j)}
 \cnum_{\psi^{(i)}_{j,\ell}\leq t<\phi^{(i)}_{j,\ell}}.
\]
\end{itemize}
For any relatively compact subset $V\subset \nbigc$,
the restriction of $|g^{(i)}_{j,k}|$
to $V\cap \nbigc^{(i)}_j$ are bounded.
As we have already observed,
$\psi^{(i)}_{j,\ell}$ and $\phi^{(i)}_{j,\ell}$
are mutually bounded for each $(i,j,\ell)$.

Let $\gbigf_1$ and $\gbigf_2$ denote the image
and the cokernel of the naturally induced morphism
$\pi^{-1}\pi_{\ast}(\gbigf)\lrarr\gbigf$.
It is easy to see that
\[
 \pi^{-1}(\cnum_{\nbigc^{(i)}_j})
 \otimes \gbigf_1
\simeq
 \bigoplus_{k\in\Gamma(i,j)}
 \cnum_{t\geq g^{(i)}_{j,k}},
\quad
 \pi^{-1}(\cnum_{\nbigc^{(i)}_j})
 \otimes \gbigf_2
\simeq
 \bigoplus_{\ell\in\Lambda(i,j)}
 \cnum_{\psi^{(i)}_{j,\ell}\leq t<\phi^{(i)}_{j,\ell}}.
\]
Then, we obtain
$\cnum^{\Ecat}_{\vecnbigc}\overset{+}{\otimes}\gbigf_2\simeq 0$
and 
$K\simeq 
 \cnum^{\Ecat}_{\vecnbigc}
 \overset{+}{\otimes}
 \gbigf
\simeq
  \cnum^{\Ecat}_{\vecnbigc}
 \overset{+}{\otimes}
 \gbigf_1$.
We may assume that $\Lambda(i,j)=\emptyset$
from the beginning,
i.e.,
\begin{equation}
\label{eq;16.7.15.10}
 \pi^{-1}(\cnum_{\nbigc^{(i)}_j})\otimes
 \gbigf
\simeq
 \bigoplus_{k\in\Gamma(i,j)}
 \cnum_{t\geq g^{(i)}_{j,k}}.
\end{equation}

By using Lemma \ref{lem;16.7.23.1}
and Lemma \ref{lem;16.7.23.2},
there exists a continuous subanalytic function $G\in \Sub(\nbigc,M)$
such that 
$G_{|\nbigc^{(i)}_j}<g^{(i)}_{j,k}$
for any $i,j,k$.
By (\ref{eq;16.7.15.10}),
we obtain
$\cnum_{t<G}\otimes\gbigf=0$.

There exists the natural isomorphism
$\Phi:
 \cnum_{\nbigc}^{\Ecat}
 \overset{+}{\otimes}
 \bigl(
 \pi^{-1}(L)\otimes\cnum_{t\geq 0}
 \bigr)
\simeq
 \bigl(
 \cnum_M^{\Ecat}
 \overset{+}{\otimes}
 \gbigf\bigr)_{|\nbigc}$
in 
$\Ecat^b_{\realc}(\IC_{\nbigc})$.
Let $V\subset \nbigc$ be 
any relatively compact open subset.
Then, $\Phi_{|V}$ corresponds to a morphism
\[
 \Psi_{V,a(V)}:
 \pi^{-1}(\cnum_{V})\otimes
  \pi^{-1}(L)\otimes\cnum_{t\geq -a(V)}
\lrarr
 \pi^{-1}(\cnum_{V})\otimes
 \gbigf
\]
for a sufficiently large  number $a(V)>0$.
Because $\cnum_{t<G}\otimes\gbigf=0$,
the following morphism is induced by $\Psi_{V,a(V)}$:
\[
 \Psi_{V,G}:
 \pi^{-1}(\cnum_V)\otimes
  \pi^{-1}(L)\otimes\cnum_{t\geq G}
\lrarr
 \pi^{-1}(\cnum_V)\otimes
 \gbigf.
\]
By enlarging $V$,
we obtain the morphism
$\Psi_{G}:
  \pi^{-1}(L)\otimes\cnum_{t\geq G}
\lrarr
 \gbigf$
on $\nbigc\times\real_{\infty}$.
It induces the following morphism
in $\Pro_{tf}(\nbigc,M)$:
\[
K_{G}=
 \Ecat\iota_{\vecnbigc}^{-1}
 \cnum_M^{\Ecat}
 \overset{+}{\otimes}
 \bigl(
 \pi^{-1}(L_M)\otimes
 \cnum_{t\geq G}
 \bigr)
\stackrel{\rho}{\lrarr}
 \Ecat\iota_{\vecnbigc}^{-1}
\cnum_M^{\Ecat}
 \overset{+}{\otimes}
 \gbigf
=K.
\]
Note that
$\Ecat i_{\vecnbigc!!}K$ and
$\Ecat i_{\vecnbigc!!}K_G$
are objects in
$\nbigh^0
 \Ecat^b_{\realc}(\IC_M)$.
We set
\[
 \nbigf^{(i)}
 \Ecat i_{\vecnbigc!!}K:=
 \pi^{-1}(\cnum_{\nbigc^{(i)}})
 \otimes
 \Ecat i_{\vecnbigc!!}K,
\quad
 \nbigf^{(i)}
 \Ecat i_{\vecnbigc!!}K_G:=
 \pi^{-1}(\cnum_{\nbigc^{(i)}})
 \otimes
 \Ecat i_{\vecnbigc!!}K_G.
\]
By Lemma \ref{lem;16.7.23.20},
they induce filtrations
in the abelian category
$\nbigh^0\Ecat^b_{\realc}(\IC_M)$.
By construction,
the morphisms
\[
 \Gr_{\nbigf}^{(i)}(\rho):
 \Gr_{\nbigf}^{(i)}
 \Ecat i_{\vecnbigc!!}K_G
\lrarr
 \Gr_{\nbigf}^{(i)}
 \Ecat i_{\vecnbigc!!}K
\]
are equal to the morphisms
induced by the epimorphisms
$\cnum_{t\geq G^{(i)}_{j}}
\lrarr
 \cnum_{t\geq g^{(i)}_{j,k}}$,
where
$G^{(i)}_{j}$
denotes the restriction of $G$ to
$\nbigc^{(i)}_{j}$.
Because 
the functor
$\gbigg\longmapsto
 \cnum_M^E\overset{+}{\otimes}
 \gbigg$
is exact
and preserves the $t$-structure,
we obtain that
$\Gr_{\nbigf}^{(i)}(\rho)$
are epimorphisms.
Hence, we obtain that
$\rho$ is an epimorphism.

Let $\rho':K_G\lrarr K$ be any morphism
in $\Pro_{tf}(\nbigc,M)$.
It is induced by a morphism
$\pi^{-1}(L)\otimes
 \cnum_{t\geq G-a}
\lrarr\gbigf$
in $\Mod(\cnum_{\nbigc\times\real})$
for some $a>0$.
We obtain the induced morphism
$\Psi_G':\pi^{-1}(L)\otimes
 \cnum_{t\geq G}
\lrarr\gbigf$,
which also induces $\rho'$.
Note that the induced morphisms
$\cnum_{\nbigc}^{\Ecat}
 \overset{+}{\otimes}\Psi'_G$
and 
$\cnum_{\nbigc}^{\Ecat}
 \overset{+}{\otimes}\Psi_G$
are the same in $\Ecat^b(\IC_{\nbigc})$.
Hence, we obtain that
$\Psi_G=\Psi_G'$.
\hfill\qed

\subsubsection{Prolongations of isomorphisms}

Let $(K_i,\iota_{K_i})\in \Pro_{tf}(\nbigc,M)$ $(i=1,2)$
be prolongations of a local system $L$ on $\nbigc$.

\begin{prop}
\label{prop;16.6.23.21}
The identity of $L$ is extended to
an isomorphism
$(K_1,\iota_{K_1})
\simeq
 (K_2,\iota_{K_2})$
in $\Pro_{tf}(\nbigc,M)$
if and only if
the following holds.
\begin{itemize}
\item
For any real analytic map
  $\gamma:(\II^{\circ},\II)\lrarr (\nbigc,M)$,
 the identity of 
 $\gamma_{|\II^{\circ}}^{-1}(L)$ is extended to
 an isomorphism
 $\Ecat\gamma^{-1}(K_1,\iota_1)
\simeq
 \Ecat\gamma^{-1}(K_2,\iota_2)$
 in $\Pro_{tf}(\II^{\circ},\II)$.
\end{itemize}
If such an isomorphism exists,
it is unique.
\end{prop}
\pf
It is enough to consider the case where $\nbigc$ is relatively compact.
The ``only if'' part is clear.
Let us prove the ``if'' part.
Let $L$ denote the local system on $\nbigc$
such that $K_a$ are prolongations of $L$.
There exist a continuous subanalytic function $G$
and morphisms
$K_G:=
 \cnum_M^{\Ecat}
 \overset{+}{\otimes}
 \bigl(
 \pi^{-1}(L)
 \otimes\cnum_{t\geq G}
 \bigr)
\stackrel{\rho_a}{\lrarr}
 K_a$
in $\Pro_{tf}(\nbigc,M)$.
We obtain the objects
$\Ker(\Ecat i_{\vecnbigc!!}(\rho_i))$
in $\Ecat^b_{\realc}(\IC_{M})$.
It is enough to prove
$\Ker(\Ecat i_{\vecnbigc!!}(\rho_1))
=\Ker(\Ecat i_{\vecnbigc!!}(\rho_2))$.

There exists a filtration 
$\nbigc=\nbigc^{(0)}
 \supset
 \nbigc^{(1)}\supset\cdots$
for both $K_1$ and $K_2$,
by closed subanalytic subsets such that
$\nbigc^{(i)}\setminus \nbigc^{(i+1)}$
are locally closed subanalytic submanifolds
of $M$ of codimension $i$.
Let $\nbigc^{(i)}\setminus \nbigc^{(i+1)}=
 \coprod_{j\in\Lambda(i)}\nbigc^{(i)}_j$
be the decomposition into connected components.
For $a=1,2$,
there exist subanalytic functions
$g^{(i)}_{a,j,p}$ $(p\in\Gamma(a,i,j))$
such that
\[
 \pi^{-1}(\cnum_{\nbigc^{(i)}_j})
 \otimes
 \Ecat i_{\vecnbigc!!}K_a
\simeq
 \bigoplus_{p\in\Gamma(a,i,j)}
 \cnum_M^{\Ecat}
 \overset{+}{\otimes}
 \cnum_{t\geq g^{(i)}_{a,j,p}}.
\]
By Lemma \ref{lem;16.6.22.11},
the identities of
$L_{|\nbigc^{(i)}_j}$
uniquely extend to isomorphisms
$\pi^{-1}(\cnum_{\nbigc^{(i)}_j})
 \otimes
 \Ecat i_{\vecnbigc!!}K_1
\simeq
 \pi^{-1}(\cnum_{\nbigc^{(i)}_j})
 \otimes
 \Ecat i_{\vecnbigc!!}K_2$.

We set
$\nbigd^{(i)}:=\nbigc\setminus\nbigc^{(\dim M+1-i)}$.
Then, $\nbigd^{(i)}$ gives a filtration by
open subanalytic subsets of $\nbigc$
such that
$\nbigd^{(i)}\setminus\nbigd^{(i+1)}$
is $i$-dimensional locally closed subanalytic 
submanifolds of $M$.
We set
$\nbigf^{(i)}N:=
 \pi^{-1}(\cnum_{\nbigd^{(i)}})
 \otimes N$
for $N\in 
\nbigh^0\Ecat^b_{\realc}(\IC_M)$.
They induce filtrations $\nbigf$
on $N$
in the category 
$\nbigh^0\Ecat^b_{\realc}(\IC_M)$.
We also set
$\Gr_{\nbigf}^{(i)}(N):=
 \nbigf^{(i)}N/\nbigf^{(i+1)}N$.
Note that
$\Ker\Gr^{(i)}_{\nbigf}(\Ecat i_{\vecnbigc!!}(\rho_1))
=\Ker\Gr^{(i)}_{\nbigf}(\Ecat i_{\vecnbigc!!}(\rho_2))$
in $\Gr_{\nbigf}^{(i)}\Ecat i_{\vecnbigc!!}K_G$.
By Lemma \ref{lem;16.7.23.21}
and Lemma \ref{lem;16.7.23.20},
we obtain
$\Gr^{(i)}_{\nbigf}\Ker(\Ecat i_{\vecnbigc!!}(\rho_1))
=\Gr^{(i)}_{\nbigf}\Ker(\Ecat i_{\vecnbigc!!}(\rho_2))$
in $\Gr_{\nbigf}^{(i)}\Ecat i_{\vecnbigc!!}K_G$.
Hence, the claim of the ``if'' part follows 
from Lemma \ref{lem;16.7.15.20} below.
The uniqueness of isomorphisms
follows from the uniqueness 
in Lemma \ref{lem;16.7.23.10}.
\hfill\qed

\subsubsection{Appendix}

Let $\nbigc$ be any locally closed relatively compact 
subanalytic subset in $M$.
Let $\nbigz\subset \nbigcbar$ be a subanalytic closed subset.
We set $\nbigc_1:=\nbigc\setminus \nbigz$
and $\nbigz_1:=\nbigz\cap\nbigc$.

Let $N_0$ be an object in
$\nbigh^0\Ecat^b_{\realc}(\IC_{M})$
obtained as 
$\iota_{\vecnbigc!!}N_0'$
for an object
$N_0'\in \Ecat^b_{\realc}(\IC_{\vecnbigc})$.
Let $N_a$ $(a=1,2)$ be subobjects of $N_0$
in $\nbigh^0\Ecat^b_{\realc}(\IC_M)$
such that
$\pi^{-1}(\cnum_{\nbigc})\otimes N_a\simeq N_a$.
There exist the following exact sequences 
$(a=0,1,2)$:
\[
 0\lrarr
 \pi^{-1}(\cnum_{\nbigc_1})\otimes N_a
\lrarr
 N_a
\lrarr
 \pi^{-1}(\cnum_{\nbigz_1})\otimes N_a
\lrarr 0.
\]
Suppose that 
$\pi^{-1}(\cnum_{\nbigz_1})\otimes N_1
=\pi^{-1}(\cnum_{\nbigz_1})\otimes N_2$
in 
$\pi^{-1}(\cnum_{\nbigz_1})\otimes N_0$,
and 
$\pi^{-1}(\cnum_{\nbigc_1})\otimes N_1
=\pi^{-1}(\cnum_{\nbigc_1})\otimes N_2$
in 
$\pi^{-1}(\cnum_{\nbigc_1})\otimes N_0$.
\begin{lem}
\label{lem;16.7.15.20}
Under the assumption,
$N_1=N_2$ holds.
\end{lem}
\pf
It is enough to prove that the induced morphism 
$N_1\lrarr N_0/N_2$ is $0$.
It is enough to prove that the induced morphism
$\pi^{-1}(\cnum_{\nbigz_1})\otimes
 N_1
\lrarr
 \pi^{-1}(\cnum_{\nbigc_1})\otimes
 (N_0/N_2)$
is $0$.

Note that
$N_1=\cnum_M^{\Ecat}\overset{+}{\otimes}\gbigf$
and 
$N_0/N_2=
 \cnum_M^{\Ecat}\overset{+}{\otimes}\gbigg$
for some 
$\real$-constructible sheaves
$\gbigf,\gbigg$
on $\vecnbigc\times\real_{\infty}$
such that
$\cnum_{t\geq 0}\overset{+}{\otimes}\gbigf\simeq\gbigf$
and 
$\cnum_{t\geq 0}\overset{+}{\otimes}\gbigg\simeq\gbigg$.
A morphism
$\pi^{-1}(\cnum_{\nbigz_1})\otimes
 N_1
\lrarr
 \pi^{-1}(\cnum_{\nbigc_1})\otimes
 (N_0/N_2)$
corresponds to
$\cnum_{t\geq -a}
\overset{+}{\otimes}
 \bigl(
 \pi^{-1}(\cnum_{\nbigz_1})\otimes\gbigf
 \bigr)
\lrarr
 \pi^{-1}(\cnum_{\nbigc_1})
 \otimes
 \gbigg$
for some $a>0$,
which has to be $0$
by the support condition.
Hence, we obtain the desired vanishing.
\hfill\qed

\subsection{A sufficient condition for the existence of global filtrations}
\label{subsection;18.11.16.113}

Let $U_1$ be an open ball in $\real^{n-1}$.
We put $U:=\II^{\circ}\times U_1$
and $\Ubar:=\II\times U_1$.
We use the coordinate system
$(r,\vecx)$ 
of $\real\times\real^{n-1}$.
For any ramified real analytic function
$f=\sum f_{\gminiy}(\vecx)r^{\gminiy}$,
set $\ord_r(f):=\min\{\gminiy\,|\,f_{\gminiy}\neq 0\}$.

Let $\nbigj=\{f^{(1)},\ldots,f^{(m)}\}$ 
be a finite tuple of ramified real analytic functions on $U$ 
of the form
$f^{(j)}=\sum f^{(j)}_{\gminiy} r^{\gminiy}$,
where $r$ is the coordinate of $I$,
and $f^{(j)}_{\gminiy}$ are real analytic functions on $U_1$
such that the following holds:
\begin{itemize}
\item
 If $j\neq k$,
 $a(j,k):=\ord_r(f^{(j)}-f^{(k)})$ is negative,
 and 
 $f^{(j)}_{a(j,k)}-f^{(k)}_{a(j,k)}$ 
 is a nowhere vanishing function on $U_1$.
\end{itemize}
Let $\gminim(1),\ldots,\gminim(m)$
be non-negative integers.
Let $K\in\Pro_{tf}(U,\Ubar)$ 
be a prolongation of $L$
such that the following holds:
\begin{itemize}
\item 
 Let $\gamma:(\II^{\circ},\II)\lrarr (U,\Ubar)$
 be any real analytic map
 such that 
$\gamma(0)\in \Ubar\setminus U$.
Then,
$\Ecat\gamma^{-1}(K)$
is a stably free enhanced ind-sheaf
induced by
the local system
$\gamma_{|\II^{\circ}}^{-1}(L)$ 
with a filtration $\nbigf$
indexed by 
 $\gamma^{\ast}\nbigj$
such that 
$\rank\Gr^{\nbigf}_{\gamma^{-1}f^{(j)}}
 \gamma^{-1}(L)=\gminim(j)$.
\end{itemize}
Let us prove the following proposition.
\begin{prop}
\label{prop;16.8.25.2}
There exists a global filtration $\nbigf$ on $L$
indexed by $\nbigj$ such that 
$K$ is induced by $(L,\nbigf)$.
\end{prop}

\subsubsection{Preliminary}

Let $\nbigc\subset\real^n$
be a locally closed subanalytic subset.
We consider tuples
$\nbigj_1=\{f_1,\ldots,f_m\},
\nbigj_2=\{g_1,\ldots,g_{\ell}\}
\subset\Sub(\nbigc,\real^n)$.
\begin{lem}
\label{lem;16.2.5.100}
We assume the following.
\begin{itemize}
\item
For any real analytic map
$\gamma:(\II^{\circ},\II)\lrarr (\nbigc,\real^n)$,
$\overline{\gamma^{\ast}\nbigj_1}
=\overline{\gamma^{\ast}\nbigj_2}$
holds
in $\Subbar(\II^{\circ},\II)$.
\end{itemize}
Then, there exists a stratification
$\nbigc=\coprod\nbigc_i$
such that 
$\overline{\nbigj_{1|\nbigc_i}}
=\overline{\nbigj_{2|\nbigc_i}}$
in $\Subbar(\nbigc_i,\real^n)$
for any $i$.
\end{lem}
\pf
We may assume that $\nbigc$ is relatively compact
of pure dimension $m$.
Let $\overline{\nbigc}$ be the closure of $\nbigc$ in $\real^n$.
According to 
\cite{Bierstone-Milman, Hironaka-Introduction-real-analytic},
there exists
a uniformization $F:M\lrarr \real^n$ of $\overline{\nbigc}$,
i.e.,
$M$ is an $m$-dimensional real analytic manifold,
and $F$ is a real analytic map
such that $F(M)=\overline{\nbigc}$.
We obtain the functions
$F^{\ast}(f_i)$ and $F^{\ast}(g_j)$ on 
the subanalytic open subsets $\nbigc':=F^{-1}(\nbigc)$,
which are subanalytic as functions 
on $(\nbigc',M)$.

According to Corollary \ref{cor;18.11.10.1},
there exists a rectilinearization
$\lambda_{\alpha}:W_{\alpha}\lrarr M$ $(\alpha\in\Lambda)$
for $\nbigc'$ and the functions
$F^{\ast}(f_i)$, $F^{\ast}(g_j)$, 
$F^{\ast}(f_i-g_j)$, $F^{\ast}(f_i-f_j)$ and 
$F^{\ast}(g_i-g_j)$
for any $i$ and $j$.
Set $\phi_{\alpha}:=F\circ\lambda_{\alpha}$.
Let $\gbigq$ be a quadrant of $W_{\alpha}$
contained in $\phi_{\alpha}^{-1}(\nbigc')$.
Let us observe that
$\overline{\{\phi_{\alpha}^{\ast}(f_i)_{|\gbigq}\}}
=\overline{\{\phi_{\alpha}^{\ast}(g_j)_{|\gbigq}\}}$
holds in $\Subbar(\gbigq,W_{\alpha})$.
Let $(x_1,\ldots,x_n)$ be a standard coordinate system
of $W_{\alpha}$.
We may assume 
$\gbigq=\{x_i>0\,(i\leq m_1),\,\,x_i=0\,(i> m_1)\}$
for some $0\leq m_1\leq m$.

If $\phi_{\alpha}^{\ast}(f_i)_{|\gbigq}-\phi_{\alpha}^{\ast}(f_j)_{|\gbigq}$
is not constantly $0$,
we obtain the expression
\[
 \phi_{\alpha}^{\ast}(f_i)_{|\gbigq}
-\phi_{\alpha}^{\ast}(f_j)_{|\gbigq}=a^{(0)}_{i,j}\cdot
 \prod_{p\leq m_1} x_p^{\alpha(i,j)^{(0)}_p}, 
\]
where
$a^{(0)}_{i,j}$ is nowhere vanishing,
and $(\alpha(i,j)^{(0)}_p)\in
 (\rnum_{\geq 0})^{m_1}
\cup
 (\rnum_{\leq 0})^{m_1}$.
Hence,
$\overline{\{\phi_{\alpha}^{\ast}(f_j)_{|\gbigq}\}}
\subset
 \Subbar(\gbigq,W_{\alpha})$
is totally ordered.
By a similar argument,
we obtain that
$\overline{\{\phi_{\alpha}^{\ast}(f_i)_{|\gbigq}\}}
\cup
 \overline{\{\phi_{\alpha}^{\ast}(g_j)_{|\gbigq}\}}$
is totally ordered.

Suppose that 
$\overline{\phi_{\alpha}^{\ast}(f_i)_{|\gbigq}}$
is not contained in 
$\overline{\{\phi_{\alpha}^{\ast}(g_j)_{|\gbigq}\}}$.
For any $j$,
we express
\[
 \phi_{\alpha}^{\ast}(f_i)_{|\gbigq}
-\phi_{\alpha}^{\ast}(g_j)_{|\gbigq}=a_{i,j}\cdot
 \prod_{p\leq m_1} x_{p}^{\alpha(i,j)_{p}},
\]
where
$a_{i,j}$ are  nowhere vanishing,
and $(\alpha(i,j)_{p})\in(\rnum_{\leq 0})^{m_1}\setminus\{(0,\ldots,0)\}$.

Let $\phi_{\alpha}^{\ast}(g_{j_1})_{|\gbigq}$
be a representative of the maximum of
$\overline{\{\phi_{\alpha}^{\ast}(g_j)_{|\gbigq}
\,|\,
 \phi_{\alpha}^{\ast}(g_j)_{|\gbigq}
\prec
  \phi_{\alpha}^{\ast}(f_i)_{|\gbigq}
 \}}$.
Let $\phi_{\alpha}^{\ast}(g_{j_2})_{|\gbigq}$
be a representative of the minimum of
$\overline{\{\phi_{\alpha}^{\ast}(g_j)_{|\gbigq}
\,|\,
  \phi_{\alpha}^{\ast}(f_i)_{|\gbigq}
\prec
 \phi_{\alpha}^{\ast}(g_j)_{|\gbigq}\}}$.
Note that 
either one of the following holds:
(1) $\alpha(i,j_1)_{p}\leq \alpha(i,j_2)_{p}$ 
for any $p$,
or
(2) $\alpha(i,j_1)_{p}\geq \alpha(i,j_2)_{p}$
for any $p$.
Hence, there exists $p$ such that
$\alpha(i,j_1)_p<0$
and $\alpha(i,j_2)_p<0$.
Let $(c_1,\ldots,c_{p-1},c_{p+1},\ldots,c_{m_i})$
be any point of $\real_{>0}^{m_i-1}$.
Let $\gamma_p:\II\lrarr \gbigq$
be the real analytic map
such that
$\gamma_p(t)
=(c_1,\ldots,c_{p-1},\epsilon t,c_{p+1},\ldots,c_{m_i},0,\ldots,0)$
for a sufficiently small $\epsilon>0$.
Then, we obtain 
$\gamma_p^{\ast}\phi_{\alpha}^{\ast}(g_{j})
\neq
 \gamma_p^{\ast}\phi_{\alpha}^{\ast}(f_i)$
in $\Subbar(\II^{\circ},\II)$
for any $j$,
which contradicts the assumption.
Hence, we obtain 
$\phi_{\alpha}^{\ast}(f_i)_{|\gbigq}
\in
 \overline{\{\phi_{\alpha}^{\ast}(g_j)_{|\gbigq}\}}$.
Similarly,
we obtain
$\phi_{\alpha}^{\ast}(g_j)_{|\gbigq}
\in 
 \overline{\{\phi_{\alpha}^{\ast}(f_i)_{|\gbigq}\}}$,
and thus
$\overline{\{\phi_{\alpha}^{\ast}(f_i)_{|\gbigq}\}}
=\overline{\{\phi_{\alpha}^{\ast}(g_j)_{|\gbigq}\}}$.

There exist subanalytic compact subsets
$N_{\alpha}\subset W_{\alpha}$
such that 
$\bigcup_{\alpha\in\Lambda}
 \lambda_{\alpha}(N_{\alpha})=M$.
We can describe $\nbigc$
as the union of
the subanalytic subsets
$\nbign(\alpha,\gbigq):=\phi_{\alpha}(\gbigq\cap N_{\alpha})$,
where $\alpha$ runs through $\Lambda$,
and $\gbigq$ runs through the set of quadrants
in $W_{\alpha}$ such that
$\gbigq\subset\phi_{\alpha}^{-1}(\nbigc)$.
Note that
$\overline{\nbigj_{1|\nbign(\alpha,\gbigq)}}
=\overline{\nbigj_{2|\nbign(\alpha,\gbigq)}}$
in $\Subbar(\nbign(\alpha,\gbigq),\real^n)$.
By using \cite[Lemma 8.3.21]{Kashiwara-Schapira},
we can construct a stratification
with the desired property
which is finer than the covering
$\nbigc=\bigcup \nbign(\alpha,\gbigq)$.
\hfill\qed

\subsubsection{Proof of Proposition \ref{prop;16.8.25.2}}

By Lemma \ref{lem;16.2.5.100},
there exists a subanalytic stratification
$U=\coprod\nbigc$
with an isomorphism
\[
 \pi^{-1}(\cnum_{\nbigc})
 \otimes K
\simeq
 \bigoplus
 \Bigl(
 \cnum^{\Ecat}
 \overset{+}{\otimes}
 \cnum_{t\geq f^{(j)}_{|\nbigc}}
 \Bigr)
 \otimes\cnum^{\gminim(j)}.
\]
Let $\nbigjbar^{\nbigc}$ denote the image of
$\nbigj$ in $\Subbar(\nbigc,\Ubar)$.
We obtain the filtration
$\nbigf^{\nbigc}$ of $L_{|\nbigc}$
indexed by $\nbigjbar^{\nbigc}$
such that 
$\rank\Gr^{\nbigf^{\nbigc}}_{f}(L_{|\nbigc})
=\sum_{f^{(j)}\equiv f}\gminim(j)$,
where $f^{(j)}$ runs through $\nbigj$
which are equivalent to $f$
in $\Subbar(\nbigc,\Ubar)$.

\begin{lem}
\label{lem;18.11.7.2}
There exists a filtration $\nbigf$ on $L$
such that for any $\nbigc$
the filtration of $L_{|\nbigc}$ induced by $\nbigf$
indexed by $\nbigjbar^{\nbigc}$
is equal to the filtration $\nbigf^{\nbigc}$.
\end{lem}
\pf
Let $\nbigc_1$ be an $n$-dimensional stratum.
Let $\nbigc_2$ be a stratum contained in
the closure of $\nbigc_1$.
Assume that the closure of $\nbigc_2$
intersects with $\Ubar\setminus U$.
We have already obtained the filtrations
$\nbigf^{\nbigc_i}$ of $L_{|\nbigc_i}$ $(i=1,2)$.
By using a path
$\gamma:(\II^{\circ},\II)\lrarr (\nbigc_1,U)$
such that $\gamma(0)\in\nbigc_2$,
we obtain a filtration on $L_{|\nbigc_2}$
induced by the filtration $\nbigf^{C_1}$,
which is also denoted by $\nbigf^{\nbigc_1}$.
It is enough to compare
the filtrations $\nbigf^{\nbigc_1}$ and $\nbigf^{\nbigc_2}$
of $L_{|\nbigc_2}$.

By using a rectilinearization of
the subanalytic sets
$\nbigc_1$ and $\nbigc_2$,
there exists a real analytic map
$\phi:[0,1]\times[0,1]\lrarr \Ubar$
such that the following holds:
\begin{itemize}
\item
$\phi(a,b)\in \nbigc_1$
for any $a>0$ and $b>0$.
\item
 $\phi(a,0)\in\nbigc_2$ for any $a>0$.
\item
$\phi(0,b)\in \Ubar\setminus U$
 for any $b$.
\end{itemize}

Set $Z:=[0,1]\times[0,1]$
and $Z^{\circ}:=Z\setminus(\{0\}\times[0,1])$.
We obtain the bordered space
$\vecZ=(Z^{\circ},Z)$.
Let us study
$\Ecat\phi^{-1}K$
in $\Ecat^b_{\realc}(\vecZ)$.

Set 
$\nbigu:=Z^{\circ}\setminus
 ([0,1]\times\{0\})$
and $Y^{\circ}:=\openclosed{0}{1}\times\{0\}$.
We set $h^{(j)}:=\phi^{\ast}(f^{(j)})$.
The restrictions to $\nbigu$ and $Y^{\circ}$
are denoted by
$h^{(j)}_{\nbigu}$
and $h^{(j)}_{Y^{\circ}}$, respectively.

Let $\nbigf^{\nbigu}$ denote the filtration
of $\phi^{-1}(L)_{|\nbigu}$
indexed by $\{h^{(j)}_{\nbigu}\}$
which is induced by $\nbigf^{\nbigc_1}$.
It naturally extends to a filtration
$\nbigf^{Z^{\circ}}$
of $\phi^{-1}(L)$
indexed by $\{h^{(j)}\}$.
By taking the restriction to $Y^{\circ}$,
we obtain the filtration
$\nbigf^{Z^{\circ}}_{|Y^{\circ}}$
of $\phi^{-1}(L)_{|Y^{\circ}}$.
Let $\nbigf^{Y^{\circ}}$ denote the filtration
of $\phi^{-1}(L)_{|Y^{\circ}}$
indexed by $\{h^{(j)}_{Y^{\circ}}\}$
which is induced by $\nbigf^{\nbigc_2}$.
It is enough to prove 
$\nbigf^{Z^{\circ}}_{|Y^{\circ}}
=\nbigf^{Y^{\circ}}$.

There exists an $\real$-constructible sheaf
$N$ on $Z\times\real_{\infty}$
such that
$\cnum^{\Ecat}\overset{+}{\otimes} N
\simeq
 \Ecat\phi^{-1}(K)$
and $\cnum_{t\geq 0}\overset{+}{\otimes}N\simeq N$.
There exist the following natural morphisms:
\[
 \cnum^{\Ecat}\overset{+}{\otimes}
 \bigl(\pi^{-1}(\cnum_{\nbigu})\otimes N\bigr)
\simeq
 \bigoplus
 \cnum^{\Ecat}\overset{+}{\otimes}
 \bigl(
 \cnum_{t\geq h^{(j)}_{\nbigu}}
 \otimes\cnum^{\gminim(j)}
 \bigr)
\lrarr
  \bigoplus
 \cnum^{\Ecat}\overset{+}{\otimes}
 \bigl(
 \cnum_{t\geq h^{(j)}}
 \otimes\cnum^{\gminim(j)}
 \bigr).
\]
We may assume that it is induced by a morphism
$\pi^{-1}(\cnum_{\nbigu})\otimes N
\lrarr
 \bigoplus \cnum_{t\geq h^{(j)}}
 \otimes\cnum^{\gminim(j)}$,
which extends to
$N\lrarr
 \bigoplus \cnum_{t\geq h^{(j)}}
 \otimes\cnum^{\gminim(j)}$.
Hence, there exists a morphism
\[
 \Ecat\phi^{-1}(K)
\lrarr
   \bigoplus
 \cnum^{\Ecat}\overset{+}{\otimes}
 \bigl(
 \cnum_{t\geq h^{(j)}}
 \otimes\cnum^{\gminim(j)}
 \bigr)
\]
whose restriction to
$Z^{\circ}$ is the identity of 
$\phi^{-1}(L)$.
Set $Y:=[0,1]\times\{0\}$,
and $\vecY=(Y^{\circ},Y)$.
Let $\iota:\vecY\lrarr \vecZ$ be the inclusion
of bordered spaces.
We obtain a morphism
\[
 \Ecat(\phi\circ\iota)^{-1}K
\lrarr
 \bigoplus \cnum^{\Ecat}
 \overset{+}{\otimes}
 \bigl(\cnum_{t\geq h^{(j)}_{Y^{\circ}}}
 \otimes\cnum^{\gminim(j)}\bigr)
\]
whose restriction to $Y^{\circ}$
is the identity of 
$\phi^{-1}(L)_{|Y^{\circ}}$.
It implies
$\nbigf^{Y^{\circ}}_{f^{(j)}}
\subset
 \bigl(
 \nbigf^{\nbigu}_{f^{(j)}}
 \bigr)_{|Y^{\circ}}$
 for any $f^{(j)}$.
By comparing the rank,
we obtain 
$\nbigf^{Y^{\circ}}_{f^{(j)}}
=
 \bigl(
 \nbigf^{\nbigu}_{f^{(j)}}
 \bigr)_{|Y^{\circ}}$.
\hfill\qed

\vspace{.1in}

Then, the claim of Proposition \ref{prop;16.8.25.2}
follows from Proposition \ref{prop;16.6.23.21}.
\hfill\qed

\part{Subanalytic functions and complex blowings up}
\label{part;18.11.16.31}

\section{Preliminary}
\label{section;18.11.15.10}

\subsection{Notation for infinite sequences of blowings up}
\label{subsection;18.11.25.20}

We prepare some notation
to study sequences of complex blowings up at points on surfaces
in a way convenient to us.

\subsubsection{The basic case}

Let $\KK$ denote $\real$ or $\cnum$.
Let $(x,y)$ be the standard coordinate system of $\KK^2$.
Let $p:\Bl_{(0,0)}\KK^2\lrarr \KK^2$
be the blowing up at $(0,0)$,
i.e.,
$\Bl_{(0,0)}\KK^2:=
\bigl\{
 \bigl(
 (x,y),[x':y']
 \bigr)\in
 \KK^2\times\proj^1(\KK)\,\big|\,
 xy'-x'y=0 \bigr\}$
and $p((x,y),[x':y']):=(x,y)$.
We set
$P_+:=\bigl((0,0),[1:0]\bigr)$ 
and
$P_-:=\bigl((0,0),[0:1]\bigr)$.
We also set
$U_+:=\{((x,y),[x':y']\in\Bl_{(0,0)}\KK^2\,|\,
 a\neq 0\}$
and
$U_-:=\{((x,y),[x':y']\in\Bl_{(0,0)}\KK^2\,|\,
 b\neq 0\}$.
We obtain the coordinate neighbourhood
$(U_+,u_+,v_+)$ induced by
$(u_+,v_+)=(x,y/x)$ around $P_+$,
and $(U_-,u_-,v_-)$ induced by 
$(u_-,v_-)=(x/y,y)$ around $P_-$.
We identify 
$U_{\pm}\simeq\KK^2$ by the coordinate systems.
The restriction of $p$ to $U_{\pm}$
are denoted by $p_{\pm}$.
The morphisms
$p_{\pm}$ are described as
$p_+(x_1,y_1)=(x_1,x_1y_1)$
and 
$p_-(x_1,y_1)=(x_1y_1,y_1)$
with respect to the standard coordinate system
$(x_1,y_1)$ on $\KK^2$.

For any $\alpha\neq 0$,
we put $P_{\alpha}:=(\alpha,0)\in U_-$.
We obtain the coordinate neighbourhood
$(U_{\alpha},u_{\alpha},v_{\alpha})$
of $P_{\alpha}$
defined by
$U_{\alpha}:=U_-$
and 
$(u_{\alpha},v_{\alpha}):=(u_--\alpha,v_-)$.
For the coordinate system
$(u_{\alpha},v_{\alpha})$,
$P_{\alpha}$ is denoted as $(0,0)$.
For the coordinate system $(U_+,u_+,v_+)$,
$P_{\alpha}$ is denoted as $(0,\alpha^{-1})$.

\subsubsection{Sequences of blowings up}
\label{subsection;18.11.14.2}

We set $\gbigp(\KK):=(\KK\setminus\{0\})\cup\{+,-\}$.
Let $\veceta=(\eta_1,\ldots,\eta_{\ell})$
be an element in $\gbigp(\KK)^{\ell}$
for some $\ell\geq 1$.
We shall construct a sequence of 
spaces
$X_{\veceta}^{(i)}$ $(i=0,\ldots,\ell)$
with points
$P_{\veceta}^{(i)}$,
maps 
$p^{(i)}:(X_{\veceta}^{(i)},P_{\veceta}^{(i)})
 \lrarr (X_{\veceta}^{(i-1)},P_{\veceta}^{(i-1)})$,
and coordinate neighbourhoods
$(U_{\veceta}^{(i)},u^{(i)}_{\veceta},v^{(i)}_{\veceta})$
around $P_{\veceta}^{(i)}$.

We set $X_{\veceta}^{(0)}:=\KK^2$
and $P^{(0)}_{\veceta}=(0,0)$.
The coordinate neighbourhood
$(U_{\veceta}^{(0)},u^{(0)}_{\veceta},v^{(0)}_{\veceta})$
is defined to be $(\KK^2,x,y)$.
Let $p_{\veceta}^{(1)}:X_{\veceta}^{(1)}\lrarr X_{\veceta}^{(0)}$
be the blowing up at $P^{(0)}_{\veceta}$.
There exists the natural isomorphism
$X_{\veceta}^{(1)}\simeq \Bl_{(0,0)}\KK^2$.
Let $P^{(1)}_{\veceta}\in X_{\veceta}^{(1)}$ be the point
corresponding to
$P_{\eta_1}\in\Bl_{(0,0)}\KK^2$.
The coordinate neighbourhood
$(U^{(1)}_{\veceta},u^{(1)}_{\veceta},v^{(1)}_{\veceta})$
is defined to be
$(U_{\eta_1},u_{\eta_1},v_{\eta_1})$.

Suppose that we have already obtained
$X_{\veceta}^{(i)}$, $P^{(i)}_{\veceta}$
and $(U^{(i)}_{\veceta},u^{(i)}_{\veceta},v^{(i)}_{\veceta})$.
Let $p^{(i+1)}_{\veceta}:X_{\veceta}^{(i+1)}\lrarr 
 X_{\veceta}^{(i)}$
be the blowing up at $P^{(i)}_{\veceta}$.
The inclusion
$U^{(i)}_{\veceta}\lrarr X_{\veceta}^{(i)}$
induces
$\Bl_{P_{\veceta}^{(i)}} U^{(i)}_{\veceta}
\lrarr
 X_{\veceta}^{(i+1)}$.
There exists the isomorphism
$\Bl_{P_{\veceta}^{(i)}} U^{(i)}_{\veceta}
\simeq
 \Bl_{(0,0)}\KK^2$
induced by the coordinate system
$(u^{(i)}_{\veceta},v^{(i)}_{\veceta})$.
Let $P^{(i+1)}_{\veceta}$
be the point corresponding to
$P_{\eta_{i+1}}$.
Let $(U^{(i+1)}_{\veceta},u^{(i+1)}_{\veceta},v^{(i+1)}_{\veceta})$
be the coordinate neighbourhood around 
$P^{(i+1)}_{\veceta}$ corresponding to
$(U_{\eta_{i+1}},u_{\eta_{i+1}},v_{\eta_{i+1}})$.

We set
$\Bl_{\veceta}\KK^2:=X_{\veceta}^{(\ell)}$,
$P_{\veceta}:=P^{(\ell)}_{\veceta}$
and
$(U_{\veceta},u_{\veceta},v_{\veceta}):=
 (U^{(\ell)}_{\veceta},u^{(\ell)}_{\veceta},v^{(\ell)}_{\veceta})$.

\subsubsection{Infinite sequences}

For any $\veceta=(\eta_1,\eta_2,\ldots)\in \gbigp(\KK)^{\infty}$,
we set $\veceta_k:=(\eta_1,\ldots,\eta_k)$.
We obtain the sequence of  the spaces
$\Bl_{\veceta_k}\KK^2$
with the base point $P_{\veceta_k}$
and the coordinate neighbourhood
$(U_{\veceta_k},u_{\veceta_k},v_{\veceta_k})$.
Let 
$p_{\veceta}^{(k)}:
 (\Bl_{\veceta_k}\KK^2,P_{\veceta_k})
\lrarr
 (\Bl_{\veceta_{k-1}}\KK^2,P_{\veceta_{k-1}})$
denote the naturally induced morphisms.

If there exists a sequence $k_j\to\infty$
such that $\eta_{k_j}\in\KK\setminus\{0\}$,
then the sequences $\veceta$
and $(\Bl_{\veceta_{k}}\KK_2,P_{\veceta_{k}})$
$(k=1,2,\ldots,)$
are called of mixed type in this paper.

\begin{rem}
Such sequences $\veceta$ correspond to 
sequence of ``infinitely near points''
in {\rm\cite[\S6.2]{Favre-Jonsson}}.
``Mixed type'' in this paper
corresponds to 
the union of Type {\rm 1} and Type {\rm 4}.
\hfill\qed
\end{rem}

\subsection{Positively linear subsets}

\subsubsection{Subsets of $S^1\times S^1$}

For $(p,q)\in\seisuu^2$
and 
$\vecphi=(\phi_1,\phi_2)
 \in S^1\times S^1$,
we set
$H(p,q,\vecphi):=
\bigl\{
 (\theta_1,\theta_2)\in S^1\times S^1\,|\,
 p(\theta_1-\phi_1)
=q(\theta_2-\phi_2)
 \bigr\}$.
Let $A$ be a closed subanalytic subset in
$S^1\times S^1$
such that $A$ is purely $1$ dimensional,
i.e.,
any open subset of $A$ is one dimensional.
We say that $A$ is positively linear
if there exist 
$(p_i,q_i,\vecphi_i)\in
 \seisuu_{>0}^2\times (S^1)^2$ 
$(i=1,\ldots,N)$
such that
$A$ is contained in
$\bigcup_{i=1}^N
 H(p_i,q_i,\vecphi_i)$.

\begin{lem}
\label{lem;16.2.4.30}
Let $B\subset S^1\times S^1$
be a closed subanalytic subset
such that $B$ is purely $1$ dimensional.
Then, there exists a decomposition
$B=B_1\cup B_2$
such that
(i) $B_1$ is positively linear,
(ii) $\dim(B_1\cap B_2)=0$,
(iii) any open subset of $B_2$ is 
 not positively linear.
\end{lem}
\pf
Let $B_1$ be the union of
the positively linear subsets in $B$.
It is enough to prove that $B_1$ is subanalytic.
Let $P$ be any point of $B$.
Because $\dim B=1$,
$B$ is semianalytic.
(See \cite{Bierstone-Milman}.)
Let $U_P$ be a small neighbourhood of $P$
in $(S^1)^2$.
If $U_P$ is sufficiently small,
there exist real analytic maps
$\kappa_i:\II_{\epsilon}\lrarr U_P$
$(i=1,\ldots,m)$
such that
$U_P\cap B=\bigcup_{i=1}^m\kappa_i(\II_{\epsilon})$.
Then, the claim of the lemma is clear.
\hfill\qed

\subsubsection{Subsets in the oriented real blowing up 
of complex surfaces}

Let $X$ be any complex surface
with a normal crossing hypersurface $H$.
Let $\varpi:\Xtilde(H)\lrarr X$
be the oriented blowing up of $X$ along $H$.
Let $P$ be any cross point of $H$.
Let $(X_P,z_1,z_2)$ be a holomorphic coordinate neighbourhood
around $P$
such that $H=\bigcup_{i=1,2}\{z_i=0\}$.
Let $z_i=r_i\exp(\sqrt{-1}\theta_i)$ be the polar decomposition.
It induces an isomorphism
$\Phi_{\vecz}:\varpi^{-1}(P)\simeq S^1\times S^1$.
We say that a locally closed purely $1$-dimensional
subanalytic subset $Z$ in $\varpi^{-1}(P)$
is positively linear
if $\Phi_{\vecz}(Z)\subset S^1\times S^1$
is positively linear.

For another choice of holomorphic coordinate system
$(w_1,w_2)$
such that $\{w_i=0\}=\{z_i=0\}$,
the isomorphism
$\Phi_{\vecw}\circ\Phi_{\vecz}^{-1}$
is a translation
$(\theta_1,\theta_2)
 \longmapsto
 (\theta_1+\theta_1^{(0)},\theta_2+\theta_2^{(0)})$.
Hence, 
the family of positively linear subsets
in $\varpi^{-1}(P)$
is defined independently from the choice of 
the holomorphic coordinate system.

\vspace{.1in}

Let $f:\Delta\lrarr X$ be a holomorphic map
such that $f(\Delta\setminus\{0\})\subset X\setminus H$
and that $f(0)=P$ is a cross point.
Let $\varpi_{\Delta}:\Deltatilde(0)\lrarr \Delta$
be the oriented real blowing up along $0$.
Let $\ftilde:\Deltatilde(0)\lrarr \Xtilde(H)$ denote
the induced morphism.
The following obvious lemma is 
the motivation to consider positively linear subsets.
\begin{lem}
$\ftilde(\varpi_{\Delta}^{-1}(0))$ is a positively linear subset
of $\varpi^{-1}(P)$.
\hfill\qed
\end{lem}

\section{Sequences of complex blowings up at cross points}
\label{section;16.9.4.1}

In this section,
we study sequences of complex blowing ups of a complex surface
at cross points of the exceptional divisors.
We shall explain the following;
for any subanalytic function on the original surface,
if the length of the sequence is large enough,
the pull back of the function
has nice property around the centers of blowings up.
In \S\ref{subsection;18.11.15.30},
we state related results concerning a sequence of
complex blowings up at cross points.
In \S\ref{subsection;18.11.15.31},
we state a more technical result 
which is useful for comparison of
sequences of complex blowings up
at cross points
and sequences of real blowings up.
The rest of the section is devoted 
to the proof of the statements.

\subsection{Statements}

\subsubsection{Resolutions by sequences of complex blowings up at cross points}
\label{subsection;18.11.15.30}

We set $H:=\{x=0\}\cup\{y=0\}\subset\cnum^2$
and $O=(0,0)$.
Let $\veceta\in\{+,-\}^{\ell}$.
We have constructed the morphism
$\psi_{\veceta}:
 (\Bl_{\veceta}\cnum^2,P_{\veceta})
\lrarr
 (\cnum^2,O)$
in \S\ref{subsection;18.11.25.20}.
We set $H_{\veceta}:=\psi_{\veceta}^{-1}(H)
\subset \Bl_{\veceta}\cnum^2$.
Let $\varpi:\widetilde{\cnum^2}(H)\lrarr \cnum^2$
and 
$\varpi_{\veceta}:
 \widetilde{\Bl_{\veceta}\cnum^2}(H_{\veceta})\lrarr 
 \Bl_{\veceta}\cnum^2$
be the oriented real blowing up.
Let 
$\widetilde{\psi_{\veceta}}:
 \widetilde{\Bl_{\veceta}\cnum^2}(H_{\veceta})
\lrarr
 \widetilde{\cnum^2}(H)$
denote the induced morphism.
We obtain the morphism
$\Psi_{\veceta}:
 \varpi_{\veceta}^{-1}(P_{\veceta})
\lrarr
 \varpi^{-1}(O)$
as the restriction of 
$\widetilde{\psi_{\veceta}}$.

We shall prove the following lemma
in \S\ref{subsection;18.11.12.2}.

\begin{lem}
\label{lem;16.4.27.30}
For any positively linear subset
$\gbiga\subset \varpi^{-1}(O)$,
there exists $\ell_0$ such that 
$\Psi_{\veceta}^{-1}(\gbiga)$
does not contain a positively linear subset
for any $\ell\geq \ell_0$ and
$\veceta\in\{+,-\}^{\ell}$.
\end{lem}

We shall prove the following proposition
in \S\ref{subsection;18.11.12.1}.

\begin{prop}
\label{prop;16.4.26.11}
Let $\gbigg$ be a closed subanalytic subset 
in $\widetilde{\cnum^2}(H)$ with $\dim_{\real}\gbigg\leq 3$.
Then, there exists $\ell_0$
such that the following holds
for any $\ell\geq \ell_0$
and $\veceta\in \{+,-\}^{\ell}$:
\begin{itemize}
\item
Let $\gbigg_{\veceta}$ denote the strict transform of
$\gbigg$ with respect to
$\widetilde{\psi_{\veceta}}$,
i.e.,
$\gbigg_{\veceta}$ is the closure of
$\psi_{\veceta}^{-1}(\gbigg\setminus \varpi^{-1}(H))$
in $\widetilde{\Bl_{\veceta}\cnum^2}(H_{\veceta})$.
Then, 
$\dim_{\real}\bigl(
 \gbigg_{\veceta}\cap
 \varpi_{\veceta}^{-1}(P_{\veceta})
 \bigr)\leq 1$
holds.
\end{itemize}
\end{prop}

Let $\gbigh$ be a closed subanalytic subset
in $\widetilde{\cnum^2}(H)$
with $\dim_{\real} \gbigh\leq 3$
such that the following holds.
\begin{itemize}
\item
We set $\gbigh^{\circ}:=\gbigh\cap(\cnum^2\setminus H)$,
and let $\overline{\gbigh^{\circ}}$
denote the closure of $\gbigh^{\circ}$
in $\cnumtilde^2(H)$.
 Then,
 $\dim\bigl(
 \varpi^{-1}(O)
 \cap
 \overline{\gbigh^{\circ}}
 \bigr)\leq 1$.
\end{itemize}
Let $f_i$ $(i=1\ldots,n)$ be real analytic functions
on $\widetilde{\cnum^2}(H)\setminus 
 (\gbigh\cup\varpi^{-1}(H))$
which is subanalytic on
$\bigl(
 \widetilde{\cnum^2}(H)\setminus (\gbigh\cup\varpi^{-1}(H)),
 \widetilde{\cnum^2}(H)\bigr)$.
We shall prove the following proposition
in \S\ref{subsection;18.11.12.3}.
\begin{prop}
\label{prop;16.4.27.20}
There exist a positive number $\ell_0$
and a closed subanalytic subset
$\gbigz\subset\varpi^{-1}(O)$
with $\dim_{\real}\gbigz\leq 1$
such that the following holds
for any 
$\ell\geq \ell_0$
and $\veceta\in\{+,-\}^{\ell}$:
\begin{itemize}
\item
 $\gbigz$ contains
 $\overline{\gbigh^{\circ}}\cap
 \varpi^{-1}(O)$.
\item
 For any 
 $Q\in 
 \varpi_{\veceta}^{-1}(P_{\veceta})
 \setminus
 \Psi_{\veceta}^{-1}(\gbigz)$,
 the functions 
 $\widetilde{\psi_{\veceta}}^{\ast}(f_i)$
 are ramified real analytic around $Q$.
\end{itemize}
\end{prop}

\subsubsection{Lifting with respect to sequences of 
local real blowings up}
\label{subsection;18.11.15.31}

Let $M$ be any $2$-dimensional real analytic manifold.
Although we are interested in the case $M=S^1\times S^1$ in this paper,
we study the problem in a slightly general way.
Let $\phi:W\lrarr \real^2\times M$
be the composition of local real blowings up
i.e.,
there exists a factorization of $\phi$ as follows:
\[
 W=W_k
 \stackrel{\phi_k}{\lrarr}
 W_{k-1}
 \stackrel{\phi_{k-1}}{\lrarr}
 \cdots
 \stackrel{\phi_2}{\lrarr}
 W_1
 \stackrel{\phi_1}{\lrarr}
 W_0=\real^2\times M.
\]
Moreover, there exist subanalytic open subsets
$U_i\subset W_i$
and closed real analytic submanifolds $C_i\subset U_i$
such that the morphisms $W_{i+1}\lrarr W_i$
are the real blowing up of $U_i$ along $C_i$.
For simplicity,
we assume that $C_i$ are subanalytic in $W_i$.

Let $\ell$ be any integer larger than $k$,
and let 
$\veceta=(\eta_1,\ldots,\eta_{\ell})$
be any element of $\{+,-\}^{\ell}$.
For each $m\leq \ell$,
we set 
$\veceta_m:=(\eta_1,\ldots,\eta_m)$.
Let $\Crit(\phi)\subset \real^2\times M$ 
be the set of the critical values of $\phi$.
Note $\dim_{\real}\Crit(\phi)\leq 2$.
Let $\gbigb(\veceta)\subset M$ be the set of
$Q\in M$ such that 
there exists a real analytic path 
$\gamma:\II\lrarr \Bl_{\veceta}\real^2\times M$
satisfying the following conditions:
\begin{itemize}
\item
$\gamma(0)=(P_{\veceta},Q)$,
and 
$\psi_{\veceta}\circ\gamma(t)\in
 (\real^2\times M)\setminus
 \Crit(\phi)$
for $t>0$.
\item
There exists a real analytic map
$\gamma':\II\lrarr W$
such that 
$\phi\circ\gamma'=(\psi_{\veceta}\times\id)\circ\gamma$.
Note that such $\gamma'$ is uniquely determined.
\end{itemize}
Note that
$\gbigb(\veceta)$ may be empty set.

We shall prove the following proposition
in \S\ref{subsection;16.5.31.3}.

\begin{prop}
\label{prop;16.4.27.10}
For any $(W,\phi)$,
there exists a $1$-dimensional closed 
subanalytic subset $Z\subset M$
such that for any $\ell\geq k$,
any $\veceta\in\{\pm\}^{\ell}$
and any $Q\in \gbigb(\veceta)\setminus Z$,
there exist a small neighbourhood $\nbigu$ of 
$(P_{\veceta},Q)$ in
$\Bl_{\veceta}\real^2\times M$,
a small neighbourhood
$W_{\gamma'(0)}$ of $\gamma'(0)$
in $W$,
a non-negative integer $i(Q)\leq \ell$,
an open embedding
$\iota_Q:W_{\gamma'(0)}\lrarr
 \Bl_{\veceta_{i(Q)}}\real^2\times M$,
such that the following holds:
\begin{itemize}
\item
$\phi_{|W_{\gamma'(0)}}=(\psi_{\veceta_{i(Q)}}\times\id)\circ\iota_Q$.
\item 
The morphism
$\Bl_{\veceta}\real^2\times M
\lrarr
 \Bl_{\veceta_{i(Q)}}\real^2\times M$
induces the real analytic map
$\nbigu\lrarr 
 \iota_Q(W_{\gamma'(0)})$.
\end{itemize}
\end{prop}

\subsection{Preliminary}

\subsubsection{Explicit descriptions}
\label{subsection;18.11.14.10}

Let $\KK$ denote $\real$ or $\cnum$.
For any 
$\veceta\in\{+,-\}^{\ell}$,
we constructed 
$\psi^{(i)}_{\veceta}:
 X^{(i)}_{\veceta}\lrarr\KK^2$
in \S\ref{subsection;18.11.14.2}.
Let us describe the restriction
$(\psi^{(i)}_{\veceta})_{|U^{(i)}_{\veceta}}:
 U^{(i)}_{\veceta}
 \lrarr\KK^2$
under the identification
$\KK^2\simeq U^{(i)}_{\veceta}$
by
$(x_1,y_1)=(u_{\veceta}^{(i)},v_{\veceta}^{(i)})$.
We set
\begin{equation}
\label{eq;16.4.26.1}
 A_+:=\left(
 \begin{array}{cc}
 1 & 0 \\ 1 & 1
 \end{array}
 \right),
\quad\quad
 A_-:=\left(
 \begin{array}{cc}
 1 & 1 \\ 0 & 1
 \end{array}
 \right).
\end{equation}
We set
$A^{(i)}_{\veceta}
=A_{\eta_1}A_{\eta_2}
\cdot \cdots\cdot
 A_{\eta_{i}}$.
Let
$\alpha^{(i)}(\veceta)$,
$\beta^{(i)}(\veceta)$,
$\gamma^{(i)}(\veceta)$
and $\delta^{(i)}(\veceta)$
be the components of
$A^{(i)}_{\veceta}$:
\[
A^{(i)}_{\veceta}
=\left(
 \begin{array}{cc}
 \alpha^{(i)}(\veceta)
 &
 \beta^{(i)}(\veceta) \\
 \gamma^{(i)}(\veceta)
 &
 \delta^{(i)}(\veceta)
 \end{array}
 \right).
\]
Then,
$\psi^{(i)}_{\veceta}(x_1,y_1)
=(x_1^{\alpha^{(i)}(\veceta)}y_1^{\beta^{(i)}(\veceta)},
 x_1^{\gamma^{(i)}(\veceta)}y_1^{\delta^{(i)}(\veceta)})$
holds.

We set
$A_{\veceta}:=A^{(\ell)}_{\veceta}$,
$\alpha(\veceta):=\alpha^{(\ell)}(\veceta)$,
$\beta(\veceta):=\beta^{(\ell)}(\veceta)$,
$\gamma(\veceta):=\gamma^{(\ell)}(\veceta)$,
and
$\delta(\veceta):=\delta^{(\ell)}(\veceta)$.
The restriction
$(\psi_{\veceta})_{|U_{\veceta}}:U_{\veceta}\lrarr \KK^2$
is denoted as
\[
 \psi_{\veceta}\bigl(
 u_{\veceta},v_{\veceta}
 \bigr)
=\bigl(
 u_{\veceta}^{\alpha(\veceta)}v_{\veceta}^{\beta(\veceta)},
\,\,
 u_{\veceta}^{\gamma(\veceta)}v_{\veceta}^{\delta(\veceta)}
 \bigr).
\]

\subsubsection{Polar decompositions}
\label{subsection;16.4.26.3}

There exists the natural inclusion
$P_{\veceta}\in\Bl_{\veceta}\real^2\subset
 \Bl_{\veceta}\cnum^2$.
For $\KK=\real$ or $\cnum$,
there exist the coordinate neighbourhoods
$(U_{\veceta,\KK},u_{\veceta,\KK},v_{\veceta,\KK})$
of $P_{\veceta}$ in $\Bl_{\veceta}\KK^2$,
where the subscript $\KK$ is included
for distinction.
The coordinate system
$(u_{\veceta,\real},v_{\veceta,\real})$
is the restriction of 
$(u_{\veceta,\cnum},v_{\veceta,\cnum})$.

Let $\Bl_{\veceta}(\real_{\geq 0}^2)$
denote the strict transform of 
$\real_{\geq 0}^2\subset\real^2$
with respect to the morphism
$\psi_{\veceta}:\Bl_{\veceta}(\real^2)
 \lrarr \real^2$.
Namely,
$\Bl_{\veceta}(\real_{\geq 0}^2)$
is the closure of
$\psi_{\veceta}^{-1}(\real_{>0}^2)$
in $\Bl_{\veceta}(\real^2)$.
Let 
$\psi_{\veceta,\real_{\geq 0}}:\Bl_{\veceta}(\real^2_{\geq 0})
\lrarr\real^2_{\geq 0}$
denote the induced morphism.
Note that
$U_{\veceta,\real}\cap
 \Bl_{\veceta}(\real^2_{\geq 0})
=\{u_{\veceta,\real}\geq 0,v_{\veceta,\real}\geq 0\}$.

\vspace{.1in}

Let $H:=\{x=0\}\cup\{y=0\}$ in $\cnum^2$.
We set 
$H_{\veceta}:=
 \psi_{\veceta}^{-1}(H)\subset\Bl_{\veceta}\cnum^2$.
Let 
$\varpi:\widetilde{\cnum^2}(H)\lrarr \cnum^2$
and 
$\varpi_{\veceta}:
 \widetilde{\Bl_{\veceta}\cnum^2}(H_{\veceta})
\lrarr
 \Bl_{\veceta}\cnum^2$
denote the oriented real blowings up.
We obtain the induced morphism
$\widetilde{\psi_{\veceta,\cnum}}:
 \widetilde{\Bl_{\veceta}\cnum^2}(H_{\veceta})
\lrarr
 \widetilde{\cnum^2}(H)$.

By the polar decompositions of the coordinate functions
$u_{\veceta,\cnum}
 =r_{\veceta,1}e^{\sqrt{-1}\theta_{\veceta,1}}$
and 
$v_{\veceta,\cnum}
 =r_{\veceta,2}e^{\sqrt{-1}\theta_{\veceta,2}}$,
we obtain
$\varpi_{\veceta}^{-1}(U_{\veceta,\cnum})
\simeq
 \real^2_{\geq 0}\times (S^1)^2$.
We also obtain the natural identification
$\widetilde{\cnum^2}(H)
\simeq
 \real_{\geq 0}^2\times (S^1)^2$
induced by the coordinate system $(x,y)$.
The following is clear by the expression of
the morphisms $\psi_{\veceta,\KK}$.

\begin{lem}
The restriction of 
$\widetilde{\psi_{\veceta,\cnum}}$
to $\varpi_{\veceta}^{-1}(U_{\veceta,\cnum})$
is identified with the product of the following morphisms:
\begin{itemize}
\item
the induced morphism
$\psi_{\veceta,\real_{\geq 0}}:
 \Bl_{\veceta}(\real^2_{\geq 0})
\cap
 U_{\veceta,\real}
\lrarr\real^2_{\geq 0}$,
\item
the isomorphism
$\Psi_{\veceta}:(S^1)^2\simeq(S^1)^2$
induced by
$\bigl(\theta_{\veceta,1},\theta_{\veceta,2}\bigr)
\longmapsto
\bigl(
 \alpha(\veceta)\theta_{\veceta,1}
+\beta(\veceta)\theta_{\veceta,2},
 \gamma(\veceta)\theta_{\veceta,1}
+\delta(\veceta)\theta_{\veceta,2}
\bigr)$.

\hfill\qed
\end{itemize}
\end{lem}

\subsection{Proof of Lemma \ref{lem;16.4.27.30}}
\label{subsection;18.11.12.2}

First, we remark the following obvious lemma.

\begin{lem}
\label{lem;16.4.26.2}
Suppose that $pq\leq 0$.
Then, for any $\veceta\in\{+,-\}^{\ell}$,
there exists 
$(p',q',\vecphi')\in\seisuu^2\times(S^1)^2$
such that 
$p'q'\leq 0$
and that
$\Psi_{\veceta}^{-1}H(p,q,\vecphi)
\subset
 H(p',q',\vecphi')$.
\end{lem}
\pf
It is enough to check the claim
in the case
$\veceta=(\eta_1)\in\{+,-\}$,
which can be checked directly.
\hfill\qed

\vspace{.1in}

We easily obtain Lemma \ref{lem;16.4.27.30}
from the following lemma
and the description of the morphism $\psi_{\veceta}$
in \S\ref{subsection;16.4.26.3}.

\begin{lem}
\label{lem;16.4.26.4}
For any positively linear subset
$\gbiga\subset S^1\times S^1$,
there exists $\ell_0$ such that
the following holds 
for any $\ell\geq \ell_0$
and $\veceta\in\{+,-\}^{\ell}$:
\begin{itemize}
\item
 $\Psi_{\veceta}^{-1}(\gbiga)$
 does not contain 
 positively linear subsets.
\end{itemize}
\end{lem}
\pf
It is enough to consider the case where
$\gbiga$ is an open subset of $H(p,q,\vecphi)$
for $(p,q)\in\seisuu_{>0}^2$
and $\vecphi=(\phi_1,\phi_2)\in S^1\times S^1$.
We use an induction on $p+q$.
We may assume that $p\geq q$.

We consider the isomorphism
$\Psi_{\eta_1}^{-1}:S^1\times S^1
\lrarr S^1\times S^1$.
If $\eta_1=+$,
$\Psi_{\eta_1}^{-1}(\gbiga)$
is an open subset of 
$H\bigl(p-q,q,(0,\phi_2-pq^{-1}\phi_1)\bigr)$.
Because $p-q+q<p+q$,
we can apply the assumption of the induction
to $\Psi_{\eta_1}^{-1}(\gbiga)$.

Let us consider the case $\eta_1=-$.
Then, $\Psi_{\eta_1}^{-1}(\gbiga)$
is an open subset of
$H\bigl(
p,q-p,(\phi_1-qp^{-1}\phi_2,0)
 \bigr)$.
Because $q-p$ is non-positive,
we obtain the claim in this case
by Lemma \ref{lem;16.4.26.2}.
\hfill\qed

\subsection{Proof of Proposition \ref{prop;16.4.26.11}}
\label{subsection;18.11.12.1}

We set $\realbar_{\geq 0}:=
 \real_{\geq 0}\cup\{\infty\}$.
Set $Y:=\realbar_{\geq 0}^2\times (S^1)^2$.
Let $\gbigg\subset Y$ be a closed subanalytic subset
with $\dim_{\real}\gbigg=3$.

\begin{prop}
\label{prop;16.4.26.12}
There exist a non-negative integer $\ell_0$
and 
a closed subanalytic subset $\gbigz\subset (S^1)^2$
with $\dim_{\real} \gbigz\leq 1$
such that the following holds
for any 
$\ell\geq \ell_0$
and 
$\veceta\in\{+,-\}^{\ell}$:
\begin{itemize}
\item
Let $\gbigg_{\veceta}$ denote the strict transform
of $\gbigg$ with respect to
$\psi_{\veceta,\realbar_{\geq 0}}\times\id:
\Bl_{\veceta}\realbar_{\geq 0}^2\times(S^1)^2
\lrarr
\realbar_{\geq 0}^2\times (S^1)^2$,
i.e.,
$\gbigg_{\veceta}$ denote the closure of 
$\gbigg\cap\bigl(\real_{>0}^2\times(S^1)^2\bigr)$
in $\Bl_{\veceta}\realbar^2_{\geq 0}\times (S^1)^2$.
Then,
$\gbigg_{\veceta}$ does not intersect with
$\{P_{\veceta}\}\times
\bigl(
 (S^1)^2\setminus \gbigz
\bigr)$.
\end{itemize}
\end{prop}
\pf
Let $p:Y\lrarr \realbar_{\geq 0}\times (S^1)^2$
be the projection
$p(t_1,t_2,\theta_1,\theta_2)
=(t_2,\theta_1,\theta_2)$.

\begin{lem}
There exists a closed subanalytic subset
$Z_0\subset\realbar_{\geq 0}\times (S^1)^2$
with $\dim Z_0\leq 2$ such that
the following holds:
\begin{itemize}
\item
 $p^{-1}(Z_0)$ contains the singular locus of $\gbigg$.
\item
 For any $Q\in \gbigg\setminus p^{-1}(Z_0)$,
 the derivative of $p_{|\gbigg}$ at $Q$ 
 is an isomorphism.
\end{itemize}
\end{lem}
\pf
There exists a proper morphism of real analytic manifolds
$\rho:\gbiggtilde\lrarr 
 \realbar_{\geq 0}^2\times (S^1)^2$
such that
(i) $\rho(\gbiggtilde)=\gbigg$,
(ii) $\dim_{\real} \gbiggtilde=3$.
Let $Z'_{0}\subset \gbiggtilde$
be the set of critical points of $p\circ\rho$.
Note 
$\dim (p\circ\rho)(Z'_0)\leq 2$.
Let $Z_0$ be the union of $p(\Sing(\gbigg))$
and $(p\circ\rho)(Z'_0)$.
Then, the claim of the lemma is clear.
\hfill\qed

\vspace{.1in}
By enlarging $Z_0$,
we may assume that 
each connected component 
$\nbigc$ 
of $\bigl(\real_{\geq 0}\times (S^1)^2 \bigr)\setminus Z_0$
is simply connected.
The map $\gbigg\cap p^{-1}(\nbigc)\lrarr \nbigc$
is proper and a local diffeomorphism.
There exists the decomposition
$\gbigg\cap p^{-1}(\nbigc)=\coprod\nbigc_i$
such that the induced maps
$\nbigc_i\lrarr \nbigc$ are isomorphisms.
We obtain the subanalytic functions $f^{\nbigc}_1,\ldots,f^{\nbigc}_m$
on $(\nbigc,\realbar_{\geq 0}\times (S^1)^2)$
such that 
$\gbigg\cap p^{-1}(\nbigc)$ is the union 
of the graph 
$\Gamma(f^{\nbigc}_i)$ of $f^{\nbigc}_i$.
We may assume that $f^{\nbigc}_i$
are real analytic on $\nbigc$.

We set $V_{\nbigc}:=\overline{\nbigc}\cap(\{0\}\times (S^1)^2)$.
Let $V_{\nbigc}^{\circ}$ denote the interior part of 
$V_{\nbigc}\subset (S^1)^2$.
Let $\overline{V_{\nbigc}^{\circ}}$ be the closure of
$V_{\nbigc}^{\circ}$ in $(S^1)^2$.

\begin{lem}
There exists a closed subanalytic subset
$Z_{1,\nbigc}\subset \overline{V_{\nbigc}^{\circ}}$
with $\dim_{\real}Z_{1,\nbigc}\leq 1$
such that the following holds.
\begin{itemize}
\item
The functions $f_j^{\nbigc}$ are
 ramified real analytic around 
 any point of $V_{\nbigc}^{\circ}\setminus Z_{1,\nbigc}$.
\item
For each $Q\in V_{\nbigc}^{\circ}\setminus Z_{1,\nbigc}$,
there exist rational numbers $b_j(Q)$ and
 positive rational numbers $c_j(Q)$
 such that
\[
 f^{\nbigc}_j=t^{b_j(Q)}\cdot
 \bigl(
 g^Q_{0,j}(\theta_1,\theta_2)
+t^{c_j(Q)}g^Q_{1,j}(\theta_1,\theta_2,t)
 \bigr).
\]
Here,
$g^Q_{0,j}$
are nowhere vanishing real analytic functions
on a neighbourhood of $Q$ in $V_{\nbigc}^{\circ}$,
and $g^Q_{1,j}$
are bounded ramified real analytic functions 
on a neighbourhood of $Q$
in $\real_{\geq 0}\times (S^1)^2$.
The numbers $b_j(Q)$ and $c_j(Q)$
are constants on the connected components of
$V_{\nbigc}^{\circ}\setminus Z_1$.
\end{itemize}
\end{lem}
\pf 
By Lemma \ref{lem;16.7.20.1},
there exists a closed subanalytic subset
$Z_1'\subset \overline{V_{\nbigc}^{\circ}}$
for which the first claim holds.
Around $Q$,
there exists the expression
\[
 f^{\nbigc}_j=
 t^{b_j(Q)}\cdot
 \bigl(g_{0,j}(\theta_1,\theta_2)
+t^{c_j(Q)}g_{1,j}(\theta_1,\theta_2,t)\bigr).
\]
By enlarging $Z_1'$,
we may assume that $g_{0,j}(\theta_1,\theta_2)$
is nowhere vanishing,
and the second claim follows.
\hfill\qed

\vspace{.1in}

We set
$Z_1:=\bigcup_{\nbigc} Z_{1,\nbigc}$.
It is a closed subanalytic subset
in $\{0\}\times (S^1)^2\subset
 \real_{\geq 0}\times (S^1)^2$.
Let $\nbigb$ denote the set of 
rational numbers of the form
$b_j(Q)$ for some $Q\in V_{\nbigc}$.
It is a finite set.

For any $\veceta\in\{+,-\}^{\ell}$,
we have constructed the space
$\Bl_{\veceta}\real^2_{\geq 0}$
with the base point $P_{\veceta}$
and the morphism
$\psi_{\veceta,\realbar_{\geq 0}}:
 \Bl_{\veceta}\realbar^2_{\geq 0}
\lrarr
 \realbar^2_{\geq 0}$.
The following lemma is well known and easy to see.
\begin{lem}
\label{lem;16.4.26.10}
Let $b$ be any positive rational number.
Let $p$ and $q$ be the positive integers
such that $b=p/q$ and $\gcd(p,q)=1$.
For any positive numbers $\alpha$ and $\beta$,
we define $\rho_{\alpha,\beta}(t):\II_{\epsilon}\lrarr 
 \realbar_{\geq 0}^2$
by $\rho_{\alpha,\beta}(t)=(\alpha t^p,\beta t^q)$.
Then, there exists $\ell_1$ such that
the following holds
for any 
$\ell\geq \ell_1$
and $\veceta\in\{+,-\}^{\ell}$:
\begin{itemize}
\item
 The strict transform of
 $\rho_{\alpha,\beta}$
 with respect to
 $\psi_{\veceta,\real_{\geq 0}}$
 does not intersect with
 $P_{\veceta}$.
\hfill\qed
\end{itemize}
\end{lem}

Let us finish the proof of Proposition \ref{prop;16.4.26.12}.
Let $Q$ be any point of 
$\bigl(\{0\}\times (S^1)^2\bigr)\setminus Z_1$.
Locally around $Q$,
$\gbigg$ is described as the union of
the closure 
$\overline{\Gamma(f^{\nbigc}_j)}$
of the graph $\Gamma(f^{\nbigc}_j)$.
If $b_j(Q)\leq 0$,
then 
$\overline{\Gamma(f^{\nbigc}_j)}$
does not intersect with 
$(0,0)\times (S^1)^2$.
By applying Lemma \ref{lem;16.4.26.10}
to each $b_j(Q)>0$,
we obtain the claim of Proposition \ref{prop;16.4.26.12}.
\hfill\qed

\vspace{.1in}

Note that
$\realbar_{\geq 0}^2\times (S^1)^2$
is a compactification of $\widetilde{\cnum^2}(H)$.
Proposition \ref{prop;16.4.26.11}
follows from Proposition \ref{prop;16.4.26.12}
and the expression of 
$\psi_{\veceta}$ in \S\ref{subsection;16.4.26.3}.
\hfill\qed

\subsection{Proof of Proposition \ref{prop;16.4.27.10}}
\label{subsection;16.5.31.3}

\subsubsection{A refined statement}

Let $q_i:W_i\lrarr M$ denote the morphisms
obtained as the composite of
the induced morphisms $W_i\lrarr W_0$
and the projection $W_0\lrarr M$.

\begin{prop}
\label{prop;16.4.27.11}
There exists a closed subanalytic subset $Z\subset M$
with $\dim Z=1$
satisfying the following conditions:
\begin{itemize}
\item
 Each connected component of
 $M\setminus Z$ is simply connected.
\item
 The restriction of $q_i$
 to $C_i\setminus q_i^{-1}(Z)$
 is proper and a local diffeomorphism to $M\setminus Z$.
\item
 Let $\nbigd$ be any connected component
 of $M\setminus Z$.
 For each $i\leq k$,
 there exist a non-negative number 
 $\ell(\nbigd,i)\leq \ell$,
a closed subanalytic subset
 $L_{\nbigd,i}\subset q_i^{-1}(\nbigd)$
 with $\dim L_{\nbigd,i}\leq 3$,
 and an open immersion 
 $\iota_{\nbigd,i}:
 q_i^{-1}(\nbigd)\setminus L_{\nbigd,i}
\lrarr
 \Bl_{\veceta_{\ell(\nbigd,i)}}\real^2\times \nbigd$.
\item
 $L_{\nbigd,0}=\emptyset$,
 and
 $L_{\nbigd,i}\supset
 \phi_i^{-1}(L_{\nbigd,i-1})$.
\item
The following diagram is commutative:
\[
 \begin{CD}
 q_{i+1}^{-1}(\nbigd)\setminus L_{\nbigd,i+1}
 @>{\iota_{\nbigd,i+1}}>>
 \Bl_{\veceta_{\ell(\nbigd,i+1)}}\real^2\times\nbigd\\
 @VV{\phi_{i+1}}V @VVV \\
 q_{i}^{-1}(\nbigd)\setminus L_{\nbigd,i}
 @>{\iota_{\nbigd,i}}>>
 \Bl_{\veceta_{\ell(\nbigd,i)}}\real^2\times\nbigd\\
 \end{CD}
\]
\item
 Set $U_{i|\nbigd}:=U_i\cap q_i^{-1}(\nbigd)$.
Then, either one of the following conditions is satisfied.
\begin{itemize}
\item
$\iota_{\nbigd,i}(U_{i|\nbigd}\setminus L_{\nbigd,i})
\supset 
 P_{\veceta_{\ell(\nbigd,i)}}\times\nbigd$.
\item
$\iota_{\nbigd,i}(U_{i|\nbigd}\setminus L_{\nbigd,i})
\cap
 \bigl(
 P_{\veceta_{\ell(\nbigd,i)}}\times\nbigd
 \bigr)
=\emptyset$.
\end{itemize}
\item
If $\iota_{\nbigd,i}(U_{i|\nbigd}\setminus L_{\nbigd,i})
 \supset
 P_{\veceta_{\ell(\nbigd,i)}}\times\nbigd$,
then 
 $\iota_{\nbigd,i+1}(q_{i+1}^{-1}(\nbigd)\setminus L_{\nbigd,i+1})
\supset
 P_{\veceta_{\ell(\nbigd,i+1)}}\times\nbigd$.
\end{itemize}
\end{prop}
\pf
First, there exists a closed subanalytic subset 
$Z\subset M$ 
with $\dim_{\real} Z\leq 1$
satisfying the following conditions.
\begin{itemize}
\item
 $M\setminus Z$ is simply connected.
\item
 The restriction of $q_i$
 to $C_i\setminus q_i^{-1}(Z)$
 are proper and local diffeomorphisms
 to $M\setminus Z$ for any $i$.
\end{itemize}
In the following, we shall enlarge $Z$.
We note that $\dim C_i\leq 2$.

Let $\nbigd$ be a connected component of
$M\setminus Z$.
Let $C_{0|\nbigd}:=C_0\cap q_0^{-1}(\nbigd)$.
Let 
$C_{0|\nbigd}
=\coprod_{j\in\Lambda(0,\nbigd)}
 C_{0,\nbigd,j}$
denote the decomposition
into connected components.
After enlarging $Z$,
we may assume either one of the following
for each $j\in\Lambda(0,\nbigd)$;
(i) $C_{0,\nbigd,j}\cap((0,0)\times\nbigd)=\emptyset$,
(ii) $C_{0,\nbigd,j}= (0,0)\times\nbigd$.
We may also assume that 
either one of the following holds;
(i) $U_0\cap((0,0)\times\nbigd)=\emptyset$,
(ii) $U_0\supset (0,0)\times\nbigd$.
We set
$\Lambda^{\circ}(0,\nbigd):=
 \bigl\{j\,\big|\,
 C_{0,\nbigd,j}=(0,0)\times\nbigd\bigr\}$
and 
$\Lambda^{\bot}(0,\nbigd):=
 \Lambda(0,\nbigd)\setminus
 \Lambda^{\circ}(0,\nbigd)$.
We set
$C^{\bot}_{0,\nbigd}:=
 \coprod_{j\in\Lambda^{\bot}(0,\nbigd)}C_{0,\nbigd,j}$.

We set $\ell(1,\nbigd):=1$ if
$\Lambda^{\circ}(0,\nbigd)\neq\emptyset$,
or 
$\ell(1,\nbigd):=0$ if 
$\Lambda^{\circ}(0,\nbigd)=\emptyset$.
We set
$L_{1,\nbigd}:=
 \phi_1^{-1}(C^{\bot}_{0,\nbigd})
\subset
 q_1^{-1}(\nbigd)$.
Then, by the construction,
there exists the natural open immersion
$\iota_{1,\nbigd}:
 q_1^{-1}(\nbigd)\setminus L_{1,\nbigd}
\lrarr
 \Bl_{\veceta_{\ell(1,\nbigd)}}
 \times\nbigd$.
If 
$(0,0)\times\nbigd\subset
 U_{0|\nbigd}\times\nbigd$,
then we obtain
$P_{\veceta_{\ell(1,\nbigd)}}\times\nbigd
\subset
 \iota_{1,\nbigd}\bigl(
 q_1^{-1}(\nbigd)\setminus L_{1,\nbigd}
 \bigr)$.

Let $C_{1|\nbigd}
=\coprod_{j\in\Lambda(1,\nbigd)}
 C_{1,\nbigd,j}$
denote the decomposition
into connected components.
After enlarging $Z$,
we may assume that
either one of the following holds for each
$j\in \Lambda(1,\nbigd)$;
(i)
 $C_{1,\nbigd,j}\cap
 L_{1,\nbigd}=\emptyset$,
(ii)
 $C_{1,\nbigd,j}\subset L_{1,\nbigd}$.
We may assume that
either one of the following holds for each
$j\in\Lambda(1,\nbigd)$;
(i)
$\iota_{1,\nbigd}(C_{1,\nbigd,j}\setminus L_{1,\nbigd})
\cap
 \bigl(
 P_{\veceta_{\ell(1,\nbigd)}}
 \times\nbigd
 \bigr)=\emptyset$,
(ii)
$\iota_{1,\nbigd}(C_{1,\nbigd,j}\setminus L_{1,\nbigd})
=P_{\veceta_{\ell(1,\nbigd)}}
 \times\nbigd$.
We may also assume that
either one of the following holds:
(i)
$\iota_{1\nbigd}(U_{1|\nbigd})
\supset
 \bigl(
 P_{\veceta_{\ell(1,\nbigd)}}\times \nbigd
 \bigr)$,
(ii)
$\iota_{1\nbigd}(U_{1|\nbigd})
\cap
 \bigl(
 P_{\veceta_{\ell(1,\nbigd)}}\times \nbigd
 \bigr)=\emptyset$.

We set 
$\Lambda^{\circ}(1,\nbigd):=
\bigl\{
 j\in\Lambda(1,\nbigd)\,\big|\,
 C_{1,\nbigd,j}\cap L_{1,\nbigd}=\emptyset,\,
 \iota_{1,\nbigd}(C_{1,\nbigd,j})=
 P_{\veceta_{\ell(1,\nbigd)}}\times\nbigd
\bigr\}$
and 
$\Lambda^{\bot}(1,\nbigd):=
 \Lambda(1,\nbigd)\setminus
 \Lambda^{\circ}(1,\nbigd)$.
We set 
$C^{\bot}_{1,\nbigd}:=
 L_{1,\nbigd}\cup
 \coprod_{j\in\Lambda^{\bot}(1,\nbigd)}
 C_{1,\nbigd,j}$.
We set
\[
 \ell(2,\nbigd):=
 \left\{
\begin{array}{ll}
 \ell(1,\nbigd)+1
 & (\Lambda^{\circ}(1,\nbigd)\neq\emptyset),
 \\
 \ell(1,\nbigd)
 & (\Lambda^{\circ}(1,\nbigd)=\emptyset).
\end{array}
 \right.
\]
We set
$L_{2,\nbigd}:=
 \phi_2^{-1}(C_{1,\nbigd}^{\bot})$.
By the construction,
there exists
the naturally defined open immersion
$q_2^{-1}(\nbigd)\setminus
 L_{2,\nbigd}
\lrarr
 \Bl_{\veceta_{\ell(2,\nbigd)}}\real^2_{\geq 0}\times \nbigd$.
If $P_{\veceta_{\ell(1,\nbigd)}}\times\nbigd
\subset
 U_{1|\nbigd}\setminus L_{1,\nbigd}$,
then we obtain
$P_{\veceta_{\ell(2,\nbigd)}}\times\nbigd
\subset
 q_2^{-1}(\nbigd)\setminus L_{2,\nbigd}$
by construction.

Inductively, we construct the desired data.
Suppose that we have already obtained
$\ell(\nbigd,i)$,
$L_{\nbigd,i}\subset q_i^{-1}(\nbigd)$,
and an open immersion
$\iota_{\nbigd,i}:
 q_{i}^{-1}(\nbigd)\setminus L_{\nbigd,i}
\lrarr
 \Bl_{\veceta_{\ell(\nbigd,i)}}\real^2\times\nbigd$.
Let 
$C_{i|\nbigd}=
 \coprod_{j\in\Lambda(i,\nbigd)}
 \nbigc_{i,\nbigd,j}$
denote the decomposition
into the connected components.
After enlarging $Z$,
we may assume that
either one of the following holds
for each $j\in \Lambda(i,\nbigd)$;
(i) $\nbigc_{i,\nbigd,j}\cap L_{i,\nbigd}=\emptyset$,
or 
(ii) $\nbigc_{i,\nbigd,j}\subset L_{i,\nbigd}$.
We may assume that either 
one of the following for each $j\in\Lambda(i,\nbigd)$;
(i) $\iota_{i,\nbigd}(C_{i,\nbigd,j}\setminus L_{i,\nbigd})
\cap
 (P_{\veceta_{\ell}(i,\nbigd)}\times\nbigd)
=\emptyset$,
(ii) 
$\iota_{i,\nbigd}(C_{i,\nbigd,j}\setminus L_{i,\nbigd})
=P_{\veceta_{\ell}(i,\nbigd)}\times\nbigd$.
We may also assume that
either one of the following holds:
(i) 
$\iota_{i,\nbigd}(U_{i|\nbigd})
\supset
 P_{\veceta_{\ell(i,\nbigd)}}\times\nbigd$,
(ii)
$\iota_{i,\nbigd}(U_{i|\nbigd})
\cap
 (P_{\veceta_{\ell(i,\nbigd)}}\times\nbigd)
=\emptyset$.

We set 
$\Lambda^{\circ}(i,\nbigd):=
\bigl\{
 j\in\Lambda(i,\nbigd)\,\big|\,
 \iota_{i,\nbigd}(C_{i,\nbigd,j})
=P_{\veceta(i,\nbigd)}\times\nbigd
\bigr\}$
and
$\Lambda^{\bot}(i,\nbigd):=
 \Lambda(i,\nbigd)\setminus
 \Lambda^{\circ}(i,\nbigd)$.
We set
$C_{2,\nbigd}^{\bot}:=
 L_{i,\nbigd}\cup
 \coprod_{j\in\Lambda^{\bot}(i,\nbigd)}
 C_{i,\nbigd,j}$.
We define
\[
 \ell(i+1,\nbigd):=
 \left\{
\begin{array}{ll}
 \ell(i,\nbigd)+1
 &
 (\Lambda^{\circ}(i,\nbigd)\neq\emptyset),\\
 \ell(i,\nbigd)
 & (\Lambda^{\circ}(i,\nbigd)=\emptyset).
\end{array}
 \right.
\]
We set
$L_{i+1,\nbigd}:=\phi_{i+1}^{-1}(C_{i,\nbigd}^{\bot})$.
By the construction,
there exists the naturally defined open immersion
$q_{i+1}^{-1}(\nbigd)\setminus
 L_{i+1,\nbigd}
\lrarr
 \Bl_{\veceta_{\ell(i+1,\nbigd)}}\real^2_{\geq 0}\times \nbigd$.
If $P_{\veceta_{\ell(i,\nbigd)}}\times\nbigd
\subset
 U_{i|\nbigd}\setminus L_{i,\nbigd}$,
then we obtain
$P_{\veceta_{\ell(i+1,\nbigd)}}\times\nbigd
\subset
 q_{i+1}^{-1}(\nbigd)\setminus L_{i+1,\nbigd}$
by construction.
In this way,
the inductive construction can proceed.
\hfill\qed

\subsubsection{Proof of Proposition \ref{prop;16.4.27.10}}

There exists the connected component $\nbigd$ of $M\setminus Z$
such that $Q\in\nbigd$.
Because $(O,Q)$ is contained in the image of $\phi$,
$(O,Q)$ is contained in $U_{0|\nbigd}$.
We set $\gamma_0:=(\psi_{\veceta}\times\id)\circ\gamma$.
There exists a lift $\gamma_1$ of $\gamma_0$ to $W_1$,
and 
$\iota_{\nbigd,1}(\gamma_1(0))=
 (P_{\veceta_{\ell(\nbigd,1)}},Q)
 \in
 P_{\veceta_{\ell(\nbigd,1)}}\times\nbigd$.
Because the image of $\gamma_1$
is contained in the image of 
$W\lrarr W_1$,
we obtain
$\bigl(
 P_{\veceta_{\ell(\nbigd,1)}}\times \nbigd
 \bigr)
\cap
 \iota_{\nbigd,1}(U_{1|\nbigd}\setminus L_{\nbigd,1})
\neq\emptyset$.
Hence,
we obtain
$\bigl(
 P_{\veceta_{\ell(\nbigd,1)}}\times \nbigd
 \bigr)
\subset
 \iota_{\nbigd,1}(U_{1|\nbigd}\setminus L_{\nbigd,1})$.
Inductively,
we obtain
$P_{\veceta_{\ell(\nbigd,i)}}\times\nbigd
\subset
 \iota_{\nbigd,i}(U_{i|\nbigd}\setminus L_{\nbigd,i})$.
In particular,
$\iota_{\nbigd,k}\bigl(
 q_k^{-1}(\nbigd)\setminus L_{\nbigd,k}
 \bigr)$
is an open neighbourhood of
$P_{\veceta_{\ell(\nbigd,k)}}\times\nbigd$.
Then, the claim of Proposition \ref{prop;16.4.27.10}
follows.
\hfill\qed

\subsection{Proof of Proposition \ref{prop;16.4.27.20}}
\label{subsection;18.11.12.3}

We naturally regard
$\widetilde{\cnum^2}(H)
=(\real_{\geq 0}\times S^1)^2$
which is naturally a closed subanalytic subset
of the real analytic manifold
$(\real\times S^1)^2$.
There exists  the rectilinearization
$(W_{\alpha},\phi_{\alpha})$ $(\alpha\in\Lambda)$
for the subanalytic functions $\{f_i\}$.
There exist compact subsets $K_{\alpha}\subset W_{\alpha}$
such that 
$\bigcup_{\alpha}\phi_{\alpha}(K_{\alpha})$
contains a neighbourhood of
$\widetilde{\cnum^2}(H)$
in $(\real\times S^1)^2$.
Each $\phi_{\alpha}:W_{\alpha}\lrarr (\real\times S^1)^2$
is factorized into local real blowings up 
as follows:
\[
 W_{\alpha}
=W_{\alpha,k(\alpha)}
\stackrel{\phi^{(k(\alpha))}_{\alpha}}{\lrarr}
 W_{\alpha,k(\alpha)-1}
\stackrel{\phi^{(k(\alpha)-1)}_{\alpha}}{\lrarr}
\cdots
\lrarr
W_{\alpha,2}
\stackrel{\phi_{\alpha}^{(2)}}{\lrarr}
 W_{\alpha,1}
\stackrel{\phi_{\alpha}^{(1)}}{\lrarr}
W_{\alpha,0}
=(\real\times S^1)^2
=\real^2\times\varpi^{-1}(O).
\]
Let $Z^{(1)}_{\alpha}\subset\varpi^{-1}(O)$
be as in Proposition \ref{prop;16.4.27.10}.
Let $Z^{(2)}_{\alpha}=\bigcup \phi_{\alpha}(\gbigq)$,
where $\gbigq$ runs through the quadrants of $W_{\alpha}$
such that $\dim\phi_{\alpha}(\gbigq)\leq 1$.
We set 
$Z:=
 \overline{W^{\circ}}\cup
 \bigcup_{\alpha\in\Lambda}
 (Z^{(1)}_{\alpha}\cup Z^{(2)}_{\alpha})$.
Let $\Crit(\phi_{\alpha})\subset 
 (\real\times S^1)^2$
denote the set of the critical values of $\phi_{\alpha}$.

Let $\ell_0$ be any integer 
strictly larger than $k(\alpha)$
for any $\alpha$
such that $\phi_{\alpha}(K_{\alpha})\cap 
\varpi^{-1}(O)\neq\emptyset$.
Let $\ell$ be any integer larger than $\ell_0$,
and let $\veceta$ be any element of
$\{+,-\}^{\ell}$.
Let $Q$ be any point of
$\varpi_{\veceta}^{-1}(P_{\veceta})
 \setminus
 \Psi_{\veceta}^{-1}(Z)$.
There exists a path 
$\gamma:\II\lrarr 
 \widetilde{\Bl_{\veceta}\cnum^2}(H_{\veceta})$
such that
 $\gamma(0)=Q$
and
\[
 \widetilde{\psi_{\veceta}}\circ
 \gamma(t)\in
 \bigl(
 \real_{>0}\times S^1
 \bigr)^2
 \setminus
\bigcup_{\alpha\in\Lambda}
 \Crit(\phi_{\alpha})
=
 \bigl(
 \cnum^2\setminus H
 \bigr)
 \setminus
\bigcup_{\alpha\in\Lambda}
 \Crit(\phi_{\alpha})
\quad (t>0).
\]
Because 
$\bigcup_{\alpha\in\Lambda}
 \phi_{\alpha}(K_{\alpha})$ contains
a neighbourhood of
$\widetilde{\cnum^2}(H)$,
there exist $\alpha_0\in\Lambda$
and $\epsilon>0$
such that 
$\widetilde{\psi_{\veceta}}\circ\gamma(\II_{\epsilon})
\subset
 \phi_{\alpha_0}(K_{\alpha_0})$.
There exists
$\gamma':\II_{\epsilon}\lrarr 
 W_{\alpha_0}$
such that
$\phi_{\alpha_0}\circ\gamma'(t)
=\widetilde{\psi_{\veceta}}\circ\gamma(t)$
for any $t\in\II_{\epsilon}$.
There exist a small neighbourhood
$W_{\alpha_0,\gamma'(0)}$ of $\gamma'(0)$
in $W_{\alpha_0}$,
an integer $i(Q)\leq \ell_0$,
and an embedding
$\iota_0:W_{\alpha_0,\gamma'(0)}\lrarr
 \Bl_{\veceta_{i(Q)}}\real^2\times (S^1)^2$
as in Proposition \ref{prop;16.4.27.10}.
In particular,
$\phi_{\alpha_0|W_{\alpha_0,\gamma'(0)}}
=(\psi_{\veceta_{i(Q)}}\times\id)\circ\iota_0$
holds.
Note that
$\nbigv:=\phi_{\alpha_0}^{-1}
 \bigl(
 \cnumtilde^2(H)\setminus(\gbigh\cup\varpi^{-1}(H))
 \bigr)$
is rectilinearized
with respect to the coordinate system
of $W_{\alpha_0}$.
Moreover, $\gamma'(0)$ is not contained
in any quadrant whose dimension is strictly smaller than $2$.
Hence,
$\iota_{0}\bigl(
 W_{\alpha_0,\gamma'(0)}
\cap\overline{\nbigv}
 \bigr)$
is equal to
$\iota_0\bigl(W_{\alpha_0,\gamma'(0)}\bigr)
\cap
 \Bigl(
 \Bl_{\veceta_{i(Q)}}\real^2_{\geq 0}
 \times (S^1)^2
\Bigr)$,
and the pull back
$(\iota_0^{-1})^{\ast}
 \phi_{\alpha_0}^{\ast}(f_i)$
is expressed as 
the power series of
$(u_{\veceta_{i(Q),\real}}^{1/\rho},v_{\veceta_{i(Q)},\real}^{1/\rho},
 \theta_1,\theta_2)$
for a positive integer $\rho$.
Here,
$(\theta_1,\theta_2)$ is a natural coordinate system
of
$(S^1)^2$,
and 
$(u_{\veceta_{i(Q),\real}}^{1/\rho},v_{\veceta_{i(Q)},\real}^{1/\rho})$
is as in \S\ref{subsection;16.4.26.3}.

Let $\nbigu_Q$ be a small neighbourhood of $Q$
in $\widetilde{\Bl_{\veceta}\cnum^2}(H_{\veceta})$,
and then there exists the factorization of
$\psitilde_{\veceta|\nbigu_Q}$
as follows:
\[
\begin{CD}
\nbigu_Q
@>{\Phi}>>
\iota_0(W_{\alpha_0,\gamma'(0)}\cap\nbigv)
@>{\phi_{\alpha_0}\circ\iota_0^{-1}}>>
 \real^2\times (S^1)^2.
\end{CD}
\]
Here,
$\Phi$ is induced by
$\psi_{\veceta_{i(Q)},\veceta}:
 \Bl_{\veceta}\real^2
\lrarr
 \Bl_{\veceta_{i(Q)}}\real^2$
and 
$\Psi_{\veceta}$.
It is easy to see that the pull back  of 
$(\iota_0^{-1})^{\ast}
 \phi_{\alpha_0}^{\ast}(f_i)$
via $\Phi$ is ramified real analytic.
Thus, we obtain the claim of Proposition 
\ref{prop;16.4.27.20}.
\hfill\qed

\section{Sequences of complex blowings up of mixed type}
\label{section;16.9.4.3}

In this section,
we study sequences of complex blowings up
of mixed type.
As in \S\ref{section;16.9.4.1},
we explain the following;
for any subanalytic function on the original surface,
if the length of the sequence is large enough,
the pull back of the function
has nice property around the centers of the blowings up.
For that purpose,
we study comparison of
sequences of complex blowings up
and sequences of real blowings up.

In \S\ref{subsection;16.6.1.40},
as a preliminary,
we introduce the families of curves and the limit curve
associated to a sequence of blowings up of mixed type.
The main results of this section are
collected in \S\ref{subsection;18.11.15.40},
which will be proved in the rest of this section.

\subsection{The induced families of curves and limit curves}
\label{subsection;16.6.1.40}

\subsubsection{The induced family of curves}
\label{subsection;16.5.31.1}

Let
$\vecY=(\veceta_1,\omega_1,\ldots,\veceta_k,\omega_k)
 \in \coprod_{\ell\geq 0} \gbigp^{\ell}$.
Here,
$\veceta_j\in \{+,-\}^{\ell(j)}$
and $\omega_j\in\cnum^{\ast}$.
For any $m\leq k$,
we set
$\vecY_m:=(\veceta_1,\omega_1,\ldots,\veceta_m,\omega_m)$.
We obtain the morphism
$\psi_{\vecY_m}:
 \Bl_{\vecY_m}\cnum^2
 \lrarr\cnum^2$,
the point
$P_{\vecY_m}\in\Bl_{\vecY_m}\cnum^2$
and the coordinate neighbourhood
$(U_{\vecY_m},u_{\vecY_m},v_{\vecY_m})$.
Instead of an explicit description of
$\psi_{\vecY_m}:
 U_{\vecY_m}\lrarr \cnum^2$
with respect to the coordinate systems
$(x_1,y_1):=(u_{\vecY_m},v_{\vecY_m})$
and $(x,y)$,
we study a description of
an induced family of curves.

\vspace{.1in}

For any small positive number $\epsilon>0$,
we obtain the family of holomorphic curves
\[
 \varphi_{m}:
(\cnum\setminus\{-\omega_m\})\times
 \Delta_{\epsilon}\lrarr U_{\vecY_m}
\]
by $\varphi_{m}(a,\zeta):=(a,\zeta)$.
As the composite with $\psi_{\vecY_m}$,
we obtain a family of holomorphic curves
\[
 \psi_{\vecY_m}\circ\varphi_{m}:
 (\cnum\setminus\{-\omega_m\})\times
 \Delta_{\epsilon}\lrarr \cnum^2.
\]
Let $(\psi_{\vecY_m}\circ\varphi_{m})_x(a,\zeta)$
and $(\psi_{\vecY_m}\circ\varphi_{m})_y(a,\zeta)$
denote the $x$-component and $y$-component
of $\psi_{\vecY_m}\circ\varphi_{m}(a,\zeta)$:
\[
 \psi_{\vecY_m}\circ\varphi_{m}(a,\zeta)
=\bigl(
 (\psi_{\vecY_m}\circ\varphi_{m})_x(a,\zeta),
 (\psi_{\vecY_m}\circ\varphi_{m})_y(a,\zeta)
 \bigr).
\]

\subsubsection{Change of parametrization}
\label{subsection;16.6.2.1}

To obtain a more convenient description
of the family of curves
$\psi_{\vecY_m}\circ\varphi_{m}$,
we change a parametrization of the curves.
Let $A_{\veceta_i}$ be as in \S\ref{subsection;18.11.14.10}.
For $i=1,\ldots,m$,
let $(\alpha_i,\beta_i,\gamma_i,\delta_i)\in\seisuu^4$
be determined as follows:
\[
 \left(
 \begin{array}{cc}
 \alpha_i & \beta_i\\
 \gamma_i & \delta_i
 \end{array}
 \right)
=A_{\veceta_i}\cdot
 \left(
 \begin{array}{cc}
 1 & 1 \\
 0 & 1
 \end{array}
 \right).
\]
We set $\vecdelta(i):=\prod_{j=1}^i\delta_j$.
We can check the following lemma by a direct computation.
\begin{lem}
The function
\[
G_m(a,\zeta):=
 \zeta^{-\vecdelta(m)}
 \cdot
 \prod_{i=1}^m
 \omega_i^{-\gamma_i\vecdelta(i-1)}
 \cdot
 \bigl(
 \psi_{\vecY_m}\circ\varphi_m
 \bigr)_y
\]
is holomorphic on 
$\bigl(\cnum\setminus\{-\omega_m\}\bigr)
 \times
 \Delta_{\epsilon}$,
and 
$G_m(a,0)=1$.
\hfill\qed
\end{lem}

We define the map $\cnum_u\lrarr\cnum_y$
by $u\longmapsto u^{\vecdelta(m)}$.
We denote the variable $u$ by $y^{1/\vecdelta(m)}$.
We fix 
a $\vecdelta(m)/\vecdelta(i-1)$-root of $\omega_i^{\gamma_i}$
for each $i$.
They determine a holomorphic map
\[
\bigl(
 \psi_{\vecY_m}\circ\varphi_m
\bigr)_y^{1/\vecdelta(m)}:
 \nbigu_{m,1}
\lrarr
 \cnum_{y^{1/\vecdelta(m)}}
\]
such that 
$\Bigl(
 \bigl(
 \psi_{\vecY_m}\circ\varphi_m
\bigr)_y^{1/\vecdelta(m)}
 \Bigr)^{\vecdelta(m)}
=(\psi_{\vecY_m}\circ\varphi_m)_y$,
where $\nbigu_{m,1}$
is a small neighbourhood of 
$(\cnum\setminus\{-\omega_m\})\times\{0\}$
in $(\cnum\setminus\{-\omega_m\})\times\Delta_{\epsilon}$.

We define the map
$\Phi_m:
 \nbigu_{m,1}\lrarr
 \cnum_a\times\cnum_{y^{1/\vecdelta(m)}}$
by
$\Phi_m(a,\zeta)=
 \Bigl(a,
 \bigl(
 \psi_{\vecY_m}\circ\varphi_m
 \bigr)_y^{1/\vecdelta(m)}(a,\zeta)
 \Bigr)$.
By shrinking $\nbigu_{m,1}$,
we may assume that 
$\Phi_m$ induces an isomorphism of
$\nbigu_{m,1}$
and $\nbigu_{m,2}:=\Phi_m(\nbigu_{m,1})$.
We obtain the following holomorphic function
on $\nbigu_{m,2}$:
\[
 g_{\vecY_m}:=
 \bigl(
 \psi_{\vecY_m}\circ\varphi_m
 \bigr)_x\circ
 \Phi_m^{-1}.
\]
We obtain the family of holomorphic curves
$F_{\vecY_m}:=\psi_{\vecY_m}\circ\varphi_m\circ\Phi_m^{-1}:
 \nbigu_{m,2}\lrarr
 \cnum_x\times\cnum_y$:
\[
 F_{\vecY_m}(a,y^{1/\vecdelta(m)})
=\bigl(
 g_{\vecY_m}(a,y^{1/\vecdelta(m)}),
 y
 \bigr).
\]
In this way,
we obtain a collection of families of holomorphic curves
$F_{\vecY_m}$ $(m=1,\ldots,k)$.
For the expansion
\[
 g_{\vecY_m}(a,y^{1/\vecdelta(m)})
=\sum g_{\vecY_m,\gminiy}(a) y^{\gminiy},
\]
we set $\nbigs(\vecY_m):=
 \bigl\{\gminiy\in\rnum
 \,\big|\,
 g_{\vecY_m,\gminiy}\neq 0\bigr\}$.

\subsubsection{Special case: sequences of blowings up at smooth points}

Let $\vecY=(\veceta_1,\omega_1,\ldots,\veceta_k,\omega_k)$
be as in \S\ref{subsection;16.5.31.1}.
In the case $\veceta_i\in \{-\}^{\ell(i)}$ for any $i$,
we can easily obtain a simple description of the family of curves.
In this case, it is convenient to  regard $-\in\gbigp$ as $0$,
i.e.,
we identify 
$\{-\}\sqcup\cnum^{\ast}$ with $\cnum$
in $\gbigp$.
Let $\vecomega=(\omega_1,\ldots,\omega_m)
\in
\cnum^m\subset\gbigp^m$.
The induced map
$\Bl_{\vecomega}\cnum^2\lrarr \cnum^2$
is denoted by 
$\psi_{\vecomega}$.
By a direct computation,
the map
$\psi_{\vecomega}:
(U_{\vecomega},u_{\vecomega},v_{\vecomega})
\lrarr
(\cnum^2,x,y)$
is described as follows:
\[
 \psi_{\vecomega}(u_{\vecomega},v_{\vecomega})
 =\Bigl(
 \sum_{i=1}^{m-1} \omega_iv_{\vecomega}^i
 +(\omega_m+u_{\vecomega})v_{\vecomega}^{m},\,\,
 v_{\vecomega}
 \Bigr).
\]
Hence, the family
$F_{\vecY}:\cnum\setminus\{-\omega_m\}
 \times\Delta
 \lrarr
 \cnum^2$
is described as follows:
\[
 F_{\vecY}(a,y)=\Bigl(
 \sum_{i=1}^{m-1}\omega_i y^i+(\omega_m+a)y^{m},
 y
 \Bigr).
\]

\subsubsection{Statements}

We obtain the following rational numbers
\[
 \kappa(\vecY_m):=
 \sum_{j=1}^m \frac{\beta_j}{\vecdelta(j)}.
\]
Clearly,
$\kappa(\vecY_{m'})<\kappa(\vecY_{m})$ holds
for $m'<m$.
We shall prove the following proposition
in \S\ref{subsection;16.2.19.10}--\ref{subsection;16.5.1.1}.
\begin{prop}
\label{prop;16.4.29.2}
If $\gminiy<\kappa(\vecY_m)$,
then $g_{\vecY_m,\gminiy}(a)$
are independent of $a$.
The coefficient
$g_{\vecY_m,\kappa(\vecY_m)}(a)$
is an affine function of $a$.
Moreover, for $m<m_1$, the following holds:
\begin{itemize}
\item
 $g_{\vecY_{m_1},\gminiy}=g_{\vecY_{m},\gminiy}$
 for $\gminiy<\kappa(\vecY_m)$,
and
 $g_{\vecY_{m_1},\kappa(\vecY_m)}(a)
 =g_{\vecY_m,\kappa(\vecY_m)}(0)$.
\end{itemize}
\end{prop}

For any positive rational number $\alpha$,
let $\alpha=\prod \gminip^{\ell_{\gminip}}$
be the prime factorization,
where $\gminip$ runs through prime numbers,
and $\ell_{\gminip}$ are integers.
We set $\ord_{\gminip}(\alpha):=\ell_{\gminip}$.
Under some additional assumption,
we obtain the following,
which will also be proved in
\S\ref{subsection;16.2.19.10}--\S\ref{subsection;16.5.1.1}.
\begin{prop}
\label{prop;16.2.19.12}
Suppose that 
there exists a prime $\gminip_0$
 such that 
 $\ord_{\gminip_0}(\kappa(\vecY_m))<0$ 
and
 $\ord_{\gminip_0}(\kappa(\vecY_m))
 <\ord_{\gminip_0}(\kappa(\vecY_{m'}))$
 for any $m'<m$.
Then, the following holds:
\begin{itemize}
\item
 $g_{\vecY_m,\kappa(\vecY_m)}(0)\neq 0$.
\item
For any $\gminiy\in\nbigs(\vecY_m)$ such that
 $\gminiy<\kappa(\vecY_m)$,
we obtain
$\ord_{\gminip_0}(\kappa(\vecY_m))
 <\ord_{\gminip_0}(\gminiy)$.
\end{itemize}
\end{prop}

We give a remark.
\begin{lem}
\label{lem;16.5.11.2}
 Let $N_0:=
 \min\bigl\{N\,\big|\,
 \kappa(\vecY_{i})\in \frac{1}{N}\seisuu\,\,
 (i=1,\ldots,m-1)
 \bigr\}$.
If $\kappa(\vecY_m)\in\frac{1}{N_0}\seisuu$,
then $\delta_m=1$.
\end{lem}
\pf
By the choice of $N_0$,
$a_i:=N_0\beta_i/\vecdelta(i)$
are integers for any $i\leq m-1$,
and 
$\gcd(a_1,\ldots,a_{m-1})=1$
holds.
Let us observe that
$N_0^{-1}\vecdelta(m-1)\in\seisuu$.
Indeed, if not,
there exists a prime $\gminip$ such that
$\ord_{\gminip}(N_0^{-1}\vecdelta(m-1))<0$.
Then, 
$\ord_{\gminip}(N_0^{-1}\vecdelta(i))<0$
holds for any $i\leq m-1$.
Because
$a_i N_0^{-1}\vecdelta(i)\in\seisuu$,
we obtain
$\ord_{\gminip}(a_i)>0$ for $i=1,\ldots,m-1$,
which contradicts
$\gcd(a_1,\ldots,a_{m-1})=1$.
Hence, we can conclude that
$A:=N_0^{-1}\vecdelta(m-1)$
is an integer.

Note that
$\kappa(\vecY_m)
=\kappa(\vecY_{m-1})+
 N_0^{-1}A^{-1}\delta_m^{-1}\beta_m$.
If $\kappa(\vecY_m)\in \frac{1}{N_0}\seisuu$,
then $A^{-1}\delta_m^{-1}\beta_m$ is an integer.
Because $\gcd(\delta_m,\beta_m)=1$,
we obtain $\delta_m=1$.
\hfill\qed

\subsubsection{Limit curves}
\label{subsection;16.5.31.33}

Let $\vecY\in\gbigp^{\infty}$.
There exists the limit
$\kappahat(\vecY):=\lim_{m\to\infty} \kappa(\vecY_m)$
in $\real\cup\{\infty\}$.
For any $\gminiy<\kappahat(\vecY)$,
there exists the limit
$g_{\vecY,\gminiy}:=\lim_{m\to\infty}g_{\vecY_m,\gminiy}$
by Proposition \ref{prop;16.4.29.2}.
When $\gminiy\geq \kappahat(\vecY)$,
we formally set
$g_{\vecY,\gminiy}:=0$.
We set
$g_{\vecY}:=\sum g_{\vecY,\gminiy}y^{\gminiy}$
in $\prod_{\gminiy\in\rnum_{\geq 0}} \cnum y^{\gminiy}$.
We set
$\nbigs(g_{\vecY}):=
 \{\gminiy\in\rnum_{\geq 0}\,|\,g_{\vecY,\gminiy}\neq 0\}$.
For any $\kappa<\kappahat(\vecY)$,
the set
$\{0\leq \gminiy\leq\kappa\}\cap \nbigs(g_{\vecY})$
is finite.
We call $g_{\vecY}(y)$ the limit curve.
(See also \cite[Proposition 4.1]{Favre-Jonsson}
which explains the correspondence between valuations
and such formal curves.
The two constructions of formal curves should be related
although the author has not yet checked it.)

We introduce a condition for infinite sequences $\vecY$.
\begin{itemize}
 \item 
 We say that the limit curve of $\vecY$ is convergent
 if $\nbigs(g_{\vecY})$ is contained in $\frac{1}{N}\seisuu$ 
for some $N\in\seisuu_{>0}$,
 and the power series
 $g_{\vecY}$ is convergent.
 Note that the condition implies
the existence of $i_0$ such that
 $\veceta_i=(-,\ldots,-)^{\ell(i)}$
 for any $i\geq i_0$,
 which follows from Lemma \ref{lem;16.5.11.2}.
It also implies that
$\kappahat(\vecY)=\infty$.
\end{itemize}

Later, we shall divide $\gbigp^{\infty}$
into the following three classes.
\begin{description}
\item[Class (i):]
 $\kappahat(\vecY)<\infty$.
\item[Class (ii):]
 $\kappahat(\vecY)=\infty$, but
 the limit curve of $\vecY$ is not convergent.
\item[Class (iii):]
 The limit curve of $\vecY$ is convergent.
\end{description}

\subsubsection{Inductive step}
\label{subsection;16.2.19.10}

Let us start the proof of Proposition \ref{prop;16.4.29.2}.
For $\veceta\in\{+,-\}^{\ell}$,
we set $\veceta-=(\veceta,-)\in\{+,-\}^{\ell+1}$.
Let $\alpha,\beta,\gamma,\delta$
be determined by
\[
 \left(
 \begin{array}{cc}
 \alpha &\beta \\
 \gamma &\delta
 \end{array}
 \right)
=A_{\veceta-}=A_{\veceta}
\cdot
 \left(
 \begin{array}{cc}
 1 & 1 \\ 0 & 1
 \end{array}
 \right).
\]
Note that $\beta>0$ and $\delta>0$.
We define
$\psi_{\veceta-}:
 (U_{\veceta-},u_{\veceta-},v_{\veceta-})
\lrarr
 (\cnum^2,x,y)$
by
\[
 \psi_{\veceta-}(u_{\veceta-},v_{\veceta-})
=\bigl(
 u_{\veceta-}^{\alpha}v_{\veceta-}^{\beta},
 u_{\veceta-}^{\gamma}v_{\veceta-}^{\delta}
 \bigr).
\]

\vspace{.1in}

Let $\nbigu$ be a neighbourhood of $0$.
Let $N$ be a positive integer.
Let $\kappa\in \frac{1}{N}\seisuu_{>0}$.
Let
$P:\nbigu\times\Delta_{\epsilon}
\lrarr \cnum$
be a holomorphic function
satisfying the following condition:
\begin{itemize}
\item
 For the expansion
 $P(a,\sigma)=\sum_{j\geq 0} P_{j/N}(a)\sigma^j$,
 if $j/N<\kappa$, 
 the coefficients
 $P_{j/N}(a)$ are independent of $a$.
 Moreover,
 $P_0\neq 0$.
\end{itemize}

We consider a family of curves
$\varphi:
\nbigu\times\Delta_{\epsilon}
\lrarr
 (U_{\veceta-},u_{\veceta-},v_{\veceta-})$
defined by
\[
 \varphi(a,\sigma)
=\bigl(
 P(a,\sigma),
 \sigma^N
 \bigr).
\]
We obtain the induced family of curves
$\psi_{\veceta-}\circ\varphi:
 \nbigu\times\Delta_{\epsilon}
\lrarr
 (\cnum^2,x,y)$:
\begin{equation}
\label{eq;16.4.29.1}
 \psi_{\veceta-}\circ\varphi(a,\sigma)
=\bigl(
 P(a,\sigma)^{\alpha}\sigma^{N\beta}, 
 P(a,\sigma)^{\gamma}\sigma^{N\delta}
\bigr).
\end{equation}
We choose $P_0^{\gamma/(N\delta)}\in\cnum$,
and then
$P(a,\sigma)^{\gamma/(N\delta)}$ is naturally defined.
Set $\upsilon:=P(a,\sigma)^{\gamma/(N\delta)}\sigma$.

We define the map
$\Phi:
 \nbigu\times\Delta_{\epsilon}
\lrarr
 \nbigu\times \cnum$
by
$\Phi(a,\sigma):=
 \bigl(a,\upsilon(a,\sigma)\bigr)$.
The restriction of $\Phi$
to a small neighbourhood of
$\nbigu\times\{0\}$ induces a bi-holomorphic map
into $\nbigu\times\cnum$.
Let $\nbigu'$ be a sufficiently smaller neighbourhood 
of $0$ in $\nbigu$,
and let $\epsilon'>0$ be sufficiently smaller than $\epsilon$.
Then, we obtain a holomorphic map
$Q:\nbigu\times\Delta_{\epsilon'}
\lrarr \cnum$
such that the following conditions are satisfied.
\begin{itemize}
\item
$\sigma=Q(a,\upsilon)\cdot\upsilon$,
i.e.,
$\Phi^{-1}(a,\upsilon)
=\bigl(a,Q(a,\upsilon)\cdot\upsilon\bigr)$.
\item
For the expansion
$Q(a,\upsilon)=\sum_{j\geq 0}Q_{j/N}(a)\upsilon^j$,
the coefficients
$Q_{j/N}(a)$ $(j/N<\kappa)$
are determined by $P_{\gminiy}$ $(\gminiy<\kappa)$.
In particular,
they are independent of $a$.
\item
$Q_0=P_0^{-\gamma/(N\delta)}\neq 0$.
\item
The function
$Q_{\kappa}(a)
+\gamma(N\delta)^{-1}
 P_0^{-1-\gamma\kappa/\delta-\gamma/(N\delta)}
 P_{\kappa}(a)$
is constant.
\end{itemize}

We can check the following lemma
by direct computations.

\begin{lem}
\label{lem;16.4.28.30}
Suppose that there exists a prime $\gminip_0$
such that the following conditions are satisfied.
\begin{itemize}
\item
$\ord_{\gminip_0}(\kappa)<0$.
\item
If $P_{\gminiy}\neq 0$ and $\gminiy<\kappa$,
then 
$\ord_{\gminip_0}(\gminiy)>\ord_{\gminip_0}(\kappa)$. 
\end{itemize}
Then, the following conditions also satisfied:
\begin{itemize}
\item
If $Q_{\gminiy}\neq 0$ and $\gminiy<\kappa$,
then 
$\ord_{\gminip_0}(\gminiy)
 >\ord_{\gminip_0}(\kappa)$.
\item
$Q_{\kappa}(a)=
-\gamma(N\delta)^{-1}
 P_0^{-1-\gamma\kappa/\delta-\gamma/(N\delta)}
 P_{\kappa}(a)$.
\hfill\qed
\end{itemize}
\end{lem}

\vspace{.1in}

We define the holomorphic map
$R:\nbigu\times\Delta_{\epsilon'}
\lrarr \cnum$
by 
\[
R(a,\upsilon):=
P\bigl(a,\upsilon Q(a,\upsilon)\bigr)
=\sum_{j\geq 0}R_{j/N}(a)\upsilon^j.
\]
The following lemma is clear.
\begin{lem}
$R_0=P_0\neq 0$ holds,
and $R_{j/N}(a)$ $(j/N<\kappa)$ are 
determined by $P_{\gminiy}$ $(\gminiy<\kappa)$.
In particular, they are independent of $a$.

Moreover, if the assumption of Lemma
{\rm\ref{lem;16.4.28.30}} is satisfied,
the following holds:
\begin{itemize}
\item
If $R_{\gminiy}\neq 0$ and $\gminiy<\kappa$,
then $\ord_{\gminip_0}(\gminiy)>\ord_{\gminip_0}(\kappa)$.
\item
$R_{\kappa}(a)=
 P_0^{-\gamma\kappa/\delta}P_{\kappa}(a)$.
\hfill\qed
\end{itemize}
\end{lem}

We set as follows:
\[
 S(a,\upsilon):=
 R(a,\upsilon)^{\alpha}
 Q(a,\upsilon)^{N\beta}
=\sum S_{j/N}(a)\upsilon^{j}.
\]
Then, we obtain the following description of
(\ref{eq;16.4.29.1}):
\[
 \psi_{\veceta-}\circ\varphi\circ\Phi^{-1}(a,\upsilon)
=\bigl(
 \upsilon^{N\beta}
 S(a,\upsilon),\,
 \upsilon^{N\delta}
 \bigr).
\]
The following lemma is clear.
\begin{lem}
\label{lem;16.4.29.3}
$S_0=P_0^{1/\delta}\neq 0$ holds,
and $S_{j/N}(a)$ $(j/N<\kappa)$
depend only on $\{P_{\gminiy}\,|\,\gminiy<\kappa\}$.
The function
$S_{\kappa}(a)
-\delta^{-1}
 P_0^{\alpha-1-\gamma\kappa/\delta-\gamma\beta/\delta}
 P_{\kappa}(a)$ is constant.

Moreover, if the assumption of Lemma
{\rm\ref{lem;16.4.28.30}} is satisfied,
the following holds:
\begin{itemize}
 \item
 If $S_{\gminiy}\neq 0$ and $\gminiy<\kappa$,
 then $\ord_{\gminip_0}(\gminiy)>\ord_{\gminip_0}(\kappa)$ holds.
 \item
 $S_{\kappa}(a)=
 \delta^{-1}P_0^{\alpha-1-\gamma\kappa/\delta-\gamma\beta/\delta}
 P_{\kappa}(a)$.
\hfill\qed
\end{itemize}
\end{lem}

\subsubsection{Proof of 
Proposition \ref{prop;16.4.29.2}
and 
Proposition \ref{prop;16.2.19.12}}
\label{subsection;16.5.1.1}

For $1\leq i\leq j\leq m$,
we define
$\vecY_{i,j}:=
 (\veceta_i,\omega_i,\ldots,\veceta_{j},\omega_j)$.
Note
$\vecY_{1,\ell}=\vecY_{\ell}$.
We also set 
$\vecdelta(i,j):=\vecdelta(j)/\vecdelta(i-1)$.
The following equality holds:
\[
 \kappa(\vecY_{i,j})=\frac{\beta_i}{\vecdelta(i,i)}
+\frac{\beta_{i+1}}{\vecdelta(i,i+1)}
+\cdots+
 \frac{\beta_j}{\vecdelta(i,j)}.
\]
By the construction,
$\kappa(\vecY_{k,k})=\beta_k\cdot \delta_k^{-1}$
and 
$\kappa(\vecY_{1,\ell})=\kappa(\vecY_{\ell})$
hold.
The following equalities also hold:
\[
 \kappa(\vecY_{i,j})
=\frac{\beta_i}{\delta_i}
+\frac{\kappa(\vecY_{i+1,j})}{\delta_i}.
\]
We consider the family of curves
$F_{\vecY_{i+1,j}}(a,y^{1/\vecdelta(i+1,j)}):
 \nbigu\times\Delta_{\epsilon}
\lrarr
 (\cnum^2,x,y)$:
\[
 F_{\vecY_{i+1,j}}(a,y^{1/\vecdelta(i+1,j)})
=\bigl(
 g_{\vecY_{i+1,j}}(a,y^{1/\vecdelta(i+1,j)}),y
 \bigr).
\]

There exists the isomorphism
$(U_{\vecY_{i,i}},u_{\vecY_{i,i}},v_{\vecY_{i,i}})
\simeq
 (\cnum^2,x,y)$
induced by the coordinate systems.
There exists the isomorphism
$(U_{\vecY_{i,i}},u_{\vecY_{i,i}},v_{\vecY_{i,i}})
\simeq
 (U_{\veceta_i-},u_{\veceta_i-},v_{\veceta_i-})$
defined by
$u_{\vecY_{i,i}}=u_{\veceta_i-}-\omega_i$
and 
$v_{\vecY_{i,i}}=v_{\veceta_i-}$.
Hence, with respect to the coordinate system
$(u_{\veceta_i-},v_{\veceta_i-})$,
$F_{\vecY_{i+1,j}}$ is described as follows:
\[
 F_{\vecY_{i+1,j}}(a,y^{1/\vecdelta(i+1,j)})
=\bigl(
 \omega_i+g_{\vecY_{i+1,j}}(a,y^{1/\vecdelta(i+1,j)}),y
 \bigr).
\]
Note $\vecdelta(i+1,j)\delta_i=\vecdelta(i,j)$.
We choose a $\vecdelta(i,j)$-root of $\omega_i^{\gamma_i}$.
Then, applying the construction in 
\S\ref{subsection;16.2.19.10}
to the composite of $\psi_{\veceta_i-}$
and $F_{\vecY_{i+1,j}}$,
we obtain the induced family
\[
 G_{\vecY_{i,j}}(a,y^{1/\vecdelta(i+1,j)\delta_i})
=\bigl(
 g'_{\vecY_{i,j}}(a,y^{1/\vecdelta(i+1,j)\delta_i}),
 y
 \bigr).
\]
By the construction,
it is equal to the induced family of curves
$F_{\vecY_{i,j}}$.
Then, the claim of Proposition \ref{prop;16.4.29.2}
follows from an induction and Lemma \ref{lem;16.4.29.3}.

\vspace{.1in}

We assume that
there exists a prime $\gminip_0$
such that $\ord_{\gminip_0}(\kappa(\vecY_m))<0$
and that
$\ord_{\gminip_0}(\kappa(\vecY_m))
<\ord_{\gminip_0}(\kappa(\vecY_{m'}))$
for any $m'<m$
as in Proposition \ref{prop;16.2.19.12}.

\begin{lem}
Let $m'<m$.
Then,
$\ord_{\gminip_0}(\kappa(\vecY_{m'+1,m}))<0$
and 
$\ord_{\gminip_0}(\kappa(\vecY_{m'+1,m}))
<\ord_{\gminip_0}(\kappa(\vecY_{m'+1,\ell}))$
hold for any $m'+1\leq \ell<m$.
\end{lem}
\pf
Because 
$\kappa(\vecY_{m'})
+ \kappa(\vecY_{m'+1,m})
 \prod_{p\leq m'}\delta_p^{-1}
=\kappa(\vecY_m)$,
the following holds:
\begin{equation}
\label{eq;16.3.3.1}
 \ord_{\gminip_0}\Bigl(
 \kappa(\vecY_{m'+1,m})
 \prod_{p\leq m'}\delta_p^{-1}
 \Bigr)
=\ord_{\gminip_0}(\kappa(\vecY_m)).
\end{equation}
Because 
$\kappa(\vecY_{m'})+
 \kappa(\vecY_{m'+1,\ell})
 \prod_{p\leq m'}\delta_p^{-1}
=\kappa(\vecY_{\ell})$,
the following holds:
\[
 \ord_{\gminip_0}\Bigl(
 \kappa(\vecY_{m'+1,\ell})
 \prod_{p\leq m'}\delta_p^{-1}
 \Bigr)
>\ord_{\gminip_0}(\kappa(\vecY_m)).
\]
Hence,
we obtain
\[
 \ord_{\gminip_0}\Bigl(
 \kappa(\vecY_{m'+1,m})
 \prod_{p\leq m'}\delta_p^{-1}
 \Bigr)
< \ord_{\gminip_0}\Bigl(
  \kappa(\vecY_{m'+1,\ell})
 \prod_{p\leq m'}\delta_p^{-1}
 \Bigr).
\]
It implies
$ \ord_{\gminip_0}\bigl(
 \kappa(\vecY_{m'+1,m})
 \bigr)
< \ord_{\gminip_0}\bigl(
 \kappa(\vecY_{m'+1,\ell})
 \bigr)$.

Suppose that
$\ord_{\gminip_0}(\kappa(\vecY_{m'+1,m}))\geq 0$
for some $m'$.
Because of (\ref{eq;16.3.3.1}),
we obtain
\[
 \ord_{\gminip_0}\Bigl(
 \prod_{p\leq m'}\delta_p
 \Bigr)
\geq
 \bigl|
 \ord_{\gminip_0}\bigl(
 \kappa(\vecY_m)
 \bigr)
 \bigr|.
\]
Put
$\ell_0:=
 \min\bigl\{m'
 \,\big|\,
 \ord_{\gminip_0}(\prod_{p\leq m'}\delta_p)
 \geq
  \bigl|
 \ord_{\gminip_0}\bigl(
 \kappa(\vecY_m)
 \bigr)
 \bigr|
 \bigr\}$.
Because
$\gcd(\gminip_0,\delta_{\ell_0})=\gminip_0$ holds,
we obtain
$\gcd(\gminip_0,\beta_{\ell_0})=1$.
Because 
$\kappa(\vecY_{\ell_0})
=\kappa(\vecY_{\ell_0-1})
+\beta_{\ell_0}\prod_{i=1}^{\ell_0}\delta_i^{-1}$,
we obtain
\[
\min\bigl\{
 \ord_{\gminip_0}(\kappa(\vecY_{\ell_0})),
 \ord_{\gminip_0}(\kappa(\vecY_{\ell_0-1}))
\bigr\}
\leq
 \ord_{\gminip_0}\bigl(
 \prod_{i=1}^{\ell_0}\delta_i^{-1}\beta_{\ell_0}
 \bigr)
=-\ord_{\gminip_0}(\delta_1\cdots\delta_{\ell_0})
\leq
 \ord_{\gminip_0}(\kappa(\vecY_{m})).
\]
It contradicts the assumption
$\ord_{\gminip_0}\kappa(\vecY_{m'})
>\ord_{\gminip_0}\kappa(\vecY_m)$
for $m'<m$.
Hence, we obtain
$\ord_{\gminip_0}(\kappa(\vecY_{m'+1,m}))<0$.
\hfill\qed

\vspace{.1in}
Then, we obtain the claim of Proposition \ref{prop;16.2.19.12}
by using an induction and Lemma \ref{lem;16.4.29.3}.
\hfill\qed

\subsection{Statements}
\label{subsection;18.11.15.40}

\subsubsection{Infinite sequences of complex blowings up}
\label{subsection;16.6.3.30}

We use the notation in \S\ref{subsection;16.6.1.40}.
Let $\vecY=(\veceta_1,\omega_1,\veceta_2,\omega_2,\ldots)
 \in\gbigp^{\infty}$,
where 
$\veceta_i\in \{+,-\}^{\ell(i)}$
and $\omega_i\in\cnum^{\ast}$.
We assume the following.
\begin{itemize}
 \item 
 The limit curve for $\vecY$ is not convergent.
\end{itemize}
We put 
$\vecY_m:=(\veceta_1,\omega_1,\ldots,\veceta_m,\omega_m)$.
We set $X_m:=\Bl_{\vecY_m}\cnum^2$.
We have constructed the sequence of morphisms of complex manifolds
\[
 \cdots\lrarr X_{m+1}
 \lrarr
 X_m
 \lrarr
 X_{m-1}
\lrarr\cdots\lrarr X_0=\cnum^2.
\]
We set the points
$P_m:=P_{\vecY_m}\in X_m$
and the coordinate neighbourhood
$(U_m,u_m,v_m):=(U_{\vecY_m},u_{\vecY_m},v_{\vecY_m})$
around $P_{m}$.
For $m'\geq m$,
the induced maps
$\psi_{\vecY_m,\vecY_{m'}}:
\Bl_{\vecY_{m'}}\cnum^2
\lrarr
 \Bl_{\vecY_{m}}\cnum^2$
are denoted by
$\psi_{m,m'}:X_{m'}\lrarr X_m$.
The induced morphisms
$\psi_{\vecY_m}:
\Bl_{\vecY_m}\cnum^2\lrarr\cnum^2$
are denoted by
$\psi_{m}:X_{m}\lrarr X_0$.

Let $H_0:=\{y=0\}\subset X_0$.
We set
$H_{m}:=\psi_{m}^{-1}(H_0)$.
Let $\varpi_{m}:
 \Xtilde_m(H_m)
\lrarr X_m$
be the oriented real blowing up
along $H_{m}$.
Let 
$\psitilde_{m}:
 \Xtilde_m(H_m)
\lrarr
 \Xtilde_0(H_0)$
and 
$\psitilde_{m,m'}:
 \Xtilde_{m'}(H_{m'})
\lrarr
 \Xtilde_m(H_m)$
denote the induced morphisms.

Note that $P_{m}$ are smooth points of 
normal crossing hypersurfaces $H_{m}$.
Let $\nbigu_{m}$ denote a small neighbourhood
of $0$ in $\{u_{m}\in\cnum\}$.
Let $\Delta_{m,\epsilon}:=\{|v_{m}|<\epsilon\}$.
We can naturally regard
$\nbigu_m\times\Delta_{m,\epsilon}$ 
as a neighbourhood of 
$P_{m}$ in $U_{m}$.
Let $\Deltatilde_{m,\epsilon}(0)\lrarr \Delta_{m,\epsilon}$
denote the oriented real blowing up at $0$.
We can naturally regard
$\nbigu_m\times\Deltatilde_{m,\epsilon}(0)$
as a neighbourhood of 
$\varpi_{m}^{-1}(P_{m})$
in $\Xtilde_m(H_m)$.
We can naturally identify
$\Deltatilde_{m,\epsilon}(0)$
with 
$\del\Deltatilde_{m,\epsilon}(0)\times \closedopen{0}{\epsilon}$,
and hence
$\nbigu_m\times \Deltatilde_{m,\epsilon}(0)
\simeq
 \nbigu_m\times
 \del\Deltatilde_{m,\epsilon}(0)
 \times
 \closedopen{0}{\epsilon}$.
For any an open subset
$\nbigi\subset
 \varpi_{m}^{-1}(P_{m})\simeq
 \del\Deltatilde_{m,\epsilon}(0)$
and a positive number $\epsilon'>0$,
we naturally regard
$\nbigi\times\openopen{0}{\epsilon'}$
and 
$\nbigi\times\closedopen{0}{\epsilon'}$
as subsets in
$\Xtilde_m(H_m)$.

\subsubsection{Lifting with respect to 
the composition of local real blowings up}
\label{subsection;16.7.26.30}

We identify $\Xtilde_0(H_0)=\cnum_x\times\widetilde{\cnum}_y(0)$
with $\nbigy_+:=\cnum\times \real_{\geq 0}\times S^1$.
We naturally embed it into
$\nbigy:=\cnum\times \real\times S^1$.
In this way, we regard $\Xtilde_0(H_0)$
as a closed subset in $\nbigy$.

Let $\phi:W\lrarr \nbigy$ be a morphism
obtained as the composite of local real blowings up:
\[
W=W^{(k)}
 \stackrel{\phi^{(k)}}{\lrarr}
 W^{(k-1)}
 \stackrel{\phi^{(k-1)}}{\lrarr}
\cdots
 \stackrel{\phi^{(2)}}{\lrarr}
 W^{(1)}
 \stackrel{\phi^{(1)}}{\lrarr}
 W^{(0)}=\nbigy.
\]
Moreover,
there exist subanalytic open subsets
$U^{(p)}\subset W^{(p)}$
and closed real analytic submanifolds
$C^{(p)}$ of $U^{(p)}$,
and the morphisms $\phi^{(p+1)}:W^{(p+1)}\lrarr W^{(p)}$
are obtained as the real blowing up of $U^{(p)}$
along $C^{(p)}$.
We may assume that 
$C^{(p)}$ are the complete intersection
of the real analytic functions
$g^{(p)}_1,\ldots,g^{(p)}_{r(p)}$ on $U^{(p)}$
such that 
$dg^{(p)}_1,\ldots,dg^{(p)}_{r(p)}$
are linearly independent
at each point of $C^{(p)}$.
Let $\xi^{(p)}:W^{(p)}\lrarr \nbigy$
denote the induced map.
We introduce a condition for $(W,\phi)$.

\begin{condition}
\label{condition;18.11.25.100}
There exist
a sequence of numbers
$m_i\to\infty$,
a sequence of points
$Q_{m_i}\in \varpi_{m_i}^{-1}(P_{m_i})$,
and
a sequence of neighbourhoods
$\nbigi_{m_i}$
of $Q_{m_i}$ in $\varpi_{m_i}^{-1}(P_{m_i})$,
such that the following holds.
\begin{description}
\item[E1]
$\psitilde_{m_i,m_j}(Q_{m_j})
=Q_{m_i}$
for $m_j>m_i$.
\item[E2]
For small positive numbers $r_{m_i}$,
the images
$\widetilde{\psi_{m_i}}(\nbigi_{m_i}
 \times\openopen{0}{r_{m_i}})$
do not intersect with 
the set of the critical values of $\phi$.
Moreover,
there exist real analytic maps
$\psihat_{m_i,0}:
 \nbigi_{m_i}\times\closedopen{0}{r_{m_i}}
\lrarr W$
such that 
$\phi\circ\psihat_{m_i,0}$
are equal to
the restriction of $\psitilde_{m_i}$
to $\nbigi_{m_i}\times\closedopen{0}{r_{m_i}}$.
Note that such $\psihat_{m_i,0}$
are uniquely determined by the property.
\end{description}
\end{condition}

We shall prove the following theorem later
(\S\ref{subsection;16.10.11.10} and \S\ref{subsection;16.10.11.11}).
\begin{thm}
\label{thm;16.5.4.1}
Suppose that Condition {\rm\ref{condition;18.11.25.100}}
is satisfied for $(W,\phi)$.
There exists $i_1$
such that the following holds for each $i\geq i_1$:
\begin{itemize}
\item
There exist a non-empty open subset
$\nbigi'_{m_i}\subset
 \nbigi_{m_i}$,
a neighbourhood $\nbigv_{m_i}$ of $\nbigi'_{m_i}$
in $\Xtilde_{m_i}(H_{m_i})$
and a real analytic map
\[
 \psihat_{m_i}:
 \nbigv_{m_i}
\lrarr
 W
\]
such that
$\phi\circ\psihat_{m_i}$
is equal to the restriction of
$\psitilde_{m_i}$
to $\nbigv_{m_i}$.
\end{itemize}
\end{thm}

\subsubsection{Lifting with respect to 
covering by the composition of local blowings up}
\label{subsection;18.11.16.101}

Let $(W_{\lambda},\phi_{\lambda})$ $(\lambda\in\Lambda)$
be a finite family of analytic maps
$\phi_{\lambda}:W_{\lambda}\lrarr \nbigy$
such that 
(i) $\phi_{\lambda}$ are the composition of local real blowings up,
(ii) there exist subanalytic compact subsets 
 $K_{\lambda}\subset W_{\lambda}$
 such that
 $\bigcup \phi_{\lambda}(K_{\lambda})$
contains the neighbourhood 
of $\varpi_0^{-1}(0,0)$.

\begin{prop}
\label{prop;16.5.6.10}
There exists $\lambda_0\in\Lambda$,
a sequence $m_i\to\infty$,
a sequence of points
$Q_{m_i}\in \varpi_{m_i}^{-1}(P_{m_i})$,
and a sequence of neighbourhoods
$\nbigi_{m_i}$
of $Q_{m_i}$ in $\varpi_{m_i}^{-1}(P_{m_i})$,
such that Condition {\rm\ref{condition;18.11.25.100}}
is satisfied 
for $(W_{\lambda_0},\phi_{\lambda_0})$.
In particular,
we may apply Theorem {\rm\ref{thm;16.5.4.1}}
to $\phi_{\lambda_0}$.
\end{prop}
\pf
Let $\Crit(\phi_{\lambda})$ denote the set of the critical values
of $\phi_{\lambda}$.
Note $\dim_{\real} \Crit(\phi_{\lambda})\leq 2$.
Let $q:\nbigy=\cnum\times(S^1\times\real)
 \lrarr S^1\times\real$
be the projection.
By construction,
$\Xtilde_0(H_0)$
is naturally identified with
$q^{-1}(S^1\times\real_{\geq 0})$.
It is standard that
there exists a $0$-dimensional subanalytic subset
$Z_0\subset S^1\times\{0\}$
with the following property.
\begin{itemize}
 \item
 For any $P\in (S^1\times\{0\})\setminus Z_0$,
 there exist a neighbourhood $\nbigu_P$ of $P$
 in $S^1\times \real_{\geq 0}$
 and continuous subanalytic functions
 $f^P_{1},\ldots,f^P_{\ell}$ on $(\nbigu_P,S^1\times\real)$
 to $\cnum$ such that 
 $\bigcup_{\lambda}\Crit(\phi_{\lambda})\cap q^{-1}(\nbigu_P)$
 is the union of the graph of
 $f^P_p$. 
\end{itemize}

We regard $Z_0\subset \varpi_0^{-1}(0,0)$
by the natural isomorphism
$\varpi_0^{-1}(0,0)\simeq S^1\times\{0\}$.

\begin{lem}
If $m_0$ is sufficiently large,
for each $m\geq m_0$,
there exists a finite subset
$Z_m\subset \varpi_m^{-1}(P_m)$
with the following property.
\begin{itemize}
\item 
 $\psitilde_m^{-1}(Z_0)\cap\varpi_m^{-1}(P_m)
 \subset Z_m$.
\item
 For any point
 $P\in \varpi_m^{-1}(P_m)\setminus Z_m$,
 there exists a neighbourhood  $\nbigu_P$
 of $P$ in $\{0\}\times\Deltatilde_{m,\epsilon}(0)$,
 a positive number $r(P,m)$,
 $\lambda(P)\in\Lambda$
 and a morphism
 $\psihat_{m,\lambda(P),P}:\nbigu_P\times\closedopen{0}{r(P,m)}
 \lrarr W_{\lambda(P)}$
 satisfying
 $\phi_{\lambda(P)}\circ\psihat_{m,\lambda(P),P}
 =\psitilde_{m|\nbigu_P}$.
\end{itemize}
\end{lem}
\pf
Note that 
$\psitilde_m(\{0\}\times\Deltatilde_{m,\epsilon}(0))$
is described as the graph of
the multivalued functions
$S^1\times\closedopen{0}{\rho_m}\lrarr \cnum$
induced by
$g_{\vecY_m|\{0\}\times\Delta_y}$,
where $g_{\vecY_m}$ are defined
as in \S\ref{subsection;16.6.1.40}.
The limit curve is assumed to be non-convergent.
Hence, if $m_1$ is sufficiently large,
for each $m\geq m_1$,
the intersection of 
$\psitilde_m(\{0\}\times\Deltatilde_{m,\epsilon}(0))$
and
$\bigcup \Crit(\phi_{\lambda})$
is at most one dimensional.

For $m\geq m_1$,
there exist the subanalytic subsets
$M_{m,\lambda}:=
 \psitilde_m^{-1}\bigl(
 \phi_{\lambda}(K_{\lambda})\bigr)
\cap
 (\{0\}\times\Deltatilde_{m,\epsilon}(0))$.
We obtain the following decomposition into connected components:
\[
 M_{m,\lambda}
\setminus
\Bigl(
 \bigl(
 \{0\}\times
 \del\Deltatilde_{m,\epsilon}(0)
 \bigr)
\cup
 \psitilde_m^{-1}\Bigl(
 \bigcup\Crit(\phi_{\lambda})
 \Bigr)
\Bigr)
=\coprod \nbigc_{m,\lambda,j}.
\]
There exist the unique morphisms
$g_{m,\lambda,j}:
 \nbigc_{m,\lambda,j}
\lrarr
 W_{\lambda}$
such that
$\phi_{\lambda}\circ g_{m,\lambda,j}$
is equal to the restriction of
$\psitilde_{m}$ to $\nbigc_{m,\lambda,j}$.
According to Lemma \ref{lem;16.7.26.20},
for each $(\lambda,j)$ such that
$\dim_{\real} \nbigc_{m,\lambda,j}=2$,
there exist a finite subset 
$N_{m,\lambda,j}\subset \overline{\nbigc}_{m,\lambda,j}$
and a real analytic map
$\overline{g}_{m,\lambda,j}:
 \overline{\nbigc}_{m,\lambda,j}\setminus N_{m,\lambda,j}
\lrarr
 W_{\lambda}$
such that 
(i) $N_{m,\lambda,j}$
 contains the singular locus of 
$\del\overline{\nbigc}_{m,\lambda,j}$,
(ii)
$\phi_{\lambda}\circ \overline{g}_{m,\lambda,j}$
is equal to
the restriction of $\psitilde_{m}$
to $\overline{\nbigc}_{m,\lambda,j}
\setminus N_{m,\lambda,j}$.

Because
$\Deltatilde_{m,\epsilon}(0)
=\bigcup_{\lambda,j} \overline{\nbigc}_{m,\lambda,j}$,
there exists a finite set
$Z_{m}\subset\del\Deltatilde_{m,\epsilon}(0)$
such that
(i) $Z_{m}\supset\psitilde_{m}^{-1}(Z_{0})\cap\varpi_m^{-1}(P_m)$,
(ii) for any $P\in\varpi_m^{-1}(P_m)\setminus Z_{m}$,
 there exists $(\lambda(P),j(P))$ such that
 $P\in \overline{\nbigc}_{m,\lambda(P),j(P)}
 \setminus N_{m,\lambda(P),j(P)}$.
Then, $Z_m$ has the desired property.
\hfill\qed

\vspace{.1in}

The union
$Z'_0:=\bigcup_m \psitilde_m(Z_m)$ is 
a countable subset in $\varpi_0^{-1}(0,0)$.
Note $Q_0\in \varpi_0^{-1}(0,0)\setminus Z_0'$.
There exists a sequence $Q_m\in \varpi_m^{-1}(P_m)$
such that $\psitilde_{m,m'}(Q_{m'})=Q_m$
and that $\psitilde_m(Q_m)=Q_0$.
There exist neighbourhoods $\nbigi_m$ of $Q_m$
in $\varpi_m^{-1}(P_m)$ such that 
(i) $\psitilde_m(\nbigi_m)
 \supset
 \psitilde_{m'}(\nbigi_{m'})$
for $m\leq m'$,
(ii)
$\nbigi_m\cap
 \bigcup_{j\leq m}
 \nutilde_{j,m}^{-1}(Z_j)
=\emptyset$.
Hence, 
there exist $\lambda_0\in\Lambda$,
a sequence $m_i\to\infty$,
a sequence of points
$Q_{m_i}\in\varpi_{m_i}^{-1}(P_{m_i})$,
a sequence of neighbourhoods
$\nbigi_{m_i}$
of $Q_{m_i}$ in $\varpi_{m_i}^{-1}(P_{m_i})$
such that
Condition \ref{condition;18.11.25.100}
is satisfied 
for $(W_{\lambda_0},\phi_{\lambda_0})$.
\hfill\qed

\subsubsection{Pull back of subanalytic functions
 in the case $\kappahat(\vecY)<\infty$}

Let $\gbigh$ be a closed subanalytic subset in 
$\Xtilde_0(H_0)$ with $\dim_{\real}\gbigh\leq 3$.
Let $f$ be a continuous subanalytic function on
$(\Xtilde_0(H_0) \setminus \gbigh,
 \Xtilde_0(H_0))$.
We assume that 
$f$ is bounded around 
any point of $\gbigh\setminus\del\Xtilde_0(H_0)$.
We shall prove the following proposition 
in \S\ref{subsection;16.10.11.1}.
\begin{thm}
\label{thm;16.5.6.1}
If $\kappahat(\vecY)<\infty$,
there exists $m_0$ such that 
the following holds for each $m\geq m_0$.
\begin{itemize}
\item
There exist a non-empty connected open subset
$\nbigi_m\subset\varpi_m^{-1}(P_m)$
and a small open neighbourhood 
$\nbigv_m$ of $\nbigi_m$
in $\Xtilde_m(H_m)$
such that 
$\psitilde_m
 \bigl(\nbigv_m\setminus\varpi_m^{-1}(H_m)\bigr)
\subset
 \Xtilde_0(H_0)\setminus \gbigh$.
\item
The function
$\psitilde_m^{\ast}(f)$ on $\nbigv_m\setminus\varpi_m^{-1}(H_m)$
is ramified real analytic along $\nbigv_m\cap\varpi_m^{-1}(H_m)$.
\item
If $\psitilde_m^{\ast}(f)$ is not constantly $0$,
the order
$\ord_{\rho_m}\psitilde_m^{\ast}(f)(u_m,\theta_m,\rho_m)$
are independent of $u_m$,
where $(\theta_m,\rho_m)$ is a polar coordinate system of 
$\cnumtilde_{v_m}(0)$
defined by $v_m=\rho_me^{\sqrt{-1}\theta_m}$.
\end{itemize}
\end{thm}

\subsubsection{Pull back of subanalytic functions
in the case $\kappahat(\vecY)=\infty$}
\label{subsection;16.10.11.20}

In this subsection,
we assume that $\vecY$ is not convergent 
and that $\kappahat(\vecY)=\infty$.
We shall prove the propositions
in \S\ref{subsection;16.10.11.2}.

Any small neighbourhood $\nbigu_m$ of $0$
in $\cnum_{u_m}=\{u_m\in\cnum\}$
naturally induces a subset
$\nbigu_m\times\{0\}$
of $U_{m}=
 (\cnum_{u_m}\setminus\{-\omega_m\})
 \times \cnum_{v_m}
\subset \Xtilde_m(H_m)$.
Let $u_m=a_m+\sqrt{-1}b_m$
denote the decomposition
into the real part and the imaginary part.
We consider open sets $\nbigu_m$
of the form 
$\bigl\{(a_m,b_m)\,\big|\,
 |a_m|<\delta_{1,m},\,
 |b_m|<\delta_{2,m}
 \bigr\}$.

Let $\gbigh$ be a closed subanalytic subset in $\Xtilde_0(H_0)$
with $\dim_{\real}\gbigh=3$.
Let $\gbigh_1$ be a closed subanalytic subset in $\gbigh$
with $\dim_{\real}\gbigh_1=2$.
\begin{thm}
\label{thm;16.10.11.30}
There exists $m_0$ such that 
for any $m\geq m_0$
there exist a non-empty connected open subset 
$\nbigi_m\subset\varpi_m^{-1}(P_m)$
and a small neighbourhood
$\nbigv_m=
\nbigu_m\times \nbigi_m\times \{0\leq \rho_m<\epsilon_m\}$
of $\nbigi_m$ in $\Xtilde_m(H_m)$
with the following property.
\begin{itemize}
 \item
Set 
$\nbigz_m:=
 (\nbigv_m\cap\psitilde_m^{-1}(\gbigh))
 \setminus\varpi_m^{-1}(H_m)$.
Then,  either one of the following holds:
(i) $\nbigz_m$ is empty,
or (ii) $\nbigz_m$ is a smooth connected hypersurface.
\item
In the case (ii),
one of the projections
$u_m\longmapsto a_m$
or $u_m\longmapsto b_m$
induces an isomorphism
\begin{equation}
 \label{eq;16.8.28.40}
\nbigz_m
\simeq
 \{|a_m|<\delta_{m,1}\}
 \times\nbigi_m\times\{0<\rho_m<\epsilon_m\}
=:\nbigv_{m,1},
\end{equation}
\begin{equation}
\label{eq;16.8.28.41}
\mbox{\rm or}\quad
\nbigz_m
\simeq
 \{|b_m|<\delta_{m,2}\}
 \times\nbigi_m\times\{0<\rho_m<\epsilon_m\}
=:\nbigv_{m,2}.
\end{equation}
\item
$(\nbigv_m\cap\psitilde_m^{-1}(\gbigh_1))\setminus
 \varpi_m^{-1}(H_m) =\emptyset$. 
\end{itemize}
\end{thm}

Let $f$ be a continuous subanalytic function
on $(\Xtilde_0(H_0)\setminus \gbigh,\Xtilde_0(H_0))$,
which is bounded around any point of
$\gbigh\setminus \del\Xtilde_0(H_0)$.
Let $f_{\gbigh}$ be a continuous subanalytic function
on $(\gbigh\setminus \gbigh_1,\Xtilde_0(H_0))$,
which is bounded around any point of
$\gbigh_1\setminus\del\Xtilde_0(H_0)$.

\begin{thm}
 \label{thm;16.10.11.31}
There exists $m_1\geq m_0$ such that 
the following holds
for any $m\geq m_1$:
\begin{itemize}
 \item
 Suppose $\nbigz_m\neq\emptyset$.
 Then, 
the restriction of 
$\psitilde_m^{\ast}(f)$
to each connected component of
$\nbigv_m\cap
 \psitilde_m^{-1}(\Xtilde_0(H_0)\setminus \gbigh)$
is described as the sum of
a bounded function
and a ramified real analytic function 
of the form $\sum_{\gminiy<0}f_{m,\gminiy}(\theta_m)\rho_m^{\gminiy}$.
Moreover, the function 
$\psitilde_m^{\ast}(f_{\gbigh})$
on $\nbigz_m$ is 
described as the sum of
a bounded function and
a ramified real analytic function
of the form
$\sum_{\gminiy<0}f_{\gbigh,m,\gminiy}(\theta_m)\rho_m^{\gminiy}$
on $\nbigv_{m,1}$ or $\nbigv_{m,2}$
under the isomorphism
{\rm(\ref{eq;16.8.28.40})}
or {\rm(\ref{eq;16.8.28.41})}.
\item
Suppose that $\nbigz_m=\emptyset$.
Then, the following holds for 
a (possibly empty) closed subanalytic subset 
$\nbiga_m\subset \nbigv_m$:
\begin{itemize}
 \item
 Set $\nbiga_m^{\circ}:=\nbiga_m\setminus\varpi_m^{-1}(H_m)$.
 If $\nbiga^{\circ}_m\neq\emptyset$,
one of the projections
$u_m\longmapsto a_m$
 or $u_m\longmapsto b_m$
induces an isomorphism
 $\nbiga_m^{\circ}\simeq \nbigv_{m,1}$ or 
 $\nbiga_m^{\circ}\simeq\nbigv_{m,2}$.
\item
The restriction of 
$\psitilde_m^{\ast}(f)$
to each connected component of
$\nbigv_m
 \setminus(\nbiga_m\cup\varpi_m^{-1}(H_m))$
is described as a sum of
a bounded function
and a ramified real analytic function 
of the form $\sum_{\gminiy<0}f_{m,\gminiy}(\theta_m)\rho_m^{\gminiy}$.
Moreover, the restriction of 
$\psitilde_m^{\ast}(f)$
to $\nbiga_m^{\circ}$
is described as a sum of
a bounded function and
a ramified real analytic function
of the form
$\sum_{\gminiy<0}f_{\nbiga,m,\gminiy}(\theta_m)\rho_m^{\gminiy}$
on $\nbigv_{m,1}$ or $\nbigv_{m,2}$
under the isomorphism
$\nbiga_m^{\circ}\simeq\nbigv_{m,1}$
or 
$\nbiga_m^{\circ}\simeq\nbigv_{m,2}$.
\end{itemize}
\end{itemize}
\end{thm}

\subsection{The case $\kappahat(\vecY)<\infty$}
\label{subsection;18.11.13.1}

\subsubsection{Proof of Theorem \ref{thm;16.5.4.1}}
\label{subsection;16.10.11.10}

In the following,
we shall make
$\nbigi_{m}$ and $r_{m}$ smaller.
Let $\nbigv_{m}$ denote a small neighbourhood
of $\nbigi_{m}\times\closedopen{0}{r_{m}}$
in $\Xtilde_{m}(H_{m})$
of the form
$\nbigu_{m}\times
 \{\theta_1<\theta<\theta_2\}
 \times
 \closedopen{0}{r_m}$.
We denote $m_i$ by $m(i)$.

Because of the existence of $\psihat_{m(i),0}$,
if $r_{m(i)}$ are sufficiently small,
we obtain
$\psitilde_{m(i)}\bigl(
 \nbigi_{m(i)}
 \times\closedopen{0}{r_{m(i)}}\bigr)
\subset
 U^{(0)}$.
Hence, we may assume that 
$\psitilde_{m(i)}\bigl(
\nbigv_{m(i)} \bigr)
\subset
 U^{(0)}$.
Recall that $C^{(0)}\subset U^{(0)}$
is described as
$\{g^{(0)}_1=\cdots =g^{(0)}_{r(0)}=0\}$.
We use the change of parametrization in \S\ref{subsection;16.6.2.1},
and we apply Theorem \ref{thm;16.3.8.20}.
Then, 
by shrinking $\nbigi_{m(i)}$, $r_{m(i)}$ and $\nbigv_{m(i)}$,
there exists $m^{(0)}$ such that 
for any $m(i)\geq m^{(0)}$
the following holds:
(i) $\psitilde_{m(i)}^{\ast}(g^{(0)}_j)$ 
are not constantly $0$,
(ii) the order
$\ord_{\rho_{m(i)}}
 \psitilde_{m(i)}^{\ast}
 g^{(0)}_j(u_{m(i)},\theta_{m(i)},\rho_{m(i)})$
are constant with respect to 
$(u_{m(i)},\theta_{m(i)})
 \in \varpi_{m(i)}^{-1}(H_{m(i)})\cap\nbigv_{m(i)}$.
Then, there exists a unique morphism
$\psitilde_{m(i)}^{(1)}:
\nbigv_{m(i)}
\lrarr
 W^{(1)}$
such that 
$\phi^{(1)}\circ\psitilde_{m(i)}^{(1)}$
is equal to
$\psitilde_{m(i)}$ for any $m(i)\geq m^{(0)}$.
By taking a subsequence,
we may assume that such morphisms
$\psitilde_{m(i)}^{(1)}$ exist
for any $m(i)$.

\vspace{.1in}

Suppose that for $\ell\geq 1$
we have already constructed
real analytic morphisms
$\psitilde_{m(i)}^{(\ell)}:
 \nbigv_{m(i)}
\lrarr
 W^{(\ell)}$
such that
$\xi^{(\ell)}\circ\psitilde_{m(i)}^{(\ell)}
=\psitilde_{m(i)}$.
Because of the existence of $\psihat_{m(i),0}$,
we may assume that 
$\Image\psitilde_{m(i)}^{(\ell)}$
is contained in $U^{(\ell)}$,
after making $\nbigv_{m(i)}$ smaller.
We obtain real analytic functions
$(\psitilde_{m(1)}^{(\ell)})^{\ast}(g^{(\ell)}_j)$
on $\nbigv_{m(1)}$,
and the real analytic functions
$(\psitilde_{m(i)}^{(\ell)})^{\ast}(g^{(\ell)}_j)$ $(i\geq 1)$
are obtained as the pull back of 
$(\psitilde_{m(1)}^{(\ell)})^{\ast}(g^{(\ell)}_j)$
by $\psitilde_{m(1),m(i)}$.
Hence, by using the change of parametrization
as in \S\ref{subsection;16.6.2.1}
and Theorem \ref{thm;16.3.8.20},
and by shrinking $\nbigi_{m(i)}$, $r_{m(i)}$ and $\nbigv_{m(i)}$,
there exist $m^{(\ell)}$
such that for any $m(i)\geq m^{(\ell)}$
the following holds:
(i) $\psitilde_{m(i)}^{\ast}(g^{(\ell)}_j)$ are not constantly $0$,
(ii) the orders
$\ord_{\rho_{m(i)}}
 \psitilde_{m(i)}^{\ast}g^{(\ell)}_j(u_{m(i)},\theta_{m(i)},\rho_{m(i)})$
are constant
with respect to
$(u_{m(i)},\theta_{m(i)})
\in \varpi_{m(i)}^{-1}(H_{m(i)})\cap \nbigv_{m(i)}$.
Hence, there exists a unique morphism
$\psitilde_{m(i)}^{(\ell+1)}:
 \nbigv_{m(i)}
\lrarr
 W^{(\ell+1)}$
such that
$\phi^{(\ell+1)}\circ\psitilde_{m(i)}^{(\ell+1)}
=\psitilde_{m(i)}^{(\ell)}$.
In this way,
the inductive construction can proceed,
and we obtain Theorem \ref{thm;16.5.4.1}
in the case $\kappahat(\vecY)<\infty$.
\hfill\qed

\subsubsection{Proof of Theorem \ref{thm;16.5.6.1}}
\label{subsection;16.10.11.1}

Let 
$(W_{\lambda},\phi_{\lambda})$ $(\lambda\in\Lambda)$
be any rectilinearization for $f$.
By Theorem \ref{thm;16.5.4.1}
and Proposition \ref{prop;16.5.6.10},
there exist $\lambda_0\in\Lambda$,
a positive number $m(1)$,
an open subset $\nbigv_{m(1)}$
in $\Xtilde_{m(1)}(H_{m(1)})$
such that
$\nbigv_{m(1)}\cap
 \varpi_{m(1)}^{-1}(P_{m(1)})
\neq\emptyset$,
and a real analytic morphism
$\psihat_{m(1)}:
 \nbigv_{m(1)}\lrarr W_{\lambda_0}$
such that 
$\phi_{\lambda_0}\circ\psihat_{m(1)}
=\psitilde_{m(1)|\nbigv_{m(1)}}$.
There exists a sequence of open subsets
$\nbigv_{m}\subset \Xtilde_{m}(H_m)$
$(m\geq m(1))$
such that
(i) $\nbigv_{m}\cap\varpi_m^{-1}(P_m)\neq\emptyset$,
(ii) $\psitilde_{m',m}(\nbigv_m)\subset \nbigv_{m'}$
for any $m\geq m'\geq m(1)$.

The set
$\phi_{\lambda_0}^{-1}(\gbigh)$
is rectilinearized,
and hence expressed as the $0$-set of a real analytic
function $h_{\lambda_0}$ on $W_{\lambda_0}$.
The set $\psitilde_{m(1)}^{-1}(\gbigh)$
is contained in the $0$-set of
$h_{m(1),\lambda_0}:=
 \psihat_{m(1)}^{\ast}(h_{\lambda_0})$.
If $m$ is sufficiently larger than $m(1)$,
after shrinking $\nbigv_m$,
we obtain that the $0$-set of 
$\psitilde_{m(1),m}^{\ast}(h_{m(1),\lambda_0})$
is contained in 
$\nbigv_m\cap\varpi_m^{-1}(H_m)$
by Corollary \ref{cor;16.5.6.20}.
It implies the first claim of Theorem \ref{thm;16.5.6.1}
in the case $\kappahat(\vecY)<\infty$.

Let $(x_1,x_2,x_3,x_4)$ be the coordinate system 
on $W_{\lambda_0}$.
We obtain the real analytic functions
$\psihat_{m(1)}^{\ast}(x_i)$.
After replacing $m(1)$ with a larger number,
we may assume that 
$x_{i,m(1)}:=\psihat_{m(1)}^{\ast}(x_i)>0$
on $\nbigv_{m(1)}$.
If $m$ is sufficiently larger than $m(1)$,
after shrinking $\nbigv_m$,
we obtain that the order
$\ord_{\rho_m}
 \psitilde_{m(1),m}^{\ast}(x_{i,m(1)})(u_m,\theta_m,\rho_m)$
are constant with respect to $(u_m,\theta_m)$.
Then, 
$\psitilde_{m(1),m}^{\ast}(x_{i,m(1)})^{1/e}$
are ramified real analytic functions on $\nbigv_{m}$.

Note that
$\psitilde_m^{\ast}(f)
=\psitilde_{m(1),m}^{\ast}
 \psihat_{m(1)}^{\ast}\bigl(
 \phi_{\lambda_0}^{\ast}(f)
 \bigr)$,
and that the function
$\phi_{\lambda_0}^{\ast}(f)$
is expressed as real analytic functions
of $(x_1^{1/e},x_2^{1/e},x_3^{1/e},x_4^{1/e})$
for some $e\in\seisuu_{>0}$.
Hence, 
if $m$ is sufficiently larger than $m(1)$,
we obtain that
$\psitilde_m^{\ast}(f)$
is a ramified real analytic function
on $\nbigv_m$,
i.e.,
we obtain the second claim of
Theorem \ref{thm;16.5.6.1}
in the case $\kappahat(\vecY)<\infty$.

By replacing $m(1)$ with a larger number,
we may assume that
$\psitilde_m^{\ast}(f)$ is a ramified real analytic function.
Then,
by Corollary \ref{cor;16.5.6.20}
we obtain the third claim of
Theorem \ref{thm;16.5.6.1}
in the case $\kappahat(\vecY)<\infty$.
Thus, the proof of Theorem \ref{thm;16.5.6.1}
is finished
in the case $\kappahat(\vecY)<\infty$.
\hfill\qed

\subsection{The case $\kappahat(\vecY)=\infty$}
\label{subsection;18.11.13.2}

\subsubsection{Change of parametrization}

As in \S\ref{subsection;16.6.2.1},
we change the parametrization of curves.
Let $\nbigu_m$ denote a small neighbourhood of
$0$ in $\cnum_{u_m}$.
We have constructed the holomorphic embedding
$\Phi_m^{-1}:
 \nbigu_m\times\Delta_{y^{1/\vecdelta(m)},\epsilon_m}
\lrarr 
 U_m$,
and the composition
$F_m=\psi_m\circ\Phi_m^{-1}$
is described as
\[
 F_m(u_m,y^{1/\vecdelta(m)})
=(g_{\vecY_m}(u_m,y^{1/\vecdelta(m)}),y).
\]
Let 
$\Ftilde_m:
 \nbigu_m\times
 \Deltatilde_{y^{1/\vecdelta(m)},\epsilon_m}(0)
\lrarr
 \Xtilde_0(H_0)$
denote the induced morphism.
It is enough to study 
$\Ftilde_m$
instead of $\psitilde_m$.

We use the polar coordinate system
$y=te^{\sqrt{-1}\phi}$,
i.e.,
$y^{1/\vecdelta(m)}
=t^{1/\vecdelta(m)}
 e^{\sqrt{-1}\phi/\vecdelta(m)}$.
Set $\gbigi:=\{0\leq \phi\leq 2\pi\}$.
We define the map
$\nbigu_m\times
 \gbigi\times
 \{0\leq t^{1/\vecdelta(m)}\leq \epsilon\}
\lrarr
 \Xtilde_0(H_0)$
by
$(u_m,\phi,t^{1/\vecdelta(m)})\longmapsto
\Ftilde_m(u_m,e^{\sqrt{-1}\phi/\vecdelta(m)}t^{1/\vecdelta(m)})$.
We may regard $\nbigi_m\subset\gbigi$.
The induced $\nbign_{\nbigu_m\times\gbigi}$-paths
in $\Xtilde_0(H_0)$
are also denoted by $\Ftilde_m$.
(See \S\ref{subsection;16.5.31.32}
for $\nbign$-paths.)

\subsubsection{Limit}

There exist the expansions
$g_{\vecY_m}(u_m,t^{1/\vecdelta(m)}e^{\sqrt{-1}\phi/\vecdelta(m)})
=\sum_{\gminiy}g_{\vecY_m,\gminiy}(u_m)e^{\sqrt{-1}\phi\gminiy}t^{\gminiy}$.
Recall that
for $\gminiy<\kappa(\vecY_m)$,
the coefficients
$g_{\vecY_{m'},\gminiy}(u_m)$ are independent of $m'>m$,
and constant with respect to $u_m$.
We define $g_{\vecY,\gminiy}$
as $g_{\vecY_m,\gminiy}\in\cnum$ $(\gminiy<\kappa(\vecY_m))$.
We set
$g_{\vecY}:=\sum g_{\vecY,\gminiy}e^{\sqrt{-1}\gminiy\phi} t^{\gminiy}$
as a section of $\nbign_{\gbigi}$.
(See \S\ref{subsection;16.5.31.32}
 for the sheaf $\nbign_{\gbigi}$ on $\gbigi$.)
We obtain the non-convergent $\nbign_{\gbigi}$-path
$\Ftilde_{\infty}$
in $\Xtilde_0(H_0)$
defined by
$\Ftilde_{\infty}(\phi,t)
=(g_{\vecY},\phi,t)$.
We can naturally regard $\Ftilde_{\infty}$
as an $\nbign_{\nbigu_m\times\gbigi}$-path
in $\widetilde{\cnum^2}(H)$ for each $m$.
By the construction,
$\ord_t(\Ftilde_{\infty},\Ftilde_{m})\geq \kappa(\vecY_m)$
holds.

\subsubsection{Proof of Theorem \ref{thm;16.5.4.1}}
\label{subsection;16.10.11.11}

Recall that for any $\nbign_Y$-path $F$ in a space $B$
the underlying map $Y\lrarr B$
is denoted by $F_0$.

We set $\gbigi^{(0)\prime}:=
\bigl((\Ftilde_{\infty})_0\bigr)^{-1}(U^{(0)})$.
The 
$\nbign_{\gbigi^{(0)\prime}}$-path
$\Ftilde_{\infty|\gbigi^{(0)\prime}}$
is induced in $U^{(0)}$.
Note that
$\Ftilde_{\infty|\gbigi^{(0)\prime}}$
does not factor through $C^{(0)}$
by Lemma \ref{lem;16.8.28.20}.
There exists the $0$-dimensional closed subset
$Z^{(0)}\subset \gbigi^{(0)\prime}$
such that
for any $P\in \gbigi^{(0)\prime}\setminus Z^{(0)}$
the following holds:
\[
\ord_{P,t}(\Ftilde_{\infty|\gbigi^{(0)\prime}},C^{(0)})
=
\ord_{P,t}(\Ftilde_{\infty|P},C^{(0)})
=:\mu^{(0)}.
\]
We set $\gbigi^{(1)}:=\gbigi^{(0)\prime}\setminus Z^{(0)}$.
As explained in \S\ref{subsection;16.6.3.1},
there exists the $\nbign_{\gbigi^{(1)}}$-path 
$\Ftilde^{(1)}_{\infty,\gbigi^{(1)}}$ in $W^{(1)}$
which is the lift of
$\Ftilde_{\infty|\gbigi^{(1)}}$,
i.e.,
$\phi^{(1)}\circ\Ftilde^{(1)}_{\infty,\gbigi^{(1)}}
=\Ftilde_{\infty|\gbigi^{(1)}}$.

We continue such a process inductively, as possible.
Suppose that we have already constructed
an $\nbign_{\gbigi^{(\ell)}}$-path
$\Ftilde_{\infty,\gbigi^{(\ell)}}^{(\ell)}$
in $W^{(\ell)}$
for $\ell\geq 1$.
We set 
$\gbigi^{(\ell)\prime}:=
 \Bigl(
 (\Ftilde_{\infty,\gbigi^{(\ell)}}^{(\ell)})_{0}
 \Bigr)^{-1}(U^{(\ell)})$.
If $\gbigi^{(\ell)\prime}\neq\emptyset$,
we obtain the induced 
$\nbign_{\gbigi^{(\ell)\prime}}$-path
$\Ftilde^{(\ell)}_{\infty,\gbigi^{(\ell)\prime}}$
in $U^{(\ell)}$.
Note that
$\Ftilde^{(\ell)}_{\infty|\gbigi^{(\ell)\prime}}$
does not factor through $C^{(\ell)}$.
There exists a $0$-dimensional subset
$Z^{(\ell)}\subset \gbigi^{(\ell)\prime}$
such that
for any $P\in \gbigi^{(\ell)\prime}\setminus Z^{(\ell)}$
the following holds:
\[
 \ord_{P,t}\Bigl(
 \Ftilde^{(\ell)}_{\infty|\gbigi^{(\ell)\prime}},
 C^{(\ell)}
 \Bigr)
=
 \ord_{P,t}\Bigl(
 (\Ftilde^{(\ell)}_{\infty|\gbigi^{(\ell)\prime}})_{|P},
 C^{(\ell)}
 \Bigr)
=:\mu^{(\ell)}.
\]
We set 
$\gbigi^{(\ell+1)}:=
 \gbigi^{(\ell)\prime}\setminus
 Z^{(\ell)}$.
There exists the $\nbign_{\gbigi^{(\ell)}}$-path
$\Ftilde^{(\ell)}_{\infty,\gbigi^{(\ell+1)}}$
in $W^{(\ell+1)}$
which is the lift of
$\Ftilde^{(\ell)}_{\infty|\gbigi^{(\ell+1)}}$.

We continue the process until either 
$\ell=k$
or $\gbigi^{(\ell)\prime}=\emptyset$
holds.

\vspace{.1in}

We set 
$N_0:=
 \sum_{\ell} \mu^{(\ell)}
<\infty$.
We can easily obtain
the claim of Theorem \ref{thm;16.5.4.1}
from the following lemma,
in the case $\kappahat(\vecY)=\infty$ for non-convergent $\vecY$.

\begin{lem}
\mbox{{}}\label{lem;16.6.2.10}
\begin{itemize}
\item
There exist
an open subset $\gbigi_0\subset\gbigi$
with $\gbigi_0\cap\nbigi_{m_i}\neq \emptyset$ $(i\geq 1)$
and an $\nbign_{\gbigi_0}$-path
$\Fhat_{\infty,\gbigi_0}$
in $W$
which is a lift of $\Ftilde_{\infty|\gbigi_0}$
with respect to $\phi$.
\item
Let $N_1$ be any positive integer.
Let $A$ be any complex manifold.
Let $\nu$ be any $\nbign_{\gbigi_0\times A}$-path
in $\Xtilde_0(H_0)$
such that
\[
 \ord_{t,P}\bigl(
 \nu_{|\gbigi_0\times\{a\}},
 \Ftilde_{\infty|\gbigi_0}
 \bigr)>N_0+N_1
\]
 for any $P\in\gbigi_0$
 and for any $a\in A$.
Then, there exists
the $\nbign_{\gbigi_0\times A}$-path 
$\nutilde$ in $W$ such that 
(i) it is the lift of $\nu$,
(ii) $\ord_{t,P}\bigl(
 \nutilde_{\gbigi\times\{a\}},
 \Fhat_{\infty,\gbigi_0}
 \bigr)>N_1$
for any $P\in\gbigi_0$
and $a\in A$.
\item
Let $N_1$ be any positive integer.
If $m_i$ is sufficiently large,
there exists
the $\nbign_{\nbigu_m\times(\gbigi_{0}\cap\nbigi_m)}$-path
$\Fhat_{m_i}$ in $W$
such that
(i) $\Fhat_{m_i}$ is the lift of 
$\Ftilde_{m|\nbigu_m\times(\gbigi_0\cap\nbigi_m)}$
with respect to $\phi$,
(ii) the following holds
at any $P\in \gbigi_0\cap\nbigi_m$ 
and any $R\in \nbigu_{m_i}$:
\[
 \ord_{t,P}\Bigl(
 \Fhat_{\infty,\gbigi_0},
 \Fhat_{m_i|\{R\}\times(\gbigi_0\cap\nbigi_m)}
 \Bigr)
>N_1.
\]
\end{itemize}
\end{lem}
\pf
If $m_i$ is sufficiently large,
$\ord(\Ftilde_{m_i,\gbigi},
 \Ftilde_{\infty,\gbigi})
>N_0$
holds.
Then, by an inductive argument on $\ell$,
we can prove that
(i) $\gbigi^{(\ell)}\neq \emptyset$,
(ii) $\gbigi^{(\ell)}_{m_i}:=
 \nbigi_{m_i}\cap\gbigi^{(\ell)}\neq \emptyset$,
(iii) the following holds:
\[
 \ord_{t,P}\Bigl(
 \Ftilde^{(\ell)}_{m_i,\gbigi^{(\ell)}_{m_i}},
  \Ftilde^{(\ell)}_{\infty,\gbigi^{(\ell)}_{m_i}}
 \Bigr)
>N_0-\sum_{p\leq\ell} \mu^{(p)}.
\]
Hence, we obtain the first claim.
We can also prove the second claim
similarly by an easy induction.
The third claim follows from the second.
Thus, we obtain Lemma \ref{lem;16.6.2.10}
and Theorem \ref{thm;16.5.4.1}.
\hfill\qed

\subsubsection{Pull back of real analytic functions}

We make a preliminary for the proof of 
Theorems \ref{thm;16.10.11.30}
and \ref{thm;16.10.11.31}.
Recall that $(x,y)$ is the coordinate system of $X_0=\cnum^2$
such that $H_0=\{y=0\}$.
Let $y=re^{\sqrt{-1}\theta}$ be the polar decomposition.
Let 
$\gbigj:=\{\theta_1<\theta<\theta_2\}\subset
 \gbigi$
be any interval.
We consider an open subset
$\nbigb:=
 \bigl\{
 (x,\theta,r)\,\big|\,
 |x|<\epsilon_1,\,\,
 \theta_1<\theta<\theta_2,\,\,
 0\leq r<\epsilon_2
 \bigr\}$
in $\Xtilde_0(H_0)$.
Let $h$ be a non-constant real analytic function on $\nbigb$.
We assume the following.
\begin{itemize}
\item
The exterior derivative of $h$ is not $0$
at every point of $h^{-1}(0)\cap\{r>0\}$.
\item
$\Ftilde_{\infty,\gbigj}^{\ast}(h)=0$
as a section of $\nbign_{\gbigj}$.
\end{itemize}
Let $(a,b)$ be the real coordinate system
determined by $x=a+\sqrt{-1}b$.

\begin{lem}
One of 
$\Ftilde_{\infty,\gbigj}^{\ast}(\del_ah)$
or $\Ftilde_{\infty,\gbigj}^{\ast}(\del_bh)$
is not $0$.
\end{lem}
\pf
Because 
$\Ftilde_{\infty,\gbigj}^{\ast}(h)=0$,
we obtain
$\Ftilde_{\infty,\gbigj}^{\ast}(dh)=0$,
and hence the following equalities:
\begin{equation}
\label{eq;16.8.28.10}
 \Ftilde_{\infty,\gbigj}^{\ast}(\del_ah)
 \del_{t}\Ftilde_{\infty,\gbigj}^{\ast}(a)
+\Ftilde_{\infty,\gbigj}^{\ast}(\del_bh)
 \del_{t}\Ftilde_{\infty,\gbigj}^{\ast}(b)
+
 \Ftilde_{\infty,\gbigj}^{\ast}(\del_rh)
=0,
\end{equation}
\begin{equation}
\label{eq;16.8.28.11}
 \Ftilde_{\infty,\gbigj}^{\ast}(\del_ah)
 \del_{\phi}\Ftilde_{\infty,\gbigj}^{\ast}(a)
+\Ftilde_{\infty,\gbigj}^{\ast}(\del_bh)
 \del_{\phi}\Ftilde_{\infty,\gbigj}^{\ast}(b)
+
 \Ftilde_{\infty,\gbigj}^{\ast}(\del_{\theta}h)
=0.
\end{equation}
By the assumption,
at least one of 
$\Ftilde_{\infty,\gbigj}^{\ast}(\del_{\kappa}h)$
$(\kappa=a,b,\theta,r)$ is not $0$.
If both
$\Ftilde_{\infty,\gbigj}^{\ast}(\del_ah)=0$
and 
$\Ftilde_{\infty,\gbigj}^{\ast}(\del_bh)=0$ hold,
then we obtain 
$\Ftilde_{\infty,\gbigj}^{\ast}(\del_{\kappa}h)=0$
$(\kappa=\theta,r)$
by the above equalities
(\ref{eq;16.8.28.10})
and 
(\ref{eq;16.8.28.11}).
Hence, we obtain the claim of the lemma.
\hfill\qed

\vspace{.1in}

Set $\gminik$ denote the minimum of
$\ord_t\Ftilde_{\infty,\gbigj}^{\ast}(\del_ah)$
and 
$\ord_t\Ftilde_{\infty,\gbigj}^{\ast}(\del_bh)$.
Note that
for the expansion
$g_{\vecY_m}=
 \sum_{\gminiy}g_{\vecY_m,\gminiy}(u_m)
 e^{\sqrt{-1}\phi\gminiy}t^{\gminiy}$,
the coefficients
$g_{\vecY_m,\gminiy}$ $(\gminiy<\kappa(\vecY_m))$
are constants and independent of $m$.

\begin{lem}
Let $m$ be any positive integer
such that $\kappa(\vecY_m)>\gminik$.
Then, there exist a non-empty interval
$\gbigj_m\subset\gbigj$,
a small neighbourhood $\nbigu_m$
and $\epsilon_m>0$
such that 
either one of 
$\Ftilde_{m,\,\nbigu_m\times\gbigj_m}^{\ast}(\del_ah)$
or 
$\Ftilde_{m,\,\nbigu_m\times\gbigj_m}^{\ast}(\del_bh)$
is nowhere vanishing on
$\nbigu_m\times\gbigj_m\times\{0<t<\epsilon_m\}$.
\end{lem}
\pf
Note that
$\Ftilde_{m,\nbigu_m\times\gbigj_m}^{\ast}(\del_{\kappa}h)
\equiv
\Ftilde_{\infty,\gbigj_m}^{\ast}(\del_{\kappa}h)$
$(\kappa=a,b)$
modulo $t^{\kappa(\vecY_m)}$.
Hence, we obtain either one of
$\ord_t\Ftilde_{m,\nbigu_m\times\gbigj_m}^{\ast}(\del_ah)
=\ord_t\Ftilde_{\infty,\gbigj_m}^{\ast}(\del_ah)=\gminik$
or 
$\ord_t\Ftilde_{m,\nbigu_m\times\gbigj_m}^{\ast}(\del_bh)
=\ord_t\Ftilde_{\infty,\gbigj_m}^{\ast}(\del_bh)=\gminik$.
Then, we obtain the claim of the lemma.
\hfill\qed

\vspace{.1in}

Note that 
$g_{\vecY_m,\kappa(\vecY_m)}(u_m)-g_{\vecY_m,\kappa(\vecY_m)}(0)$
are $\cnum$-linear functions of $u_m$.
Let $(a_m,b_m)$ be as in \S\ref{subsection;16.10.11.20}.
We obtain the following.
\begin{lem}
\label{lem;16.8.28.30}
Let $m$ be a large integer 
such that $\kappa(\vecY_m)>\gminik$.
Then, there exist a non-empty interval
$\gbigj_m\subset\gbigj$,
a small neighbourhood $\nbigu_m$
and $\epsilon_m>0$
such that 
either $\del_{a_m}\Ftilde_{m,\nbigu_m\times\gbigj_m}^{\ast}(h)$
or $\del_{b_m}\Ftilde_{m,\nbigu_m\times\gbigj_m}^{\ast}(h)$
is nowhere vanishing on
$\nbigu_m\times\gbigj_m\times\{0<t<\epsilon_m\}$.
As a result,
either one of the following holds.
\begin{itemize}
\item
$\Ftilde_{m,\nbigu_m\times\gbigj_m}^{\ast}(h)^{-1}(0)
 \cap \{t>0\}
\simeq
 \{|a_m|<\delta_{1,m}\}
 \times
 \gbigj_m
 \times\{0<t<\epsilon_m\}$
 induced by
 $u_m\longmapsto a_m$.
\item
$\Ftilde_{m,\nbigu_m\times\gbigj_m}^{\ast}(h)^{-1}(0)
 \cap \{t>0\}
\simeq
 \{|b_m|<\delta_{2,m}\}
 \times
 \gbigj_m
 \times\{0<t<\epsilon_m\}$
 induced by
 $u_m\longmapsto b_m$.
\hfill\qed
\end{itemize}
\end{lem}

\subsubsection{Proof of Theorems \ref{thm;16.10.11.30}
and \ref{thm;16.10.11.31}}
\label{subsection;16.10.11.2}

Let $\phi_{\lambda}:
 W_{\lambda}\lrarr \nbigy$
$(\lambda\in\Lambda)$
be a rectilinearization of $f$.
(See \S\ref{subsection;16.7.26.30}
 for $\nbigy$.)
We may assume that 
the sets $\phi_{\lambda}^{-1}(\gbigh)$
are rectilinearized,
i.e., they are the union of some tuples of 
quadrants $\gbigq\subset W_{\lambda}$
with $\dim \gbigq\leq 3$.
Let $\gbigq_{\lambda,q}\subset W_{\lambda}$
be a $4$-dimensional quadrant such that
$\gbigq_{\lambda,q}
 \subset
 \phi_{\lambda}^{-1}(\Xtilde_0(H_0)\setminus \gbigh)$.
We obtain the expression
\[
 \phi_{\lambda}^{\ast}(f)_{|\gbigq_{\lambda,q}}=
 a(\lambda,q)\cdot \prod_{p=1}^4
 (\pm x_p)^{\ell_p(\lambda,q)/\rho(\lambda,q)},
\]
where 
$a(\lambda,q)$ is a nowhere vanishing ramified real analytic function,
$\rho(\lambda,q)$ is a positive integer,
and the signature of $\pm x_p$
are chosen so that
$\pm x_p> 0$ on $\gbigq_{\lambda,q}$.
Either one of the following holds:
(i)
$\ell_p(\lambda,q)\geq 0$ for any $p$,
or 
(ii) 
$\ell_p(\lambda,q)\leq 0$ for any $p$.
If $\del \gbigq_{\lambda,q}\cap
 \phi_{\lambda}^{-1}(X_0\setminus H_0)
 \neq\emptyset$,
then
$\phi_{\lambda}^{\ast}(f)_{|\gbigq_{\lambda,q}}$
is bounded around any point of 
$\del \gbigq_{\lambda,q}\cap 
 \phi_{\lambda}^{-1}(X_0\setminus H_0)$.
Hence, 
$\phi_{\lambda}^{\ast}(f)_{|\gbigq_{\lambda,q}}$
extends to a continuous function
on 
$\overline{\gbigq_{\lambda,q}}\setminus
 \phi_{\lambda}^{-1}(X_0\setminus H_0)$,
and there exists the restriction of
$\phi_{\lambda}^{\ast}(f)_{|\gbigq_{\lambda,q}}$
to any quadrant contained in
$\del \gbigq_{\lambda,q}\cap
 \phi_{\lambda}^{-1}(X_0\setminus H_0)$.

By Theorem \ref{thm;16.5.4.1},
there exist $\lambda_0$,
an interval $\gbigj_1\subset \gbigi$,
and $\nbign_{\nbigu_m\times\gbigj_1}$-paths
$\Fhat_m$ in $W_{\lambda_0}$
for any large $m$,
which are lifts of the restriction of $\Ftilde_m$.
There also exists 
the $\nbign_{\gbigj_1}$-path
$\Fhat_{\infty,\gbigj_1}$ in $W_{\lambda_0}$
which is the lift of the restriction of
$\Ftilde_{\infty}$.
For any $N_1>0$,
there exists $m_0$
such that the following holds for any $m\geq m_0$
and for any $(R_m,P)\in\nbigu_m\times\gbigj_1$:
\begin{equation}
\label{eq;16.6.3.10}
 \ord_{P,t}\Bigl(
 \Fhat_{\infty,\gbigj_1|P},
 \Fhat_{m|(R_m,P)}
 \Bigr)
>N_1.
\end{equation}

We set 
$I_{t^{1/\vecdelta(m)}}=\{0\leq t^{1/\vecdelta(m)}\leq \epsilon_m\}$.
We regard
the $\nbign_{\nbigu_m\times\gbigj_1}$-path
$\Ftilde_{m|\nbigu_m\times\gbigj_1}$
as a real analytic map
\[
 \Ftilde_{m,\nbigu_m\times\gbigj_1}:
  \nbigu_m\times \gbigj_1\times I_{t^{1/\vecdelta(m)}}
\lrarr
 \Xtilde_0(H_0).
\]
Let $\gbigh^{(m)}
 \subset
 \nbigu_m\times \gbigj_1\times I_{t^{1/\vecdelta(m)}}$
be the pull back of 
$\gbigh$ by $\Ftilde_{m,\nbigu_m\times\gbigj_1}$.
Note $\dim_{\real} \gbigh^{(m)}\leq 3$.

Let $W_{\lambda_0}^{\langle \ell\rangle}$
denote the union of the quadrants $\gbigq$ 
of $W_{\lambda_0}$ with $\dim \gbigq\leq \ell$.
Let us consider the following three cases.
\begin{description}
 \item[(a1)] 
 $\Fhat_{\infty,\gbigj_1}$ factors through
 $W_{\lambda_0}^{\langle 2\rangle}$.
 \item[(a2)]
 $\Fhat_{\infty,\gbigj_1}$ factors through
 $W_{\lambda_0}^{\langle 3\rangle}$,
 but does not factor through
 $W_{\lambda_0}^{\langle 2\rangle}$.
 \item[(a3)]
 $\Fhat_{\infty,\gbigj_1}$ does not factor through
 $W_{\lambda_0}^{\langle 3\rangle}$.
\end{description}

Let us observe that the case {\bf (a1)} does not occur.
Because the sheaf of rings $\nbign_{\gbigj_1}$ is integral,
we obtain that
$\Fhat_{\infty,\gbigj_1}$ factors through
a $2$-dimensional linear subspace of $W_{\lambda_0}$.
Then $\Fhat_{\infty,\gbigj_1}$ is convergent by 
Lemma \ref{lem;16.8.28.20},
which contradicts our assumption.
Hence, 
 $\Fhat_{\infty,\gbigj_1}$ does not factor through
$W_{\lambda_0}^{\langle 2\rangle}$.

\vspace{.1in}
Let us consider the case {\bf (a3)}.
Note that
$\Fhat_{\infty,\gbigj_1}^{\ast}(x_p)\neq 0$ for $p=1,2,3,4$.
By shrinking $\gbigj_1$,
we may assume that
$\ord_t\Fhat_{\infty,\gbigj_1|P}^{\ast}(x_p)$
are independent of $P\in\gbigj_1$.
Then, if $m$ is sufficiently large,
$\ord_t\Fhat_{m|(u_m,P)}^{\ast}(x_p)=
 \ord_t\Fhat_{\infty,\gbigj_1|P}^{\ast}(x_p)$
holds
for any $(u_m,P)\in\nbigu_m\times\gbigj_1$.
We obtain that
$\Fhat_m^{-1}(W_{\lambda_0}^{\langle 3\rangle})
\subset
 \nbigu_m\times\gbigj_1\times\{0\}$
for any sufficiently large $m$.
It implies that
$\gbigh^{(m)}\setminus
 (\nbigu_m\times\gbigj_1\times\{0\})
=\emptyset$.
Hence, the claims of Theorem \ref{thm;16.10.11.30}
holds in the case {\bf (a3)}.
Let us study the claim of Theorem \ref{thm;16.10.11.31}
in the case {\bf(a3)}.
We may assume that 
$(\Fhat_{m})^{\ast}(x_p)>0$
$(p=1,2,3,4)$
on $\nbigu_m\times\gbigj_1\times\{t^{1/\vecdelta(m)}>0\}$.
Because of (\ref{eq;16.6.3.10}),
the orders
$\ord_{t}(\Fhat_{m})^{\ast}(x_i)(u_m,\phi,t)$
are constant with respect to $(u_m,\phi)$.
Hence, 
$(\Fhat_{m})^{\ast}(x_p^{1/e})$
are ramified real analytic on
$\nbigu_m\times\gbigj_1
 \times I_{t^{1/\vecdelta(m)}}$.
We consider the expansion
\[
 \bigl(
 \Fhat_{m}
 \bigr)^{\ast}(f)
=\sum_{\gminiy\geq -M} f^{(m)}_{\gminiy}(u_m,\phi)t^{\gminiy}.
\]
By Lemma \ref{lem;16.6.3.20} below,
if $m$ is sufficiently large,
$f^{(m)}_{\gminiy}(u_m,\phi)$ 
$(\gminiy<0)$ are  independent of $u_m$.
Hence, the claim of Theorem \ref{thm;16.10.11.31}
holds with $\nbigz_m=\emptyset$
and $\nbiga_m=\emptyset$.

\vspace{.1in}

Let us study the case {\bf (a2)}.
We may assume that 
$\Fhat_{\infty,\gbigj_1}$ factors through
$\{x_1=0\}$.
Note that
$\Fhat_{\infty,\gbigj_1}^{\ast}(x_p)\neq 0$ $(p=2,3,4)$.
By shrinking $\gbigj_1$,
we may assume that 
$\gminiy(p):=\ord_t\Fhat_{\infty,\gbigj_1|P}^{\ast}(x_p)$ $(p=2,3,4)$
are independent of $P$.
We may assume that
$\Fhat_{\infty,\gbigj_1}^{\ast}(x_p)_{\gminiy(p)}>0$ $(p=2,3,4)$
for the expansion
$\Fhat_{\infty,\gbigj_1}^{\ast}(x_p)
=\sum \Fhat_{\infty,\gbigj_1}^{\ast}(x_p)_{\gminiy}t^{\gminiy}$.
There exists $m_{10}$
such that 
$\Fhat_{m}^{\ast}(x_p)>0$ $(p=2,3,4)$
for any $m\geq m_{10}$
on $\nbigu_m\times\gbigj_1
 \times(I_{t^{1/\vecdelta(m)}}\setminus\{0\})$.
Applying Lemma \ref{lem;16.8.28.30}
to the pull back of 
$\Fhat_{m_{10}}^{\ast}(x_1)$
via the induced morphisms
$\nbigu_m\times\gbigj_1\times I_{t^{1/\vecdelta(m)}}
\lrarr
\nbigu_{m_{10}}\times\gbigj_1\times 
 I_{t^{1/\vecdelta(m_{10})}}$
for any sufficiently large $m$,
we obtain that 
$\Fhat_{m}^{\ast}(x_1)^{-1}(0)
 \setminus
 (\nbigu_m\times\gbigj_1\times\{0\})
\simeq 
 \nbigv_{m,1}$
or 
$\Fhat_{m}^{\ast}(x_1)^{-1}(0)
 \setminus
 (\nbigu_m\times\gbigj_1\times\{0\})
\simeq 
 \nbigv_{m,2}$
by the map
$u_m\longmapsto a_m$
or $u_m\longmapsto b_m$.
If $\gbigh^{(m)}\setminus
 (\nbigu_{m}\times\gbigj_1\times\{0\})\neq\emptyset$,
it is equal to
$\Fhat_m^{\ast}(x_1)^{-1}(0)
\setminus
 (\nbigu_m\times\gbigj_1\times\{0\})$.
We also note that
$\Fhat_m^{-1}(W_{\lambda_0}^{\langle 2\rangle})
\subset
 \nbigu_m\times\gbigj_1\times\{0\}$.
Hence, we obtain that
the claims of Theorem \ref{thm;16.10.11.30}
hold in the case {\bf (a2)}.

Let us study the claims of Theorem \ref{thm;16.10.11.31}
in the case {\bf (a2)}.
We set
$\gbigq:=\{x_1=0,x_2>0,x_3>0,x_4>0\}$,
$\gbigq_+:=\{x_1>0,x_2>0,x_3>0,x_4>0\}$
and 
$\gbigq_-:=\{x_1<0,x_2>0,x_3>0,x_4>0\}$.
We set 
$f_{\pm}:=\phi_{\lambda_0}^{\ast}(f)_{|\gbigq_{\pm}}$.
We may naturally regard $f_{\pm}$
as the functions on $\gbigq_{\pm}\cup \gbigq$.
In particular, the restrictions $f_{\pm|\gbigq}$ are well defined.
Let $q_{\pm}:\gbigq_{\pm}\lrarr \gbigq$ be the projection
forgetting $x_1$.
Then, we obtain the functions
$q_{\pm}^{\ast}(f_{\pm|\gbigq})$.
We express
$f_{\pm}-q_{\pm}^{\ast}(f_{\pm|\gbigq})
=x_1^{\gamma_{\pm}}\cdot
 h_{\pm}
 \cdot \prod_{p=2}^4x_p^{\ell_p/\rho}$
for some $\gamma_{\pm}>0$,
where $h_{\pm}$
is a ramified real analytic function on 
the closure of $\gbigq_{\pm}$.
We may assume that
$\ord_t \Fhat_{m}^{\ast}(x_1)$ 
is sufficiently large.
Hence, we obtain that
$\Fhat_{m}^{\ast}(f_{\pm})
-\Fhat_m^{\ast}(q_{\pm}^{\ast}f_{\pm|\gbigq})$
are bounded.
By an argument as in the case {\bf (a3)},
we obtain that 
$\Fhat_m^{\ast}(q_{\pm}^{\ast}f_{\pm|\gbigq})$
is ramified real analytic,
and the coefficients of $t^{\gminiy}$ $(\gminiy<0)$
depend only on $\theta_m$.

Suppose
$\gbigh^{(m)}\setminus
 (\nbigu_m\times\gbigj_1\times\{0\})\neq\emptyset$.
The quadrant
$\gbigq=\{x_1=0\}\cap\bigcap_{i=2,3,4}\{x_i>0\}$ of $W_{\lambda_0}$
is contained in $\phi_{\lambda_0}^{-1}(\gbigh)$.
We obtain the continuous subanalytic function
$\phi_{\lambda_0}^{\ast}(f_{\gbigh})$ on 
$(\gbigq,W_{\lambda_0})$.
Note $W_{\lambda_0}=\real_{x_1}\times\real^3_{x_2,x_3,x_4}$.
By applying Theorem \ref{thm;16.7.20.2}
to $\gbigq\subset\real^3_{x_2,x_3,x_4}$ with 
$f_{\gbigq}:=\phi_{\lambda_0}^{\ast}(f_{\gbigh})_{|\gbigq}$,
we can obtain a locally finite family of 
real analytic morphisms
$\psi_{1,\beta}:V_{1,\beta}\lrarr \real^3_{x_2,x_3,x_4}$
$(\beta\in\Gamma)$
which is a rectilinearization of $\gbigq$,
such that the following holds:
\begin{itemize}
 \item Let $\gbigq'\subset V_{1,\beta}$ 
 be any quadrant contained in
 $\psi_{1,\beta}^{-1}(\gbigq)$.
 Then, $\psi_{1,\beta}^{\ast}(f_{\gbigq})_{|\gbigq'}\in \RNC(\gbigq')$.
\end{itemize}
We set $V_{\beta}:=\real_{x_1}\times V_{1,\beta}$
and $\psi_{\beta}:=\id\times \psi_{1,\beta}$.
We obtain a rectilinearization
$\{(V_{\beta},\psi_{\beta})\,|\,\beta\in\Gamma\}$
of $\gbigq\subset\real^4$
such that the following holds:
\begin{itemize}
 \item Let $\gbigq'\subset V_{\beta}$ 
 be any quadrant contained in
 $\psi_{\beta}^{-1}(\gbigq)$.
 Then, $\psi_{\beta}^{\ast}(f_{\gbigq})_{|\gbigq'}\in \RNC(\gbigq')$.
\end{itemize}
By using the arguments in the proof of
Theorem \ref{thm;16.5.4.1}
and Proposition \ref{prop;16.5.6.10},
after making $m$ larger,  shrinking $\gbigi_1$ and $\nbigu_m$,
and making $\epsilon_m$ smaller,
we can obtain $\beta_0\in\Gamma$
and morphisms
$\Fhat'_{m}:
 \nbigu_m\times\gbigi_1\times I_{t^{1/\vecdelta(m)}}
\lrarr
 V_{\beta_0}$ for any sufficiently large $m$
such that 
$\psi_{\beta_0}\circ\Fhat_{m}'=\Fhat_m$.
On $\gbigh^{(m)}\setminus(\nbigu_m\times\gbigi_1\times\{0\})$,
$\Ftilde_m^{\ast}(f_{\gbigh})
=(\Fhat'_m)^{\ast}(\psi_{\beta_0}(f_{\gbigq}))$
holds.
By the argument in the case of {\bf (a3)},
we obtain that it is ramified real analytic
as a function on $\nbigv_{1,m}$ or $\nbigv_{2,m}$,
and the coefficient of $t^{\gminiy}$ $(\gminiy<0)$
depend only on $\theta_m$.

Suppose
$\gbigh^{(m)}\setminus(\nbigu_m\times\gbigj_1\times\{0\})
=\emptyset$.
Set 
$\nbiga_m:=\Fhat_m^{-1}(\overline{\gbigq})$.
By an argument as in the case {\bf(a3)},
we obtain that 
the restriction of 
$ \Fhat_m^{\ast}\bigl(
 f_{|\gbigq}
 \bigr)$
on 
$\nbiga^{\circ}_m:=
\Fhat_m^{-1}(\gbigqbar)\setminus
 (\nbigu_m\times\gbigj_1\times\{0\})$
is ramified real analytic
as a function on $\nbigv_{1,m}$ or $\nbigv_{2,m}$,
and that the coefficients of $t^{\gminiy}$ $(\gminiy<0)$
depend only on $\theta_m$.
Hence, the claims of Theorem \ref{thm;16.10.11.31}
holds in the case {\bf (a2)}.
\hfill\qed

\subsubsection{Appendix}
\label{subsection;16.4.11.2}

Let $U$ be an open subanalytic subset 
in a real analytic manifold $M$.
Let $\nbigh$ be a hypersurface of $M$ 
defined as $\{h=0\}$ for a real analytic function $h$ on $M$.
Let $F$ be a subanalytic function on 
$(U\setminus \nbigh,M)$.
Let $Y$ be a real analytic manifold 
equipped with a global coordinate system $(x_1,\ldots,x_n)$,
and let $\Phi:Y\lrarr U$ be a real analytic morphism
such that
$\Phi^{\ast}(h)=a\cdot\prod_{i=1}^k x_i^{m_i}$
and 
$\Phi^{\ast}(F)=b\cdot\prod_{i=1}^k x_i^{\ell_i}$,
where $a$ is a nowhere vanishing real analytic function,
$b$ is a nowhere vanishing ramified real analytic function,
$k$ is an integer such that $1\leq k\leq n$,
and $(m_i)\in\seisuu_{>0}^k$
and $(\ell_i)\in(\rnum_{\geq 0})^k\cup(\rnum_{\leq 0})^k$.

Let $P_0$ denote a one point set.
Let $\gamma'$ be an $\nbign_{P_0}$-path in $Y$
such that $(\gamma')_0(P_0)=(0,\ldots,0)$.
Let $\gamma=\Phi\circ\gamma'$ be 
the induced $\nbign_{P_0}$-path in $U$.

\begin{lem}
\label{lem;16.6.3.20}
There exists an integer $N_0>0$
depending only on $(\ell_i)$
and $\ord\gamma^{\ast}(h)$
such that the following holds.
\begin{itemize}
\item
 Let $\sigma(t)$ be an $\nbign_{P_0}$-path
in $Y$ such that
$\sigma_0(P_0)=(0,\ldots,0)$
and that 
 $\ord_t(\sigma,\gamma')\geq N_0$.
Then, the polar parts of
 $(\gamma')^{\ast}\Phi^{\ast}F$
and $\sigma^{\ast}\Phi^{\ast}F$
are the same.
\end{itemize}
\end{lem}
\pf
It is enough to consider the case
$(\ell_i)\in\rnum_{\leq 0}^n$.
We describe
$\gamma'(t)=
 \bigl(\gamma'_1(t),\ldots,\gamma'_n(t)\bigr)$.
Because
$\gamma^{\ast}(h)=
 (\gamma')^{\ast}\Phi^{\ast}(h)$,
$\ord_t(\gamma_i')\leq
 \ord_t(\gamma^{\ast}(h))$
holds.
We set $d_i:=\ord_t(\gamma_i')$
and $d:=\ord_t(\gamma^{\ast}h)$.
Let $N_0$ be any integer larger than
$\bigl(n+\sum_i |\ell_i|\bigr)d$.
For the description
$\sigma(t)=(\sigma_1(t),\ldots,\sigma_n(t))$,
$\sigma_i(t)=
 \gamma'_i(t)+t^{N_0}b_i(t)$ holds.
We obtain
\[
 \prod_i
 \Bigl(
 \gamma'_i(t)+t^{N_0}b_i(t)
 \Bigr)^{\ell_i}
=\prod_i\gamma'_i(t)^{\ell_i}
 (1+t^{N_0-d_i}c_i(t))^{\ell_i}
=\prod_i
 \gamma'_i(t)^{\ell_i}
 (1+t^{N_0-d_i}c_i'(t)),
\]
for some power series
$c_i(t)$ and $c_i'(t)$.
Then, the claim of the lemma is clear.
\hfill\qed

\part{Enhanced ind-sheaves and holonomic $\nbigd$-modules}
\label{part;18.11.16.30}

\section{Meromorphic flat bundles and enhanced ind-sheaves}
\label{section;18.11.15.50}

Let $X$ be a complex manifold
with a complex hypersurface $H$.
We define the bordered space
$\vecX(H):=(X\setminus H,X)$.
The dimension of $X$ is denoted by $d_X$
unless otherwise specified.

In \S\ref{subsection;18.11.16.120},
we introduce two subcategories 
$\Ecat^b_{\mero}(\IC_{\vecX(H)})$ and 
$\Ecat^b_{\circledcirc}(\IC_{\vecX(H)})$
of $\Ecat^b_{\realc}(\IC_{\vecX(H)})$.
The subcategory 
$\Ecat^b_{\mero}(\IC_{\vecX(H)})$
is defined as the essential image of
the category of meromorphic flat bundles
by the enhanced de Rham functor.
The subcategory 
$\Ecat^b_{\circledcirc}(\IC_{\vecX(H)})$
is defined by the curve test.
Theorem \ref{thm;18.11.12.10} states
that they are actually the same,
that is one of the main theorems in this paper.
The proof will be given in 
\S\ref{section;18.11.15.20}--\ref{section;18.11.15.21}.

In \S\ref{subsection;18.11.16.121},
we shall explain another way 
to define the subcategory $\Ecat^b_{\circledcirc}(\vecX(H))$.
In \S\ref{subsection;16.7.25.30},
we shall give a preliminary for the proof of
Theorem \ref{thm;18.11.12.10}.

Let us introduce some notation.
For any $\epsilon>0$,
we set $\Delta_{\epsilon}:=
 \bigl\{z\in\cnum\,\big|\,|z|<\epsilon\bigr\}$
and $\Delta^{\ast}_{\epsilon}:=\Delta\setminus\{0\}$.
Let $\vecDelta_{\epsilon}(0)$
denote the bordered space
$(\Delta_{\epsilon}^{\ast},\Delta_{\epsilon})$.
Let $\varpi_0:\Deltatilde_{\epsilon}(0)\lrarr\Deltatilde_{\epsilon}$
denote the oriented real blowing up.
For any subsets $S_i\subset X$ $(i=1,2)$,
a holomorphic map
$\varphi:(\Delta_{\epsilon}^{\ast},0,\Delta_{\epsilon})
 \lrarr (S_1,S_2,X)$
means 
a holomorphic map
$\varphi:\Delta_{\epsilon}\lrarr X$
such that
$\varphi(\Delta_{\epsilon}^{\ast})\subset S_1$
and 
$\varphi(0)\in S_2$.
If $\epsilon=1$,
we omit to denote $\Delta$ instead of $\Delta_{1}$.
We use similar notation $\Delta^{\ast}$, $\vecDelta(0)$, etc.

\subsection{Curve test for enhanced ind-sheaves and 
meromorphic flat bundles}
\label{subsection;18.11.16.120}

Let $\Hol(\nbigd_X)$ denote the category of
holonomic $\nbigd_X$-modules.
Let $\Mero_{\ast}(X,H)\subset \Hol(\nbigd_X)$ 
be the full subcategory of 
holonomic $\nbigd_X$-modules $M$
which are coherent over $\nbigo_X(\ast H)$.
According to 
\cite{DAgnolo-Kashiwara1},
the enhanced de Rham functor
$\DR^{\Ecat}_{\vecX(H)}:
 \Mero_{\ast}(X,H)\lrarr
 \Ecat^b_{\realc}(\IC_{\vecX(H)})$
is a fully faithful functor.
\begin{df}
Let $\Ecat^b_{\mero}(\IC_{\vecX(H)})
 \subset \Ecat^b_{\realc}(\IC_{\vecX(H)})$
denote the essential image of
$V \longmapsto
 \DR^{\Ecat}_{\vecX(H)}(V)[-d_X]$.
\hfill\qed
\end{df}

Let $\Mero_!(X,H)\subset\Hol(\nbigd_X)$
denote the essential image of
$\DDD:\Mero_{\ast}(X,H)^{\op}\lrarr
 \Hol(\nbigd_X)$,
where $\DDD$ denotes the duality functor
on $\Hol(\nbigd_X)$,
and $\Mero_{\ast}(X,H)^{\op}$
denotes the opposite category of
$\Mero_{\ast}(X,H)$.
We may also regard
$\Ecat^b_{\mero}(\IC_{\vecX(H)})$
as the essential image of the functor
$\Mero_!(X,H)
\lrarr
 \Ecat^b_{\realc}(\IC_{\vecX(H)})$
defined by
$M\longmapsto 
 \DR^{\Ecat}_{\vecX(H)}(M)[-d_X]$.

Let $j:(X\setminus H,X)\lrarr (X,X)$ be the natural inclusion
of the bordered spaces.
There exists the fully faithful functor
$\Ecat j_{!!}:
\Ecat^b_{\realc}(\IC_{\vecX(H)})
\lrarr \Ecat^b_{\realc}(\IC_{X})$.
The essential image of the induced functor
$\Ecat j_{!!}:
 \Ecat^b_{\mero}(\IC_{\vecX(H)})
\lrarr
 \Ecat^b_{\realc}(\IC_{X})$
is equal to the essential image of
the functor
$\Mero_!(X,H)
\lrarr
 \Ecat^b_{\realc}(\IC_{X})$
defined by 
$M\longmapsto 
\DR^{\Ecat}_{X}(M)[-d_X]$.

We introduce another subcategory of
$\Ecat^b_{\realc}(\IC_{\vecX(H)})$.

\begin{df}
Let
$\Ecat^b_{\circledcirc}(\IC_{\vecX(H)})\subset
 \Ecat^b_{\realc}(\IC_{\vecX(H)})$
be the full subcategory
of objects $K$ satisfying the following conditions.
\begin{itemize}
\item
 $K_{|X\setminus H}$ is 
 a locally free $\cnum_{X\setminus H}$-module.
\item
For any holomorphic map
$\varphi:(\Delta^{\ast},0,\Delta)\lrarr (X\setminus H,H,X)$,
 $\Ecat\varphi^{-1}K$
is an object in
$\Ecat^b_{\mero}(\IC_{\vecDelta(0)})$.
\hfill\qed
\end{itemize}
\end{df}

By the compatibility of the enhanced de Rham functor
and the $6$-functors for 
holonomic $\nbigd$-modules and
$\real$-constructible enhanced ind-sheaves
in \cite{DAgnolo-Kashiwara1},
$\Ecat^b_{\mero}(\IC_{\vecX(H)})$
is a full subcategory of
$\Ecat^b_{\circledcirc}(\IC_{\vecX(H)})$.
If $\dim_{\cnum} X=1$,
$\Ecat^b_{\mero}(\IC_{\vecX(H)})
=\Ecat^b_{\circledcirc}(\IC_{\vecX(H)})$
clearly holds.

We shall prove the following theorem 
in \S\ref{section;18.11.15.20}--\S\ref{section;18.11.15.21}.

\begin{thm}
\label{thm;18.11.12.10}
\label{thm;16.5.11.20}
$\Ecat^b_{\mero}(\IC_{\vecX(H)})
=\Ecat^b_{\circledcirc}(\IC_{\vecX(H)})$.
\end{thm}

\subsection{Another description}
\label{subsection;18.11.16.121}

\subsubsection{Subcategory determined by sector test}

Let $S_{\epsilon,\delta}\subset\Deltatilde_{\epsilon}(0)$ 
denote the sector
$\{(r,\theta)\in\Deltatilde_{\epsilon}(0)\,|\,
 -\delta<\theta<\delta \}$.
Set $Z_{\delta}:=\{(0,\theta)\,|\,-\delta<\theta<\delta\}
\subset S_{\epsilon,\delta}$,
and 
$S^{\circ}_{\epsilon,\delta}:=S_{\epsilon,\delta}\setminus Z_{\delta}$.
Note that
$S^{\circ}_{\epsilon,\delta}=
 \{z\,|\,0<|z|<\epsilon,
 -\delta<\arg(z)<\delta\}
\subset\cnum$.
Let $\vecS^{\circ}_{\epsilon,\delta}$
denote the bordered space
$(S^{\circ}_{\epsilon,\delta},S_{\epsilon,\delta})$.

For any real analytic map
$\gamma:(\II^{\circ},0,\II)\lrarr (X\setminus H,H,X)$,
if $\epsilon$ is sufficiently small,
we obtain the induced holomorphic map
$\gamma_{\cnum}:
 \Delta_{\epsilon}\lrarr X$.
It induces the real analytic map of bordered spaces
$\gammatilde_{\epsilon,\delta}:
 \vecS^{\circ}_{\epsilon,\delta}
\lrarr \vecX(H)$.

\begin{df}
Let $\Ecat^b_{\circledcirc'}(\IC_{\vecX(H)})\subset
 \Ecat^b_{\realc}(\IC_{\vecX(H)})$
be the full subcategory of the objects $K$
which have the following property:
\begin{itemize}
\item
$K_{|X\setminus H}$ is 
 a locally free $\cnum_{X\setminus H}$-module.
\item
Let $\gamma:(\II^{\circ},0,\II)\lrarr (X\setminus H,H,X)$
be any real analytic map.
Then, there exist 
small positive numbers $\delta$ and $\epsilon$,
and $\cnum$-valued subanalytic functions
$g_{1},\ldots,g_m$ on 
$(S^{\circ}_{\epsilon,\delta},\Deltatilde(0))$
which are holomorphic as functions on $S^{\circ}_{\epsilon,\delta}$,
such that
\[
 \Ecat\gammatilde_{\epsilon,\delta}^{-1}K
\simeq
 \bigoplus_{i=1}^m
 \cnum_{S_{\epsilon,\delta}}^{\Ecat}\otimes
 \cnum_{t\geq \Re(g_i)}
\]
in
$\Ecat^b_{\realc}(\IC_{\vecS^{\circ}_{\epsilon,\delta}})$.
\hfill\qed
\end{itemize}
\end{df}

The following is clear by the definitions.
\begin{lem}
$\Ecat^b_{\circledcirc'}(\IC_{\vecX(H)})$
is a full subcategory of
$\Pro_{tf}(\vecX(H))$.
\hfill\qed
\end{lem}

\subsubsection{Basic isomorphisms}

Let $f$ be a meromorphic function on $(X,H)$.
We set
$\nbige^f_X(\ast H):=
 \bigl(\nbigo_X(\ast H),d+df\bigr)$
in $\Mero_{\ast}(X,H)$,
and
$\nbige^f_X(!H):=
 \DDD_X
 \bigl(
 \nbige^{-f}_X(\ast H)
 \bigr)$,
where $\DDD_X$ denotes 
the duality functor on $\Hol(X)$.
According to \cite{DAgnolo-Kashiwara1},
there exist natural isomorphisms
in $\Ecat^b_{\realc}(\IC_X)$:
\begin{equation}
\label{eq;16.6.10.1}
 \DR^{\Ecat}_X\bigl(
 \nbige^f_X(\ast H)
 \bigr)[-d_X]
\simeq
 \cnum_X^{\Ecat}
 \overset{+}{\otimes}
 R\nihom\bigl(
 \pi^{-1}(\cnum_{X\setminus H}),
 \cnum_{t=\Re(f)}
 \bigr),
\end{equation}
\begin{equation}
\label{eq;16.6.10.2}
 \DR_X^{\Ecat}\bigl(
 \nbige^f_X(! H)
 \bigr)[-d_X]
\simeq
 \cnum_X^{\Ecat}
 \overset{+}{\otimes}
 \cnum_{t\geq \Re(f)}.
\end{equation}

\subsubsection{Comparison of curve test and sector test}

\begin{prop}
\label{prop;16.6.10.3}
The subcategories
$\Ecat^b_{\circledcirc}(\IC_{\vecX(H)})$ and 
$\Ecat^b_{\circledcirc'}(\IC_{\vecX(H)})$
are the same.
In particular,
if $\dim X=1$,
$\Ecat^b_{\mero}(\IC_{\vecX(H)})
=\Ecat^b_{\circledcirc'}(\IC_{\vecX(H)})$
holds.
\end{prop}
\pf
We begin with a lemma for functions
on one dimensional sectors.

\begin{lem}
\label{lem;16.1.25.10}
Let $\nbigs:=
 \bigl\{z\in\cnum\setminus\{0\}
 \,\big|\,
 0<|z|<\epsilon,\,\,
 |\arg(z)|<\delta\bigr\}$.
Let $f$ be a subanalytic $\cnum$-valued function on
$(\nbigs,\cnum)$.
Suppose that $f$ is holomorphic on $\nbigs$.
Then, there exist
a neighbourhood $U$ of $0$,
a positive integer $\rho$,
and a meromorphic function $g$ on $(U,0)$
such that
$f_{|U\cap \nbigs}
=g(z^{1/\rho})_{|U\cap\nbigs}$.
\end{lem}
\pf
Let $\varpi:\cnumtilde(0)\lrarr \cnum$
be the oriented real blowing up.
We may regard $\nbigs$
as an open subset in $\cnumtilde(0)$.
We may regard $f$ as a subanalytic function
on $(\nbigs,\cnumtilde(0))$.
Let $\overline{\nbigs}$ denote 
the closure of
$\nbigs$ in $\cnumtilde(0)$.
We set $Z=\overline{\nbigs}\cap\varpi^{-1}(0)$.
Let $Z^{\circ}$ be the interior part of
$Z\subset\varpi^{-1}(0)$.
There exists a sufficiently large integer $N$
such that
$z^{N}f$ extends to a continuous function 
on $\nbigs\cup Z^{\circ}$,
and 
$(z^{N}f)_{|Z^{\circ}}=0$.
By Lemma \ref{lem;15.12.28.10},
there exists a discrete subset
$Z'\subset Z^{\circ}$
and a positive integer $\rho_0>0$
such that any $\theta_0\in Z\setminus Z'$
has a neighbourhood $\nbigu_0$
of $(0,\theta_0)$ in $\cnumtilde(0)$
on which
$(z^Nf)_{|\nbigu_0}$
is expressed by a convergent power series
for some $a_{i,j}\in\cnum$:
\[
(z^Nf)_{|\nbigu_0}
=\sum_{j>0}\sum_{i\geq 0}
 a_{i,j}(\theta-\theta_0)^i\cdot r^{j/\rho_0}.
\]
Because $(z^Nf)_{|\nbigu_0\setminus\varpi^{-1}(0)}$ is holomorphic,
$(z^Nf)_{|\nbigu_0}
=\sum_{j>0} b_jz^{j/\rho_0}$
holds for some $b_j\in\cnum$.
Then, the claim of the lemma follows.
\hfill\qed

\vspace{.1in}

To prove Proposition \ref{prop;16.6.10.3},
it is enough to consider the case
$X=\Delta$ and $H=\{0\}$.
It follows from the following lemma.

\begin{lem}
\label{lem;16.1.25.100}
An object $K$ of $\Ecat^b_{\realc}(\IC_{\vecDelta(0)})$
is contained in $\Ecat^b_{\mero}(\IC_{\vecDelta(0)})$
if and only if there exist the following:
\begin{itemize}
\item
a finite covering
$\Delta\setminus \{0\}=\bigcup_{i=1}^N\nbigs_i$,
where
$\nbigs_i=
 \bigl\{z\in\cnum\setminus\{0\}\,\big|\,
 0<|z|<1,\,\,
 \theta_{i,1}<\arg(z)<\theta_{i,2}
 \bigr\}$,
\item
subanalytic $\cnum$-valued functions
$g_{ij}$ $(j=1,\ldots,m)$ on $(\nbigs_i,\cnum)$
such that 
$g_{ij}$ are holomorphic on $\nbigs_i$,
\item
isomorphisms
$\pi^{-1}(\cnum_{\nbigs_i})\otimes K
\simeq
 \bigoplus_{j=1}^m
 \pi^{-1}(\cnum_{\nbigs_i})\otimes
 \Bigl(
 \cnum_{X}^{\Ecat}\overset{+}{\otimes}
 \cnum_{t\geq \Re(g_{ij})}
 \Bigr)$
in $\Ecat^b_{\realc}(\IC_{\vecDelta(0)})$.
\end{itemize}
\end{lem}
\pf
The ``only if'' part follows from
the basic isomorphisms
(\ref{eq;16.6.10.1}, \ref{eq;16.6.10.2})
and the classical asymptotic analysis for meromorphic 
flat bundles on curves.
Let us study the ``if'' part.
By shrinking $\nbigs_i$,
we may assume that
$g_{ij}$ are bounded
on $\{\epsilon'<|z|<\epsilon\,,\,\theta_{i,1}<\theta<\theta_{i,2}\}$
for any $\epsilon'>0$.
By pulling back via an appropriate ramified covering
$\psi:(\Delta,0)\lrarr(\Delta,0)$,
we may assume that
$\psi^{\ast}g_{ij}$ extend to meromorphic functions on
$(\Delta,0)$,
by Lemma \ref{lem;16.1.25.10}.
We may assume that
they are elements in $z^{-1}\cnum[z^{-1}]$.

By the condition,
$K_{|\Delta^{\ast}}$ induces
a $\cnum_{\Delta^{\ast}}$-locally free sheaf
on $\Delta^{\ast}$,
which is denoted by $L$.
Let $\Ltilde$ denote the local system 
on $\Deltatilde(0)$ induced by $L$.

We set $I_{i}:=\{\Re(g_{ij})\,|\,j=1,\ldots,m\}$.
Let $\nbigsbar_i$ denote the closure of $\nbigs_i$
in $\Deltatilde(0)$.
Let $\Ibar_i$ denote the image of
$I_i\lrarr \Subbar(\nbigs_i,\nbigsbar_i)$.

As remarked in \S\ref{subsection;16.1.25.20},
there exist  the induced filtrations $\nbigf^{\nbigs_i}$
on $L_{|\nbigs_i}$
indexed by
$(\Ibar_i,\prec)$.
If $\nbigs_{i,j}:=\nbigs_i\cap\nbigs_j\neq \emptyset$,
$\Ibar_{i|\nbigs_{i,j}}
=\Ibar_{j|\nbigs_{i,j}}$ holds
in $\Subbar(\nbigs_{i,j},\nbigsbar_{i,j})$.
There also exist the induced filtrations
$\nbigf^{\nbigs_i\cap\nbigs_j}$
on $L_{|\nbigs_i\cap\nbigs_j}$
which are equal to the filtrations
induced by $\nbigf^{\nbigs_k}$ $(k=i,j)$.
Hence, 
$\Ltilde$ is equipped with a Stokes structure
$\vecnbigf$,
i.e., a family of Stokes filtrations.

Let $(V,\nabla)$ be the meromorphic flat bundle
corresponding to $(\Ltilde,\vecnbigf)$.
The isomorphisms
$\DR^{\Ecat}_{\nbigs_i}(V,\nabla)[-1]
\simeq
 K_{|\nbigs_i}$
extend to  isomorphisms
$\pi^{-1}(\cnum_{\nbigs_i})\otimes
 \DR^{\Ecat}_{\vecDelta(0)}(V,\nabla)
\simeq
 \pi^{-1}(\cnum_{\nbigs_i})
 \otimes
 K$
in $\Ecat^b_{\realc}(\IC_{\vecDelta(0)})$.
Such isomorphisms are unique 
by Lemma \ref{lem;16.1.25.32}.
Hence, we obtain an isomorphism
$\DR^{\Ecat}_{\vecDelta(0)}(V,\nabla)
 \simeq
 K$.
Thus,
Lemma \ref{lem;16.1.25.100},
and Proposition \ref{prop;16.6.10.3}
are proved.
\hfill\qed

\subsubsection{Prolongations of isomorphisms of enhanced ind-sheaves}

Let $(V,\nabla)\in\Mero_{\ast}(X,H)$.
Let $L_V$ be the local system on $X\setminus H$
corresponding to $(V,\nabla)_{|X\setminus H}$.
Let $K\in \Ecat^b_{\circledcirc}(\IC_{\vecX(H)})$
such that
$K_{|X\setminus H}$ comes from
a locally free $\cnum_{X\setminus H}$-module $L_K$.
Let $\Phi:L_V\simeq L_K$
be an isomorphism 
 of $\cnum_{X\setminus H}$-modules
such that the following holds:
\begin{itemize}
\item
For any holomorphic map
$\varphi:(\Delta^{\ast},0,\Delta)\lrarr (X\setminus H,H,X)$,
 $\varphi_{|\Delta^{\ast}}^{-1}\Phi$
 extends to an isomorphism
\[
 \DR^{\Ecat}_{\vecDelta(0)}(\varphi^{\ast}(V))[-1]
\simeq
 \Ecat\varphi^{-1}K
\]
 in
 $\Ecat^b_{\realc}(\IC_{\vecDelta(0)})$.
\end{itemize}
We obtain the following proposition
as a consequence of Proposition \ref{prop;16.6.23.21}
and Proposition \ref{prop;16.6.10.3}.

\begin{prop}
\label{prop;16.6.23.20}
$\Phi$ extends to
an isomorphism
$\DR^{\Ecat}_{\vecX(H)}(V)[-d_X]
\simeq
 K$
in $\Ecat^b_{\realc}(\IC_{\vecX(H)})$.
\hfill\qed
\end{prop}

\subsection{Auxiliary conditions}
\label{subsection;16.7.25.30}
\subsubsection{Set of ramified irregular values}

Let $P$ be any point of $H$.
Let $(X_P,z_1,\ldots,z_{d_X})$ denote 
a small coordinate neighbourhood of $P$
in $X$ such that
$H_P:=H\cap X_P=\bigcup_{i=1}^{\ell}\{z_i=0\}$
and $P=(0,\ldots,0)$.
For any positive integer $e$,
let $X_P^{(e)}$ be an open neighbourhood of $(0,\ldots,0)$
in $\cnum^n=\{(\zeta_1,\ldots,\zeta_n)\}$,
and let $\psi_{e}:X_P^{(e)}\lrarr X_P$
be the ramified covering defined by
\[
 \psi_e(\zeta_1,\ldots,\zeta_n)
=(\zeta_1^e,\ldots,\zeta_{\ell}^e,\zeta_{\ell+1},\ldots,\zeta_n).
\]
Set $H_P^{(e)}:=\psi_e^{-1}(H_P)$.
A set of ramified irregular values at $P$
is a tuple 
$(g_1,\ldots,g_m)$ 
in $\nbigo_{X^{(e)}_P}(\ast H_P^{(e)})_P\big/\nbigo_{X^{(e)}_P,P}$
for some positive integer $e$,
which is invariant
under the action of the Galois group of $\psi_e$.
If $\{g_1,\ldots,g_m\}$ is a good set of irregular 
values on $(X_P^{(e)},H_P^{(e)})$,
then it is called a good set of ramified irregular
values at $P$.
A multiplicity function on $\{g_1,\ldots,g_m\}$
is a function
$\gminim:\{g_1,\ldots,g_m\}
\lrarr \seisuu_{\geq 0}$,
which is invariant under the action of the Galois group.

Let $\varphi:(\Delta^{\ast},0,\Delta)\lrarr 
(X\setminus H,H_P,X)$ 
be any holomorphic map.
There exists $\epsilon>0$
such that 
the restriction of $\varphi$ to $\Delta_{\epsilon}$
induces 
$\varphi_{\epsilon}:
 (\Delta^{\ast}_{\epsilon},0,\Delta_{\epsilon})
\lrarr
 (X_P\setminus H_P,H_P,X)$.
For any positive integer $e$,
let $\Delta^{(e)}$ be a small neighbourhood of $0$
in $\cnum=\{\zeta\}$.
We define $\psi_{\Delta,e}:\Delta^{(e)}\lrarr\Delta_{\epsilon}$
by $\psi_{\Delta,e}(\zeta)=\zeta^e$.
Then, there exists a holomorphic map
$\varphi^{(e)}:\Delta^{(e)}\lrarr X^{(e)}_P$
such that 
$\psi_e\circ \varphi^{(e)}
=\varphi\circ\psi_{e,\Delta}$.
For any set of ramified irregular values $\nbigi$ at $P$,
we define the pull back
$\varphi^{\ast}\nbigi
 \subset\nbigo_{\Delta^{(e)}}(\ast 0)_0\big/
 \nbigo_{\Delta^{(e)},0}$ by using 
$\varphi^{(e)}$ for some $e$.
Moreover,
for any multiplicity function $\gminim:\nbigi\lrarr \seisuu_{\geq 0}$,
we define
$\varphi^{\ast}(\gminim)(f):=
 \sum_{\varphi^{(e)\ast}(g)=f}\gminim(g)$.
By the Galois invariance,
they are independent of the choice of $\varphi^{(e)}$.

\subsubsection{Some conditions}

We prepare a notation.
\begin{notation}
For any 
$K=\DR^{\Ecat}_{\vecDelta(0)}(V,\nabla)[-1]
\in\Ecat^b_{\mero}(\IC_{\vecDelta(0)})$,
we set
$\Irr(K):=\Irr(V,\nabla)$.
(See {\rm\S\ref{subsection;18.11.26.1}}
for $\Irr(V,\nabla)$.)
\hfill\qed
\end{notation}

We introduce a condition
for objects $K\in \Ecat^{b}_{\circledcirc}(\IC_{\vecX(H)})$
and $P\in H$.
\begin{condition}
\label{condition;18.11.26.2}
There exist a set of ramified irregular values
 $\nbigi_P$ at $P$,
 a multiplicity function
 $\gminim_P:\nbigi_P\lrarr
 \seisuu_{>0}$,
 and a neighbourhood $H_P$ of $P$ in $H$,
such that the following holds:
\begin{itemize}
\item
For any holomorphic map
 $\varphi:(\Delta^{\ast},0,\Delta)\lrarr 
 (X\setminus H,H_P,X)$,
 $\Irr(\Ecat\varphi^{-1}(K))
 =\varphi^{\ast}\nbigi_P$
 holds,
and the multiplicity functions are the equal.
\end{itemize}
Note that $\nbigi_P$ is not assumed to be
a good set of ramified irregular values.
\hfill\qed
\end{condition}

\begin{lem}
If Condition {\rm\ref{condition;18.11.26.2}}
is satisfied for $K\in\Ecat^b_{\circledcirc}(\IC_{\vecX(H)})$
and $P\in H$,
then the set of ramified irregular values $\nbigi_P$
and the multiplicity function $\gminim_P$
are uniquely determined.
\end{lem}
\pf
Suppose that Condition \ref{condition;18.11.26.2}
for $(K,P)$ is satisfied for $(\nbigi_P,\gminim_P)$
and $(\nbigi_P',\gminim_P')$.
There exists a positive integer $e$
such that
$\nbigi_P$ and $\nbigi_{P}'$
are contained in
$\nbigo_{X_P^{(e)}}(\ast H_P^{(e)})/\nbigo_{X_P^{(e)}}$.
If $\nbigi_P\neq\nbigi_{P'}$,
there exists a holomorphic map
$\varphi:(\Delta^{\circ},0,\Delta)\lrarr 
 (X_P^{(e)}\setminus H_P^{(e)},H_P^{(e)},X_P^{(e)})$
such that
$\varphi^{\ast}\nbigi_P\neq\varphi^{\ast}\nbigi_P'$.
It contradicts 
$\varphi^{\ast}(\nbigi_P)
=\Irr(\Ecat\varphi^{-1}\Ecat\psi_e^{-1}K)
=\varphi^{\ast}(\nbigi_P')$.
The multiplicity functions
$\gminim_P$ and $\gminim_P'$ are compared similarly.
\hfill\qed

\vspace{.1in}

Let $P$ be any point of $H$.
For a positive integer $e$,
let $X_P^{(e)}$, $H_P^{(e)}$ and $\psi_e$
be as above.
Set $P^{(e)}:=\psi_e^{-1}(P)$.

\begin{lem}
\label{lem;18.11.6.1}
Let $K\in\Ecat^b_{\circledcirc}(\IC_{\vecX(H)})$.
Then, Condition {\rm\ref{condition;18.11.26.2}}
is satisfied for $(K,P)$
if and only if
Condition {\rm\ref{condition;18.11.26.2}}
is satisfied for
$(\Ecat\psi_e^{-1}(K),P^{(e)})$.
\end{lem}
\pf
Set $X_P:=\psi_e(X_P^{(e)})$
and $H_P:=\psi_e(H_P^{(e)})$.
Let $\varphi:(\Delta^{\ast},0,\Delta)\lrarr 
 (X_P\setminus H_P,H_P,X_P)$ be any holomorphic map.
There exists a holomorphic map
$\varphi^{(e)}:\Delta^{(e)}\lrarr X^{(e)}_P$
such that
$\psi_e\circ\varphi^{(e)}=\varphi\circ\psi_{e,\Delta}$.
Then, 
$\Ecat\varphi^{-1}(K)$
is contained in
$\Ecat^b_{\mero}(\IC_{\vecDelta(0)})$
if and only if
$\Ecat\psi_{e,\Delta}^{-1}\Ecat\varphi^{-1}(K)$
is contained in
$\Ecat^b_{\mero}(\IC_{\vecDelta^{(e)}(0)})$.
Then, the claim is clear.
\hfill\qed

\vspace{.1in}
The following is a stronger condition
for $K\in \Ecat^b_{\circledcirc}(\IC_{\vecX(H)})$
and $P\in H$.
\begin{condition}
\label{condition;18.11.26.3}
Condition {\rm\ref{condition;18.11.26.2}}
is satisfied for $(K,P)$,
and moreover,
$\nbigi_P$ is a good set of ramified irregular values.
\hfill\qed
\end{condition}

\begin{df}
Let $K\in \Ecat^b_{\realc}(\IC_{\vecX(H)})$.
We say that 
Condition {\rm\ref{condition;18.11.26.2}}
is satisfied for $K$,
(resp. Condition {\rm\ref{condition;18.11.26.3}})
if 
Condition {\rm\ref{condition;18.11.26.2}}
(resp. Condition {\rm\ref{condition;18.11.26.3}})
is satisfied for $K$ and any $P\in H$.
\hfill\qed
\end{df}

\subsubsection{Good smooth case}

We assume that $H$ is smooth.

\begin{prop}
\label{prop;16.7.25.40}
If $K$ satisfies Condition {\rm\ref{condition;18.11.26.3}},
then there exists a good meromorphic flat bundle
$(V,\nabla)$ on $(X,H)$
with an isomorphism
\[
 \DR^{\Ecat}_{\vecX(H)}(V)[-d_X]
\simeq
 K.
\]
\end{prop}
\pf
Let $P$ be any point of $H$.
There exists the good set of ramified irregular values
$\nbigi_P$
with the multiplicity function 
$\gminim_P:\nbigi_P\lrarr\seisuu_{> 0}$
as in Condition \ref{condition;18.11.26.2}.
We may assume that
$\nbigi_P$ is unramified,
i.e.,
elements of $\nbigi_P$
are meromorphic functions on $(X_P,H_P)$,
where $X_P$ is a small neighbourhood of $P$ in $X$,
and $H_P:=X_P\cap H$.
Set $K_P:=K_{|X_P}$.
Let $L_P$ be the local system on $X_P\setminus H_P$
corresponding to
$K_{|X_P\setminus H_P}$.
Set $n:=d_X$.
We may assume
$X_P=\Delta^{n}=\{(z_1,\ldots,z_{n})\,|\,|z_i|<1\}$,
$H_P=\{z_1=0\}$
and $P=(0,\ldots,0)$.

Let $\varphi_P:\Delta\lrarr X_P$
be the map defined by 
$\varphi_P(\zeta)=(\zeta,0,\ldots,0)$.
There exists a meromorphic flat bundle $V_{P,0}$
with an isomorphism
$\DR^{\Ecat}_{\vecDelta(0)}V_{P,0}[-1]
 \simeq
 \Ecat\varphi_P^{-1}(K)$.
Note that
$\Irr(V_{P,0})=
 \varphi_P^{\ast}\nbigi_P$ holds
and that the multiplicity functions are equal.
By using Theorem \ref{thm;17.12.27.1},
we obtain an unramifiedly good meromorphic flat bundle
$V_P$ on $(X_P,H_P)$
such that $\Irr(V_P)=\nbigi_P$
with an isomorphism
$\varphi_P^{\ast}V_P\simeq V_{P,0}$.
Such $V_P$ is uniquely determined
up to canonical isomorphisms.

Let $\varpi_P:\Xtilde_P(H_P)\lrarr X_P$
be the oriented real blowing up.
For $f,g\in\nbigi_P$ with $f\neq g$,
we obtain the $C^{\infty}$-function
$F_{f,g}:=\Re(f-g)|z_1|^{-\ord(f-g)}$
on $\Xtilde_P(H_P)$.
Set $Z$
as the union of
$F_{f,g}^{-1}(0)\cap\varpi_P^{-1}(H_P)$ for such $(f,g)$.
Note that $\dim_{\real}Z=2n-2$.

The object
$\Ecat\varpi_P^{-1}(K)$
on the bordered space
$\vecXtilde_P(H_P)=
(X_P\setminus H_P,\Xtilde_P(H_P))$
is a prolongation of $L_P$.
Let $Q$ be any point of $\varpi_P^{-1}(H_P)\setminus Z$.
Let $\nbigu_Q$ be a small neighbourhood
of $Q$ in $\Xtilde_P(H_P)$
such that
$\nbigu_Q\cap Z=\emptyset$.
We set $\nbigu_Q^{\circ}:=
 \nbigu_Q\setminus\varpi_P^{-1}(H_P)$.
By Proposition \ref{prop;16.8.25.2},
there exists a filtration $\nbigf^Q$
of $L_{P|\nbigu_Q^{\circ}}$
such that
$\Ecat\varpi_P^{-1}(K_P)_{|\nbigu_Q}$
is the stably free enhanced ind-sheaf
induced by the local system $L_{P|\nbigu_Q^{\circ}}$
with the filtration $\nbigf^Q$.
By considering the restriction to
$\varpi^{-1}_P(\Delta\times\{(0,\ldots,0)\})$,
we obtain that there exists an isomorphism
\[
 \Ecat\varpi_P^{-1}\DR^{\Ecat}_{\vecX_P(H_P)}(V_P)[-n]_{|\nbigu_Q}
\simeq
 \Ecat\varpi_P^{-1}(K)_{|\nbigu_Q}
\]
extending the isomorphism
of the local systems on $\nbigu_Q^{\circ}$.

Let $\varpi_P^{-1}(P)\setminus Z=\coprod \gbigi_j$
be the decomposition into the connected components.
For each $j$, let $Q_j$ be any point of $\gbigi_j$.
Let $\nbigu_{Q_j}$ be a neighbourhood of $Q_j$ as above.
By shrinking $X_P$,
we assume $\varpi_P(\nbigu_{Q_j})=H_P$ for each $j$.
Let $P'=(0,z^{(0)}_2,\ldots,z^{(0)}_n)$ be any point of $H_P$,
and let $\varphi_{P'}:\Delta\lrarr X_P$ be the holomorphic map
defined by
$\varphi_{P'}(\zeta)=(\zeta,z_2^{(0)},\ldots,z_n^{(0)})$.
Then, by the comparison of the Stokes filtrations
and by using Lemma \ref{lem;16.2.7.10},
we obtain an isomorphism
$\DR^{\Ecat}_{\vecDelta(0)}\varphi_{P'}^{\ast}(V_P)[-1]
\simeq
 \Ecat\varphi_{P'}^{-1}(K_P)$
extending the isomorphism of the local systems
on $\Delta\setminus\{0\}$.

Let $\varphi:(\Delta^{\ast},0,\Delta)\lrarr 
 (X_P\setminus H_P,0,X_P)$ be any holomorphic map.
Similarly, by comparing the Stokes filtrations,
we obtain an isomorphism 
$\DR^{\Ecat}_{\vecDelta(0)}\varphi^{\ast}(V_P)[-1]
\simeq
 \Ecat\varphi^{-1}(K_P)$
extending the isomorphism of the local systems
on $\Delta\setminus\{0\}$.

By varying $P\in H$,
and by gluing $V_P$,
we obtain a good meromorphic flat bundle $V$ on $(X,H)$.
By Proposition \ref{prop;16.6.23.20},
there exists an isomorphism
$\DR^{\Ecat}_{\vecX(H)}(K)[-n]
\simeq
 K$
extending the isomorphism
of the local systems on $X\setminus H$.
\hfill\qed

\section{Generic part of normal crossing hypersurfaces}
\label{section;18.11.15.20}

Let $X$ be an $n$-dimensional complex manifold.
Let $H$ be a simple normal crossing hypersurface of $X$.
In this section,
we shall prove the following proposition
as a preliminary of the proof of Theorem \ref{thm;18.11.12.10}.
\begin{prop}
\label{prop;16.6.26.10}
For any $K\in\Ecat^b_{\circledcirc}(\IC_{\vecX(H)})$,
there exist a closed subanalytic subset
$Z\subset H$
with $\dim_{\real}Z\leq 2n-4$,
a good meromorphic flat bundle
$V$ on $(X',H'):=(X\setminus Z,H\setminus Z)$
and an isomorphism
\[
 \DR^{\Ecat}_{\vecX'(H')}(V)[-n]
\simeq
 K_{|X'\setminus H'}
\]
in $\Ecat^b(\IC_{\vecX'(H')})$.
\end{prop}

\subsection{Ramified irregular values outside of 
subset of real codimension one}
\label{subsection;16.7.25.10}

Let $n>1$ be an integer.
Let $H$ be an $(n-1)$-dimensional complex manifold.
We set $X:=\Delta\times H$.
We identify $H$ and $\{0\}\times H$.
Let $\varpi:\Xtilde(H)\lrarr X$ be the oriented real blowing up.
Let $K\in \Ecat^b_{\circledcirc}(\IC_{\vecX(H)})$.
Let $P_0$ be a point of $H$,
and let $Q_0$ be a point of $\varpi^{-1}(P_0)$.

\begin{lem}
\label{lem;16.2.1.3}
We assume the following.
\begin{itemize}
\item
There exist a neighbourhood $\nbigu_{Q_0}$
of $Q_0$ in $\Xtilde(H)$
and real analytic functions
$h^{Q_0}_1,\ldots,h^{Q_0}_m$ on 
$\nbigu_{Q_0}^{\circ}=\nbigu_{Q_0}\setminus\varpi^{-1}(H)$
such that 
(i) they are ramified real analytic around $Q_0$
(see {\rm\S\ref{subsection;18.11.16.40}} for ramified real analyticity),
(ii) they control the growth order of
$\pi^{-1}(\cnum_{\nbigu^{\circ}_{Q_0}})\otimes K$,
i.e.,
$\pi^{-1}(\cnum_{\nbigu^{\circ}_{Q_0}})\otimes K
=\bigoplus_{i=1}^m
 \cnum^{\Ecat}_{\Xtilde(H)}
 \otimes
 \cnum_{t\geq h^{Q_0}_i}$.
\end{itemize}
Then, 
Condition {\rm\ref{condition;18.11.26.2}}
is satisfied for $(K,P_0)$.
\end{lem}
\pf
We may assume that 
$H$ is equipped with a holomorphic coordinate system
$\vecw=(w_1,\ldots,w_{n-1})$.
For any $\vecw\in H$,
let $\varphi_{\vecw}:\Delta\lrarr X$ be 
the holomorphic map
defined by
$\varphi_{\vecw}(z)=(z,\vecw)$.
Because
$\Ecat\varphi_{\vecw}^{-1}(K)\in
 \Ecat_{\mero}^b(\IC_{\vecDelta(0)})$,
we obtain the set of ramified irregular values
$\Irr(\Ecat\varphi_{\vecw}^{-1}K)$.
We set $H_{P_0}:=\varpi(\nbigu_{Q_0})$.
We may assume that
$\nbigu_{Q_0}$ 
is the product of $H_{P_0}$
and a closed sector in $\Deltatilde(0)$.

Let $(z,w_1,\ldots,w_{n-1})$ be 
the holomorphic coordinate system
of $X=\Delta\times H$.
Let $(r,\theta)$ be the polar coordinate system
defined by $z=re^{\sqrt{-1}\theta}$.
By considering a ramified covering of $(X,H)$,
we may assume that $h^{Q_0}_i$
are expressed as follows:
\[
 h^{Q_0}_i=
 \sum_{-N_1\leq j<0}
 \alpha_{i,j}(\theta,\vecw)\cdot r^{j}.
\]
Here,
$\alpha_{i,j}$ are real analytic functions
on a neighbourhood of $Q_0$ in 
$\del\Xtilde(H)=S^1\times H$.
For each $\vecw\in H_{P_0}$,
because
\[
 \varphi_{\vecw}^{\ast}(h^{Q_0}_i)
=\sum_{-N_1\leq j<0}
\alpha_{i,j}(\theta,\vecw)\cdot r^j
\]
are the real part of 
ramified meromorphic functions on $(\Delta,0)$
up to bounded functions,
there exist 
$g_{i,\vecw}(z)=\sum \beta_{i,j,\vecw}z^j
 \in z^{-1}\cnum[z^{-1}]$
such that 
$\varphi_{\vecw}^{\ast}(h^{Q_0}_i)
=\Re\Bigl(
 g_{i,\vecw}(z)
 \Bigr)$.
Because
$\Re(\beta_{i,j,\vecw}e^{\sqrt{-1}j\theta})
=\alpha_{i,j}(\theta,\vecw)$,
we can deduce that
$\beta_{i,j,\vecw}$ are real analytic functions
of $\vecw$.
We set 
$g_{i}(z,\vecw):=g_{i,\vecw}(z)$
and $\beta_{i,j}(\vecw):=\beta_{i,j,\vecw}$.
Then,
$h^{Q_0}_i=\varpi^{\ast}(g_i)$ holds
on $\nbigu_{Q_0}$,
and the coefficients $\beta_{i,j}$ are 
$\cnum$-valued real analytic functions.

Let
$(\vecw,\vecb,M)$
be any element of
$H_{P_0}\times\cnum^{n-1}\times\seisuu_{>0}$.
If $\epsilon>0$ is sufficiently small,
the holomorphic map 
$\varphi_{M,\vecw,\vecb}:
 \Delta_{\epsilon}\lrarr X$
is defined by
$\varphi_{M,\vecw,\vecb}(\zeta)
=(\zeta^M,\vecw+\zeta\vecb)$.
Let $\varphitilde_{M,\vecw,\vecb}:
 \Deltatilde_{\epsilon}(0)\lrarr \Xtilde(H)$
be the induced map.
There exists a small sector $S$ in $\Deltatilde_{\epsilon}(0)$
such that 
 $\varphitilde_{M,\vecw,\vecb}(\Stilde)
 \subset \nbigu_{Q_0}$.
The growth order of
$\pi^{-1}(\cnum_{S})\otimes
 \Ecat(\varphitilde_{M,\vecw,\vecb}\circ\varpi)^{-1}K$
is controlled by
the functions
$\varphitilde^{\ast}_{M,\vecw,\vecb}h^{Q_0}_{i}
=\varpi_0^{\ast}\varphi^{\ast}_{M,\vecw,\vecb}
 \Re(g_i)$.
We use the polar coordinate system
$(r,\theta)$ determined by
$\zeta=re^{\sqrt{-1}\theta}$ on $\Deltatilde_{\epsilon}(0)$.
As the Taylor expansion at $r=0$,
we obtain the following:
\begin{multline}
\varpi_0^{\ast}\varphi^{\ast}_{M,\vecw,\vecb}(\Re g_i)
=\Re\Bigl(
 \sum_{-N_1\leq j<0}\beta_{ij}(\vecw+\zeta\vecb)\zeta^{Mj}
 \Bigr)
 \\
=\sum_{-N_1M\leq k<0}
 \gamma_{i,k,(M,\vecw,\vecb)}(\theta)r^{k}
+\sum_{k\geq 0}
 \gamma_{i,k,(M,\vecw,\vecb)}(\theta)r^{k}
=:A_{i,(M,\vecw,\vecb),-}
+A_{i,(M,\vecw,\vecb),+}.
\end{multline}
There exists
$f\in\Irr(\Ecat\varphi_{M,\vecw,\vecb}^{-1}K)$
such that
$\Re(f)-A_{i,(M,\vecw,\vecb),-}$ is bounded on $S$.
Hence, for $k<0$,
there exists $\delta_{i,k,(M,\vecw,\vecb)}\in\cnum$
such that
$\gamma_{i,k,(M,\vecw,\vecb)}(\theta)r^k$
is of the form
$\Re\bigl(
 \delta_{i,k,(M,\vecw,\vecb)}\zeta^k
 \bigr)$.

For any fixed $i$, let $j_0=j_0(i)$ be the maximum
of $\{j\,|\,\beta_{i,-j}\neq 0\}$.
Then, we obtain
\begin{multline}
 \Re\Bigl(
 \beta_{i,-j_0}(\vecw+\zeta\vecb)\zeta^{-Nj_0}
 \Bigr)
= \\
\Re\Bigl(
 \beta_{i,-j_0}(\vecw)\zeta^{-Nj_0}
+\sum_k\overline{b}_k\del_{\wbar_k}\beta_{i,-j_0}(\vecw)\zetabar\zeta^{-Nj_0}
+\sum_k b_k\del_{w_k}\beta_{i,-j_0}(\vecw)\zeta\zeta^{-Nj_0}
+\cdots
 \Bigr).
\end{multline}
By looking at the term
$\gamma_{i,-Nj_0+1,(M,\vecw,\vecb)}(\theta,\vecw)r^{-Nj_0+1}$,
we can deduce
$\del_{\wbar_k}\beta_{i,-j_0}(\vecw)=0$
for any $\vecw\in H_{P_0}$.
Hence, we obtain that
$\beta_{i,-j_0}$ is holomorphic.
By a descending inductive argument,
we obtain that 
$\beta_{i,-j}$ are holomorphic
for any $j>0$.

Let $\varphi:(\Delta^{\ast},0,\Delta)\lrarr (X\setminus H,H,X)$
be any holomorphic map.
Let $\varphitilde:\Deltatilde(0)\lrarr \Xtilde(H)$
denote the induced map.
There exists a sector
$S=\bigl\{(\theta,r)\in\Deltatilde(0)\,\big|\,
 \theta_1<\theta<\theta_2,\,\,0\leq r<\epsilon'\bigr\}$
such that
$\varphitilde(S)\subset\nbigu_{Q_0}$.
The growth order of 
$\pi^{-1}(\cnum_S)\otimes
 \Ecat(\varpi\circ\varphitilde)^{-1}(K)$
is controlled by
$\varphitilde^{\ast}h_i^{Q_0}$.
Hence, we obtain that
$\Irr\bigl(\Ecat\varphi^{-1}K\bigr)
=\{\varphi^{\ast}g_i\}$,
which is clearly compatible with the multiplicity.
Thus, we obtain that  
Condition \ref{condition;18.11.26.2}
is satisfied for $K$ and $P_0$
with the index set $\{g_i\}$.
\hfill\qed

\subsection{Meromorphic functions and subanalytic functions}

\subsubsection{Extension of meromorphic functions}
\label{subsection;16.8.26.10}

Let $X,$ $H$ and $\varpi:\Xtilde(H)\lrarr X$
be as in \S\ref{subsection;16.7.25.10}.
Let $\nu$ be a real analytic function on $H$
such that $d\nu\neq 0$ at 
$H_0:=\nu^{-1}(0)$.
We set $H_{\pm}:=\{P\in H\,|\,\pm \nu(P)>0\}$.
We assume that 
$H=H_0\times \openopen{-1}{1}$
as a real analytic manifold,
where the projection to $\openopen{-1}{1}$
is induced by $\nu$.

Let $f$ be a real analytic function defined on
$X_+:=\Delta\times H_+$
of the form
\[
 f(z,P)=\Re\Bigl(
 \sum_{0<j\leq N_1}
 \alpha_j(P)z^{-j}
 \Bigr),
\]
where $\alpha_j$ are holomorphic functions
defined on $H_+$.

We use the polar decomposition $z=re^{\sqrt{-1}\theta}$.
For any positive integer $m$ and positive numbers
$\theta_1$, $c_1$ and $\epsilon_1$,
and for any subset $A\subset H_0$,
we set 
\[
B(A,m,c_1,\theta_1,\epsilon_1):=
\bigl\{
 (r,\theta,\nu)\,\big|\,
 \nu>0,\,\,
 0<r<c_1\nu^m,\,\,
 |\theta|<\theta_1
\bigr\}
\times A
\subset
 X\setminus H.
\]
Its closure in $\Xtilde(H)$
is denoted by 
$\overline{B}(A,m,c_1,\theta_1,\epsilon_1)$.
Let $m$, $c_1$, $\theta_1$
and $\epsilon_1$ be as above.
Let $g$ be a function 
on $B(H_0,m,c_1,\theta_1,\epsilon_1)$
defined by a convergent power series
\[
 g=
 \sum_{i\geq -N_1}
\sum_{j\geq -N_2}
  (r\nu^{-m})^{i/\rho}\nu^{j/\rho}\beta_{ij}(\theta,v),
\]
where $v$ varies on $H_0$.

\begin{lem}
\label{lem;16.6.25.1}
We assume the following.
\begin{itemize}
\item
$f-g$ is bounded
around any 
$(0,\theta,\nu,v)
\in \Bbar(H_0,m,c_1,\theta_1,\epsilon_1)$
such that $\nu>0$.
\end{itemize}
Then, the function
$\sum_{0<j\leq N_1} \alpha_j z^{-j}$
extends to a meromorphic function
on a neighbourhood of
$\Delta\times\{\nu\geq 0\}$.
\end{lem}
\pf
Fix a point $P_0\in H_0$.
By shrinking $H_0$ around $P_0$,
we may assume that there exist
holomorphic functions
$w_1,\ldots,w_{n-2}$ on $H$
such that
$F_{|H_0}:H_0\lrarr \cnum^{n-2}$
is submersive,
where $F:H\lrarr \cnum^{n-2}$ is induced by
$(w_1,\ldots,w_{n-2})$.
By shrinking $H_0$ around $P_0$,
we may assume that there exists
a $\cnum$-valued real analytic function
$\eta$ on $H$
such that
(i) $\eta_{|F^{-1}(P')}$
are holomorphic
for any $P'\in \cnum^{n-2}$,
(ii) $H_0=\{\Re\eta=0\}$.
We may assume that the coordinate system
$(\eta,w_1,\ldots,w_{n-2})$
induces an open embedding of $H$
to $\cnum^{n-1}$.
We may assume that
$\nu=\Re(\eta)$.
We set $\mu:=\Image(\eta)$.
Tuples $(w_1,\ldots,w_{n-2})$
are denoted by $\vecw$.

By comparing the power expansions with respect to $r$,
we obtain the following equality:
\begin{equation}
\label{eq;16.6.28.30}
 f=\sum_{-N_1\leq i<0}\sum_{j\geq -N_2}
 (r\nu^{-m})^{i/\rho}\nu^{j/\rho}
 \beta_{ij}(\theta,\mu,\vecw).
\end{equation}
Moreover,
for any fixed $(\theta,\vecw)$,
the function
$\sum_j\nu^{-mi/\rho}\nu^{j/\rho}\beta_{ij}(\theta,\mu,\vecw)$
is harmonic with respect to
$\del_{\nu}^2+\del_{\mu}^2$.
So, we obtain that 
the only non-negative integer powers of $\nu$ appear
in (\ref{eq;16.6.28.30}),
and that $\alpha_{j}$ extend
to holomorphic functions
on a neighbourhood of $\{\nu\geq 0\}$.
\hfill\qed

\subsubsection{Boundedness}

We obtain the boundedness of $f-g$ under
an additional assumption.

\begin{lem}
\label{lem;16.8.25.1}
In addition to the assumption in Lemma {\rm\ref{lem;16.6.25.1}},
we assume the following to $g$.
\begin{description}
\item[(P)]
Let $\varphi:(\Delta^{\ast},0,\Delta)\lrarr (X\setminus H,H_0,X)$ 
be any holomorphic map.
Then, 
for each connected component $\nbigc$
of $\varphi^{-1}B(H_0,m,c_1,\theta_1,\epsilon_1)$,
there exist $\kappa\in\seisuu_{>0}$
and $h\in\cnum[\zeta^{-1/\kappa}]$
such that
$\varphi^{\ast}g-\Re(h)$ is bounded on $\nbigc$.
\end{description}
Let $H_1\subset H_0$ be any relatively compact open subset.
Then,
$f-g$ is bounded on $B(H_1,m,c_1,\theta_1,\epsilon_1)$.
\end{lem}
\pf
We use the notation in the proof of 
Lemma \ref{lem;16.6.25.1}.
It is enough to prove 
$\beta_{i,j}(\theta,\mu,\vecw)=0$ for $i\geq 0$ and $j<0$.
We assume that
$\sum_{i\geq 0}\sum_{0>j\geq -N_2}
 (r\nu^{-m})^{i/\rho}
 \nu^{j/\rho}
 \beta_{ij}(\theta,\mu,\vecw)\neq 0$,
and we shall derive a contradiction.
We may assume that 
there exist $i_0\geq 0$
and $(\mu_0,\vecw_0)$
such that 
$\beta_{i_0,-N_2}(\theta,\mu_0,\vecw_0)$
is not constantly $0$.
We may assume $\mu_0=0$.
We obtain the expansion
$\beta_{i,j}(\theta,\mu,\vecw_0)=\sum_{k\geq 0}
 \beta_{i,j,k}(\theta,\vecw_0)\mu^k$.

We set
\[
 g_0:=
 \sum_{i\geq 0}\sum_{j\geq -N_2}
 (r\nu^{-m})^{i/\rho}
 \nu^{j/\rho}
 \beta_{ij}(\theta,\mu,\vecw).
\]
The condition (P) is satisfied for $g_0$
by Lemma \ref{lem;16.6.25.1}.

We define the holomorphic map
$\varphi:\{(\zeta,a)\,|\,|\zeta|<\epsilon,\,\,|1-a|<\epsilon\}
 \lrarr X$
by 
$\varphi(\zeta,a)=(\zeta^m,a\zeta,\vecw_0)$,
where we use the coordinate system
$(z,\eta,\vecw)$ on $X$.
Note that $\eta$ is holomorphic on $F^{-1}(\vecw_0)$.
On the domain
$\varphi^{-1}(B(H_1,m,c_1,\theta_1,\epsilon_1))$,
we obtain
\begin{multline}
 \varphi^{\ast}(g_0)(a,\zeta)
=
\sum_{i\geq 0}
 \sum_{j\geq -N_2}
 \sum_{k\geq 0}
 |\zeta|^{mi/\rho}
 \Re(a\zeta)^{(-mi+j)/\rho}
 \Image(a\zeta)^k
 \beta_{i,j,k}(m\arg(\zeta),\vecw_0)
 \\
=
 \sum_{i\geq 0}
 \sum_{j\geq -N_2}
 \sum_{k\geq 0}
 |\zeta|^{j/\rho+k}
 |a|^{(-mi+j)/\rho+k}
 \cos(\arg(a\zeta))^{(-mi+j)/\rho}
 \sin(\arg(a\zeta))^{k}
 \beta_{i,j,k}(m\arg(\zeta),\vecw_0).
\end{multline}
By the argument in the proof of Lemma \ref{lem;16.2.1.3},
we obtain that
the coefficient of $|\zeta|^{-N_2/\rho}$
is harmonic with respect to $a$,
i.e.,
the function
\begin{multline}
 \sum_{i\geq 0}
  |a|^{(-mi-N_2)/\rho}
 \cos(\arg(a\zeta))^{(-mi-N_2)/\rho}
 \beta_{i,-N_2,0}(m\arg(\zeta),\vecw_0)
= \\
 \sum_{i\geq 0}
 \Re(a\zeta)^{(-mi-N_2)/\rho}
 |\zeta|^{(mi+N_2)/\rho}
 \beta_{i,-N_2,0}(m\arg(\zeta),\vecw_0)
\end{multline}
is harmonic with respect to $a$.
It implies 
$\beta_{i,-N_2,0}(m\arg(\zeta),\vecw_0)=0$
for any $i\geq 0$,
which contradicts the assumption.
Thus, we obtain Lemma \ref{lem;16.8.25.1}.
\hfill\qed

\subsection{Extension along subsets of real codimension three}

\subsubsection{Preliminary}
\label{subsection;17.12.28.1}

Let $X$, $H$, $\nu$ and $H_0$, $H_{\pm}$
be as in \S\ref{subsection;16.8.26.10}.
For any $P\in H_0$,
let $X_P$ be a small neighbourhood of $P$ in $X$,
and we set $H_P:=H\cap X_P$,
$H_{0,P}:=H_0\cap H_P$
and $H_{\pm,P}:=H_{\pm}\cap H_P$.
For a positive integer $m$
and positive numbers $c_1$, $\theta_1$ and $\epsilon_1$,
we set
\[
B_{\pm}(H_P,m,c_1,\theta_1,\epsilon_1)
:=
 \bigl\{(r,\theta,\nu)\,\big|\,
 \pm\nu>0,\,\,
 0<r<c_1(\pm\nu)^{m},\,\,
 |\theta|<\theta_1
 \bigr\}
 \times H_{0,P}
\subset X\setminus H.
\]

Let $K\in\Ecat^b_{\circledcirc}(\IC_{\vecX(H)})$.
Suppose that 
Condition \ref{condition;18.11.26.3}
is satisfied for $K$ and 
any $P\in H\setminus H_0$.
For simplicity, we assume that 
the ramified sets of irregular values $\nbigi_P$
are contained in $z^{-1}\nbigo_{H,P}[z^{-1}]$,
where $\nbigo_{H,P}$ denotes the local ring at $P$.
We also impose the following condition along $H_0$.

\begin{assumption}
\label{assumption;17.12.27.10}
For each $P\in H_0$,
there exist good sets of irregular values 
$\nbigi_{P,\kappa}
 \subset z^{-1}\nbigo_{H,P}[z^{-1}]$  $(\kappa=\pm)$
with a multiplicity function
$\gminim_{P,\kappa}:\nbigi_{P,\kappa}\lrarr\seisuu_{>0}$
such that the following holds
for a neighbourhood $H_P$
and some positive constants $m(P),c_1(P),\theta_1(P),\epsilon_1(P)$:
\begin{itemize}
\item
 Let $\varphi:\Delta\lrarr X$ be any holomorphic map
 such that 
 $\varphi^{-1}(B_{\kappa}(H_P,m(P),c_1(P),\theta_1(P),\epsilon_1(P)))$
 contains a sector of $\Delta^{\ast}$.
 Then,
 $\Irr(\Ecat\varphi^{-1}K)
=\varphi^{\ast}\nbigi_{P,\kappa}$ holds,
 and the multiplicity functions are equal.
\end{itemize}
\end{assumption}

\begin{lem}
We obtain that
 $\nbigi_{P,+}=\nbigi_{P,-}$
 and $\gminim_{P,+}=\gminim_{P,-}$.
\end{lem}
\pf
It is enough to consider the case
$\dim_{\cnum}H=1$.
We may assume that there exists 
a holomorphic coordinate $w$
on $H$ such that $w(P)=0$.
Suppose that
$f\in\nbigi_{P,+}\setminus\nbigi_{P,-}$,
and we shall derive a contradiction.
Let $\nbigi_{P,-}=\{g_1,,\ldots,g_{\ell}\}$.
We consider the expansion
\[
 (f-g_j)(z,w)
=\sum_{i=1}^{N(j)}a_{j,i}(w)z^{-i},
\]
where $a_{j,N(j)}$ are not constantly $0$.
If $a_{j,i}$ is not constantly $0$,
there exists a holomorphic function $b_{j,i}(w)$
with $b_{j,i}(0)\neq 0$
such that 
$a_{j,i}(w)=w^{\ell(j,i)}\cdot b_{j,i}(w)$.
For each $j$, there exists $m(j)\in\seisuu_{>0}$ such that 
$-m(j)N(j)+\ell(j,N(j))<-m(j)i+\ell(j,i)$
for any $i$ such that $a_{j,i}\neq 0$.
We set $m_0:=\max\{m(j)\}$.
Note that 
$-mN(j)+\ell(j,N(j))<-mi+\ell(j,i)$
for any $m\geq m_0$
and for any $i$ and $j$ such that $a_{i,j}\neq 0$.

Let $m$ be any integer larger than $m_0$.
We define $\varphi_1:\Delta_{\epsilon}\lrarr X$
by $\varphi_1(\zeta)=(\zeta^m,\zeta)$.
Then, by our choice of $m$,
we obtain
$\varphi_1^{\ast}(f-g_j)\neq 0$
in 
$\nbigo_{\Delta}(\ast 0)/\nbigo_{\Delta}$.
If $m$ is sufficiently large,
then both 
\[
\varphi^{-1}(B_{\kappa}(H_P),m(P),c_1(P),\theta_1(P),\epsilon_1(P))
\quad
(\kappa=+,-)
\]
contain sectors of $\Delta^{\ast}$.
Hence, we have arrived at a contradiction,
i.e.,
$\nbigi_{P,+}\subset\nbigi_{P,-}$.

Similarly, we obtain that
$\nbigi_{P,-}\subset\nbigi_{P,+}$.
We can also compare the multiplicity functions
$\gminim_{P,+}$ and $\gminim_{P,-}$ in a similar way.
\hfill\qed

\vspace{.1in}

We denote $\nbigi_{P,\pm}$ and $\gminim_{P,\pm}$
by $\nbigi_P$ and $\gminim_P$,
respectively.

\subsubsection{Extension}

We continue to use the notation in \S\ref{subsection;17.12.28.1}.
We assume that there exists
a global holomorphic coordinate system
$(x_1,\ldots,x_{n-1})$ of $H$,
by which we regard $H$
as an open subset of $\cnum^{n-1}$.
We also regard $X$ as an open subset of $\cnum^n$.

Let $P\in H_0$
and $\veca_0\in \cnum^{n-1}$.
Let $m$ be a positive integer strictly larger than $m(P)$.
Let $A(\veca_0)$ be a neighbourhood of
a point $\veca_0$ in $\cnum^{n-1}$.
We define $\Phi_{P,\veca_0}:
 \Delta_{\epsilon}\times A(\veca_0)
 \times H_{0,P}\lrarr \cnum^n$
by
$\Phi_{P,\veca_0}(\xi,\veca,\vecx)=(\xi^m,\vecx+\xi\veca)$.
We set 
$Y_{P,\veca_0}:=\Phi_{P,\veca_0}^{-1}(X)$
and $J_{P,\veca_0}:=
 \Phi_{P,\veca_0}^{-1}(H)=\{0\}\times(A(\veca_0)\times H_{0,P})$.
We set 
 $\vecY_{P,\veca_0}(J_{P,\veca_0}):=
 (Y_{P,\veca_0}\setminus J_{P,\veca_0},Y_{P,\veca_0})$.
We obtain
$\Ecat\Phi^{-1}_{P,\veca_0}(K)$
in $\Ecat^b_{\realc}(\IC_{\vecY_{P,\veca_0}(J_{P,\veca_0})})$.
Let $\varpi_{Y_{P,\veca_0}}:\Ytilde_{P,\veca_0}(J_{P,\veca_0})\lrarr Y_{P,\veca_0}$
denote the oriented real blowing up
along $J_{P,\veca_0}$.

\begin{assumption}
\label{assumption;16.8.27.10}
For any $P\in H_0$,
there exist $\veca(P)\in\cnum^{n-1}$,
a neighbourhood $A(\veca(P))$,
and a subanalytic subset
$\nbigz\subset\varpi_{Y_{P,\veca(P)}}^{-1}(J_{P,\veca(P)})$
with $\dim_{\real}\nbigz\leq \dim_{\real}J_{P,\veca(P)}$
such that the following holds:
\begin{itemize}
\item
 $\nbigz$ is horizontal with respect to
 $\varpi_{Y_{P,\veca(P)}}$.
\item
Any $Q\in \varpi_{Y_{P,\veca(P)}}^{-1}(J_{P,\veca(P)})
 \setminus\nbigz$,
has a neighbourhood $\nbigu_Q$
in $\Ytilde_{P,\veca(P)}(J_{P,\veca(P)})$
such that 
\[
 \pi^{-1}(\cnum_{\nbigu_Q})\otimes
 \Ecat\Phi_{P,\veca(P)}^{-1}(K)
=\bigoplus_{i=1}^m
 \cnum^{\Ecat}
 \overset{+}{\otimes}
 \cnum_{t\geq h_{Q,i}}
\]
for a tuple of ramified real analytic functions $h_{Q,i}$
on $\nbigu_Q$.
\end{itemize}
\end{assumption}
We remark the following.
\begin{lem}
\label{lem;16.10.12.1}
$\{h_{Q,i}\}$ in Assumption {\rm \ref{assumption;16.8.27.10}}
is equal to
$\{\Re\Phi_{P,\veca(P)}^{\ast}(f)\,|\,f\in\nbigi_P\}$
in $\Subbar(\nbigu_Q^{\circ},\nbigu_Q)$
compatible with the multiplicity. 
\end{lem}
\pf
For each $\vecx\in H_{0,P}$,
there exists a real subspace
$T(\vecx)\subset\cnum^{n-1}$
with $\dim_{\real}T(\vecx)=2n-3$
such that if $\veca\not\in T(\vecx)$
then $|\nu(\veca\xi+\vecx)|\sim b(\veca,\vecx)|\xi|$
for $b(\veca,\vecx)>0$.
For such $(\veca,\vecx)\in J_{P,\veca(P)}$,
let us consider the holomorphic map
$\varphi:\Delta\lrarr\cnum^n$
defined by
$\varphi(\zeta)=(\zeta^m,\veca\zeta+\vecx)$.
Then,
$\varphi^{-1}(B_{\pm}(H_P,m(P),c_1(P),\theta_1(P),\epsilon_1(P)))$
contains sectors.
Hence, we obtain that 
$\{h_{Q,i}\}$ is equal to
$\{\Re\Phi_{P,\veca(P)}^{\ast}(f)\,|\,f\in\nbigi_P\}$
in $\Subbar(\nbigu_Q^{\circ},\nbigu_Q)$
compatible with the multiplicity
for $Q\in \varpi_{Y_{P,\veca(P)}}^{-1}(\veca,\vecx)$.
We obtain the claim for general $Q$
by using the continuity
because the functions $h_{Q,i}$ are ramified real analytic functions.
\hfill\qed

\begin{prop}
\label{prop;18.11.26.20}
Condition {\rm\ref{condition;18.11.26.3}}
is satisfied for $K$.
\end{prop}
\pf
It is enough to prove the following.
\begin{itemize}
\item
For any holomorphic map
$\varphi:(\Delta^{\ast},0,\Delta)\lrarr 
 (X\setminus H,H_0,X)$,
$\Irr(\Ecat\varphi^{-1}(K))=\varphi^{\ast}\nbigi_{\varphi(0)}$
holds,
and the multiplicity functions are the same.
\end{itemize}

Let us begin with some general lemmas.

\begin{lem}
\label{lem;16.8.30.10}
Let $g:\Delta\lrarr \real_{\geq 0}$ 
be a non-constant real analytic function.
Let $m_1<m_2$ be positive integers.
Let $\theta_1<\theta_2$.
Let $\ell$ be any positive number.
We set
\[
 T_1:=
 \bigl\{\zeta\,\big|\,
 \theta_1<\arg(\zeta)<\theta_2,\,\,
 |\zeta|^{\ell}<g(\zeta)^{m_1}\bigr\},
\quad
 T_2:=
 \bigl\{\zeta\,\big|\,
 g(\zeta)^{m_2}<|\zeta|^{\ell}
 \bigr\}.
\]
Then, 
either one of the following holds:
(i) $T_1$ contains a sector of $(\Delta,0)$,
or (ii) $T_2$ contains a neighbourhood of $0$.
\end{lem}
\pf
Let $(x,y)$ be the real coordinate system
obtained as $\zeta=x+\sqrt{-1}y$.
There exists a positive integer $k$
such that 
$g=\sum_{i+j\geq k} a_{i,j}x^iy^j$,
and that at least one of $a_{i,j}$ $(i+j=k)$ is not $0$.

If $k>\ell/m_2$,
then $T_2$ clearly contains a neighbourhood of $0$.
Let us consider the case $k\leq \ell/m_2<\ell/m_1$.
By considering a coordinate change
from $\zeta$ to $\beta\zeta$,
we may assume $a_{k,0}\neq 0$,
and $\theta_1<0<\theta_2$.
If $\epsilon_i>0$  $(i=1,2)$ are sufficiently small,
then $g$ and $x^k$ are mutually bounded on
$S:=\bigl\{(x,y)\,\big|\,
 0<x<\epsilon_1,\,\, |y|<\epsilon_2x\bigr\}$.
Moreover,
$S$ is contained  in
$\{\zeta\,|\,\theta_1<\arg(\zeta)<\theta_2\}$,
and 
$|\zeta|$ and $x$ are mutually bounded
on $S$.
Then, the claim is clear.
\hfill\qed

\begin{lem}
\label{lem;16.8.30.11}
Let $g_2:\Delta\lrarr \real_{\geq 0}$
be a real analytic function.
Suppose that
$g_2(\zeta)\leq |\zeta|^{\ell}$
for a positive integer $\ell$.
Let $\psi:[0,\epsilon]\lrarr \Delta$
be a real analytic curve
such that $\psi(0)=0$.
Suppose that
$g_2\circ\psi(t)=0$ for any $t$.
Then, for any $\delta>0$,
there exists a sector $S_{\delta}$ of $\Delta^{\ast}$
and $0<\epsilon'<\epsilon$
such that
(i) $S_{\delta}$ contains $\psi(\openopen{0}{\epsilon'})$,
(ii) $g_2\leq \delta|\zeta|^{\ell}$ on $S_{\delta}$.
\end{lem}
\pf
There exist $k>0$
such that 
$g_2=\sum_{i+j\geq k}a_{i,j}x^iy^j$
and that at least one of $a_{i,j}$ $(i+j=k)$
is not $0$.
Clearly, $k\geq \ell$ holds.
It is enough to study the case $k=\ell$.

Let us consider the case $\psi(t)=t$.
We obtain $a_{k,0}=0$.
Hence, if $\rho_1$ is sufficiently small,
then on $\{|y|<\rho_1x\}$,
the inequality $g_2\leq \delta|\zeta|^{\ell}$ holds.

Let us consider the general case.
Let $\psi_{\cnum}:\Delta_{\epsilon}\lrarr \Delta$
denote the holomorphic map
whose restriction to
$\{0\leq t\}$ is equal to $\psi$.
There exists a sector 
$S$ in $\Delta^{\ast}_{\epsilon}$
on which 
$g_2\circ\psi\leq \delta|\psi^{\ast}(\zeta^{\ell})|$
holds.
Because
$\psi_{\cnum}(S)$ contains a small sector,
we obtain the claim of the lemma.
\hfill\qed

\vspace{.1in}
Let $P$ be any point of $H_0$.
If $H_P$ is sufficiently small,
there exists a real analytic map
$q:H_P\lrarr H_{0,P}$
whose restriction to $H_{0,P}$ is the identity.

Let $\varphi:(\Delta^{\ast},0,\Delta)\lrarr 
 (X_P\setminus H_P,H_P,X_P)$
be any holomorphic map.
It is described as
$(\varphi_{z}(\zeta),\varphi_{\vecx}(\zeta))$.
By Lemma \ref{lem;16.8.30.10},
either one of the following holds.
\begin{itemize}
\item
$\varphi^{-1}(B_{\pm}(H_P,m(P),c_1(P),\theta_1(P),\epsilon_1(P)))$
contain sectors.
\item
 $T_3=\{|\varphi_z(\zeta)|^2>\nu(\varphi_{\vecx}(\zeta))^{2m}\}$
 contains a neighbourhood of $0$.
\end{itemize}

If the first case occurs,
$\Irr(\Ecat\varphi^{-1}K)
=\varphi^{\ast}\nbigi_{P}$
holds,
and the multiplicity functions are the same,
by Assumption \ref{assumption;17.12.27.10}.
Suppose that the second case occurs.
We may assume that
$\varphi_{z}(\zeta)=\zeta^{mk}$
for a positive integer $k$.
We obtain 
$|\nu(\varphi_{\vecx}(\zeta))|<|\zeta|^{k}$
on $T_3$.
Let $\veca(P)$ be as in
Assumption \ref{assumption;16.8.27.10}.
We define $\psi_{\vecx}:T_3\lrarr H_0$
by 
\[
 \psi_{\vecx}(\zeta):=
 q(\varphi_{\vecx}(\zeta)-\veca(P)\zeta^{k}).
\]
There exists the one dimensional real analytic subset
$\nbigc$ of $T_3$
such that $0\in \nbigc$
and that 
$\varphi_{\vecx}(\zeta)-\veca(P)\zeta^{k}
 \in H_0$ for $\zeta\in \nbigc$.
Because
$|\nu(\varphi_{\vecx}(\zeta)-\veca(P)\zeta^{k})|
<C_1|\zeta|^{k}$ for some $C_1>0$,
$\bigl|
 \varphi_{\vecx}(\zeta)-\veca(P)\zeta^k
-\psi_{\vecx}(\zeta)
 \bigr|\leq C_2|\zeta|^k$
holds
for some $C_2>0$.
Then, by Lemma \ref{lem;16.8.30.11},
for any small neighbourhood 
$A(\veca(P))$ of $\veca(P)$ in $\cnum^{n-1}$,
there exists a small sector $S$ of $(\Delta,0)$
with a real analytic function 
$\psi_{\veca}:S\lrarr A(\veca(P))$
such that 
$\varphi_{\vecx}(\zeta)
-\psi_{\vecx}(\zeta)
=\psi_{\veca}(\zeta)\zeta^{k}$.
We define a real analytic map
$\psi:S\lrarr Y_{P,\veca(P)}$ by
$\psi(\zeta)=(\zeta^{k},\psi_{\veca}(\zeta),\psi_{\vecx}(\zeta))$,
then 
$\Phi_{P,\veca(P)}\circ\psi=\varphi$ holds.
Let $\Sbar$ denote the closure of $S$
in $\Deltatilde(0)$.
The map $\psi$ extends to
the real analytic map $\Sbar\lrarr \Ytilde_{P,\veca(P)}(J_{P,\veca(P)})$.
Then,
by Lemma \ref{lem;16.10.12.1},
$\Irr(\Ecat\varphi^{-1}K)=\varphi^{\ast}\nbigi_P$
holds,
and the multiplicity functions are the same.
Thus, we obtain the claim of the proposition.
\hfill\qed

\vspace{.1in}
We obtain the following from Proposition \ref{prop;16.7.25.40}
and Proposition \ref{prop;18.11.26.20}.

\begin{cor}
\label{cor;16.8.27.12}
Under Assumption {\rm\ref{assumption;17.12.27.10}}
and Assumption {\rm\ref{assumption;16.8.27.10}},
there exist a good meromorphic flat bundle $V$
on $(X,H)$
and an isomorphism
$\DR^{\Ecat}_{\vecX(H)}(V)[-n]
\simeq
 K$.
\hfill\qed
\end{cor}

\subsubsection{Remark on Assumption \ref{assumption;16.8.27.10}}

Let $X$ be any $n$-dimensional complex manifold.
Let $H$ be a normal crossing hypersurface of $X$.
Let $H^{[2]}$ be the singular locus of $H$.
Let $C\subset H$
be a closed subanalytic subset with
$\dim_{\real}C=2n-3$
such that $H^{[2]}\subset C$.
Let $K\in \Ecat^b_{\circledcirc}(\IC_{\vecX(H)})$.

We prepare a notation.
Let $P$ be any smooth point of $C\setminus H^{[2]}$.
Let $m$ be any positive integer.
For any small holomorphic coordinate neighbourhood
$(X_P;z,x_1,\ldots,x_{n-1})$ of $X$ around $P$
such that $H_P:=H\cap X_P=\{z=0\}$,
we set $C_P:=X_P\cap C$.
We may regard $X_P$ as an open subset of $\cnum^n$
by the coordinate system.
We define
  $\Phi_P:\cnum\times\cnum^{n-1}\times C_{P}
 \lrarr \cnum^n$
by
 $\Phi_P(\xi,\veca,\vecx)
=(\xi,\vecx+\xi^m\veca)$.
 We set 
 $Y_P:=\Phi_P^{-1}(X_P)$
 and $J_P:=\Phi_P^{-1}(H_P)=\{0\}\times\cnum^{n-1}\times C_P$.
 Set $\vecY_P(J_P)=(Y_P\setminus J_P,Y_P)$.
We obtain
 $\Ecat\Phi_P^{-1}(K)
 \in \Ecat^b_{\realc}(\IC_{\vecY_P(J_P)})$.
 Set $\varpi_{Y_P}:\Ytilde_P(J_P)\lrarr Y_P$ be 
the oriented real blowing up along $J_P$.

\begin{lem}
\label{lem;16.8.27.11}
There exists a closed subanalytic subset $\gbigw\subset C$
with $\dim_{\real}\gbigw=2n-4$
such that the following holds.
\begin{itemize}
\item
$\gbigw$ contains $H^{[2]}$
and the singular locus of $C$.
\item
Let $P$ be any point of $C\setminus\gbigw$.
We use the above notation.
Then, there exist 
a closed subanalytic subset 
$\nbigz_P\subset \varpi_P^{-1}(J_P)$
such that 
(i) $\nbigz_P$ is horizontal over $J_P$,
(ii) any $Q\in \varpi_{Y_P}^{-1}(J_P)\setminus \nbigz_P$
 has a neighbourhood $\nbigu_Q$
 in $\Ytilde_P(J_P)$ 
 on which
 $\pi^{-1}(\cnum_{\nbigu_Q^{\circ}})
 \otimes \Ecat\Phi_P^{-1}(K)$ is controlled by 
 ramified real analytic functions.
In other words,
the assumption {\rm\ref{assumption;16.8.27.10}}
is satisfied around $P$.
\end{itemize}
\end{lem}
\pf
It is enough to consider the case 
where $X=\Delta_z\times\Delta_{\vecx}^{n-1}$
and $H=\{z=0\}\cup\bigcup_{j=1}^{\ell}\{x_j=0\}$.
We set $H^{(0)}:=\{z=0\}$.
It is enough to consider the case
where $C$ is contained in $H^{(0)}$.

There exists
a $(2n-3)$-dimensional real analytic manifold
$M$ with a proper map
$\rho:M\lrarr X$ such that $\rho(M)=C$.
Note that we may assume that there exists
a closed real analytic subset 
$\gbigw_M\subset M$ with $\dim \gbigw_M\leq 2n-4$
such that 
(i) $\rho^{-1}\rho(\gbigw_M)=\gbigw_M$,
(ii) $M\setminus \gbigw_M\simeq C\setminus \rho(\gbigw_M)$,
(iii) $\gbigw_M$ contains the pull back of
$C\cap H^{[2]}$ and the singular locus of $C$ by $\rho$.

We define 
$\Phi:\cnum\times\cnum^{n-1}\times M
\lrarr \cnum^n$
by 
$\Phi(\xi,\veca,y)=(\xi,\rho(y)+\xi^m\veca)$.
We set $Y:=\Phi^{-1}(X)$
and $J:=\Phi^{-1}(H^{(0)})=\{0\}\times\cnum^{n-1}\times M$.
We also set $\nbigh:=\Phi^{-1}(H)$.

Let $\varpi_{Y}:\Ytilde(J)\lrarr Y$
be the oriented real blowing up.
Let $j:(Y\setminus\nbigh,\Ytilde(J))\lrarr 
\Ytilde(J)=(\Ytilde(J),\Ytilde(J))$
be the inclusion of the bordered spaces.
Let $\Ktilde:=\Ecat j_{!!}\Ecat\Phi^{-1}K$.
There exists a filtration 
$\Ytilde(J)=\Ytilde(J)^{(0)}
\supset
 \Ytilde(J)^{(1)}\supset\cdots$
by closed subanalytic subsets
such that 
(i) $\Ytilde(J)^{(j)}\setminus \Ytilde(J)^{(j+1)}$
 are submanifolds of codimension $j$,
(ii) for each connected component $\nbigc$
 of $\Ytilde(J)^{(j)}\setminus \Ytilde(J)^{(j+1)}$,
 there exist subanalytic functions
 $h^{\nbigc}_1,\ldots,h^{\nbigc}_{k(\nbigc)}$
 on $(\nbigc,\Ytilde(J))$
 such that 
$\pi^{-1}(\cnum_{\nbigc})
\otimes \Ktilde
=\bigoplus
 \cnum^{\Ecat}\overset{+}{\otimes}
 \cnum_{t\geq h^{\nbigc}_{\ell}}$.
Applying Proposition \ref{prop;16.1.27.10}
to 
the components $\nbigc$ with maximal dimension
and the functions $h^{\nbigc}_{\ell}$,
we obtain a closed subanalytic subset
$\gbigz\subset \varpi_Y^{-1}(J)$
with 
$\dim_{\real}\gbigz\leq \dim_{\real} J$
such that the following holds.
\begin{itemize}
\item
For any $Q\in \varpi_Y^{-1}(J)\setminus \gbigz$,
there exist a neighbourhood $\nbigu_Q$ of $Q$
in $\Ytilde(J)$ and 
ramified real analytic functions $h^{Q}_i$
such that
$\pi^{-1}(\cnum_{\nbigu_Q})\otimes\Ktilde
=\bigoplus \cnum^{\Ecat}\overset{+}{\otimes}\cnum_{t\geq h^Q_i}$.
\end{itemize}
There exists a closed subanalytic subset
$\gbigr\subset J$
with $\dim_{\real}\gbigr<\dim_{\real}J$
such that
$\gbigz\setminus \varpi_Y^{-1}(\gbigr)$
is horizontal over $J$.
There exists a closed subanalytic subset
$\gbigw\subset C$
with $\dim_{\real} \gbigw\leq 2n-4$
such that 
(i) $\gbigw$ contains  $\rho(\gbigw_M)$, $H^{[2]}$
 and the singular locus of $C$,
(ii) for any $P\in C\setminus \gbigw$,
 $\dim_{\real}(\gbigr\cap \Phi^{-1}(P))
 <\dim_{\real}\Phi^{-1}(P)$ holds.
Then, the claim of the lemma follows.
\hfill\qed

\subsection{Proof of Proposition \ref{prop;16.6.26.10}}

Let $X$ be an $n$-dimensional complex manifold,
and $H$ be a simple normal crossing hypersurface of $X$.

\subsubsection{Preliminary}

Let $\varpi:\Xtilde(H)\lrarr X$ be the oriented real blowing up.
We define the bordered space
$\vecXtilde(H):=(X\setminus H,\Xtilde(H))$,
and $\varpi$ induces an isomorphism
of bordered spaces
$\vecXtilde(H)\lrarr \vecX(H)$.
We obtain
$\Ecat\varpi^{-1}K\in \Ecat^b_{\realc}(\vecXtilde(H))$.
Let 
$\Xtilde(H)
 =\Xtilde(H)^{(0)}
 \supset
 \Xtilde(H)^{(1)}\supset\cdots$
denote a filtration for $\Ecat\varpi^{-1}K$.
For any connected component $\nbigc$ 
of $\Xtilde(H)\setminus \Xtilde(H)^{(1)}$,
there exist subanalytic functions 
$h^{\nbigc}_1,\ldots,h^{\nbigc}_{\ell}$ on 
$(\nbigc,\Xtilde(H))$
such that
$\pi^{-1}(\cnum_{\nbigc})\otimes K
 \simeq
 \bigoplus_{j}
 \cnum_{\Xtilde(H)}^{\Ecat}\overset{+}{\otimes}
 \cnum_{t\geq h^{\nbigc}_j}$.
By Proposition \ref{prop;16.1.27.10}
and Lemma \ref{lem;16.8.27.11},
we obtain the following.

\begin{lem}
\label{lem;16.7.25.20}
There exist a $(2n-3)$-dimensional closed subanalytic subset
$\gbigz_{0}\subset H$
and a $(2n-2)$-dimensional closed subanalytic subset
$\gbigr_0\subset \del \Xtilde(H)$
with the following property.
\begin{itemize}
\item
$\gbigz_0$ contains the singular locus $H^{[2]}$ of $H$.
\item
 $\gbigr_0$ contains
 $\varpi^{-1}(H^{[2]})$.
 Let $\gbigw$ be the closure of $\Xtilde(H)^{(1)}\setminus\del\Xtilde(H)$
 in $\Xtilde(H)$.
 Then, the intersection 
 $\gbigw\cap\del\Xtilde(H)$
 is also contained in $\gbigr_0$.
\item
The induced map
$\gbigr_0\setminus\varpi^{-1}(\gbigz_0)\lrarr H$
is relatively $0$-dimensional.
\item
Let $Q$ be any point of 
$\del\Xtilde(H)\setminus \gbigr_0$.
Then,
there exist a neighbourhood $\nbigu_Q$
of $Q$ in $\Xtilde(H)$
and real analytic functions
$h^Q_{1},\ldots,h^Q_{\ell}$
on $\nbigu_Q^{\circ}=\nbigu_Q\setminus\varpi^{-1}(H)$
which are ramified real analytic around $Q$,
such that 
$\pi^{-1}(\cnum_{\nbigu_Q})\otimes K
=\bigoplus_j
 \cnum_{\Xtilde(H)}^{\Ecat}\overset{+}{\otimes}
 \cnum_{t\geq h^Q_j}$. 
\end{itemize}
There exist a $(2n-3)$-dimensional 
closed subanalytic subset
$\gbigr_1\subset
 \gbigz_0\times_{H}\del\Xtilde(H)$
and 
a $(2n-4)$-dimensional closed subanalytic subset
$\gbigz_1\subset\gbigz_0$
with the following property.
\begin{itemize}
\item
 $\gbigz_1$ contains the singular locus of $\gbigz_0$ 
and $H^{[2]}$.
\item
 $\gbigr_1\setminus\varpi^{-1}(\gbigz_1)\lrarr H$ 
 is relatively $0$-dimensional.
\item
 Let $P_0$ be any point of $\gbigz_0\setminus \gbigz_1$,
 and let $H_{P_0}$ be any relatively compact
 open neighbourhood of $P_0$ in $H$.
 Let $Q$ be any point of
 $\varpi^{-1}(\gbigz_0\setminus \gbigz_1)\setminus \gbigr_1$
 such that $P:=\varpi(Q)\in H_{P_0}$.
Let  $(\nbign;y_1,\ldots,y_{2n-2})$ denote 
a real analytic coordinate neighbourhood
of $H$ around $P$
 such that
 $\gbigz_0\cap \nbign=\{y_{2n-2}=0\}$.
Let  $\theta_1<\theta_2$
be real numbers such that
the interval
 $\closedclosed{\theta_1}{\theta_2}$
 is a small neighbourhood of $Q$
 in $\varpi^{-1}(P)\setminus \gbigr_1$.
Then, there exist
 a positive integer $\rho$,
 a positive integer $m$,
 a positive number $C>0$,
 connected components
 $\nbigc(Q,\pm)$ of $\Xtilde(H)\setminus\Xtilde(H)^{(1)}$
such that
\[
 \nbigu_{\pm}=
 \bigl\{(y_1,\ldots,y_{2n-2},\theta,r)\,\big|\,
 (y_1,\ldots,y_{2n-2})\in\nbign,\,\,\,
 \theta_1<\theta<\theta_2,\,\,\,
 0<r<C(\pm y_{2n-2})^{m}\,\,
 \bigr\}
\subset
 \nbigc(Q,\pm),
\]
and that the restriction of 
$h_p^{\nbigc(Q,\pm)}$ to $\nbigu_{\pm}$
are expressed as
\[
 h_p^{\nbigc(Q,\pm)}
=\sum_{i\geq -N_1}\sum_{j\geq -N_2}
 \alpha_{(\pm,p),i,j}(y_1,\ldots,y_{2n-3},\theta)\cdot
 y_{2n-2}^{i/\rho}\cdot
 (y_{2n-2}^{-m}r)^{j/\rho}.
\]
Here, $\alpha_{(\pm,p),i,j}$ are real analytic functions.
Moreover, Assumption {\rm\ref{assumption;16.8.27.10}}
is satisfied around $P_0$.
\hfill\qed
\end{itemize}
\end{lem}

\subsubsection{Proof of Proposition \ref{prop;16.6.26.10}}

Let $P\in H$.
There exists a neighbourhood $X_P$
of the form
$\Delta^{\ell}\times X_{P,0}$
such that
$H_P=H\cap X_P=\bigcup_{i=1}^{\ell}\{z_i=0\}$.
Set $\Xtilde_P(H_P):=\varpi^{-1}(X_P)$.
Let $\nbigu$ be a connected component of
$\del\Xtilde_P(H_P)\setminus \gbigr_0$.
By enlarging $\gbigr_i$ and $\gbigz_i$,
we may assume that $\nbigu$ is simply connected.
For any $Q\in \del\Xtilde_P(H_P)\setminus\gbigr_0$,
there exist functions
$h^Q_{1},\ldots,h^Q_{\ell}$
as in the condition
of Lemma \ref{lem;16.7.25.20}.
We consider the case $Q\in\varpi^{-1}(\{z_1=0\})$.
Let $z_1=r_1e^{\sqrt{-1}\theta_1}$
be the polar coordinate system.
Let 
$h^Q_{i}=\sum_{j\geq -N(i)} h^Q_{i,j}r_1^{j/\rho}$
denote the expansion of $h^Q_i$
around $Q$.
We obtain the real analytic functions
$h^Q_{i,j}$ $(j<0)$ on a neighbourhood of $Q$
in $\del \Xtilde_P(H_P)$.
By varying $Q$ in $\nbigu$,
we obtain real analytic functions
$h^{\nbigu}_{i,j}$ on $\nbigu$
for $i=1,\ldots,\ell$ and for $j=-1,\ldots,-N(i)$.
Set $h^{\nbigu}_i:=\sum_{j<0} h^{\nbigu}_{i,j}r_1^{j/\rho}$.

\begin{lem}
\label{lem;16.7.25.30}
$h^{\nbigu}_{i,j}$ are subanalytic functions
on $(\nbigu,\Xtilde(H))$.
\end{lem}
\pf
There exists a connected component 
$\nbigc$ 
of $\Xtilde(H)\setminus \Xtilde(H)^{(1)}$
such that
$\nbigu\subset\overline{\nbigc}$.
We may assume that
$h^{\nbigc}_i-h^{\nbigu}_{i}$
are bounded on any neighbourhood of
$Q\in\nbigu$ in $\Xtilde(H)$.
Hence, 
$h^{\nbigu}_{i,-N(i)}$
is described as the restriction of
$r_1^{N(i)}h^{\nbigc}_i$
to $r_1=0$.
Hence, we obtain that
$h^{\nbigc}_{i,-N(i)}$
is subanalytic 
on $(\nbigu,\del\Xtilde(H))$.
Suppose that we have already known that
$h^{\nbigc}_{i,j}$ are subanalytic
on $(\nbigu,\del\Xtilde(H))$
for $-N(i)\leq j\leq -m-1$ for $m>0$.
Then,
$h^{\nbigu}_{i,-m}$
is described as the restriction of
$r_1^{m/\rho}
\Bigl(
 h_i^{\nbigc}
-\sum_{j<-m}
 r_1^{-j/\rho}h^{\nbigu}_{i,j}
\Bigr)$
to $r_1=0$.
Hence, we obtain the claim 
Lemma \ref{lem;16.7.25.30}.
\hfill\qed

\vspace{.1in}

Let 
$P_1$ be any point of $(H_P\setminus \gbigz_0)\cap\{z_1=0\}$.
By Lemma \ref{lem;16.2.1.3},
we obtain the set of ramified irregular values  $\nbigi_{P_1}$
at $P_1$.
Each $f^{P_1}_i\in\nbigi_{P_1}$ is of the form
$f^{P_1}_i=\sum_{j=1}^{N(i)}
 f^{P_1}_{i,j}z^{-j/\rho}$.
We may assume that 
each connected component $U$ of
$(H_P\setminus \gbigz_0)\cap\{z_1=0\}$
is simply connected.
We obtain a meromorphic function
$f^{U}_i=\sum f^{U}_{i,j}z^{-j/\rho}$
whose restriction to a neighbourhood of $P_1\in U$
is equal to $f^{P_1}_{i}$.
In particular,
we obtain holomorphic functions
$f^{U}_{i,j}$ on $U$.

\begin{lem}
$\Re f^{U}_{i,j}$
and $\Image f^{U}_{i,j}$
are subanalytic functions
on $(U,H)$.
\end{lem}
\pf
Set $H_{P,j}:=\{z_j=0\}$
and 
$H_{P,1}^{\circ}:=
 H_{P,1}\setminus
 \Bigl(
 \bigcup_{j=2}^{\ell}H_{P,j}\Bigr)$.
We may assume that
$U$ is contained in
$H_{P,1}^{\circ}$.
Note that
$\varpi^{-1}(H_{P,1}^{\circ})
\simeq
 S^1\times H_{P,1}^{\circ}$.
For any $\theta\in S^1$,
the set
$\{\theta\}\times
 H_{P,1}^{\circ}$
is a locally closed subanalytic subset
in $\del\Xtilde(H)$.
There exists a $0$-dimensional closed subset
$Z_{S^1}\subset S^1$
such that
$\dim\Bigl(
 \bigl(
 \theta\times H_{P,1}^{\circ}
 \bigr)
 \cap 
 \gbigr_0
\Bigr)
\leq 2n-3$
for any $\theta\in S^1\setminus Z_{S^1}$.
Let $\nbigb_{\theta}$ be any connected component of
$(\theta\times H_{P,1}^{\circ})
\setminus \gbigr_0$.
By the previous lemma,
we obtain that
$\Re(f^U_{i,j}e^{j\sqrt{-1}\theta})$
are subanalytic functions on
$(\nbigb_{\theta}\cap (\theta\times U),\del\Xtilde(H))$.
Hence, we obtain that 
$\Re(f^U_{i,j}e^{j\sqrt{-1}\theta})$
are subanalytic functions on 
$(\theta\times U,\del\Xtilde(H))$.
Then, we obtain the claim of the lemma.
\hfill\qed

\begin{lem}
Let $P_1$ be any point of $\gbigz_0\setminus \gbigz_1$.
If $P_1\in\Ubar$,
$f^U_{i,j}$ are holomorphic
on a neighbourhood of $P_1$.
\end{lem}
\pf
It follows from 
Lemma \ref{lem;16.6.25.1}.
\hfill\qed

\vspace{.1in}
There exists a closed subanalytic subset
$Z_{\Ubar}\subset\Ubar$
with $\dim_{\real} Z_{\Ubar}\leq 2n-4$
such that the following holds.
\begin{itemize}
\item
 $Z_{\Ubar}\supset
 \gbigz_1\cap\Ubar$.
\item
 For any $P_1\in \Ubar\setminus Z_{\Ubar}$,
 $\{f^U_i\}$ is a good set of ramified irregular values
 at $P_1$.
\end{itemize}

For each 
$P_1\in U\setminus Z_{\Ubar}$,
let $X_{P_1}$ be any neighbourhood of $P_1$
in $X_P$,
and we set $H_{P_1}:=H_P\cap X_{P_1}$.
By Proposition \ref{prop;16.7.25.40},
there exist a good meromorphic flat bundle
$V_{P_1}$ on $(X_{P_1},H_{P_1})$
whose set of ramified irregular values
is $\{f^U_i\}$,
and an isomorphism
\[
 \DR^{\Ecat}_{\vecX_{P_1}(H_{P_1})}(V_{P_1})[-n]
\simeq
 K_{|X_{P_1}}.
\]
We put $C(U):=(H\setminus U)\cup Z_{\Ubar}$,
$X_U:=X_P\setminus C(U)$
and $H_U:=H_P\setminus C(U)$.
By gluing $V_{P_1}$ 
for $P_1\in U\setminus Z_{\Ubar}$,
we obtain a good meromorphic flat bundle $V_{U}$
on $(X_U,H_U)$
and an isomorphism
\[
 \DR^{\Ecat}_{\vecX_{U}(H_{U})}(V_{U})[-n]
\simeq
 K_{|X_{U}}.
\]

Let $Z^{(1)}_{P}$ denote the union of $Z_{\Ubar}$
for the connected components
$U$ of $(H_P\setminus \gbigz_0)\cap\{z_1=0\}$.
We set 
$X_P^{(1)}:=
 X_P\setminus\bigl(Z_P^{(1)}\cup\bigcup_{j=2}^{\ell}\{z_j=0\}\bigr)$,
and 
$H_P^{(1)}:=H_P\cap X_P^{(1)}$.

Let $P_1$ be any point of $H_P^{(1)}\cap \gbigz_0$.
Note that $\gbigz_0$ is smooth around $P_1$.
Let $U_i$ $(i=1,2)$ denote the connected components
of $(H_P\setminus \gbigz_0)\cap\{z_1=0\}$
such that $P_1\in\Ubar_i$.
Note that 
$\{f_i^{U_a}\}$ $(a=1,2)$ are good at $P$.
Let $X_{P_1}$ be a small neighbourhood 
of $P_1$ in $X_P^{(1)}$,
and set $H_{P_1}:=X_{P_1}\cap H_P^{(1)}$.
Then,
there exist good meromorphic flat bundles
$V_{P_1}^{(a)}$ $(a=1,2)$
on $(X_{P_1},H_{P_1})$
such that
the restriction of $V_{P_1}^{(a)}$
and $V_{U_a}$
to $X_{P_1}\cap X_{U_a}$
are isomorphic.
Note that Assumption \ref{assumption;17.12.27.10}
is satisfied around $P_1$ due to Lemma \ref{lem;16.8.25.1}.
By our choice of $\gbigz_1$,
Assumption \ref{assumption;16.8.27.10} is satisfied
around $P_1$.
By Corollary \ref{cor;16.8.27.12},
we obtain 
$V^{(1)}_{P_1}\simeq V^{(2)}_{P_1}=:V_{P_1}$,
and there exists an isomorphism
\[
 \DR^{\Ecat}_{\vecX_{P_1}(H_{P_1})}(V_{P_1})[-n]
\simeq
 K_{|X_{P_1}}.
\]

By gluing $V_U$ for connected components $U$,
and $V_{P_1}$ for $P_1\in H_P^{(1)}\cap \gbigz_0$,
we obtain a good meromorphic flat bundle
$V^{(1)}$ 
on $(X^{(1)}_P,H^{(1)}_P)$
and an isomorphism
\[
  \DR^{\Ecat}_{\vecX^{(1)}_P(H^{(1)}_{P})}(V^{(1)}_{P})[-n]
\simeq
 K_{|X^{(1)}_{P}}.
\]

Similarly,
for each $i=2,\ldots,\ell$,
we obtain the following:
\begin{itemize}
\item
 a closed subanalytic subset
 $Z_P^{(i)}$ of $H_P\cap\{z_i=0\}$
 for which
 $(\gbigz_1\cap \{z_i=0\})\subset Z_P^{(i)}$,
\item
 a good meromorphic flat bundle
 $V^{(i)}$
 on $(X_P^{(i)},H_P^{(i)})$
where
\[
  X_P^{(i)}:=X_P\setminus\Bigl(
 Z_P^{(i)}\cup\bigcup_{\substack{1\leq j\leq \ell\\ j\neq i}}
 \{z_j=0\}
 \Bigr),
\quad 
 H_P^{(i)}=X_P^{(i)}\cap H_P,
\]
\item
 an isomorphism
$\DR^{\Ecat}_{\vecX^{(i)}_P(H^{(i)}_{P})}(V^{(i)}_{P})[-n]
\simeq
 K_{|X^{(i)}_{P}}$.
\end{itemize}
We set $Z_P:=\bigcup_{i=1}^{\ell}Z_P^{(i)}$.
By gluing $V_P^{(i)}$ $(i=1,\ldots,\ell)$,
we obtain a good meromorphic flat bundle
$V_P$ on $(X_P',H_P'):=(X_P\setminus Z_P,H_P\setminus Z_P)$,
and an isomorphism
$\DR^{\Ecat}_{\vecX'_P(H'_{P})}(V_{P})[-n]
\simeq
 K_{|X'_{P}}$.

There exist a discrete subset $S\subset X$
such that $X=\bigcup_{P\in S} X_P$.
Let $Z$ denote the union of the closures of $Z_P$ $(P\in S)$.
By gluing $V_P$ $(P\in S)$,
we obtain a good meromorphic flat bundle $V$
and the desired isomorphism.
Thus, the proof of Proposition \ref{prop;16.6.26.10}
is finished.
\hfill\qed

\section{Proof of Theorem \ref{thm;18.11.12.10}}
\label{section;18.11.15.21}

We prove Theorem \ref{thm;18.11.12.10}.
The surface case is studied in \S\ref{subsection;18.11.16.130}.
The higher dimensional case is studied in 
\S\ref{subsection;18.11.16.131}.

\subsection{Surface case}
\label{subsection;18.11.16.130}

Let $X$ be a complex surface
with a normal crossing complex hypersurface $H$.
Let $K\in \Ecat^b_{\circledcirc}(\IC_{\vecX(H)})$.
In this subsection,
we shall prove the following proposition
which is the surface case of Theorem \ref{thm;18.11.12.10}.
\begin{prop}
\label{prop;16.5.10.10}
$K$ is contained in
$\Ecat^b_{\mero}(\IC_{\vecX(H)})$.
\end{prop}

\subsubsection{Preliminary}

Let $Y$ be any complex surface
with a normal crossing hypersurface $H_Y$.
Let $\varpi:\Ytilde(H_Y)\lrarr Y$
denote the oriented real blowing up of $Y$
along $H_Y$.
Let $K\in \Ecat^b_{\circledcirc}(\IC_{\vecY(H_Y)})$.
Let $P$ be a cross point of $H_Y$.
We introduce an auxiliary condition
for $K$ and $P$.

\begin{condition}
\label{condition;18.11.26.100}
There exists a filtration
 $\Ytilde(H_Y)=\Ytilde(H_Y)^{(0)}\supset
 \Ytilde(H_Y)^{(1)}\supset\cdots$ for $K$
such that the following holds.
\begin{itemize}
\item
 Let $\gbigw$ denote the closure of
 $\Ytilde(H_Y)^{(1)}\setminus\del\Ytilde(H_Y)$
 in $\Ytilde(H_Y)$.
 Then, 
 $\dim_{\real}(\gbigw\cap\varpi^{-1}(P))\leq 1$
 holds.
 Moreover, 
 $\gbigw\cap \varpi^{-1}(P)$
 does not contain any positively linear subspace. 
\hfill\qed
\end{itemize}
\end{condition}

\subsubsection{Resolution for auxiliary conditions at cross points}
\label{subsection;18.11.15.11}

Let $X:=\Delta^2$ and
$H:=\{x_1=0\}\cup\{x_2=0\}$,
or $H:=\{x_1=0\}$.
Set $O:=(0,0)$.
We use the notation in \S\ref{section;18.11.15.10}.
For any $\veceta\in\{+,-\}^{\ell}$,
let $\Bl_{\veceta}X$ 
(resp. $H_{\veceta}$)
denote the inverse image of $X$
(resp. $H$)
by the morphism
$\psi_{\veceta}:
 \Bl_{\veceta}\cnum^2\lrarr \cnum^2$.
The induced morphism
$\Bl_{\veceta}X\lrarr X$
is also denoted by $\psi_{\veceta}$.
We shall prove the following proposition.

\begin{prop}
\label{prop;16.2.4.100}
For any $K\in\Ecat^b_{\circledcirc}(\IC_{\vecX(H)})$,
there exists a positive integer
$\ell_0$
such that the following holds
for any $\ell\geq \ell_0$ and $\veceta\in\{+,-\}^{\ell}$:
\begin{itemize}
\item
Condition {\rm\ref{condition;18.11.26.2}}
and Condition {\rm\ref{condition;18.11.26.100}}
are satisfied for 
$\Ecat\psi_{\veceta}^{-1}K$ 
and $P_{\veceta}$.
\end{itemize}
\end{prop}
\pf
Let us consider the following conditions
in the case $H=\{z_1=0\}\cup\{z_2=0\}$.
\begin{description}
\item[(C1)]
 Let $\varpi:\Xtilde(H)\lrarr X$ be the oriented real blowing up.
 There exists a closed subanalytic subset $Z\subset\varpi^{-1}(O)$
 with $\dim Z\leq 1$
 such that the following holds
 for any $Q\in \varpi^{-1}(O)\setminus Z$:
 \begin{itemize}
 \item[-]
 There exist a neighbourhood $\nbigu_Q$ of $Q$ in $\Xtilde(H)$
 and ramified real analytic functions $g_1^Q,\ldots,g_m^Q$
 on $\nbigu_Q$.
 \item[-]
 The growth order of
 $\pi^{-1}(\cnum_{\nbigu_Q})\otimes K$
 is controlled by
 $g_1^Q,\ldots,g_m^Q$.
 \end{itemize}
\item[(C2)]
 $Z$ does not contain
 any positively linear subset.
\end{description}

\begin{lem}
\label{lem;16.1.16.2}
Suppose that 
the conditions {\bf (C1)} and {\bf(C2)} are satisfied.
Then, there exist a neighbourhood $X'$ of $O$ in $X$
and a set of ramified irregular values $\nbigi$ at $O$
such that the following holds.
\begin{itemize}
\item
 Let $\varphi:(\Delta^{\ast},0,\Delta)\lrarr 
 (X'\setminus H,H\cap X',X')$
be any holomorphic map.
Then,
$\Irr(\Ecat\varphi^{-1}K)
=\varphi^{\ast}\nbigi$
holds.
\end{itemize}
\end{lem}
\pf
We use the polar coordinate system
$(r_1,r_2,\theta_1,\theta_2)$
determined by $z_i=r_ie^{\sqrt{-1}\theta_i}$.
Let $Q$ be any point 
in $\varpi^{-1}(O)\setminus Z$.
Let $\nbigu_Q$ be a small neighbourhood of $Q$
in $\Xtilde(H)$ 
of the form
$\nbigu_Q=\bigl\{
 (r_1,r_2,\theta_1,\theta_2)\,\big|\,
 0\leq r_i\leq \epsilon_1,\,\,
 |\theta_i-\theta^0_{i}|<\epsilon_2
 \bigr\}$
for some positive numbers $\epsilon_i$
and $(\theta_{1}^0,\theta_2^0)\in\real^2$.
Let $g_1^Q,\ldots,g_m^Q$ 
denote the functions on $\nbigu_Q$
as above.
They are expressed as the convergent power series:
\[
 g_i^Q=\sum_{j\geq -N_1}\sum_{k\geq -N_2}
 \alpha_{i,j,k}(\theta_1,\theta_2)
 r_1^{j/\rho}r_2^{k/\rho}.
\]
By Lemma \ref{lem;18.11.6.1},
by considering the pull back via a ramified covering,
we may assume $\rho=1$.
We set
\[
 g_{i1}^Q:=\sum_{-N_1\leq j}\sum_{-N_1\leq k<0}
  \alpha_{i,j,k}(\theta_1,\theta_2)
 r_1^{j}r_2^{k},
\quad
 g_{i2}^Q:=\sum_{-N_1\leq j<0}\sum_{-N_1\leq k}
  \alpha_{i,j,k}(\theta_1,\theta_2)
 r_1^{j}r_2^{k}.
\]
Set $H_2:=\{z_2=0\}$.
Let $Q'$ be any point of $\nbigu_Q\cap
 \varphi^{-1}(H_1\setminus H_2)$.
Then, $g_i$ and $g_{i,2}$ are mutually bounded
around $Q'$.
By applying the argument in the proof of
Lemma \ref{lem;16.2.1.3} to $g_{i,2}$,
we obtain meromorphic functions
\[
 f_{i2}^Q=\sum_{-N_1\leq j<0}
 \sum_{-N_1\leq k}
 a^2_{i,j,k}z_1^{j}z_2^{k}, 
\]
such that
$\Re(f_{2a}^Q)=g_{2a}^Q$.
Similarly, 
we obtain meromorphic functions
\[
 f_{i1}^Q=\sum_{-N_1\leq j}
 \sum_{-N_1\leq k<0}
 a^1_{i,j,k}z_1^{j}z_2^{k}
\]
such that
$\Re(f_{i1}^Q)=g_{i1}^Q$.
Because
\[
 \Re\Bigl(
 \sum_{-N_1\leq j<0}
 \sum_{-N_1\leq k<0}
 a^1_{i,j,k}z_1^{j}z_2^{k}
 \Bigr)
=
\sum_{-N_1\leq j<0}
 \sum_{-N_1\leq k<0} 
\alpha_{i,j,k}(\theta_1,\theta_2)
 r_1^{j}r_2^{k}
=\Re\Bigl(
 \sum_{-N_1\leq j<0}
 \sum_{-N_1\leq k<0}
 a^2_{i,j,k}z_1^{j}z_2^{k}
 \Bigr),
\]
we obtain
\[
 \sum_{-N_1\leq j<0}
 \sum_{-N_1\leq k<0}
 a^1_{i,j,k}z_1^{j}z_2^{k}
=\sum_{-N_1\leq j<0}
 \sum_{-N_1\leq k<0}
 a^2_{i,j,k}z_1^{j}z_2^{k}.
\]
We set
\[
f^Q_i:=f_{i1}^Q+f_{i2}^Q
- \sum_{-N_1\leq j<0}
 \sum_{-N_1\leq k<0}
 a^1_{i,j,k}z_1^{j}z_2^{k}.
\]
Then, 
$\Re(f^Q_i)-g_i^Q$
are bounded
on $\nbigu_Q$.

By shrinking $X$,
we may assume that
$X=\varpi(\nbigu_Q)$.
Let $P$ be any point of 
$H\setminus \{(0,0)\}$.
Let $\varphi:\Delta\lrarr X$ be a holomorphic map
such that $\varphi(0)=P$.
Then, 
$\Irr(\Ecat\varphi^{-1}K)
=\{\varphi^{\ast}f^Q_i\}$
holds.
In particular,
$\{f^Q_i\}$ is independent of $Q$.
We denote them by
$\{f_i\}$.

Let $\varphi:\Delta_{\epsilon}\lrarr X$
be any holomorphic map.
If $\varphi(0)\neq (0,0)$,
then $\Irr(\Ecat\varphi^{-1}K)=
 \{\varphi^{\ast}f_i\}$ holds,
as already mentioned.
Suppose that $\varphi(0)=(0,0)$.
Let 
$\varphitilde:
 \Deltatilde_{\epsilon}(0)
 \lrarr \Xtilde(H)$
denote the induced map.
Note
$\dim(\varphitilde(\varpi_0^{-1}(0))\cap Z)=0$
because $Z$ does not contain
any positively linear subset.
Then, we obtain that
$\Irr(\Ecat\varphi^{-1}K)
=\{\varphi^{\ast}f_i\}$.
Thus, we obtain Lemma \ref{lem;16.1.16.2}.
\hfill\qed

\vspace{.1in}

By Lemma \ref{lem;16.4.27.30},
Proposition \ref{prop;16.4.27.11}
and Proposition \ref{prop;16.4.27.20},
we can assume that
the conditions {\bf (C1,2)}
are satisfied
for $\Ecat\psi_{\veceta}^{-1}K$
at $P_{\veceta}$.
Then, the claim of Proposition \ref{prop;16.2.4.100}
follows from Lemma \ref{lem;16.1.16.2}.
\hfill\qed

\subsubsection{Resolutions at cross points
and generic parts in the surface case}

We continue to use the notation in \S\ref{subsection;18.11.15.11}.
\begin{prop}
\label{prop;16.2.7.40}
Suppose that 
Condition {\rm\ref{condition;18.11.26.2}}
is satisfied for 
$K\in\Ecat_{\circledcirc}^b(\vecX(H))$ and $O$.
Then, there exists a meromorphic flat bundle
$(V,\nabla)$ on $(X,H)$
with an isomorphism
$\DR^{\Ecat}_{\vecX(H)}(V)[-2]
\simeq K$.
\end{prop}
\pf
Let $\nbigi$ be the set of ramified irregular values at $O$
for which Condition \ref{condition;18.11.26.2} is 
satisfied for $K$ and $O$.
After $X$ is shrunk,
it induces a system of ramified irregular values
$\vecnbigi$ on $(X,H)$.
By using \cite[Proposition 15.1.5]{Mochizuki-MTM},
there exists a projective morphism
$G:X'\lrarr X$
such that 
(i) $H'=G^{-1}(H)$ is normal crossing,
(ii) $G$ induces an isomorphism
 $X'\setminus H'\simeq X\setminus H$,
(iii) $G^{\ast}\vecnbigi$ is 
a good system of ramified irregular values.
It is enough to construct 
a meromorphic flat connection $(V'_P,\nabla_P)$
with the desired property
around any point of $P\in H'$.
By Proposition \ref{prop;16.7.25.40},
it is enough to 
study the case $P$ is a cross point of $H'$.
Hence, we may assume that $\nbigi_O$ is 
a good set of ramified irregular values
from the beginning.
By considering a ramified covering,
we may assume that $\nbigi_O$
is a good set of unramified irregular values.
By Proposition \ref{prop;16.2.4.100},
we may also assume that 
Condition \ref{condition;18.11.26.100}
is satisfied for $K$ and $O$.
Let $\Ltilde$ denote the local system
on $\Xtilde(H)$
induced by $K_{|X\setminus H}$.

Let $H_i=\{z_i=0\}$.
We may assume that 
the natural map
$\nbigi\lrarr \nbigo_{X}(\ast H)/\nbigo_X(\ast H_2)$
is injective.

According to Proposition \ref{prop;16.6.26.10},
after shrinking $X$ around $O$,
we may assume that
there exists a good meromorphic flat bundle
$(V^{\circ},\nabla^{\circ})$
on $(X\setminus O,H\setminus O)$
such that 
$\DR^{\Ecat}(V^{\circ},\nabla^{\circ})[-2]
\simeq
 K_{|X\setminus O}$.
Set
$(V_1,\nabla_1):=
 (V^{\circ},\nabla^{\circ})_{|X\setminus H_2}$.
By Proposition \ref{prop;16.2.7.21},
there exists a unique good meromorphic flat bundle
$(V,\nabla)$ on $(X,H)$
whose restriction to $X\setminus H_2$
is isomorphic to $(V_1,\nabla_1)$.

Let $\varpi:\Xtilde(H)\lrarr X$
be the oriented real blowing up of $X$ along $H$.
Let $\Xtilde(H)=\Xtilde(H)^{(0)}\supset\Xtilde(H)^{(1)}\supset
 \cdots$
be a filtration for $\Ecat\varpi^{-1}K$.
Let $\gbigw$ be the closure of 
$\Xtilde(H)^{(1)}\setminus\del\Xtilde(H)$.
By Condition \ref{condition;18.11.26.100},
we may assume 
$\dim_{\real}\bigl(\gbigw\cap\varpi^{-1}(O)\bigr)\leq 1$,
and that $\gbigw\cap\varpi^{-1}(O)$
does not contain any positively linear subset.
\begin{lem}
For any point $Q$ of $\varpi^{-1}(O)\setminus \gbigw$,
there exists a small neighbourhood $\nbigu_Q$
in $\Xtilde(H)$
such that
$\pi^{-1}(\cnum_{\nbigu^{\circ}_Q})\otimes K$
is controlled by
$\Re(f)$ $(f\in\nbigi_O)$,
i.e.,
\begin{equation}
\label{eq;16.2.7.30}
 \pi^{-1}(\cnum_{\nbigu_Q^{\circ}})\otimes K
\simeq
 \bigoplus
 \cnum_{\Xtilde(H)}^{\Ecat}\overset{+}{\otimes}
 \bigl(
 \cnum_{t\geq \Re(f)}\otimes V_f
 \bigr).
\end{equation}
Here, $\nbigu_{Q}^{\circ}=\nbigu_Q\setminus\varpi^{-1}(H)$.
\end{lem}
\pf
For a small neighbourhood $\nbigu_Q$,
there exist continuous subanalytic functions
$h_1,\ldots,h_{r}$ on $(\nbigu_Q^{\circ},\Xtilde(H))$
such that 
$\pi^{-1}(\cnum_{\nbigu_Q^{\circ}})\otimes K
\simeq
 \bigoplus
 \cnum_{\Xtilde(H)}^{\Ecat}\overset{+}{\otimes}
 \bigl(
 \cnum_{t\geq h_i}
 \bigr)$.
If $\{h_1,\ldots,h_r\}\neq \{\Re(f_i)\}$
in $\Subbar(\nbigu_Q^{\circ},\Xtilde(H))$,
there exists an analytic path
$\gamma:(\II^{\circ},0,\II)\lrarr 
(\nbigu_Q^{\circ},\varpi^{-1}(H),\Xtilde(H))$
such that
$\gamma^{\ast}\{h_1,\ldots,h_r\}
\neq
 \gamma^{\ast}\{\Re(f_i)\}$
in $\Subbar(\II^{\circ},\II)$.
For a small $\epsilon>0$,
let $\gamma_{\cnum}:\Delta_{\epsilon}\lrarr X$
denote the holomorphic map
induced by $\gamma$.
The condition
$\Irr(\Ecat\gamma_{\cnum}^{-1}K)
=\{\gamma_{\cnum}^{\ast}(f_i)\}$,
implies 
$\gamma^{\ast}\{h_1,\ldots,h_r\}
\neq
 \gamma^{\ast}\{\Re(f_i)\}$
in $\Subbar(\II^{\circ},\II)$.
Thus, we have arrived at a contradiction.
\hfill\qed

\vspace{.1in}

There exists the canonical filtration
on $\pi^{-1}(\cnum_{\nbigu_Q^{\circ}})\otimes K$
as in \S\ref{subsection;16.1.25.20},
which we denote by $\nbigf^Q$.
We may assume that
$\leq_Q$ on $\nbigi_O$
is equal to
$\prec$
on $\{\Re(f)_{|\nbigu_Q^{\circ}}\,|\,f\in\nbigi_O\}$.
The decomposition (\ref{eq;16.2.7.30})
is compatible with 
the Stokes filtration of $(V_1,\nabla_1)$
at $Q'\in \varpi^{-1}(H_1\setminus H_2)$.
Hence, by the construction in the proof of
Proposition \ref{prop;16.2.7.21},
we obtain that
$\nbigf^Q$  is equal to 
the Stokes filtration of $(V,\nabla)$ at $Q$.
Let $Q'\in \varpi^{-1}(H_2\setminus H_1)$
be contained in $\nbigu_Q$.
Then, the Stokes filtration of $(V,\nabla)$ at $Q'$
is equal to the filtration induced by the decomposition
(\ref{eq;16.2.7.30}).

Let $\varphi:(\Delta^{\ast},0,\Delta)\lrarr (X\setminus H,H,X)$
be any holomorphic map.
By the previous consideration,
we obtain that
the Stokes filtrations of
$\Ecat\varphi^{-1}K$
and $\varphi^{\ast}V$
are the same 
at general points of $\varpi_{0}^{-1}(0)$,
where we recall that $\varpi_{0}:\Deltatilde(0)\lrarr \Delta$
denotes the oriented real blowing up.
Hence, we obtain that
$\DR^{\Ecat}\varphi^{\ast}(V)[-1]
\simeq
 \Ecat\varphi^{-1}K$.
Then, the claim of Proposition \ref{prop;16.2.7.40}
follows from Proposition \ref{prop;16.6.23.20}.
\hfill\qed

\begin{cor}
\label{cor;16.5.10.1}
For any $K\in\Ecat^b_{\circledcirc}(\IC_{\vecX(H)})$,
there exists a positive integer
$\ell_0$
such that the following holds
for any 
$\ell\geq \ell_0$
and $\veceta\in\{+,-\}^{\ell}$:
\begin{itemize}
\item
There exists a meromorphic flat bundle $V$
defined on a neighbourhood $U$ around $P_{\veceta}$,
and an isomorphism
$\Ecat\psi_{\veceta}^{-1}(K)_{|U}
\simeq
 \DR^{\Ecat}(V)$.
\end{itemize}
\end{cor}
\pf
By Proposition \ref{prop;16.2.4.100},
we may assume that 
Condition \ref{condition;18.11.26.2}
and Condition \ref{condition;18.11.26.100}
are satisfied for $\Ecat\psi_{\veceta}^{-1}K$
and $P_{\veceta}$.
Then, the claim of the corollary follows from
Proposition \ref{prop;16.2.7.40}.
\hfill\qed

\subsubsection{Proof of Proposition \ref{prop;16.5.10.10}}

Let $X$ be any complex manifold
with a normal crossing hypersurface $H$.
It is enough to prove the following
for any $P\in H$:
\begin{description}
\item[$\nbigc(P)$:]
There exists a neighbourhood $U$ of $P$ in $X$
such that $K_{|U}\in \Ecat^b_{\mero}(\IC_{\vecU(H\cap U)})$.
\end{description}
Hence, we may assume that
$X$ is a relatively compact subset
in another  complex surface $X'$,
and $K$ is the restriction of
$K'\in \Ecat^b_{\circledcirc}(\IC_{\vecX'(H')})$.
By Proposition \ref{prop;16.6.26.10},
Proposition \ref{prop;16.2.7.40}
and Corollary \ref{cor;16.5.10.1},
there exists a projective morphism
$\rho:X_1\lrarr X$ such that
(i) $H_1:=\rho^{-1}(H)$ is normal crossing,
(ii) $X_1\setminus H_1\simeq X\setminus H$,
(iii) at any cross point of $P_1\in H_1$,
 the condition $\nbigc(P_1)$ for $\Ecat\rho^{-1}(K)$
 is satisfied.
Once we know that
$\Ecat\rho^{-1}(K)\in
 \Ecat^b_{\mero}(\IC_{\vecX_{1}(H_1)})$,
then we obtain that
$K\in \Ecat^b_{\mero}(\IC_{\vecX(H)})$.
Hence, we may assume that 
$\nbigc(P)$ holds for any cross point $P$ of $H$
from the beginning.

Let $P$ be a smooth point of $H$.
Suppose that $\nbigc(P)$ does not hold
for $K$.
We set $X^{(0)}:=X$ and $P^{(0)}:=P$.
Let 
$\Turn^{(1)}$ denote the set of the points
$P'\in (\rho^{(1)})^{-1}(P^{(0)})$
such that 
$\nbigc(P^{'})$ does not hold for
$K^{(1)}:=\Ecat(\rho^{(1)})^{-1}(K)$.
We have already known that $\Turn^{(1)}$ is finite
by Proposition \ref{prop;16.6.26.10}.
If $\Turn^{(1)}$ is not empty,
we choose a point 
$P^{(1)}\in \Turn^{(1)}$.
Inductively,
we construct a sequence of morphisms
\[
 X^{(\ell)}\stackrel{\rho^{(\ell)}}{\lrarr}
 X^{(\ell-1)}\stackrel{\rho^{(\ell-1)}}{\lrarr}
 \cdots
 \stackrel{\rho^{(2)}}{\lrarr}
 X^{(1)}\stackrel{\rho^{(1)}}{\lrarr}
 X^{(0)}
\]
with $P^{(i)}\subset(\rho^{(i)})^{-1}(P^{(i-1)})$
$(i=0,1,\ldots,\ell-1)$
such that 
(i) $\rho^{(i)}$ is the blowing up of $X^{(i-1)}$
 at $P^{(i-1)}$,
(ii) $\nbigc(P^{(i)})$ does not hold 
for $\Ecat(\psi^{(i)})^{-1}(K)$
at $P^{(i)}$,
where $\psi^{(i)}:X^{(i)}\lrarr X$
denotes the induced morphism.
We stop the process 
if $\nbigc(Q)$ holds for 
$\Ecat(\psi^{(\ell)})^{-1}(K)$
at any $Q\in (\rho^{(\ell)})^{-1}(P^{(\ell-1)})$.
\begin{lem}
If any such sequence is finite,
its length is bounded independently of the sequence.
\end{lem}
\pf
Suppose that there exists a sequence 
$\ell_j\to\infty$,
and sequences of morphisms 
\[
 X_j^{(\ell_j)}\stackrel{\rho_j^{(\ell_j)}}{\lrarr}
 X_j^{(\ell_j-1)}\stackrel{\rho_j^{(\ell_j-1)}}{\lrarr}
 \cdots
 \stackrel{\rho_j^{(2)}}{\lrarr}
 X_j^{(1)}\stackrel{\rho_j^{(1)}}{\lrarr}
 X_j^{(0)}=X
\]
with points $P_j^{(i)}\in X_j^{(i)}$ as above.
We shall derive a contradiction.

By the construction,
$X_j^{(0)}=X$ and $P_j^{(0)}=P$ 
are independent of $j$.
Hence,
$X_j^{(1)}$ 
and $(\rho_j^{{(1)}})^{-1}(P)$
are independent of $j$.
According to Proposition \ref{prop;16.6.26.10},
there exists a finite subset 
$D^{(1)}\subset (\rho_j^{(1)})^{-1}(P)$
such that
$\nbigc(Q)$ holds for
any $Q\in (\rho_j^{(1)})^{-1}(P)
\setminus D^{(1)}$.
Hence, by going to a sub-sequence of $(\ell_j)$,
we may assume that 
$P_j^{(1)}$ are independent of $j$.
Similarly, 
for any $k_0$, 
after going to sub-sequences,
we may assume that
$X_j^{(k)}$ and 
$P_j^{(k)}$ $(k\leq k_0)$ are independent of $j$.
In this way,
we obtain an infinite sequence of complex blowings up
\[
 \cdots\lrarr
 X^{(i)}\stackrel{\rho^{(i)}}{\lrarr}
 X^{(i-1)}\stackrel{\rho^{(i-1)}}{\lrarr}
 \cdots
 \lrarr X^{(1)}
 \stackrel{\rho^{(1)}}{\lrarr} X^{(0)}=X
\]
with points $P^{(i)}\in X^{(i)}$ as above.
It contradicts the assumption.
Thus, we obtain the claim of the lemma.
\hfill\qed

\vspace{.1in}

Let us prove that any such sequence is finite.
Suppose that there exists an infinite sequence,
and we shall deduce a contradiction.
By Corollary \ref{cor;16.5.10.1},
it is described as an infinite sequence of complex blowings up 
associated to 
$\vecY=(\veceta_1,\omega_1,\veceta_2,\omega_2,\ldots)
 \in \prod_{i=1}^{\infty}(\{\pm\}^{\ell(i)}\times \cnum^{\ast})$.
We shall use the notation in \S\ref{subsection;16.6.3.30}.

There exists a filtration
$\Xtilde(H)=\Xtilde(H)^{(0)}
 \supset\Xtilde(H)^{(1)}\supset\cdots$
for $K$.
Let
$\Xtilde(H)\setminus \Xtilde(H)^{(1)}
=\coprod \nbigc_j$
denote the decomposition into the connected components.
We may assume that there exist
subanalytic functions $g_{j,k}$
on $(\nbigc_j,\Xtilde(H))$
which control the growth order of
$\pi^{-1}(\cnum_{\nbigc_j})\otimes K$.

\vspace{.1in}
Suppose  $\kappahat(\vecY)<\infty$.
We use Theorem \ref{thm;16.5.6.1}.
If $m$ is sufficiently large,
there exists an open subset
$\nbigu$ in $\Xtilde_m(H_m)$
such that
(i) $\nbigu\cap \varpi_{m}^{-1}(P_{m})
 \neq\emptyset$,
(ii)
$\psitilde_{m}
 \bigl(\nbigu\setminus\varpi_{m}^{-1}(H_{m})
 \bigr)
\subset
\nbigc_{j_0}$ for some $j_0$.
The functions 
$\psitilde_{m}^{\ast}(g_{j_0,k})$
and 
$\psitilde_{m}^{\ast}(g_{j_0,k_1}-g_{j_0,k_2})$
are ramified real analytic,
and 
if the functions are not constantly $0$,
their $0$-sets are contained in 
$\varpi_{m}^{-1}(H_{m})$.
Moreover,
$\ord_{\rho_m}
 \psitilde_{m}^{\ast}(g_{j_0,k})(u_m,\theta_m,\rho_m)$
and 
$\ord_{\rho_m}
 \psitilde_{m}^{\ast}(g_{j_0,k_1}-g_{j_0,k_2})(u_m,\theta_m,\rho_m)$
are independent of $(u_m,\theta_m)$.
Then, 
by the argument in the proof of Lemma \ref{lem;16.2.1.3},
we obtain that 
Condition \ref{condition;18.11.26.3}
is satisfied for $K$ and $P_m$.
By Proposition \ref{prop;16.7.25.40},
we obtain that 
$\nbigc(P_m)$ holds.
But, it contradicts the assumption
on the infinite sequence.

\vspace{.1in}

Suppose that $\vecY$ is not convergent
and that $\kappahat(\vecY)=\infty$.
Let $m$ be sufficiently large.
Let $\nbigi_m$ be an interval
in $\varpi_m^{-1}(P_m)$,
and let $\nbigv_m$ be a neighbourhood
of $\nbigi_m$ in $\Xtilde_m(H_m)$
as in Theorems \ref{thm;16.10.11.30} and \ref{thm;16.10.11.31}.
Let us consider the case
$\nbigz_m:=\psitilde_m^{-1}(\Xtilde(H)^{(1)})
\cap(\nbigv_m\setminus\varpi_m^{-1}(H_m))
 \neq\emptyset$.
We explain only the case
$\psitilde_m^{-1}(\Xtilde(H)^{(1)})
\cap(\nbigv_m\setminus\varpi_m^{-1}(H_m))
\simeq
 \nbigv_{1,m}$
because the other case can be argued similarly.
Let
$\nbigv_m\setminus
(\varpi_m^{-1}(H_m)\cup\nbigz_m)
=\nbigv_{m,+}\sqcup\nbigv_{m,-}$
denote the decomposition into the connected components.
Let $v_m=\rho_m e^{\sqrt{-1}\theta_m}$
be the polar decomposition.
By Theorems \ref{thm;16.10.11.30} and \ref{thm;16.10.11.31},
there exist ramified real analytic functions
$h^{\pm}_1(\theta_m,\rho_m),\ldots,
 h^{\pm}_k(\theta_m,\rho_m)$ on $\nbigv_{m,\pm}$
such that 
$\pi^{-1}(\cnum_{\nbigv_{m,\pm}})
\otimes K
=\bigoplus \cnum^{\Ecat}
 \overset{+}{\otimes}
 \cnum_{t\geq h^{\pm}_j}$.
There also exist 
subanalytic functions
$h^{0}_1(\theta_m,\rho_m),\ldots,
 h^{0}_k(\theta_m,\rho_m)$
on $\nbigz_m$
such that 
$\pi^{-1}(\cnum_{\nbigz_m})
\otimes K
=\bigoplus \cnum^{\Ecat}
 \overset{+}{\otimes}
 \cnum_{t\geq h^{0}_j}$,
and they are ramified real analytic
as functions on $\nbigv_{1,m}$.
We may assume that 
$h^{\pm}_j$ and $h^0_j$ 
have only negative powers of $\rho_m$.
Recall that 
for a small neighbourhood $\nbigb_m$ of $P_{m}$ in $X_m$
there exists a meromorphic flat bundle
$(V,\nabla)$ on 
$(\nbigb_m\setminus P_{m},
 (\nbigb_m\cap H_{m})\setminus P_{m})$.
Let $q_m:\nbigb_m\lrarr H_m\cap \nbigb_m$
be the projection 
$q_m(u_m,v_m)=u_m$.
By comparing 
$\{h^{\pm}_j\}$, $\{h^0_j\}$ and
the irregular values of
$(V,\nabla)_{|q_m^{-1}(P')}$ $(P'\in\nbigb_m\cap H_m\setminus\{P_m\})$,
we obtain the following.
\begin{itemize}
\item $\{h^+_j\}=\{h^-_j\}=\{h^0_j\}$, 
and the multiplicity functions are the same.
\item There exists a good meromorphic flat bundle
 $(\Vtilde,\nablatilde)$ on $(\nbigb_m,\nbigb_m\cap H_m)$
 such that the restriction of
 $(\Vtilde,\nablatilde)$
 to $\nbigb_m\setminus\{P\}$ is
 $(V,\nabla)$.
\item
 $\Irr(\Vtilde,\nablatilde)$ is contained in 
 $\cnum(\!(v_m^{1/\rho})\!)/\cnum[\![v_m^{1/\rho}]\!]$
 for some $\rho\in\seisuu_{>0}$,
 i.e., it is independent of $u_m$.
 Moreover, the equalities
 $\{h^{\pm}_j\}=\{h^0_{j}\}=
 \{\Re(\gminia)\,|\,\gminia\in\Irr(\Vtilde,\nablatilde)\}$
 holds, 
 and the multiplicity functions are the same.
\end{itemize}
In particular,
we obtain that Condition \ref{condition;18.11.26.3}
is satisfied for 
$\Ecat\psi_m^{-1}K$ around $P_m$.
By Proposition \ref{prop;16.7.25.40},
we obtain an isomorphism
$\DR^{\Ecat}(\Vtilde,\nablatilde)[-2]
\simeq
 \Ecat\psi_{m}^{-1}(K)_{|\nbigb_m}$,
i.e.,
$\nbigc(P_{m})$ holds
for $\Ecat\psi_m^{-1}(K)$.
But, it contradicts the construction of
the infinite sequence.
The case
$\nbigz_m=\emptyset$
can be argued similarly.

\vspace{.1in}

Suppose that $\vecY$ is convergent.
By Lemma \ref{lem;16.5.11.2},
after a finite step,
each blowing up is taken at a smooth point.
We may assume that
each blowing up is taken at a smooth point
from the beginning.
By the convergence,
there exists a complex curve $C$ in $X$
such that 
(i) $C$ is transversal with $H$,
(ii) each blowing up is taken at
the intersection of the exceptional fiber
and the strict transform of $C$.
We may assume that
$X=\Delta^2$,
$H=\{y=0\}$
and $C=\{x=0\}$.
Then, we obtain that the process will stop 
after a finite step
by Corollary \ref{cor;16.5.10.1}.
It contradicts the construction of
the infinite sequence.
Thus, we obtain Proposition \ref{prop;16.5.10.10}.
\hfill\qed

\subsection{Higher dimensional case}
\label{subsection;18.11.16.131}

Let us prove Theorem \ref{thm;18.11.12.10}.
Let $X$ be any $n$-dimensional complex manifold
with a complex hypersurface $H$.
The case $n=2$ has already been proved.
It is enough to consider the case where
$H$ is a normal crossing hypersurface.
Let $(V_0,\nabla)$ be a flat bundle on $X\setminus H$
with an isomorphism
$\DR(V_0,\nabla)\simeq K_{|X\setminus H}$.

\subsubsection{Smooth case}
\label{subsection;16.5.11.10}

Let us consider the case
where $X=\Delta^n$ and $H=\{z_1=0\}$.
According to Proposition \ref{prop;16.6.26.10},
there exist a closed subanalytic subset $Z\subset H$
with $\dim_{\real}Z\leq 2n-4$,
a good meromorphic flat bundle $(V_1,\nabla)$
on $(X\setminus Z,H\setminus Z)$,
and an isomorphism
 $\DR^{\Ecat}(V_1,\nabla)=K_{|X\setminus Z}$.

Let us observe that 
$(V_1,\nabla)$ extends to 
a meromorphic flat connection
on $(X,H)$.
Recall that
$V_1$ has the Deligne-Malgrange lattice
$V_1^{DM}\subset V_1$
which is 
a locally free $\nbigo_{X\setminus Z}$-submodule
of $V_1$
such that
$V_1^{DM}(\ast (H\setminus Z))=V_1$.
(See \cite{malgrange}. See also \cite{mochi9}.)

There exists a decomposition 
$Z=Z'\cup Z_{1}$,
where $Z'$ is smooth of $\dim_{\real} Z'=2n-4$,
and $\dim_{\real} Z_1\leq 2n-5$.
Let $P$ be any point of $Z'$.
There exists a small holomorphic coordinate system
$(X_P,w_1,\ldots,w_n)$ around $P$
such that 
(i) $X_P\simeq \Delta^n$
by the coordinate system,
(ii) $H\cap X_P=\{w_1=0\}$,
(iii) $\max_{Q\in Z'\cap X_P}\{|w_2(Q)|\}\leq\epsilon<1$.
For any $\vecb=(w^0_3,\ldots,w^0_n)$,
set $X_{P,\vecb}=\{(w_1,w_2,w^0_3,\ldots,w^0_n)\}$
and $H_{P,\vecb}:=X_{P,\vecb}\cap H$.
The restriction of
$V_1^{DM}$
to $X_{P,\vecb}\setminus Z'$
is the Deligne-Malgrange lattice of
$V_{1|X_{P,\vecb}\setminus Z'}$.
By the result in the surface case
(Proposition \ref{prop;16.5.10.10}),
$V_{1|X_{P,\vecb}\setminus Z'}$ extends
to a meromorphic flat connection
on $(X_{P,\vecb},H_{P,\vecb})$.
Hence, 
$V^{DM}_{1|X_{P,\vecb}\setminus Z'}$
extends to a coherent $\nbigo_{X_{P,\vecb}}$-module.
According to \cite[Theorem 7.4]{siu-extension},
$V^{DM}_{1|X_P\setminus Z'}$
is uniquely extended to 
a coherent reflexive $\nbigo_{X_P}$-module.
Hence, $V_1$ extends to a coherent reflexive 
$\nbigo_{X\setminus Z_1}(\ast (H\setminus Z_1))$-module
$V_2$,
i.e.,
$(V_1,\nabla)$
extends to a meromorphic flat connection
on $(X\setminus Z_1,H\setminus Z_1)$.
There exists the Deligne-Malgrange lattice $V_2^{DM}$
of $(V_2,\nabla)$,
which is the canonical reflexive
$\nbigo_{X\setminus Z_1}$-submodule of $V_2$
such that
$V_2^{DM}(\ast (H\setminus Z_1))=V_2$.

There exists the decomposition
$Z_1=Z_1'\cup Z_2$,
where $Z_1'$ is smooth,
and $\dim_{\real}Z_2\leq 2n-6$.
Let $P$ be any point of $Z'$.
There exists a small holomorphic coordinate neighbourhood
$(X_P,w_1,\ldots,w_n)$ around $P$
such that 
(i) the coordinate system induces an isomorphism
$X_P\simeq \Delta^n$,
(ii) $H\cap X_P=\{w_1=0\}$,
(iii) $\max_{Q\in Z_1'\cap X_P}\{|w_2(Q)|\}\leq\epsilon<1$.
Let $A_P:=\{(w_1,\ldots,w_n)\in X_P\,|\,w_1=w_2=0\}$.
Let $q:X_P\lrarr A_P$ be the projection.
Because $\dim_{\real} Z_1'=2n-5$,
there exists a closed subanalytic subset
$Z_1''\subset A_P$ such that 
$q^{-1}(Q)\cap (Z_1'\cap X_P)=\emptyset$
for any $Q\in A_P\setminus Z_1''$.
Hence, by using 
\cite[Theorem 7.4]{siu-extension},
we obtain that
$V^{DM}_{2|X_P\setminus Z_1'}$
is uniquely extended to 
a reflexive $\nbigo_{X_P\setminus Z_1'}$-module.
It implies that
$(V_2,\nabla)$ extends to
a meromorphic flat connection
$(V_3,\nabla)$
on $(X\setminus Z_2,H\setminus Z_2)$.

By an induction with a similar argument,
we obtain that $(V_3,\nabla)$ extends to 
a meromorphic flat connection
$(V,\nabla)$ on $(X,H)$.

\begin{lem}
For any holomorphic map
$\varphi:(\Delta^{\ast},0,\Delta)\lrarr 
 (X\setminus H,H,X)$,
$\Ecat\varphi^{-1}K$
and $\DR^{\Ecat}\varphi^{\ast}V[-1]$
are naturally isomorphic.
\end{lem}
\pf
If $\varphi(0)\in H\setminus Z$,
the claim holds by the construction of $(V,\nabla)$.
For a general $\varphi$,
there exists a holomorphic map 
$\psi:\Delta^2\lrarr X$
such that 
(i) $\psi(w_1,0)=\varphi(w_1)$,
(ii)  $\psi(0,w_2)\in H\setminus Z$
for general $w_2$.
By applying the result in the surface case
(Proposition \ref{prop;16.5.10.10}),
there exists a meromorphic flat bundle
$(V_{100},\nabla)$
on $(\Delta^2,\{w_1=0\})$
such that
$\DR^{\Ecat}(V_{100})\simeq
 \Ecat\psi^{-1}(K)$.
There exists the isomorphism
$V_{100}\simeq \psi^{\ast}V$
outside of 
$Z_{100}\subset \{w_1=0\}$
with $\dim Z_{100}=0$.
Hence, we obtain
$V_{100}\simeq\psi^{\ast}(V)$
on $\Delta^2$.
It implies the claim of the lemma.
\hfill\qed

\vspace{.1in}

Then, we obtain a global isomorphism
$\DR^{\Ecat}(V)[-n]\simeq K$
from Proposition \ref{prop;16.6.23.20}
in this case.

\subsubsection{Normal crossing case}

Let us consider the case 
where $H=\bigcup H_j$ is normal crossing.
Set $H^{[2]}=\bigcup_{i\neq j} (H_i\cap H_j)$.
We have already extended
$(V_0,\nabla)$ on $X\setminus H$
to $(V_1,\nabla)$ on $X\setminus H^{[2]}$.
The Deligne-Malgrange lattice of $V_1$
is a coherent reflexive $\nbigo$-module
on $X\setminus H^{[2]}$.
By using the result in the surface case
(Proposition \ref{prop;16.5.10.10}),
and by using the extension theorem of Siu
\cite[Theorem 7.4]{siu-extension},
we obtain that $(V_0,\nabla)$
extends to a meromorphic flat sheaf
$(V,\nabla)$
on $(X,H)$.
Then, the claim of Theorem \ref{thm;16.5.11.20}
follows from 
Proposition \ref{prop;16.6.23.20}
and the next lemma.

\begin{lem}
\label{lem;16.6.16.2}
For any holomorphic map
$\varphi:(\Delta^{\ast},0,\Delta)\lrarr 
 (X\setminus H,H,X)$,
$\Ecat\varphi^{-1}(K)$ and 
$\DR^{\Ecat}\varphi^{\ast}V[-1]$
are naturally isomorphic.
\end{lem}
\pf
There exists a projective birational morphism
$\rho:X'\lrarr X$
such that 
(i) $H'=\rho^{-1}(H)$ is normal crossing,
(ii) $X'\setminus H'\simeq X\setminus H$,
(iii) for the strict transform $\varphi':\Delta\lrarr X'$
 of $\varphi$,
 the point $\varphi'(0)$ is contained in the smooth part 
 of the exceptional divisor of $\rho$.
Applying the consideration in the beginning of this subsection
to $\Ecat\rho^{-1}(K)$,
we obtain a meromorphic flat bundle
$(V',\nabla)$ on $(X',H')$ as the extension of
$\rho^{\ast}(V_0,\nabla)$
corresponding to $\Ecat\rho^{-1}K$.
Then, we obtain
$\DR^{\Ecat}(\varphi')^{\ast}(V')[-1]
 \simeq
 \Ecat(\varphi')^{-1}\rho^{-1}(K)$
by the result in the smooth case
\S\ref{subsection;16.5.11.10}.
Because
$\rho_{\ast}V'$ is naturally isomorphic to $V$
as meromorphic flat connections,
we obtain
$V'\simeq\rho^{\ast}V$.
Then, the claim of Lemma \ref{lem;16.6.16.2} follows.
The proof of Theorem \ref{thm;16.5.11.20}
is also completed.
\hfill\qed

\section{Holonomic $\nbigd$-modules and enhanced ind-sheaves}

We prove the main results of this paper,
that is the curve test for enhanced ind-sheaves and 
holonomic $\nbigd$-modules.
In \S\ref{subsection;18.11.27.1},
we introduce the full subcategory
$\Ecat^b_{\sankaku}(\IC_X)\subset
 \Ecat^b_{\realc}(\IC_X)$
determined by the curve test,
and we state the main result (Theorem \ref{thm;16.7.11.1}),
i.e.,
the solution complex functor induces
an equivalence
$\Sol^{\Ecat}:\Dcat^b_{\hol}(\nbigd_X)
\lrarr
 \Ecat^b_{\sankaku}(\IC_X)$.
In \S\ref{subsection;16.10.9.3},
we introduce a functor $\Upsilon^{\Ecat}$
constructing $\nbigd$-modules
from enhanced ind-sheaves
by following the reconstruction formula
in \cite{DAgnolo-Kashiwara1},
and we observe that there exists 
a natural morphism
$K\lrarr \Sol^{\Ecat}\circ\Upsilon^{\Ecat}(K)$
for any  $K\in\Ecat^b_{\realc}(\IC_{X})$.
It is enough to prove that
for any $K\in\Ecat^b_{\sankaku}(\IC_X)$,
$\Upsilon^{\Ecat}(K)$
is contained in $\Dcat^b_{\hol}(\nbigd_X)$,
and the natural morphism
$K\lrarr \Sol^{\Ecat}\circ\Upsilon^{\Ecat}(K)$
is an isomorphism.
The claim is already known for 
$K\in \Ecat^b_{\circledcirc}(\IC_{\vecX(H)})$
for any complex hypersurface $H$.
Hence, it is enough to prove that
$K$ is expressed as a kind of gluing of
objects induced by
$\Ecat^b_{\circledcirc}(\IC_{\vecZ_i(H_i)})$
for some complex manifolds $Z_i$
with complex hypersurfaces $H_i$.
For that purpose,
we study the complex analyticity of 
the support and the singular locus of $K$
in \S\ref{subsection;18.11.27.10}--\S\ref{subsection;18.11.27.11}.
Then, the proof of Theorem \ref{thm;16.7.11.1}
is obtained by a formal argument
in \S\ref{subsection;18.11.27.12}.

\subsection{Statement}
\label{subsection;18.11.27.1}

Let $X$ be a complex manifold.
Let
$\Ecat^b_{\sankaku}(\IC_X)
 \subset
 \Ecat^b_{\realc}(\IC_X)$
denote the full subcategory
of objects $K$ with the following property.
\begin{itemize}
\item
 Let $\varphi:\Delta\lrarr X$ be any holomorphic map.
 Then, 
$\Ecat\varphi^{-1}(K)$ comes from
a cohomologically holonomic $\nbigd_{\Delta}$-complex,
i.e., 
there exist an object
$M\in \Dcat^b_{\hol}(\nbigd_{\Delta})$
and an isomorphism
 $\Ecat\varphi^{-1}(K)
\simeq
 \DR^{\Ecat}_{\Delta}(M)$.
\end{itemize}

The functor
$\Sol^{\Ecat}:
 \Dcat^b_{\hol}(\nbigd_X)^{\op}
\lrarr
 \Ecat^b_{\realc}(\IC_X)$
factors through
$\Ecat^b_{\sankaku}(\IC_X)$.
(See \cite[Definition 9.1.1]{DAgnolo-Kashiwara1}
for the enhanced solution complex functor
$\Sol^{\Ecat}$.)
It is fully faithful according to 
\cite[Theorem 4.9.12, Corollary 9.4.9, Theorem 9.5.3]{DAgnolo-Kashiwara1}.
The following theorem is one of the main results
in this paper, which will be proved in the rest of this paper.

\begin{thm}
\label{thm;16.7.11.1}
$\Sol^{\Ecat}:
 \Dcat^b_{\hol}(\nbigd_X)^{\op}
\lrarr
 \Ecat^b_{\sankaku}(\IC_X)$
is an equivalence.
As a consequence,
$\DR^{\Ecat}:\Dcat^b_{\hol}(\nbigd_X)
\lrarr
 \Ecat^b_{\sankaku}(\IC_X)$
is an equivalence.
\end{thm}

\subsection{Construction of quasi-inverse}
\label{subsection;16.10.9.3}

Let $K\in \Ecat^b(\IC_X)$.
We define
$\Upsilon^{\Ecat}(K):=
 \nhom^{\Ecat}(K,\nbigo^{\Ecat})$
 in $\Dcat^b(\nbigd_X)$
by following the reconstruction formula
\cite[Proposition 9.5.1]{DAgnolo-Kashiwara1}.
Let us observe that there exists a natural morphism
$\Phi_K:K\lrarr
 \Sol^{\Ecat}(\Upsilon^{\Ecat}(K))$.
Let $\pibar:X\times\realbar\lrarr X$
denote the projection.
Let $j_X:X\times\real_{\infty}\lrarr X\times\realbar$
denote the natural morphism of bordered spaces.
We set
$\Ltilde^{\Ecat}(G):=
 Rj_{X!!}L^{\Ecat}G$
and 
$\Rtilde^{\Ecat}(G):=
 Rj_{X\ast}R^{\Ecat}G$
in $\Dcat^b(\IC_{X\times\realbar})$.
According to \cite[Definition 4.5.13]{DAgnolo-Kashiwara1},
there exists the following natural isomorphism:
\[
 \Upsilon^{\Ecat}(K)
 \simeq
 R\pibar_{\ast}
 \nrhom_{\IC_{X\times\realbar}}
 \bigl(\Ltilde^{\Ecat}K,\Rtilde^{\Ecat}\nbigo^{\Ecat}\bigr).
\]
There exists the following natural morphism:
\begin{multline}
\label{eq;16.7.29.2}
 R\nihom_{\beta\pibar^{-1}\nbigd_X}
 \bigl(
 \beta R\nhom_{\IC_{X\times\realbar}}
 (\Ltilde^{\Ecat}K,\Rtilde^{\Ecat}\nbigo^{\Ecat}),
 \Rtilde^{\Ecat}\nbigo^{\Ecat}
 \bigr)
\lrarr 
 \\
 R\nihom_{\beta\pibar^{-1}\nbigd_X}
 \bigl(
 \beta\pibar^{-1}
 R\pibar_{\ast}R\nhom_{\IC_{X\times\realbar}}
 (\Ltilde^{\Ecat}K,\Rtilde^{\Ecat}\nbigo^{\Ecat}),
 \Rtilde^{\Ecat}\nbigo^{\Ecat}
 \bigr).
\end{multline}
Let $Q:\Dcat^b(\IC_{X\times\realbar})
\lrarr \Ecat^b(\IC_X)$
denote the natural quotient functor.
We obtain the following natural isomorphism:
\begin{multline}
\label{eq;16.7.29.1}
 \Sol^{\Ecat}(\Upsilon^{\Ecat}(K))
\simeq
 R\nihom_{\beta\pi^{-1}\nbigd_X}
 \bigl(
 \beta\pi^{-1}
 R\pibar_{\ast}R\nhom_{\IC_{X\times\realbar}}
 (\Ltilde^{\Ecat}K,\Rtilde^{\Ecat}\nbigo^{\Ecat}),
 \nbigo^{\Ecat}
 \bigr)
 \\
\simeq
 QR\nihom_{\beta\pibar^{-1}\nbigd_X}
 \bigl(
 \beta\pibar^{-1}
 R\pibar_{\ast}R\nhom_{\IC_{X\times\realbar}}
 (\Ltilde^{\Ecat}K,\Rtilde^{\Ecat}\nbigo^{\Ecat}),
 \Rtilde^{\Ecat}\nbigo^{\Ecat}
 \bigr).
\end{multline}
By (\ref{eq;16.7.29.2}) and (\ref{eq;16.7.29.1}),
we obtain the following natural morphism:
\begin{equation}
\label{eq;16.7.8.1}
 QR\nihom_{\beta\pibar^{-1}\nbigd_X}
 \bigl(
 \beta R\nhom_{\IC_{X\times\realbar}}
 (\Ltilde^{\Ecat}K,\Rtilde^{\Ecat}\nbigo^{\Ecat}),
 \Rtilde^{\Ecat}\nbigo^{\Ecat}
 \bigr)
\lrarr
  \Sol^{\Ecat}(\Upsilon^{\Ecat}(K)).
\end{equation}

According to \cite[Theorem 5.4.23]{Kashiwara-Schapira-ind-sheaves},
there exists the following natural isomorphism:
\begin{multline}
 \label{eq;18.1.3.1}
R\nhom_{\IC_{X\times\realbar}}
 \Bigl(
 \Ltilde^{\Ecat}K,
 R\nihom_{\beta\pibar^{-1}\nbigd_X}
 \bigl(\beta R\nhom_{\IC_{X\times\realbar}}
 (\Ltilde^{\Ecat}K,\Rtilde^{\Ecat}\nbigo^{\Ecat}),
 \Rtilde^{\Ecat}\nbigo^{\Ecat}
 \bigr)
 \Bigr)
\simeq \\
 R\nhom_{\mathrm{I}(\beta\pibar^{-1}\nbigd_X)}
 \Bigl(
 \Ltilde^{\Ecat}K
\otimes_{\IC_{X\times\realbar}}
\beta R\nhom_{\IC_{X\times\realbar}}
 (\Ltilde^{\Ecat}K,\Rtilde^{\Ecat}\nbigo^{\Ecat}),
 \Rtilde^{\Ecat}\nbigo^{\Ecat}
 \Bigr).
\end{multline}
There exist the following natural isomorphisms:
\begin{multline}
\label{eq;16.7.8.2}
  R\nhom_{\mathrm{I}(\beta\pibar^{-1}\nbigd_X)}
 \Bigl(
 \Ltilde^{\Ecat}K
\otimes_{\IC_{X\times\realbar}}
\beta R\nhom_{\IC_{X\times\realbar}}
 (\Ltilde^{\Ecat}K,\Rtilde^{\Ecat}\nbigo^{\Ecat}),
 \Rtilde^{\Ecat}\nbigo^{\Ecat}
 \Bigr)\simeq
 \\
 R\nhom_{\mathrm{I}(\beta\pibar^{-1}\nbigd_X)}
 \Bigl(
\beta R\nhom_{\IC_{X\times\realbar}}
 (\Ltilde^{\Ecat}K,\Rtilde^{\Ecat}\nbigo^{\Ecat}),
 R\nihom_{\IC_{X\times\realbar}}
 \bigl(
 \Ltilde^{\Ecat}K,
 \Rtilde^{\Ecat}\nbigo^{\Ecat}
 \bigr)
 \Bigr)
  \simeq
 \\
 R\nhom_{\pibar^{-1}\nbigd_X}
 \Bigl(
 R\nhom_{\IC_{X\times\realbar}}
 (\Ltilde^{\Ecat}K,\Rtilde^{\Ecat}\nbigo^{\Ecat}),
 R\nhom_{\IC_{X\times\realbar}}
 \bigl(
 \Ltilde^{\Ecat}K,
 \Rtilde^{\Ecat}\nbigo^{\Ecat}
 \bigr)
 \Bigr).
\end{multline}
Here, 
we obtain the first isomorphism 
from more general isomorphisms (\ref{eq;16.10.9.2}) below
and \cite[Proposition 5.4.11]{Kashiwara-Schapira-ind-sheaves},
and we obtain the second isomorphism
by \cite[Theorem 5.6.2 (ii)]{Kashiwara-Schapira-ind-sheaves}.
By (\ref{eq;16.7.8.1}), (\ref{eq;18.1.3.1}) and (\ref{eq;16.7.8.2}),
we obtain the morphism
$\Phi_K:
 K\lrarr
 \Sol^{\Ecat}(\Upsilon^{\Ecat}(K))$
in $\Ecat^b(\IC_X)$
corresponding to the identity:
\[
 \id\in
 R\nhom_{\pibar^{-1}\nbigd_X}
 \Bigl(
 R\nhom_{\IC_{X\times\realbar}}
 (\Ltilde^{\Ecat}K,\Rtilde^{\Ecat}\nbigo^{\Ecat}),
 R\nhom_{\IC_{X\times\realbar}}
 \bigl(
 \Ltilde^{\Ecat}K,
 \Rtilde^{\Ecat}\nbigo^{\Ecat}
 \bigr)
 \Bigr).
\]

\begin{lem}
Let $\psi:K_1\lrarr K_2$ be a morphism
in $\Ecat^b(\IC_X)$.
Then, the following induced diagram is commutative:
\[
 \begin{CD}
 K_1 @>{a_1}>> 
 \Sol^{\Ecat}(\Upsilon^{\Ecat}(K_1))\\
 @V{\psi}VV @V{\psi_{\ast}}VV \\
 K_2 @>{a_2}>> 
 \Sol^{\Ecat}(\Upsilon^{\Ecat}(K_2)).
 \end{CD}
\]
Here, the right vertical arrow is the morphism
induced by $\psi$.
\end{lem}
\pf
It is easy to check that
both $\psi_{\ast}\circ a_1$
and $a_2\circ \psi$
correspond to 
\[
 \psi^{\ast}
\in
 R\nhom_{\pibar^{-1}\nbigd_X}
 \Bigl(
 R\nhom_{\IC_{X\times\realbar}}
 (\Ltilde^{\Ecat}K_2,\Rtilde^{\Ecat}\nbigo^{\Ecat}),
 R\nhom_{\IC_{X\times\realbar}}
 \bigl(
 \Ltilde^{\Ecat}K_1,
 \Rtilde^{\Ecat}\nbigo^{\Ecat}
 \bigr)
 \Bigr).
\]
\hfill\qed

\subsubsection{Appendix}

Let $Y$ be a good topological space.
Let $\nbiga$ be a ring in $\IC_Y$.
There exists the following natural isomorphism
for any
$\nbign,\nbigm\in \IA$
and $M\in \IC_Y$:
\begin{equation}
\label{eq;16.10.9.1}
 \nihom_{\nbiga}\bigl(\nbign\otimes_{\IC_Y}M,\nbigm\bigr)
\simeq
 \nihom_{\nbiga}\bigl(\nbign,\nihom_{\IC_Y}(M,\nbigm)\bigr).
\end{equation}
Indeed, 
in the case where $\nbiga=\cnum_Y$,
it follows from \cite[Corollary 4.2.9]{Kashiwara-Schapira-ind-sheaves}.
Then, we can easily obtain (\ref{eq;16.10.9.1}) for general $\nbiga$
from (\ref{eq;16.10.9.1}) for $\nbiga=\cnum_Y$
and the definition of
$\nihom_{\nbiga}(\cdot,\cdot)$
\cite[Definition 5.4.9]{Kashiwara-Schapira-ind-sheaves}.

There exists the following natural isomorphism
in $D^+(\IC_Y)$
for any 
$\nbign\in D^-(\IA)$,
$\nbigm\in D^+(\IA)$
and $M\in D^-(\IC_Y)$:
\begin{equation}
\label{eq;16.10.9.2}
 R\nihom_{\nbiga}\bigl(
 \nbign\otimes_{\IC_Y} M,\nbigm
 \bigr)
\simeq
 R\nihom_{\nbiga}\bigl(
 \nbign,
 R\nihom_{\IC_Y}(M,\nbigm)
 \bigr). 
\end{equation}
It follows from the construction of
the derived functor
$R\nihom_{\nbiga}(\cdot,\cdot)$
in \cite[\S5.1, \S5.4]{Kashiwara-Schapira-ind-sheaves}.

\subsection{Preliminary}
\label{subsection;18.11.27.10}

\subsubsection{$\real$-linear subspaces}

Let $V$ be a finite dimensional $\cnum$-vector space.
Let $H$ be an $\real$-subspace of $V$.
For any real number $a$,
we set
$[a]_+:=\min\{\ell\in\seisuu\,|\,\ell\geq a\}$
and 
$[a]_-:=\max\{\ell\in\seisuu\,|\,\ell\leq a\}$.

\begin{lem}
\label{lem;16.8.10.1}
There exists a $\cnum$-subspace $L$ of $V$
such that $L\cap H=0$
and $\dim_{\cnum}L=\dim_{\cnum}V-[\dim_{\real} H/2]_+$.
\end{lem}
\pf
We set $H':=H\cap\sqrt{-1}H$
which is a $\cnum$-subspace of $V$.
There exists an $\real$-subspace $H''\subset H$
such that $H=H'\oplus H''$.
Note
$[\dim_{\real} H/2]_+
=\dim_{\cnum}H'+[\dim_{\real}H''/2]_+$.
The $\cnum$-subspace
$H+\sqrt{-1}H$
is equal to
$H'\oplus(H''\oplus \sqrt{-1}H'')$.
Hence,
we obtain
$\dim_{\cnum}(H+\sqrt{-1}H)
=\dim_{\cnum}H'
+\dim_{\real}H''$.

For any $\ell\in\seisuu_{\geq 1}$,
we can easily construct a $\cnum$-subspace 
$L_1$ of $\cnum^{\ell}$
such that 
$L_1\cap\real^{\ell}=0$
and $\dim_{\cnum}L_1=[\ell/2]_-$.
Hence, we can easily construct a $\cnum$-subspace
$L\subset V$
such that $V\cap H=0$
and 
\[
 \dim_{\cnum}L
=\dim_{\cnum}V
-\dim_{\cnum}(H+\sqrt{-1}H)
+[\dim_{\real}H''/2]_-
=\dim_{\cnum}V
-[\dim_{\real}H/2]_+.
\]
Thus, we obtain the claim of the lemma.
\hfill\qed

\begin{lem}
\label{lem;16.8.12.1}
Suppose that $H$ is not a $\cnum$-subspace of $V$,
and that $\dim_{\real}H=2\dim_{\cnum}V-2$.
Then, for any $u\in V/H$,
there exists $v\in H$ such that
$\sqrt{-1}v$ is mapped to $u$
via the projection $V\lrarr V/H$.
\end{lem}
\pf
Note that $H+\sqrt{-1}H$ is a $\cnum$-subspace of $V$.
Because $H$ is not a $\cnum$-subspace of $V$,
$\dim_{\cnum}(H+\sqrt{-1}H)>\dim_{\cnum} V-1$,
and hence $H+\sqrt{-1}H=V$.
Set $H':=H\cap\sqrt{-1}H$.
Because $H+\sqrt{-1}H=V$,
we obtain
$\dim_{\real}H'=2\dim_{\cnum}V-4$,
and hence
$\dim_{\cnum}H'=\dim_{\cnum}V-2$.

Set $V_1:=V/H'$
and $H_1:=H/H'$.
Because
$\dim_{\cnum}H'=\dim_{\cnum}V-2$,
we obtain
$\dim_{\cnum}V_1=\dim_{\real}H_1=2$.
It implies that
$V_1=H_1\oplus\sqrt{-1}H_1$.
Then, the claim of the lemma is clear.
\hfill\qed

\vspace{.1in}
We set
$k:=\dim_{\cnum}V-[\dim_{\real}H/2]_+$.
Let $\Gr_{\cnum}(V,k)$ denote the Grassmannian variety
of $k$-dimensional $\cnum$-subspaces of $V$.
Note
$\dim_{\real}\Gr_{\cnum}(V,k)=2k(\dim_{\cnum}V-k)$.
Let $Z(V,H,k):=\{L\in \Gr_{\cnum}(V,k)\,|\,L\cap H\neq 0\}$,
which is a closed real analytic subset of $\Gr_{\cnum}(V,k)$.

\begin{lem}
$\dim_{\real}Z(V,H,k)<\dim_{\real}\Gr_{\cnum}(V,k)$
holds.
\end{lem}
\pf
Let $\proj_{\real}(H)$ denote the projective space of
the real $1$-dimensional subspaces of $H$.
Let us consider the following subspace:
\[
 \Ztilde:=\Bigl\{
 (L_1,L)\in\proj_{\real}(H)\times\Gr_{\cnum}(V,k)\,\Big|\,
 L_1\subset L
 \Bigr\}.
\]
The projection
$\Ztilde\lrarr \Gr_{\cnum}(V,k)$ induces a surjection
$\Ztilde\lrarr Z(V,H,k)$.
It is enough to prove 
$\dim_{\real}\Ztilde<\dim_{\real}\Gr_{\cnum}(V,k)$.

By the projection
$\Ztilde\lrarr \proj_{\real}(H)$,
$\Ztilde$ is a fiber bundle over $\proj_{\real}(H)$.
The fiber over $L_1$ is isomorphic to
$\Gr_{\cnum}(V/L_1\otimes\cnum,k-1)$.
Hence, we obtain
\[
 \dim_{\real} \Ztilde=
 2(k-1)(\dim_{\cnum} V-k)
+\dim_{\real}H-1
=\dim_{\real}\Gr_{\cnum}(V,k)
-2(\dim_{\cnum}V-k)+\dim_{\real}H-1.
\]
By our choice of $k$,
we obtain
$2k+\dim_{\real}H\leq 2\dim_{\cnum}V$.
Hence, we obtain
$\dim_{\real}\Ztilde<\dim_{\real}\Gr_{\cnum}(V,k)$.
\hfill\qed

\begin{cor}
\label{cor;16.9.1.1}
Let $H_1,\ldots,H_N$ be $\ell$-dimensional $\real$-subspaces of $V$.
\begin{itemize}
 \item 
There exists a $\cnum$-subspace $L$ of $V$
such that 
(i) $\dim_{\cnum}L=\dim_{\cnum}V-[\ell/2]_+$,
(ii) $L\cap H_i=0$ $(i=1,\ldots,N)$.
 \item
Suppose that $\ell=2m-1$ for a positive integer $m$.
Let $\nbigu(V,H_1,\ldots,H_N,m)\subset\Hom(V,\cnum^m)$ be 
the set of $\cnum$-linear maps
$f:V\lrarr \cnum^m$
such that $f_{|H_i}$ are injective.
Then, $\nbigu(V,H_i,m)$ is 
a non-empty open subset.
\end{itemize}
\end{cor}
\pf
Let us prove the first claim.
Set $k:=\dim_{\cnum}V-[\ell/2]_+$.
We obtain
$\Gr_{\cnum}(V,H,k)\setminus\bigcup_iZ(V,H_i,k)\neq\emptyset$
because $\dim_{\real}Z(V,H_i,k)<\dim_{\real}\Gr_{\cnum}(V,H,k)$.
Then,
any $L\in \Gr_{\cnum}(V,H,k)\setminus\bigcup_iZ(V,H_i,k)$
has the desired property.

Let us prove the second claim.
It is enough to check that $\nbigu(V,H_1,\ldots,H_N,m)$
is non-empty.
By the first claim,
there exists a $\cnum$-subspace $L\subset V$
such that $L\cap H_i=0$ $(i=1,\ldots,N)$
and $\dim_{\cnum}L=\dim_{\cnum} V-m$.
Then, we have only to compose
the projection
$V\lrarr V/L$ and a $\cnum$-isomorphism
$V/L\simeq\cnum^m$.
\hfill\qed

\subsubsection{Subanalytic subsets
which are generically complex analytic}

Let $X$ be a complex manifold.
Let $Z\subset X$ be 
a purely $k$-dimensional closed subanalytic subset.
Let $Z_0\subset Z$ be a closed subanalytic subset
with $\dim_{\real} Z_0<k$.
Suppose that $Z\setminus Z_0$ is a complex submanifold of $X$.
In particular, $k$ is even.
We obtain the following proposition
from a generalization
of the theorem of Remmert-Stein \cite{Remmert-Stein}
due to Shiffman \cite{Shiffman}.
\begin{thm}
\label{thm;16.9.2.1}
If $\dim_{\real} Z_0<k-1$,
then $Z$ is a complex analytic subvariety of $X$.
\hfill\qed
\end{thm}

Let us study a description of $Z$
around a general point of $Z_0$
in the case $\dim_{\real} Z_0=k-1$.

\begin{prop}
\label{prop;16.9.1.10}
Suppose $\dim_{\real}Z_0=k-1$.
Then, there exists a closed subanalytic subset
$Z_1\subset Z_0$ with $\dim_{\real}Z_1\leq k-2$
such that the following holds.
\begin{itemize}
\item
$Z_1$ contains the singular locus of $Z_0$.
\item
 For any $P\in Z_0\setminus Z_1$,
 there exists a neighbourhood $X_P$ of $P$ in $X$
 with a real analytic function $f_P:X_P\lrarr\real$
 and a closed $(k/2)$-dimensional
 complex submanifold $\Ztilde_P\subset X_P$
 such that
 (i) $X_P\cap Z_0=\Ztilde_P\cap f^{-1}(0)$,
 (ii) $df_{|\Ztilde_P}$ is nowhere vanishing,
(iii) $X_P\cap Z$ is equal to
 $\Ztilde_P$ or $\Ztilde_P\cap f^{-1}(\real_{\geq 0})$.
\end{itemize}
\end{prop}
\pf
We may assume that $X$ is 
an open subset of $\cnum^n$.
Let $Q$ be any point of $Z_0$.
Let $X_Q$ be any relatively compact neighbourhood 
of $Q$ in $X$.
Set $Z_Q:=Z\cap X_Q$
and $Z_{0,Q}:=Z_0\cap X_Q$.
Let $Z_{0,Q}^{\sm}$
denote the set of the smooth points of $Z_{0,Q}$.
It is decomposed into the union of
the connected components
$\nbigc_i$ $(i=1,\ldots,\ell_1(Q))$.
We decompose $Z_{Q}\setminus Z_{0,Q}$
into the connected components
$\nbigd_j$ $(j=1,\ldots,\ell_2(Q))$.
We choose any points
$R(\nbigc_i)\in \nbigc_i$
$(i=1,\ldots,\ell_1(Q))$
and
$R(\nbigd_j)\in \nbigd_j$
$(j=1,\ldots,\ell_2(Q))$.

By Corollary \ref{cor;16.9.1.1},
there exists a $\cnum$-linear map
$\phi:\cnum^n\lrarr \cnum^{k/2}$
such that 
(i) the restriction of 
$\phi$ to $T_{R(\nbigc_i)}\nbigc_i$
are injective for any $i$,
(ii) the restriction of $\phi$
to $T_{R(\nbigd_j)}\nbigd_j$ 
are isomorphisms for any $j$.
Hence, there exists a closed subanalytic subset
$\gbigw_{1}\subset Z_{0,Q}$
with $\dim_{\real}\gbigw_{1}\leq k-2$,
and $\gbigw_{2}\subset Z_Q$
with $\dim_{\real}\gbigw_{2}\leq k-1$,
such that
(i) the restriction of $\phi$
to $\nbigc_i\setminus \gbigw_{1}$
are  immersions for any $i$,
(ii) the restriction of $\phi$
to $\nbigd_j\setminus \gbigw_{2}$
 are immersions for any $j$.
We obtain the $k$-dimensional subanalytic subset
$\phi(Z_{Q})\subset \cnum^{k/2}\simeq\real^{k}$.
There exists a subanalytic closed subset
$\gbigw_{3}\subset\phi(Z_{Q})$
such that 
(i) $\dim_{\real}\gbigw_{3}\leq k-1$,
(ii) each connected component of
 $\phi(Z_Q)\setminus \gbigw_{3}$
 is simply connected,
(iii) $\gbigw_{3}$ contains
 $\del\phi(Z_Q)$
 and $\phi(Z_{0,Q})$,
(iv) the restriction of $\phi$ to
 $Z_{Q}\setminus \phi^{-1}(\gbigw_{3})
 \lrarr \cnum^{k/2}$
 is a local diffeomorphism.

There exist a $\cnum$-vector space $H$ with
$\dim_{\cnum} H=n-k/2$
and a $\cnum$-isomorphism
$\cnum^n\simeq
 \cnum^{k/2}\oplus H$
such that $\phi$ is 
identified with the projection onto the first component.
On each connected component
$\nbign$ of $\phi(Z_Q)\setminus \gbigw_{3}$,
there exist $H$-valued
subanalytic functions $h_{\nbign,p}$
$(p\in\Lambda(i,\nbign))$
on $(\nbign,\cnum^{k/2})$
such that
$\phi^{-1}(\nbign)\cap Z_Q$
is expressed as the union
of the graph $\Gamma(h_{\nbign,p})$.
By Lemma \ref{lem;15.12.28.10},
there exists a closed subanalytic subset 
$\gbigw_{4}\subset\gbigw_{3}$ 
such that
(i) $\dim_{\real}\gbigw_{4}\leq k-2$,
(ii) for each connected component $\nbign$
 of $\phi(Z_Q)\setminus \gbigw_{3}$,
 the singular locus of $\del\nbign$ 
 is contained in $\gbigw_{4}$,
 and $h_{\nbign,p}$ are ramified real analytic around any point of
 $\del\nbign\setminus \gbigw_{4}$.
We remark the following.
\begin{lem}
\label{lem;18.1.12.1}
Set $I_x:=\bigl\{x\in\real\,\big|\,|x|<\epsilon\bigr\}$
and $I_y:=\bigl\{y\in\real\,\big|\,|y|<\epsilon\bigr\}$.
Let $U_{\vecw}$ be a neighbourhood of 
$(0,\ldots,0)$ in $\cnum^m$.
There exists the embedding
$I_x\times I_y\times U_{\vecw}
\lrarr \cnum\times\cnum^m$
by $(x,y,\vecw)\longmapsto
 (x+\sqrt{-1}y,\vecw)$
with which we regard $I_x\times I_y\times U_{\vecw}$
as a complex manifold.
Let $f:I_x\times U_{\vecw}\lrarr I_y$
be a real analytic function,
and we put $A:=\bigl\{(x,y,\vecw)\,\big|\,
 y\geq f(x,\vecw)
 \bigr\}$.
Let $g$ be a $\cnum$-valued
ramified real analytic continuous function on $A$.
If $g_{|A\setminus \del A}$ is holomorphic,
there exists a holomorphic function $\gtilde$
on a neighbourhood of $A$ such that
$\gtilde_{|A}=g$.
\end{lem}
\pf
Let us consider the case $m=0$.
There exist a positive integer $\rho>0$
and the expansion
$g=\sum_{j=0}^{\infty} a_j(x)(y-f(x))^{j/\rho}$
on a neighbourhood of any point of $\del A$,
where $a_j(x)$ are $\cnum$-valued real analytic functions of $x$.
Let us observe that $a_j=0$ unless $j\in\rho\seisuu$.
Suppose 
$S:=\{j\in \seisuu\setminus\rho\seisuu\,|\,a_j\neq 0\}$
is not empty,
and we shall derive a contradiction.
Let $j_0$ be the minimum of $S$.
Both $\del_xg$ and $\del_yg$ are ramified real analytic
functions on $A\setminus\del A$,
with the expansions:
\[
 \del_xg=\sum_{j=0}^{\infty}\del_xa_j(x)(y-f(x))^{j/\rho}
-\sum_{j=0}^{\infty}a_j(x)j\rho^{-1}(y-f(x))^{(j/\rho)-1}\del_xf,
\]
\[
 \del_yg=\sum_{j=0}^{\infty}a_j(x)j\rho^{-1}(y-f(x))^{(j/\rho)-1}.
\]
Note $\del_xg-\sqrt{-1}\del_yg=0$.
By considering the coefficient of
$(y-f(x))^{(j_0/\rho)-1}$,
we obtain that $a_{j_0}(x)=0$,
which contradicts the choice of $j_0$.
Hence, we obtain that $S=\emptyset$.
It means that $g$ is the restriction of
a $\cnum$-valued real analytic function $\gtilde$
defined on a neighbourhood of $A$.
We can easily see that $\gtilde$ is holomorphic,
and  we are done in the case $m=0$.

Let us consider the general case.
By considering the restriction to
$I_x\times I_y\times\{\vecw_0\}$ for each $\vecw_0\in U_{\vecw}$,
we obtain the expression
$g=\sum_{j=0}^{\infty} a_j(x,\vecw)(y-f(x))^j$
on a neighbourhood of any point of $\del A$.
Hence, $g$ is the restriction of 
a real analytic function $\gtilde$ defined
on a neighbourhood of $A$.
We can easily see that $\gtilde$ is holomorphic,
and  we obtain the claim of Lemma \ref{lem;18.1.12.1}.
\hfill\qed

\vspace{.1in}

 Because $h_{\nbign,p}$ are holomorphic on $\nbign$,
 there exists a neighbourhood $\nbigm_R$
 for each point $R\in\del\nbign\setminus \gbigw_{4}$
 such that $h_{\nbign,p|\nbigm_R\cap\nbign}$
extend to holomorphic functions
 on $\nbigm_R$.

\begin{lem}
\label{lem;17.1.4.1}
$\dim_{\real}\bigl(
 \phi^{-1}(\gbigw_{4})
 \cap Z_{0,Q}
 \bigr)\leq k-2$
 holds.
Similarly,
$\dim_{\real}(\phi^{-1}(\gbigw_{3})\cap Z_Q)\leq k-1$
holds.
\end{lem}
\pf
Note
$\dim_{\real}
 \bigl(
 \phi^{-1}(\gbigw_{4})\cap
 \gbigw_1\bigr)
\leq
 \dim_{\real}\gbigw_1\leq k-2$.
Because 
the restriction of $\phi$
to $Z_{0,Q}\setminus \gbigw_1$
is an immersion,
$\dim_{\real}\bigl(
 \phi^{-1}(\gbigw_{4})\cap Z_{0,Q}
\setminus \gbigw_1
 \bigr)
\leq 
\dim_{\real}\gbigw_{4}\leq k-2$
holds.
Thus, we obtain the first claim of the lemma.
The second claim can be obtained similarly.
\hfill\qed

\vspace{.1in}
Let $\Gammabar(h_{\nbign,p})$
denote the closure of $\Gamma(h_{\nbign,p})$ in $Z_Q$.
We set
\[
 E((\nbign,p),(\nbign',p')):=
 \Gammabar(h_{\nbign,p})
\cap
 \Gammabar(h_{\nbign',p'}).
\]
Let $E((\nbign,p),(\nbign',p'))_{q}^{\sm}$
denote the set of the $q$-dimensional smooth points
of $E((\nbign,p),(\nbign',p'))$.
\begin{lem}
Let $P$ be any point of 
$E((\nbign,p),(\nbign',p'))_{k-1}^{\sm}
\cap
 Z_{0,Q}^{\sm}
 \setminus 
 \bigl(\gbigw_1\cup
 \phi^{-1}(\gbigw_{4})
 \bigr)$.
Suppose that $(\nbign,p)\neq(\nbign',p')$.
Then, the intersection of 
$\Gammabar(h_{\nbign,p})
\cup
 \Gammabar(h_{\nbign',p'})$
and a neighbourhood of $P$
is a complex submanifold.
\end{lem}
\pf
Let $X_P$ be a small neighbourhood of $P$ in $X$.
By construction,
there exist a complex submanifold
$Y(\nbign,p)$ and a real analytic map
$f:X_P\lrarr\real$
such that 
$\Gammabar(h_{\nbign,p})
=Y(\nbign,p)\cap f^{-1}(\real_{\geq 0})$.
There also exist a complex submanifold
$Y(\nbign',p')$ and a real analytic map
$f':X_P\lrarr\real$
such that 
$\Gammabar(h_{\nbign',p'})
=Y(\nbign',p')\cap (f')^{-1}(\real_{\geq 0})$.
Because 
$\dim_{\real}(Y(\nbign,p)\cap Y(\nbign',p'))
 \geq k-1$,
we obtain
$Y(\nbign,p)=Y(\nbign',p')$.
Then, the claim of the lemma is clear.
\hfill\qed

\vspace{.1in}
Let 
$\Sing
 \overline{
 E((\nbign,p),(\nbign',p'))_{k-1}^{\sm}}$
denote the singular locus of
$\overline{
 E((\nbign,p),(\nbign',p'))_{k-1}^{\sm}}$.
We set
\[
 \gbigw_5:=
 \bigcup_{(\nbign,p)\neq(\nbign',p')}
\Bigl(
 \bigcup_{q\leq k-2}
 \overline{
 E((\nbign,p),(\nbign',p'))_{q}^{\sm}}
\cup
 \Sing
 \overline{
 E((\nbign,p),(\nbign',p'))_{k-1}^{\sm}}
\Bigr).
\]
By construction, $\dim \gbigw_5\leq k-2$ holds.
Let $\gbigw_{6}$ denote the closure of
$\phi^{-1}(\gbigw_{3})\cap(Z_Q\setminus Z_{0,Q})$
in $Z_Q$.
Note $\dim_{\real}(\gbigw_6\cap Z_{0,Q})\leq k-2$.
Let $\gbigw_7$ denote the union of 
the singular locus of $Z_{0,Q}$
and
$\gbigw_1
\cup
 (\phi^{-1}(\gbigw_{4})\cap Z_{0,Q})
\cup
 (\gbigw_{6}\cap Z_{0,Q})
\cup \gbigw_5$.
Note $\dim_{\real} \gbigw_7\leq k-2$.

Let $P$ be any point of $Z_{0,Q}\setminus \gbigw_7$,
and let $X_P$ be a small neighbourhood
of $P$ in $X$.
There exists at least one $(\nbign,p)$ such that
$P\in \Gammabar(h_{\nbign,p})$.
If there exists another 
$(\nbign',p')$ such that
$P\in\Gammabar(h_{\nbign',p'})$,
then 
$\Gammabar(h_{\nbign,p})
\cup
\Gammabar(h_{\nbign',p'})$
is a complex submanifold.
Thus, we obtain
the claim of Proposition \ref{prop;16.9.1.10}.
\hfill\qed

\vspace{.1in}
Let $A\subset X$ be a closed complex subvariety.
Let $B\subset X$ be a closed subset
such that
(i) $B$ is a purely $2q$-dimensional subanalytic subset of $X$,
(ii) $B$ is equal to the closure of $B\setminus A$ in $X$,
(iii) $B\setminus A$ is a complex subvariety 
of $X\setminus A$.

\begin{cor}
\label{cor;16.9.2.10}
Under the assumption,
$B$ is a complex subvariety of $X$.
\end{cor}
\pf
Let $\Sing(B)$ denote the singular locus of $B$.
It is enough to prove that
$\dim_{\real}\Sing(B)<2q-1$
by Theorem \ref{thm;16.9.2.1}.
By the assumption,
$\dim_{\real}(\Sing(B)\setminus A)<2q-1$
holds.
Hence, it is enough to prove
$\dim_{\real}(\Sing(B)\cap A)<2q-1$.

We assume that 
$\dim_{\real}(\Sing(B)\cap A)=2q-1$,
and we shall derive a contradiction.
There exists a closed subanalytic subset
$\gbigw\subset \Sing(B)\cap A$
with $\dim_{\real}\gbigw<2q-1$
such that the following holds
for any $P\in (\Sing(B)\cap A)\setminus \gbigw$.
\begin{itemize}
\item
$P$ is a smooth point of $\Sing(B)$.
\item
Let $X_P$ be any small neighbourhood 
of $P$ in $X$.
Then, there exist a closed complex submanifold
$\Btilde_P$ of $X_P$
and a real analytic function $f_P:\Btilde_P\lrarr \real$
whose exterior derivative is nowhere vanishing,
such that
$B\cap X_P$ is $\Btilde_P\cap f_P^{-1}(\real_{\geq 0})$.
\end{itemize}
Because $\Btilde_P\not\subset A$,
we obtain
$\dim_{\real}(\Btilde_P\cap A)
 \leq \dim_{\real}\Btilde_P-2=2q-2$
which contradicts
$\dim_{\real}(\Sing(B)\cap A)=2q-1$.
Thus, we obtain Corollary \ref{cor;16.9.2.10}.
\hfill\qed

\subsubsection{$\real$-constructible sheaves}

Let $X$ be any complex manifold.
Let
$\Dcat^b_{\sankaku}(\cnum_X)
 \subset
 \Dcat^b_{\realc}(\cnum_X)$
denote the full subcategory
of objects $K$ with the following property.
\begin{itemize}
\item
 Let $\varphi:\Delta\lrarr X$ be any holomorphic map.
 Then,
 $\varphi^{-1}(K)$  is cohomologically $\cnum$-constructible.
\end{itemize}

\begin{lem}
\label{lem;16.8.12.20}
$\Dcat^b_{\sankaku}(\cnum_X)$
is a triangulated subcategory of
$\Dcat^b_{\realc}(\cnum_X)$.
\end{lem}
\pf
Let $K_1\lrarr K_2\lrarr K_3\lrarr K_1[1]$
be a distinguished triangle in $\Dcat^b_{\realc}(\cnum_X)$
such that $K_i$ $(i=1,2)$
are objects in $\Dcat^b_{\sankaku}(\cnum_X)$.
Let $\varphi:\Delta\lrarr X$
be any holomorphic map.
We set $K'_i:=\Ecat\varphi^{-1}(K_i)$.
We obtain a distinguished triangle
$K'_1\lrarr K'_2\lrarr K'_3\lrarr K'_1[1]$
in $\Dcat^b_{\realc}(\cnum_{\Delta})$.
For $i=1,2$,
$K'_i$ are cohomologically $\cnum$-constructible.
Then, $K_3'$ is also cohomologically $\cnum$-constructible.
Hence, 
we obtain $K_3\in \Dcat^b_{\sankaku}(\cnum_X)$.
\hfill\qed

\begin{lem}
\label{lem;16.8.9.1}
An object $K\in \Dcat^b_{\realc}(\cnum_X)$
is contained in
$\Dcat^b_{\sankaku}(\cnum_X)$
if and only if
$\nbigh^j(K)$ $(j\in\seisuu)$
are objects in $\Dcat^b_{\sankaku}(\cnum_X)$.
Here,
$\nbigh^j(K)$ denotes the $j$-th cohomology sheaf.
\end{lem}
\pf
Let us prove the ``if'' part.
For any $K$,
we set
$\ell(K):=\max\{n\,|\,\nbigh^n(K)\neq 0\}
-\min\{n\,|\,\nbigh^n(K)\neq 0\}$.
We use an induction on $\ell(K)$.
Let $K$ be an object in $\Dcat^b_{\realc}(\cnum_X)$
such that
$\nbigh^j(K)$ are objects in
$\Dcat^b_{\sankaku}(\cnum_X)$.
If $\ell(K)=0$,
$K$ is clearly contained in 
$\Dcat^b_{\sankaku}(\cnum_X)$.
There exists $n$ such that
$K=\tau^{\leq n}(K)$
and $\nbigh^n(K)\neq 0$.
There exists the distinguished triangle
$K\lrarr \nbigh^n(K)
\lrarr \tau^{\leq n-1}(K)[1]\lrarr K[1]$
in $\Dcat^b_{\realc}(\cnum_{\Delta})$.
By the assumption of the induction,
$\tau^{\leq n-1}(K)[1]$
is an object in $\Dcat^b_{\sankaku}(\cnum_{\Delta})$.
Because $\Dcat^b_{\sankaku}(\cnum_{\Delta})$
is triangulated,
we obtain that $K\in \Dcat^b_{\sankaku}(\cnum_{\Delta})$.

Let us prove the ``only if'' part.
Let $K$ be an object in $\Dcat^b_{\sankaku}(\cnum_X)$.
Let $\varphi:\Delta\lrarr X$ be any holomorphic map.
The cohomology sheaves
$\nbigh^j\varphi^{-1}(K)$ on $\Delta$
are $\cnum$-constructible.
Because 
$\nbigh^j\varphi^{-1}(K)\simeq
 \varphi^{-1}\nbigh^j(K)$,
we obtain that
$\nbigh^j(K)\in \Dcat^b_{\sankaku}(\cnum_X)$.
\hfill\qed

\vspace{.1in}

Let $Y\subset X$ be a real analytic submanifold
with $\dim_{\real}Y\leq \dim_{\real}X-2$.
Let $K$ be an $\real$-constructible sheaf on $X$
such that
$K_{|X\setminus Y}$ and $K_{|Y}$ are local systems
on $X\setminus Y$ and $Y$, respectively.
\begin{lem}
\label{lem;16.8.9.2}
Suppose that $K\in \Dcat^b_{\sankaku}(\cnum_X)$.
Then, for any $P\in Y$,
there exists a neighbourhood $X_P$ of $P$ in $X$
such that one of the following holds;
(i) $K_{|X_P}$ is a local system on $X_P$,
(ii) $Y\cap X_P$ is a complex submanifold of $X_P$.
\end{lem}
\pf
We shall shrink $X$ and $Y$ around $P$.
We may assume that $Y$ is simply connected,
and 
$X=Y\times \openopen{-1}{1}^{a}$
as a $C^{\infty}$-manifold,
where $a=\dim_{\real}X-\dim_{\real}Y$.

Let us consider the case 
$\dim_{\real}Y<\dim_{\real}X-2$.
There exists a local system $L$ on $X$
with an isomorphism 
$\kappa_0:K_{|X\setminus Y}\simeq L_{|X\setminus Y}$.
There exists a natural morphism
of $\real$-constructible sheaves $\kappa:K\lrarr L$
such that $\kappa_{|X\setminus Y}=\kappa_0$.
Note that
$L\in \Dcat^b_{\sankaku}(\cnum_X)$.
By Lemma \ref{lem;16.8.9.1},
$\Cok(\varphi)$
and $\Ker(\varphi)$
are objects in $\Dcat^b_{\sankaku}(\cnum_X)$.
Note that 
$\Cok(\varphi)$
and $\Ker(\varphi)$
are local systems on $Y$.
Hence, at least one of the following holds;
(i) $\Cok(\varphi)=\Ker(\varphi)=0$,
 i.e., $\varphi$ is an isomorphism,
or 
(ii) $Y$ is a complex submanifold.

\vspace{.1in}
Let us consider the case
$\dim_{\real} Y=\dim_{\real}X-2$.
We set $L_0:=K_{|X\setminus Y}$.
Let $F:L_0\lrarr L_0$ denote the automorphism 
obtained as the monodromy along a loop
along $Y\times\{0\}$ with counter clock-wise direction.
Let $L_0=\bigoplus_{\alpha\in\cnum} \EE_{\alpha}L_0$
denote the generalized eigen decomposition,
i.e., it is compatible with the action of $F$
and for each $\alpha$ the restriction of $F-\alpha\id$
to $\EE_{\alpha}L_0$ is nilpotent.

Let $\iota:X\setminus Y\lrarr X$
denote the inclusion.
We obtain the decomposition
$K=\bigoplus_{\alpha\in\cnum} K_{\alpha}$,
where 
$K_{\alpha}\simeq \iota_{\ast}\EE_{\alpha}L_{0}$
$(\alpha\neq 1)$,
and 
$K_{1|X\setminus Y}\simeq \EE_1L_0$.
Let us prove that
$Y$ is a complex submanifold of $X$
in the case where
$\EE_{\alpha}L_0\neq 0$ for one of $\alpha\neq 1$.
Suppose that $Y$ is not a complex submanifold.
There exists a holomorphic map
$\varphi:\Delta\lrarr X$ such that
$\varphi^{-1}(Y)=\{t\in\real\}\cap\Delta$.
Then, $\varphi^{-1}(K)$ is not $\cnum$-constructible on $\Delta$,
which contradicts our assumption
$K\in \Dcat^b_{\sankaku}(\cnum_X)$.
Hence, we obtain that $Y$ is a complex submanifold.
Let us consider the case where
$\EE_{\alpha}L_0=0$ for any $\alpha\neq 1$.
We shall prove the claim of the lemma in this case
by using an induction on $L_0$.
If $\rank L_0=0$,
$K$ is a local system on $Y$,
and hence the claim of the lemma is clear.
If $\rank L_0>0$,
there exists an epimorphism
$L_0\lrarr \cnum_{X\setminus Y}$
because the monodromy is unipotent.
We obtain the induced morphism
$\iota_{\ast}L_0\lrarr \cnum_X$.
Let $\rho:K\lrarr \cnum_X$
denote the composition of the morphisms
$K\lrarr \iota_{\ast}L_0\stackrel{\iota_{\ast}\rho}{\lrarr}\cnum_X$.
Note that $\cnum_X\in \Dcat^b_{\sankaku}(\cnum_X)$.
Hence, by Lemma \ref{lem;16.8.9.1},
we obtain that
$\Ker(\rho)$ and $\Cok(\rho)$
are objects in
$\Dcat^b_{\sankaku}(\cnum_X)$.
Note that $\Cok(\rho)$ is a local system on $Y$.
Hence, if $\Cok(\rho)\neq 0$,
we obtain that $Y$ is a complex submanifold of $X$.
Suppose that $\Cok(\rho)=0$.
Note that 
$\Ker(\rho)_{|X\setminus Y}$ and
$\Ker(\rho)_{|Y}$ are local systems.
Moreover,
the monodromy of 
$\Ker(\rho)_{|X\setminus Y}$ is unipotent,
and 
$\rank\Ker(\rho)_{|X\setminus Y}
<\rank L_0$.
Hence, by using the assumption of the induction,
we obtain the claim of the lemma in this case.
Thus, the proof of Lemma \ref{lem;16.8.9.2}
is finished.
\hfill\qed

\vspace{.1in}

Let $Y$ and $X$ be as above.
Let $f$ be a real analytic function $Y\lrarr\real$
such that $df$ is nowhere vanishing.
Set $Y_{\geq 0}:=f^{-1}(\real_{\geq 0})$
and $Y_0:=f^{-1}(0)$.
Let $K$ be an $\real$-constructible sheaf on $X$
such that 
$K_{|X\setminus Y_{\geq 0}}$,
$K_{|Y_{\geq 0}\setminus Y_0}$
and 
$K_{|Y_0}$
are local systems on
$X\setminus Y_{\geq 0}$,
$Y_{\geq 0}\setminus Y_0$
and $Y_0$,
respectively.

\begin{lem}
\label{lem;16.9.2.20}
Suppose that
$K\in \Dcat^b_{\sankaku}(\cnum_X)$.
Then, for any $P\in Y_0$,
there exists a small neighbourhood $X_P$ of $P$ in $X$
such that 
$K_{|X_P\setminus Y_0}$ is a local system.
Moreover, one of the following holds;
(i) $K_{|X_P}$ is a local system on $X_P$,
(ii) $X_P\cap Y_0$ is a complex submanifold.
\end{lem}
\pf
Let $P$ be any point of $Y_0$,
and let $X_P$ be a small neighbourhood of $P$ in $X$.
Because
$X_P\setminus Y_{\geq 0}$ is simply connected,
there exists a local system $L_P$ on $X_P$
with an isomorphism
$L_{P|X_P\setminus Y_{\geq 0}}
 \simeq 
K_{|X_P\setminus Y_{\geq 0}}$.
We obtain a natural morphism
$\varphi:K_{|X_P}\lrarr L_P$.
We set $K_1:=\Ker\varphi$
and $K_2:=\Cok\varphi$.
Then, 
we obtain $K_{i|X\setminus Y_{\geq 0}}=0$,
and 
$K_{i|Y_{\geq 0}\setminus Y_0}$
and 
$K_{i|Y_0}$
are local systems 
on $Y_{\geq 0}\setminus Y_0$
and $Y_0$, respectively.
Moreover,
$K_i\in \Dcat^b_{\sankaku}(\cnum_X)$.
Then, we can easily deduce the claim.
\hfill\qed

\subsection{Complex analyticity}
\label{subsection;18.11.27.11}

\subsubsection{Support}
Let $K\in\Ecat^b_{\circledcirc}(\IC_X)$.
Let $\nbigs(K)$ denote the collection of 
closed subsets $C$ in $X$ such that
$K_{|X\setminus C}=0$.
Let $Z=\bigcap_{C\in \nbigs(K)}C$.

\begin{lem}
$K=0$ in 
$\Ecat^b_{\realc}(\IC_{(X\setminus Z,X)})$.
\end{lem}
\pf
We obtain the claim 
by using the curve test
and an argument in the proof of Lemma \ref{lem;16.8.8.1}.
\hfill\qed

\subsubsection{Complex analyticity of the support}

Let $Z^{\sm}_k$ denote the set of $k$-dimensional smooth points
of $Z$. 
By an argument in the proof of Lemma \ref{lem;16.8.8.1},
we can prove that there exists a subanalytic closed subset
$(Z^{\sm}_k)^{(1)}
\subset
 Z^{\sm}_k$
with $\dim_{\real} (Z^{\sm}_k)^{(1)}\leq k-1$
such that the following holds.
\begin{itemize}
\item
Let $\nbigc$ be any connected component of 
$Z^{\sm}_k\setminus(Z^{\sm}_k)^{(1)}$.
Then, there exists a tuple of subanalytic functions
$\Lambda(\nbigc)$ on $(\nbigc,X)$
such that 
$\pi^{-1}(\cnum_{\nbigc})\otimes K
=\bigoplus_{g_i\in\Lambda(\nbigc)}
 \cnum^{\Ecat}_X\overset{+}{\otimes}\cnum_{t\geq g_i}$.
\end{itemize}
In particular,
$K_{|\nbigc}$ comes from a local system on $\nbigc$.

\begin{lem}
$Z^{\sm}_k$ are complex analytic submanifolds of $X$.
In particular,
$Z^{\sm}_k$ are empty
unless $k$ is even.
\end{lem}
\pf
Suppose that $Z^{\sm}_k$ is not a complex submanifold.
There exists a point 
$P\in Z^{\sm}_k\setminus (Z^{\sm}_k)^{(1)}$ such that 
$Z^{\sm}_k$ is not a complex manifold at $P$.
Let $\nbigc$ be the connected component of
$Z^{\sm}_k\setminus(Z^{\sm}_k)^{(1)}$
which contains $P$.
There exists a holomorphic map
$\varphi:\Delta\lrarr X$ such that
$\varphi(0)=P$
and that
$\varphi^{-1}(\nbigc)$ is real one dimensional.
Then, the support of $E\varphi^{-1}K$ is 
strictly real one dimensional,
which contradicts the condition
$K\in \Ecat^b_{\sankaku}(X)$.
Hence, we obtain that
$Z^{\sm}_k$ is a complex submanifold of $X$.
\hfill\qed

\begin{lem}
\label{lem;16.9.2.11}
The closure $\overline{Z^{\sm}_k}$
of $Z^{\sm}_k$ in $X$
is a complex analytic subset in $X$.
\end{lem}
\pf
We set $Y_k:=\overline{Z^{\sm}_k}$.
Let $Y_k^{\sm}$ denote the set of smooth points,
and we set $\Sing(Y_k):=Y_k\setminus Y^{\sm}_k$.
By Theorem \ref{thm;16.9.2.1},
it is enough to prove that
$\dim_{\real}\Sing(Y_k)\leq k-2$.
We assume that $\dim_{\real}\Sing(Y_k)=k-1$,
and we shall deduce a contradiction.

Let us consider the case $k=\dim_{\real}Z$.
By Proposition \ref{prop;16.9.1.10},
there exists a subanalytic closed subset
$\gbigw\subset \Sing(Y_k)$
with $\dim \gbigw\leq k-2$
such that the following holds
for any $P\in \Sing(Y_k)\setminus \gbigw$.
\begin{itemize}
\item 
$P$ is a smooth point of $\Sing(Y_k)$.
\item
$P$ is not contained in $Y_{k'}$ for any $k'<k$.
Note that $Y_{k-1}=\emptyset$.
\item
Let $X_P$ be a small neighbourhood
of $P$ in $X$.
Then, there exist closed complex submanifold
$\Ytilde_{k,P}$ of $X_P$ and 
a real analytic function 
$f:\Ytilde_{k,P}\lrarr\real$
whose exterior derivative is nowhere vanishing,
such that
$Y_k\cap X_P=\Ytilde_{k,P}\cap f^{-1}(\real_{\geq 0})$.
\end{itemize}
Then, there exists a holomorphic map
$\varphi:\Delta\lrarr X$
such that
$\varphi^{-1}(Z_k^{\sm})$
and 
$\varphi^{-1}(X\setminus Z)$
are non-empty open subsets of $\Delta$.
It contradicts the assumption that
$\Ecat\varphi^{-1}(K)$
comes from a cohomologically holonomic 
$\nbigd_{\Delta}$-complex.
Hence, we are done in the case $k=\dim_{\real}Z$.

Suppose that we have already known
that $Y_{k'}$ are complex subvarieties 
for any $k'>k$.
By using Corollary \ref{cor;16.9.2.10},
it is enough to prove that
$Y_{k}\setminus\bigcup_{k'>k}Y_{k'}$
is a complex subvariety of $X\setminus\bigcup_{k'>k}Y_{k'}$,
which can be argued
as in the case of $k=\dim_{\real}Z$.
Thus, the proof of Lemma \ref{lem;16.9.2.11}
is completed.
\hfill\qed

\begin{lem}
$Z^{\sm}_k$ is a complex analytic subset
of $X$.
\end{lem}
\pf
Note that 
$Z=\bigcup_{k\geq 0} \overline{Z^{\sm}_k}$,
and that 
$\overline{Z^{\sm}_k}\setminus
  Z^{\sm}_k
 =
  \Sing\bigl(
  \overline{Z^{\sm}_k}\bigr)
  \cup
    \bigcup_{j\neq k}
   \Bigl(
     \overline{Z^{\sm}_j}
     \cap
     \overline{Z^{\sm}_k}
   \Bigr)$
which is closed and complex analytic.
Hence, we obtain the claim of the lemma.
\hfill\qed

\subsubsection{Complex analyticity of the singular locus}
\label{subsection;18.11.16.102}

\begin{lem}
\label{lem;16.9.2.50}
There exists a closed subset
$(Z^{\sm}_k)^{(1)}_0
\subset
(Z^{\sm}_k)^{(1)}$,
which is subanalytic in $X$,
with the following property.
\begin{itemize}
\item
$\dim_{\real}(Z^{\sm}_k)^{(1)}_0
\leq
 k-2$ holds,
and 
$(Z^{\sm}_k)^{(1)}_0$
contains the singular locus of
$(Z^{\sm}_k)^{(1)}$.
\item
For any point 
 $P\in (Z^{(\sm)}_k)^{(1)}
 \setminus
 (Z^{\sm}_k)^{(1)}_0$,
and for any connected component $\nbigc$ 
 of $Z^{\sm}_k\setminus (Z^{\sm}_k)^{(1)}$
 such that $P$ is contained in the closure of $\nbigc$,
 the functions $g\in\Lambda(\nbigc)$ are bounded
 around $P$.
\item
$\pi^{-1}(\cnum_{(Z^{\sm}_k)^{(1)}})\otimes K$
is controlled by bounded functions
around any point
$P\in(Z^{\sm}_k)^{(1)}\setminus 
	(Z^{\sm}_k)^{(1)}_0$.
\end{itemize}
\end{lem}
\pf
Let $\nbigc$ be a connected component of
$Z^{\sm}_k\setminus (Z^{(\sm)}_k)^{(1)}$.
By Lemma \ref{lem;16.7.20.1},
there exists $\gbigz_1(\nbigc)\subset \del\nbigc$ such that
(i)
$\dim_{\real}\gbigz_1(\nbigc)\leq k-2$,
(ii) $\gbigz_1(\nbigc)$ contains the singular locus of $\del\nbigc$,
(iii) around any $P\in \del\nbigc\setminus \gbigz_1(\nbigc)$,
$g_i\in\Lambda(\nbigc)$ are ramified real analytic.

Let $P$ be any point of
$(Z^{\sm}_k)^{(1)}
 \setminus
\bigcup_{\nbigc} \gbigz_1(\nbigc)$.
Let $\nbigc$ be a connected component 
of $Z^{\sm}_k\setminus (Z_k^{\sm})^{(1)}$
such that the closure of $\nbigc$ contains $P$.
Suppose that $g_{i_0}\in\Lambda(\nbigc)$ 
is not bounded around $P$.
It is ramified real analytic at $P$.
Then, there exists a neighbourhood $U$ of $P$
in $(Z^{\sm}_k)^{(1)}$ 
such that $g_{i_0}$ is not bounded
around any $Q\in U$.
We can find a holomorphic map $\varphi:\Delta\lrarr X$
such that (i) $\varphi(0)=P$,
(ii) $\varphi(\Delta)\subset Z^{\sm}_k$,
(iii) $\varphi^{-1}((Z^{\sm}_k)^{(1)})$
is real $1$-dimensional.
The function
$\varphi^{-1}g_{i_0}$ is unbounded
around any point of
$\varphi^{-1}((Z^{\sm}_k)^{(1)})$.
However, it contradicts the curve test.
Hence, we obtain the claim of Lemma \ref{lem;16.9.2.50}.
\hfill\qed

\begin{lem}
\label{lem;16.8.8.2}
The cohomology sheaves of
the restriction of $K$
to $Z^{\sm}_k\setminus (Z^{\sm}_k)^{(1)}_0$
are local systems.
\end{lem}
\pf
Let $U_P$ be a neighborhood of $P$ in $Z^{\sm}_k$.
We set $W_P:=U_P\cap(Z^{\sm}_k)^{(1)}$,
and $U_P^{\circ}:=U_P\setminus W_P$.
We may assume that
$W_P$ and 
the connected components of $U_P^{\circ}$ 
are simply connected.
By construction, the restrictions of $K$ to
$U_P^{\circ}$ and $W_P$ come from
the direct sums of shifts of local systems.
By Lemma \ref{lem;16.9.2.50},
$\pi^{-1}(\cnum_{U_P^{\circ}})\otimes K_{|U_P}$
and 
$\pi^{-1}(\cnum_{W_P})\otimes K_{|U_P}$
come from cohomologically $\real$-constructible
complexes.
By \cite[Proposition 4.7.15]{DAgnolo-Kashiwara1},
we can conclude that
$K_{|U_P}$ comes from 
a cohomologically $\real$-constructible complex
on $U_P$.
We may assume
$U_P=U_{P,0}\times \Delta$
and $W_P=U_{P,0}\times\{w\in\Delta\,|\,\Image(w)=0\}$.
By the curve test,
we obtain that $K_{|U_{P,0}\times\{w\}}$
comes from a direct sum of 
shifts of local systems.
Then, the claim of the lemma follows.
\hfill\qed

\vspace{.1in}
By Lemma \ref{lem;16.8.8.2},
there exists a closed subset
$A\subset Z^{\sm}_k$
such that 
(i) $A$ is subanalytic in $X$
(ii) $\dim_{\real}A\leq k-2$,
(iii) $\nbigh^i(K_{|Z^{\sm}_k\setminus A})$
come from local systems on $Z^{\sm}_k\setminus A$.
We assume that $A$ is the minimum
among such closed subanalytic subsets.
Let $A^{\sm}_m$ denote the set of
the $m$-dimensional smooth points of $A$.

\begin{prop}
\label{prop;16.10.9.10}
$A^{\sm}_m$ is 
a complex submanifold of
$Z^{\sm}_k$.
In particular,
$A^{\sm}_m=\emptyset$
unless $m$ is even.
\end{prop}
\pf
We set 
$Y:=Z_k^{\sm}$ and $H:=A_m^{\sm}$
to simplify the notation.
Let $P$ be any point of $H$.
We assume that $T_PH$
is not a $\cnum$-subspace of $T_PY$,
and we shall derive a contradiction.
By shrinking $Y$,
we may assume that
$T_{P}H$ are not $\cnum$-subspaces of $T_{P}Y$
for any $P\in H$.

Let us consider the case where $k-m=2$.
Let $P$ be any point of $H$.
We obtain the real vector subspace
$T_PH\subset T_PY$.
There exists a complex line
$\nbigv\subset T_PY$
such that $\nbigv\cap T_PH=\{0\}$
by Lemma \ref{lem;16.8.10.1}.
Let 
$\Phi:
 \Delta_z\times\Delta_{\vecw}^{m/2}
\lrarr
 Y$
be any holomorphic map
such that 
(i) $T_{(z,\vecw)}\Phi:
 T_{(z,\vecw)}\Delta_z\times\Delta_{\vecw}
\lrarr
 T_{\Phi(z,\vecw)}Y$ 
are isomorphisms for any $(z,\vecw)$,
(ii) $T_{(0,0)}\Phi(T_0\Delta_z)=\nbigv$.
By shrinking $Y$,
we may assume that
$\Phi$ is an isomorphism,
by which we identify 
$Y$ and $\Delta_z\times\Delta^{m/2}_{\vecw}$.
We may also assume the existence of 
a real analytic function
$h:\Delta^{m/2}_{\vecw}\lrarr\Delta$ such that
$\Phi^{-1}(H)\subset\Delta_z\times\Delta^{m/2}_{\vecw}$
is the graph of $h$.
We set $\zeta(z,\vecw):=z-h(\vecw)$.
The restriction of 
$\zeta$ to $\Delta_z\times\{\vecw\}$
is holomorphic.
There exists the open embedding
$Y\lrarr \cnum_{\zeta}\times\cnum_{\vecw}^{n-1}$
defined by  $(\zeta,\vecw)$.

Let $\varpi:\Ytilde(H)\lrarr Y$
be the oriented real blowing up along $H$.
Let $\zeta=re^{\sqrt{-1}\theta}$
be the polar decomposition.
Then,
$(r,\theta,\vecw)$ is a local coordinate system
around $Q$
for any $Q\in\del\Ytilde(H)$.

There exists a closed subanalytic subset
$\Ytilde(H)^{(1)}\subset\Ytilde(H)$
with $\dim_{\real}(\Ytilde(H)^{(1)})=\dim_{\real}Y-1$
such that
(i) $\del\Ytilde(H)\subset \Ytilde(H)^{(1)}$,
(ii) there exist continuous subanalytic functions
$g_1,\ldots,g_p$
on $(\Ytilde(H)\setminus\Ytilde(H)^{(1)},\Ytilde(H))$,
and integers $m_i$,
such that
\[
 \pi^{-1}(\cnum_{\Ytilde(H)\setminus \Ytilde(H)^{(1)}})
\otimes
 \Ecat\varpi^{-1}(K)
\simeq
 \bigoplus_{i=1}^p
 \cnum^{\Ecat}\overset{+}{\otimes}
 \cnum_{t\geq g_i}[m_i].
\]
Let $\gbigw$ denote the closure of
$\Ytilde(H)^{(1)}\setminus \del\Ytilde(H)$.
There exists a closed subanalytic subset
$\gbigr\subset \del\Ytilde(H)$
and $\gbigz\subset H$ 
such that the following holds.
\begin{itemize}
\item
 $\dim_{\real} \gbigr<\dim_{\real}\del\Ytilde(H)$
and $\dim_{\real}\gbigz<\dim_{\real}(H)$.
\item
$\gbigr$ contains $\gbigw\cap\del\Ytilde(H)$.
Moreover,
$g_1,\ldots,g_p$ are ramified real analytic around 
any $Q\in \del\Ytilde(H)\setminus \gbigr$.
\item
$\gbigr\setminus\varpi^{-1}(\gbigz)\lrarr H$ is horizontal.
\end{itemize}
By enlarging $\Ytilde(H)^{(1)}$ and $\gbigz$,
we may assume that $\gbigr\setminus\varpi^{-1}(\gbigz)$ 
is contained in $\gbigw$.

Let $P_1=(z_1,\vecw_1)$ be 
any point of $H\setminus \gbigz$,
and let $Q_1$ be any point of 
$\varpi^{-1}(P_1)\setminus \gbigr$.
Let us prove that 
any $g_i$ are bounded around $Q_1$.
Suppose that 
some $g_{i_0}$ is unbounded around $Q_1$,
and we shall derive a contradiction.
We may assume that
$g_{i_0}=r^{-e/\rho}g^{(0)}_{i_0}$,
where $g^{(0)}_{i_0}$ is nowhere vanishing
analytic function of
$(r^{1/\rho},\theta,\vecw)$,
and $e$ and $\rho$ are positive integers.

Note that
$\varpi^{-1}(P_1)$
is identified with
$(T_{P_1}Y/T_{P_1}H\setminus\{0\})/\real_{>0}$.
Let $u$ be any element of 
$T_{P_1}Y/T_{P_1}H\setminus\{0\}$
such that $[u]=Q_1$.
Because it is assumed that
$T_{P_1}H$ is not a $\cnum$-subspace,
there exists $v\in T_{P_1}H$
such that 
$\sqrt{-1}v$ 
is mapped to
$u$ by the projection
$T_{P_1}Y\lrarr T_{P_1}Y/T_{P_1}H$
by Lemma \ref{lem;16.8.12.1}.
Let 
$\gamma:\openopen{-\epsilon_1}{\epsilon_1}
\lrarr
 H$
any real analytic map such that
$\gamma(0)=P_1$
and 
$\gamma'(0)=v$.
Let $\varphi:\Delta_{\epsilon_2}\lrarr Y$
be the complexification of $\gamma$.
Let $x+\sqrt{-1}y$ be the real coordinate system
on $\Delta_{\epsilon_2}$.
By construction,
$T_{(0,0)}\varphi(\del_x)=\gamma'(0)=v$
and $T_{(0,0)}\varphi(\del_y)=\sqrt{-1}v$
hold.
Hence,
$\varphi^{-1}(H)=
 \Delta_{\epsilon_2}\cap\real$.
Moreover, we obtain
$\varphi(\{(x,y)\in\Delta_{\epsilon}\,|\,y>0\})
\subset
 \Ytilde(H)\setminus \Ytilde(H)^{(1)}$,
and 
$\varphi^{\ast}(g_{i_0})$ is unbounded
around any point of
$\Delta_{\epsilon_2}\cap\real$.
It contradicts 
the assumption that
$E\varphi^{-1}K$ comes from
a holonomic $\nbigd$-module.
Hence, we obtain that 
$g_{i}$ are bounded around 
$Q_1\in \varpi^{-1}(P_1)\setminus\gbigr$.

Then, we can derive the following lemma easily.
\begin{lem}
\label{lem;16.8.12.2}
For any point $P'\in H\setminus \gbigz$,
let $\varphi_{P'}:\Delta_{\epsilon}\lrarr Y$
be the holomorphic map defined by
$\varphi_{P'}(a)=(a+h(P'),P')$.
Then, $\Ecat\varphi_{P'}^{-1}K$
comes from a regular holonomic $\nbigd$-module.
\hfill\qed
\end{lem}

We may regard
$\Ytilde(H)$ as an open subset
of $\{0\leq r\}\times S^1\times H$.
Let $\eta: \{0\leq r\}\times S^1\times H
\lrarr \{0\leq r\}\times H$
be the projection.
There exists a closed subanalytic subset
$\gbigb\subset \{0\leq r\}\times H$
with the following property.
\begin{itemize}
\item
$\dim_{\real}\gbigb=\dim_{\real}H$ holds,
and $\gbigb$ contains $\{0\}\times H$.
\item
 Each connected component of
 $(\{0\leq r\}\times H)\setminus \gbigb$
 is simply connected.
\item
 $\Ytilde(H)^{(1)}\setminus \eta^{-1}(\gbigb)\lrarr 
 \{0\leq r\}\times H$
 is proper and a local diffeomorphism.
\item
 There exist subanalytic functions
 $g^{(1)}_i$ on 
 $(\Ytilde(H)^{(1)}\setminus \eta^{-1}(\gbigb),\Ytilde(H))$,
and integers $m^{(1)}_i$,
 such that
\[
 \pi^{-1}(\cnum_{\Ytilde(H)^{(1)}\setminus \eta^{-1}(\gbigb)})
\otimes
 K
\simeq
 \bigoplus_{i=1}^m
 \cnum^{\Ecat}\overset{+}{\otimes}
 \cnum_{t\geq g^{(1)}_i}[m^{(1)}_i].
\]
\end{itemize}
On each connected component
$\nbigc$ of $(\{0\leq r\}\times H)\setminus \gbigb$,
the set
$\eta^{-1}(\nbigc)\cap \Ytilde(H)^{(1)}$
is described as the union of 
the graphs of subanalytic functions 
$h^{\nbigc}_i:\nbigc\lrarr S^1$.
Set $\gbigb_1:=\gbigb\setminus(\{0\}\times H)$.
We set 
$\gbigz_1:=\gbigz\cup (\overline{\gbigb_1}\cap (\{0\}\times H))$.

By Lemma \ref{lem;16.8.12.2},
the restriction of 
$g_{i}$ to 
$\varpi^{-1}(\Delta_{z}\times \{P_2\})
 \cap
 (\Ytilde(H)\setminus\Ytilde(H)^{(1)})$
are bounded
for any $P_2\in \Delta^{n-1}_{\vecw}\setminus \gbigz_1$.
The restrictions of  $g^{(1)}_i$ to 
$\varpi^{-1}(\Delta_{z}\times \{P_2\})
 \cap
 (\Ytilde(H)^{(1)}\setminus\del\Ytilde(H))$
are also bounded
for any $P_2\in \Delta^{n-1}_{\vecw}\setminus \gbigz_1$.
According to Lemma \ref{lem;16.8.12.11},
if $\gbigz_1$ is appropriately enlarged,
for any $P_2\in \Delta^{n-1}_{\vecw}\setminus \gbigz_1$,
there exists a neighbourhood $U$ of $P_2$
in $\Delta^{n-1}_{\vecw}$
such that the restrictions of 
$g_{i}$ to 
$\varpi^{-1}(\Delta_{z}\times U)
 \cap
 (\Ytilde(H)\setminus\Ytilde(H)^{(1)})$
are bounded,
and that the restrictions of 
$g^{(1)}_i$ to 
$\varpi^{-1}(\Delta_{z}\times U)
 \cap
 (\Ytilde(H)^{(1)}\setminus\del\Ytilde(H))$
are bounded.

Then, there exists a closed subanalytic subset
$\gbigz_2\subset H$
with $\dim_{\real}\gbigz_2<\dim_{\real} H$
such that the following holds.
\begin{itemize}
\item
 For any $P\in H\setminus \gbigz_2$,
 there exists  a neighbourhood $Y_P$ of $P$ in $Y$
 such that 
 $K_{|Y_P}$ comes from 
 an object in $\Dcat^b_{\realc}(Y_P)$.
\end{itemize}
Then, $H$ has to be a complex submanifold
by Lemma \ref{lem;16.8.12.20},
and we obtain a contradiction
in the case $k-m=2$.

\vspace{.1in}

Let us consider the case $k-m>2$.
Let $\varphi:\Delta\lrarr Y$ be any holomorphic map
such that $\varphi(0)\in H$,
and $\varphi(\Delta)\not\subset H$.
Because $\dim_{\real} H<\dim_{\real}Y-2$,
there exists a holomorphic map
$\Phi:\Delta\times \Delta\lrarr Y$
such that 
(i) $\Phi_{|\Delta\times\{0\}}=\varphi$,
(ii) $\dim_{\real}\Phi^{-1}(H)\leq 1$.
Note that $\Phi^{-1}(H)$ is real analytic.
Then, there exists a complex curve
$\nbigz\subset\Delta\times\Delta$
which contains $\Phi^{-1}(H)$.
Let $j:
\nbigy=(\Delta^2\setminus\nbigz,\Delta^2)
\lrarr
 \Delta^2$
be the inclusion.
We obtain
$\Ecat j_{!!}\Ecat j^{-1}
 \Ecat\Phi^{-1}K
\in \Ecat_{\circledcirc}(\IC_{\nbigy})$.
It comes from an object
in $\Mero_!(\Delta^2,\nbigz)$.
By considering the restriction to
a neighbourhood of any point of $\nbigz\setminus\Phi^{-1}(H)$,
we obtain that 
$\nbigv$ is regular singular.
Hence, $\Ecat\varphi^{-1}K$
comes from a regular singular holonomic $\nbigd$-complex
on $\Delta$.

There exists a closed subanalytic subset $\gbigz_{10}\subset H$
such that
(i) $\dim_{\real}\gbigz_{10}<\dim_{\real}H$,
(ii) $K_{|H\setminus \gbigz_{10}}$ comes from a local system.
Then, we obtain that
$K_{|Y\setminus \gbigz_{10}}$ comes from
an $\real$-constructible complex on $Y$.
Then, we obtain that 
$H\setminus \gbigz_{10}$ has to be a complex submanifold
by Lemma \ref{lem;16.8.9.2},
and we arrived at a contradiction.
Thus, the proof of Proposition \ref{prop;16.10.9.10}
is completed.
\hfill\qed

\begin{lem}
\label{lem;16.9.2.30}
The closure $\Abar$ of $A$ in $X$
is a complex analytic subvariety of $X$.
\end{lem}
\pf
Let $\overline{A^{\sm}_m}$
denote the closure of $A^{\sm}_m$ in $X$.
It is enough to prove that
$\overline{A^{\sm}_m}$
are complex analytic subvarieties of $X$.
By Corollary \ref{cor;16.9.2.10},
it is enough to prove that 
$A'_m:=\overline{A^{\sm}_m}\cap Z_k^{\sm}$ 
is a complex subvariety in $Z_k^{\sm}$.
Let $\Sing(A'_m)$ denote the singular locus of $A'_m$.
We have only to prove that
$\dim_{\real}\Sing A'_m<m-1$.

Let us consider the case $m=\dim_{\real} A$.
We assume that 
$\dim_{\real}\Sing A'_m=m-1$,
and we shall deduce a contradiction.
There exists a closed subanalytic subset
$\gbigw_m'\subset \Sing A'_m$
with $\dim_{\real}\gbigw_m'<m-1$
such that the following holds
for any $P\in \Sing A'_m\setminus \gbigw_m'$:
\begin{itemize}
\item
 $P$ is a smooth point of $\Sing A_m'$.
\item
 Let $X_P$ be a small neighbourhood 
 of $P$ in $X$.
 Then,  there exist a complex submanifold
 $\Atilde'_{m,P}$ of $X_P$
 and a real analytic function
 $f_P:\Atilde'_{m,P}\lrarr\real$
 whose exterior derivative is nowhere vanishing,
 such that
 $A'_m\cap X_P=\Atilde'_{m,P}\cap f_P^{-1}(\real_{\geq 0})$.
\end{itemize}
We set $Z_{k,P}:=Z_k\cap X_P$
and $A'_{m,P}:=A'_m\cap X_P$.

\begin{lem}
There exists a closed subanalytic subset
$\gbign_P\subset\Sing A_{m,P}'$
with $\dim \gbign_P<m-1$
such that the following holds.
\begin{itemize}
\item
For any $P_1\in \Sing A_{m,P}'\setminus\gbign_P$,
there exists a small neighbourhood $X_{P_1}$ of $P_1$ in $X$
such that 
$K_{|X_{P_1}}$ comes from
an object of $\Dcat^b_{\realc}(\cnum_{X_{P_1}})$.
\end{itemize}
\end{lem}
\pf
There exists a filtration
$Z_{k,P}=Z_{k,P}^{(0)}\supset Z_{k,P}^{(1)}
 \supset\cdots$
by closed subanalytic subsets
for $K_{|X_P}$,
where $Z_{k,P}^{(i)}\setminus Z_{k,P}^{(i+1)}$
are real analytic submanifolds of $Z_{k,P}$
of codimension $i$.
We may assume that
for each connected component $\nbigc$
of $Z_{k,P}^{(i)}\setminus Z_{k,P}^{(i+1)}$,
there exist continuous subanalytic functions
$h^{\nbigc}_j$ on $(\nbigc,X_P)$
such that 
$\pi^{-1}(\cnum_{\nbigc})\otimes K_{|X_P}
=\bigoplus \cnum^{\Ecat}_{X_P}\overset{+}{\otimes}
 \cnum_{t\geq h^{\nbigc}_j}$.
We may assume that
for each connected component $\nbigc$ satisfies
one of the following:
(i) $\nbigc\subset Z_{k,P}\setminus \Atilde'_{m,P}$, 
(ii)
$\nbigc\subset \Atilde'_{m,P}\setminus \Sing A'_{m,P}$,
(iii) $\nbigc\subset \Sing\Atilde'_{m,P}$.

Let $\nbigc$ be a connected component of
 $Z_{k,P}^{(i)}\setminus Z_{k,P}^{(i+1)}$
such that $\nbigc\subset Z_{k,P}\setminus \Atilde'_{m,P}$.
Let us observe that $h^{\nbigc}_j$ are bounded.
Let $\gamma:(\II^{\circ},0,\II)\lrarr 
 (\nbigc,A_m',Z_{k,P})$
be an analytic path.
Let $\gamma_{\cnum}:\Delta_{\epsilon}
 \lrarr Z_{k,P}$ be the induced holomorphic map.
Let $\nbigh_P\subset X_P$ be 
any complex hypersurface 
such that 
(i) the image of $\gamma_{\cnum}$ 
is not contained in $\nbigh_P$,
(ii) $\Atilde'_{m,P}\subset \nbigh_P$.
Note that
the set 
$(Z_{k,P}\cap \nbigh_P)\setminus A'_m$ is non-empty
because $\Atilde'_{m,P}\setminus A'_m$ is non-empty.
Let $j_{\nbigh_P}:(X_P\setminus \nbigh_P,X_P)\lrarr (X_P,X_P)$
be the inclusion of the bordered spaces.
By Theorem \ref{thm;16.5.11.20},
each cohomology of
$\Ecat j_{\nbigh_P!!}\Ecat j_{\nbigh_P}^{-1}K_{|X_P}$
comes from an object $\nbigm$
in $\Mero_!(Z_{k,P},\nbigh_P\cap Z_{k,P})$.
For each point $P'$ of $(\nbigh_P\cap Z_{k,P})\setminus A'_m$,
the restriction of $\nbigm$ to a neighbourhood of $P'$
is regular singular.
Hence, we obtain that $\nbigm$ is regular singular.
It implies that
$\gamma_{\cnum}^{\ast}\nbigm$
is regular singular.
We obtain that 
the functions $h^{\nbigc}_j$ are bounded along $\gamma$.
Hence, we obtain that $h^{\nbigc}_j$ are bounded.
Then, by using the argument 
in the proof of Lemma \ref{lem;16.9.2.50},
we obtain the existence of 
a closed subanalytic subset $\gbign_P$
with the desired property.
\hfill\qed

\vspace{.1in}

After enlarging $\gbigw_m'$,
we may assume that each cohomology 
$K^i_P$ of $K_{|X_P}$ satisfies that
(i) $K^i_{P|Z_{k,P}\setminus A'_m}$
 is a local system on $Z_{k,P}\setminus A'_m$,
(ii) $K^i_{P|X_P\cap (A'_m\setminus \Sing(A'_m))}$
 is a local system on 
 $X_P\cap (A'_m\setminus \Sing(A'_m))$,
(iii) $K^i_{P|X_P\cap \Sing(A'_m)}$
 is a local system on $X_P\cap \Sing(A'_m)$.
By Lemma \ref{lem;16.9.2.20},
we obtain that $K^i_{|Z_{k,P}\setminus\Sing(A'_m)}$
are local systems on $X_P$,
which contradicts our choice of $A$.
Hence, we obtain
$\dim_{\real}\Sing A'_m<m-1$,
and hence
$A'_m$ is a complex analytic subvariety
in the case $m=\dim_{\real}A$.

Assume that we have already known that
$A'_{m'}$ are complex analytic subvarieties
for $m'>m$.
By Corollary \ref{cor;16.9.2.10},
it is enough to prove that
$A'_m\setminus\bigcup_{m'>m}A'_{m'}$
is a complex analytic subvariety
of $X\setminus \bigcup_{m'>m}A'_{m'}$,
which can be argued
as in the case of $m=\dim A$.
Thus, we obtain Lemma \ref{lem;16.9.2.30}.
\hfill\qed

\vspace{.1in}

Let $P$ be any point of $Z$.
Let $X_P$ be a small neighbourhood of $P$,
and let $X_P\cap Z=\bigcup Z_i$ 
be the irreducible decomposition
as a germ of complex analytic sets at $P$.
For each $Z_i$,
there exists a complex hypersurface $H_i$ of $X_P$
such that 
(i) $Z_i\not\subset H_i$,
(ii) $\bigcup_{j\neq i}Z_j\subset H_i$,
(iii) $Z_i\setminus H_i$ is smooth,
(iv) $K_{|Z_i\setminus H_i}$ comes from 
 an $\real$-constructible complex
 whose cohomology sheaves are local systems.
Let $\rho_i:\Ztilde_i\lrarr Z_i$
be a projective morphism 
such that
(i) $\Ztilde_i$ is a complex manifold,
(ii) $\Htilde_i:=\rho_i^{-1}(H_i)$ is normal crossing,
(ii) $\Ztilde_i\setminus \Htilde_i\simeq Z_i\setminus H_i$.
Let $\iota_i:Z_i\lrarr X$ be the inclusion.
Let $j_i:(\Ztilde_i\setminus \Htilde_i,\Ztilde_i)\lrarr \Ztilde_i$ 
be the inclusion of the bordered spaces.
We obtain $\Ecat j_{i!!}\Ecat (\iota_i\circ\rho_i)^{-1}K$
in $\Ecat^b_{\sankaku}(\IC_{\Ztilde_i})$.

\begin{lem}
There exists a cohomologically holonomic
$\nbigd_{\Ztilde_i}$-complex $V$
with an isomorphism
\[
 \DR^{\Ecat}V\simeq
 \Ecat j_{i!!}\Ecat(\iota_i\circ\rho_i)^{-1}K.
\]
\end{lem}
\pf
The $k$-th cohomology object
of 
$\Ecat j_{i!!}\Ecat(\iota_i\circ\rho_i)^{-1}K$
comes from an object in $\Mero_!(\Ztilde_i,\Htilde_i)$
with the shift of degree.
By the fully faithfulness of
the Riemann-Hilbert correspondence
\cite{DAgnolo-Kashiwara1},
we obtain that
$\Ecat j_{i!!}\Ecat(\iota_i\circ\rho_i)^{-1}K$
comes from a cohomologically holonomic
$\nbigd_{\Ztilde_i}$-complex. 
\hfill\qed

\vspace{.1in}
Let $k_i:(X_P\setminus H_i,X_P)\lrarr X_P$
be the inclusion of the bordered spaces.
We obtain the following from the previous lemma.

\begin{lem}
\label{lem;18.1.12.2}
$\Ecat k_{i!!}\Ecat k_i^{-1}K$
comes from a cohomologically holonomic
$\nbigd_{X_P}$-complexes.
\hfill\qed
\end{lem}

\subsection{Proof of Theorem \ref{thm;16.7.11.1}}
\label{subsection;18.11.27.12}

Let us prove Theorem \ref{thm;16.7.11.1}.
It is enough to prove the essential surjectivity.
For that purpose,
it is enough to prove the following
for any $K\in \Ecat^b_{\sankaku}(\IC_X)$.
\begin{itemize}
\item
 $\Upsilon^{\Ecat}(K)$ are objects 
 in $\Dcat^b_{\hol}(\nbigd_X)$.
\item
 The natural morphisms
 $\Phi_K:K\lrarr \Sol^{\Ecat}\Upsilon^{\Ecat}(K)$
 are isomorphisms.
\end{itemize}
We have only to check these properties
locally around any point of $X$.
We use an induction on the dimension of the support of $K$.

Let $Z$ be the support of $K$.
Let $P\in Z$.
Let $X_P$ denote a small neighbourhood of $P$ in $X$.
Let $Z_1$ denote the union of 
the $\dim_{\cnum}(Z)$-dimensional components of $Z\cap X_P$.
There exists a hypersurface $H_P$ of $X_P$
such that 
(i) any irreducible component of $Z_1$ is not contained in $H_P$,
(ii) $Z_1\setminus H_P$ is smooth,
(iii) $K_{|Z_1\setminus H_P}$
 comes from an $\real$-constructible complex
 whose cohomology sheaves are local systems.
Let $j:(X_P\setminus H_P,X_P)\lrarr X_P$
be the inclusion of the bordered spaces.
As in Lemma \ref{lem;18.1.12.2},
$K_1:=\Ecat j_{!!}\Ecat j^{-1}K$
comes from a cohomologically holonomic
$\nbigd_{X_P}$-complex.
Hence, we obtain
a cohomologically holonomic
$\nbigd_{X_P}$-complex $\nbigm$
with an isomorphism
$\Sol^{\Ecat}\nbigm
\simeq
 K_1$.
Note that
$\nbigm(\ast H_P)\simeq\nbigm$.
According to \cite{DAgnolo-Kashiwara1},
there exists a canonical isomorphism
$\Upsilon^{\Ecat}\Sol^{\Ecat}(\nbigm)\simeq\nbigm$.
There exists the following morphisms:
\[
\begin{CD}
 \Sol^{\Ecat}(\nbigm)
 @>{a_1}>>
 \Sol^{\Ecat}\Upsilon^{\Ecat}\Sol^{\Ecat}(\nbigm)
 @>{a_2}>{\simeq}>
  \Sol^{\Ecat}(\nbigm)
\end{CD}
\]
Here, $a_1$ is $\Phi_{\Sol^{\Ecat}(\nbigm)}$,
and $a_2$ is the isomorphism induced by
$\Upsilon^{\Ecat}\Sol^{\Ecat}(\nbigm)
 \simeq\nbigm$.
It is easy to check that
the restriction of 
$a_2\circ a_1$ to $X_P\setminus H_P$
is an isomorphism,
and hence 
that $a_2\circ a_1$ is an isomorphism.
Hence, we obtain that
$\Upsilon^{\Ecat}(K_1)$
is a cohomologically holonomic 
$\nbigd_{X_P}$-complexes,
and 
$\Phi_{K_1}:K_{1}\lrarr
 \Sol^{\Ecat}
 \Upsilon^{\Ecat}(K_{1})$
is an isomorphism.

There exists the natural morphism 
$K_1\lrarr K_{|X_P}$
in $\Ecat^b_{\sankaku}(X_{|P})$.
There exists the distinguished triangle
$K_1\lrarr K_{|X_P}
\lrarr
 K_2\lrarr K_{1}[1]$
in $\Ecat^b_{\sankaku}(X_{|P})$.
We can apply the assumption of the induction
to $K_2$.
Then, we obtain that
$\Upsilon^{\Ecat}K_{|X_P}$
comes from a cohomologically holonomic
$\nbigd_{X_P}$-complex,
and the natural morphism
$\Phi_{K_{|X_P}}:
 K_{|X_P}\lrarr \Sol^{\Ecat}\Upsilon^{\Ecat}(K_{|X_P})$
is an isomorphism.
Thus, the proof of Theorem \ref{thm;16.7.11.1}
is completed.
\hfill\qed

\vspace{.1in}

\noindent
{\em Address\\
Research Institute for Mathematical Sciences,
Kyoto University,
Kyoto 606-8502, Japan\\
takuro@kurims.kyoto-u.ac.jp
}

\end{document}